\documentclass[12pt]{book}
\usepackage{latexsym}
\usepackage{amsfonts}
\usepackage{amssymb}
\usepackage{amsmath}
\usepackage{mathrsfs}
\input xypic
\usepackage{hyperref}
\usepackage{tikz}
\usetikzlibrary{matrix,arrows}
\allowdisplaybreaks

\newcommand{\be}{\begin{equation}}
\newcommand{\ee}{\end{equation}}
\newcommand{\bea}{\begin{eqnarray}}
\newcommand{\eea}{\end{eqnarray}}
\newcommand{\bean}{\begin{eqnarray*}}
\newcommand{\eean}{\end{eqnarray*}}
\newcommand{\brray}{\begin{array}}
\newcommand{\erray}{\end{array}}

\newtheorem{dfn}{Definition}[section]
\newtheorem{thm}[dfn]{Theorem}
\newtheorem{lmma}[dfn]{Lemma}
\newtheorem{ppsn}[dfn]{Proposition}
\newtheorem{crlre}[dfn]{Corollary}
\newtheorem{xmpl}[dfn]{Example}
\newtheorem{rmrk}[dfn]{Remark}
\newtheorem{xrcs}{Exercise}[section]

\newcommand{\bdfn}{\begin{dfn}\rm}
\newcommand{\bthm}{\begin{thm}}
\newcommand{\blmma}{\begin{lmma}}
\newcommand{\bppsn}{\begin{ppsn}}
\newcommand{\bcrlre}{\begin{crlre}}
\newcommand{\bxmpl}{\begin{xmpl}}
\newcommand{\brmrk}{\begin{rmrk}\rm}

\newcommand{\edfn}{\end{dfn}}
\newcommand{\ethm}{\end{thm}}
\newcommand{\elmma}{\end{lmma}}
\newcommand{\eppsn}{\end{ppsn}}
\newcommand{\ecrlre}{\end{crlre}}
\newcommand{\exmpl}{\end{xmpl}}
\newcommand{\ermrk}{\end{rmrk}}

\newcommand{\bbc}{\mathbb{C}}

\newcommand{\bbn}{\mathbb{N}}
\newcommand{\bbr}{\mathbb{R}}

\newcommand{\bbt}{\mathbb{T}}

\newcommand{\cla}{\mathcal{A}}

\newcommand{\clh}{\mathcal{H}}

\newcommand{\clk}{\mathcal{K}}

\newcommand{\bxrcs}{\begin{xrcs}\rm\footnotesize}

\author{S. Sundar\\}
\title{Notes on $C^{*}$-algebras}
\date{}

\begin{document}
\maketitle

\newpage
\newpage
\newpage
\newpage

\chapter*{PREFACE} 

These  notes grew out of two courses given by the author on $C^{*}$-algebras, one at IMSC, Chennai and an another  at IIT, Gandhinagar. 
The audience  were graduate students interested in pursuing research in topics such as NCG, K-theory, Quantum groups, $E_0$-semigroups etc... 
The aim of  the courses was to acquaint the students with a few basic notions in the theory of $C^{*}$-algebras which, in the author's opinion, are 
indispensable to read current papers in the above mentioned topics. The topics discussed in this notes, we merely scratch the surface as the motivation is to make the reader converse in the language as opposed to giving a complete treatment, are universal $C^{*}$-algebras, group $C^{*}$-algebras, crossed products, Hilbert $C^{*}$-modules, Morita equivalence, and K-theory. 

To make these notes accessible to those who have completed a first course in functional analysis, the basic `Gelfand theory' is developed in the first chapter of these notes. 
The reader is highly recommended to read \cite{Arveson_invitation} and \cite{Arveson_spectral} to learn the basic theory in more detail and depth. 

The organisation of this notes is as follows. 

Preliminaries on $C^{*}$-algebras such as the Gelfand-Naimark theorem for commutative $C^{*}$-algebra, the GNS construction and the double commutant result are worked out in Chapter 1. There are certain omissions. Instead of  a proper treatment of approximate identities and a discussion on ideals and quotients, we have merely collected the results, without proofs, at the end of Chapter 1.

In Chapter 2, we discuss a few examples. We start with a treatment of finite dimensional $C^{*}$-algebras that culminates with the theorem that states that a finite dimensional $C^{*}$-algebra is a direct sum of matrix algebras.  The algebra of compact operators is next discussed and is realised as a universal $C^{*}$-algebra given in terms of generators and relations. This serves
as a model for the notion of universal $C^{*}$-algebras which we take up next. This allows us to quickly define group $C^{*}$-algebras and crossed products, two classes of examples extensively studied in the literature. Several important $C^{*}$-algebras studied in the literature has this universal prescription so it is appropriate to give a rigorous treatment. An important example of an universal $C^{*}$-algebra, the Toeplitz algebra is discussed in full detail. The fundamental short exact sequence of the Toeplitz algebra and the Wold-decomposition of a single isometry are derived.   After quickly reviewing the measure theoretic preliminaries in Section 2.5, we discuss in Section 2.6 group $C^{*}$-algebras associated to  locally compact groups.

In Chapter 3, we take up crossed products of $C^{*}$-algebras and Hilbert $C^{*}$-modules. After defining the full and reduced crossed product in Section 3.1, we introduce Hilbert $C^{*}$-modules in Section 3.2 as a tool to prove that the reduced $C^{*}$-norm is independent of the choice of the representation
that one chooses. The author believes that it is an appropriate point to introduce the notion of Hilbert $C^{*}$-modules to the reader. As an application of the machinery of crossed products, Stone-von Neumann theorem regarding the uniqueness of irreducible Weyl representations, or equivalently the description of the spectrum of the Heisenberg group is worked out in Section 3.3. As another application, Cooper's theorem, that provides the Wold decomposition picture of a strongly continuous one parameter semigroup of isometries, is proved in Section 3.4.   After a short discussion on the non-commutative torus, we   discuss Mackey's imprimitivity theorem in Section 3.5 cast in Rieffel's language. The notion of Morita equivalence is introduced and a proof of the imprimitivity theorem is provided in the discrete setting. 

Chapter 4 in itself constitutes a short crash course on $K$-theory. The material in this chapter were tried out at ATM schools held at IMSc over the years. Unlike the first three chapter, we get more economical, even stingy, with details. The reader should approach this chapter only with the intention to get an exposure to $K$-theory (otherwise, he/she will be left disappointed) and should turn to the excellent books (\cite{Bla}, \cite{Rordam}, \cite{Olsen}) already available for a more detailed study. The contents of this chapter are briefly described next. After deriving the basic properties of $K_0$ and $K_1$, the chapter culminates with the proof of Bott periodicity due to Cuntz. The inspiration for the treatment on $K$-theory, and also on universal $C^{*}$-algebras  come from the three lectures given by Prof. Cuntz during the Oberwolfach conference on semigroup $C^{*}$-algebras held in Oct 2014. It also borrows material from the Master's thesis of Prakash Kumar Singh, a former student of CMI, Chennai , done under  the author's supervision. 

There are several excellent resources to read about the material covered in this notes. The bibliography contains a sample list. It is certainly not exhaustive 
and I apologise sincerely for any omission. 
The author claims no originality for the material presented.

I would  like to end this short introduction by thanking a few people who have helped me immensely so far. First, I would like to thank my advisor Prof. Partha Sarathi Chakraborty and Prof. V. S. Sunder for teaching me several aspects of mathematics (and life) which has enriched my understanding of the subject (and has made me appreciate life better). I thank Prof. Arup Kumar Pal for having several discussions on $K$-theory and operator algebras with me over the years.  I thank Prof. Renault for several discussions, in particular about groupoids, a topic closer to the author's heart but left out in these notes as the author feels it is better to learn about groupoids from Renault's books directly.  I thank Bipul Saurabh for inviting me to give a short course on $C^{*}$-algebras at IIT Gandhinagar and I thank the students of IITG, in particular, Shreema and Akshay for attending the course and the enthusiasm shown which acted as a source of encouragement to complete these notes. Finally, I thank my collaborators and students namely Murugan, Anbu, Sruthy, Prakash, Piyasa and Namitha for reading these notes and for their feedback. 

\vspace{10 mm}
\noindent
{\sc S. Sundar}\\
(\texttt{sundarsobers@gmail.com})\\
         {\footnotesize  Institute of Mathematical Sciences (HBNI),\\ CIT Campus, Taramani,\\
  Chennai, 600113,\\  Tamilnadu, INDIA.}


\tableofcontents
\chapter{Preliminaries on $C^{*}$-algebras}
\section{Gelfand-Naimark theorems}
All the algebras that we consider are over $\bbc$. Moreover, they are assumed to be associative. Let $A$ be an algebra over $\bbc$. By a $^*$-structure on $A$, we mean a map $A \ni a \to a^{*} \in A$ such that the following conditions hold. 
\begin{enumerate}
\item[(1)] For $a \in A$, $(a^*)^*=a$, 
\item[(2)] for $a,b \in A$, $(a+b)^{*}=a^{*}+b^{*}$, 
\item[(3)] for $a \in A$ and $\lambda \in \bbc$, $(\lambda a)^{*}=\overline{\lambda}a^*$, and
\item[(4)] for $a,b \in A$, $(ab)^{*}=b^{*}a^{*}$. 
\end{enumerate}
An algebra with a $^*$-structure is called a $^{*}$-algebra. 

\begin{dfn}
Let $A$ be a $^*$-algebra and $||~||$ be a norm on $A$. The pair $(A,||~||)$ is called a $C^{*}$-algebra if 
\begin{enumerate}
\item[(1)] $(A,||~||)$ is a Banach space, 
\item[(2)] for $a,b \in A$, $||ab|| \leq ||a|| ||b||$, and
\item[(3)] for $a \in A$, $||a||^{2}=||a^{*}a||$. 
\end{enumerate}
\end{dfn}
The equality $||a||^{2}=||a^*a||$ is usually called the $C^{*}$-identity. Let $A$ be a $C^{*}$-algebra. If $A$ contains an element $1_A \in A$ such that $a.1_A=1_A.a=a$, then we say that $A$ is unital. 
Moreover, the element $1_A$ which is necessarily unique is called the multiplicative identity or the unit of $A$. As usual, we abuse notation and write $1_A$ simply as $`1'$. If $A$ is not unital, we call $A$ non-unital. 

If the $^*$-structure is absent, then we call such objects as Banach algebras. More precisely, let $A$ be an algebra with a norm $||~||$. Then, $(A,||~||)$ is called a \emph{Banach algebra} if $(1)$ and $(2)$
are satisfied. 

\begin{xmpl}
\label{unital commutative}
Let $X$ be a compact Hausdorff space. Let $C(X)$ denote the algebra of continuous complex valued functions on $X$. Define the $^*$-structure and the norm $||~||$ on $C(X)$ as follows.
\begin{align*}
f^*(x):&=\overline{f(x)},\\
||f||:&=\sup\{|f(x)|:x \in X\}.
\end{align*}
Then, $C(X)$ is a unital $C^{*}$-algebra.  
\end{xmpl}

\begin{xmpl}
Let $X$ be a locally compact Hausdorff space. A continuous function $f:X \to \bbc$ is said to vanish at `$\infty$' if for every $\epsilon>0$, there exists a compact set $K \subset X$ such that 
$|f(x)|<\epsilon$ whenever $x \notin K$. Denote the algebra of continuous complex valued functions on $X$ that vanish at $\infty$ by $C_0(X)$.  It is easy to prove that if $f \in C_0(X)$, then
$f$ is bounded. 

Define the $^*$-structure and the norm $||~||$ as in the previous example. Then, $C_0(X)$ is a $C^{*}$-algebra.  
\end{xmpl}

\begin{xrcs}
Let $X$ be a locally compact Hausdorff space. Prove that $C_0(X)$ is unital if and only if $X$ is compact. 
\end{xrcs}

\begin{xmpl}
Let $\clh$ be a Hilbert space. Denote the algebra of bounded linear operators on $\clh$ by $B(\clh)$. For $T \in B(\clh)$, recall that its adjoint $T^{*}$ is the unique linear operator on $\clh$ such that 
\[
\langle T\xi|\eta \rangle=\langle \xi|T^{*}\eta \rangle\]
for every $\xi,\eta \in \clh$. The map $B(\clh) \ni T \to T^{*} \in B(\clh)$ defines a $^*$-structure on $B(\clh)$. 

For, $T \in B(\clh)$, recall that the operator norm $||T||$ is defined by 
\[
||T||:=\sup\{||T\xi||: ||\xi||=1\}.\]
Then, $B(\clh)$ with the $^*$-structure and the norm defined above is a unital $C^{*}$-algebra. We verify the $C^{*}$-identity in the next proposition. 
\end{xmpl}

\begin{ppsn}
For $T \in B(\clh)$, $||T^{*}T||=||T||^{2}$. 
\end{ppsn}
\textit{Proof.} Note that $||T||=\sup\{|\langle T\xi|\eta \rangle|: ||\xi||=||\eta||=1\}$. This implies, in particular, that the map $T \to T^{*}$ is norm preserving. 

Let $T \in B(\clh)$ be given. Clearly, $||T^{*}T|| \leq ||T^{*}||||T||=||T||^2$. Let $\xi \in \clh$ be of unit norm. Calculate as follows to observe that 
\begin{align*}
||T\xi||^{2}&=\langle T\xi|T\xi \rangle \\
 &=\langle T^{*}T\xi| \xi \rangle \\
 & \leq |\langle T^{*}T\xi|\xi \rangle |\\
 & \leq ||T^{*}T\xi||~||\xi||\\
 & \leq ||T^{*}T||.
 \end{align*}
 Consequently, 
 \[
 ||T||^{2}=\sup\{||T\xi||^2:||\xi||=1\} \leq ||T^*T||.\]
 The proof is now complete. \hfill $\Box$
 
 \begin{xmpl}
 Let $A \subset B(\clh)$ be a vector subspace. Assume that 
 \begin{enumerate}
 \item[(1)] $A$ is closed under multiplication, 
 \item[(2)] for $a \in A$, $a^{*} \in A$, i.e. $A$ is $^*$-closed, and
 \item[(3)] the subset $A$ is closed in $B(\clh)$, when $B(\clh)$ is given the norm topology.
 \end{enumerate}
Then, with the norm and the $^*$-structure `restricted' to $A$, $A$ becomes a $C^{*}$-algebra. 

Such a $C^{*}$-algebra is called a $C^{*}$-subalgebra of $B(\clh)$. If in addition, the identity operator $1_{\clh} \in A$, $A$ is said to be a unital $C^{*}$-subalgebra of $B(\clh)$. More generally, we can define $C^{*}$-subalgebras of an arbitrary $C^{*}$-algebra. 
 \end{xmpl}
 
 \begin{dfn}
 \label{iso defn}
 Let $A$ and $B$ be two $C^{*}$-algebras. A linear map $\pi:A \to B$ is said to be an isomorphism if 
 \begin{enumerate}
 \item[(1)] the map $\pi$ is a bijection, 
 \item[(2)] the map $\pi$ is multiplicative, i.e. $\pi(ab)=\pi(a)\pi(b)$ for $a,b \in A$, 
 \item[(3)] $\pi$ preserves the $^*$-structure, i.e. $\pi(a^*)=\pi(a)^{*}$ for $a \in A$, and
 \item[(4)] the map $\pi$ is an isometry, i.e. $||\pi(a)||=||a||$ for every $a \in A$.
 \end{enumerate}
 We say $A$ and $B$ are isomorphic if there exists a linear map $\pi:A \to B$ such that $\pi$ is an isomorphism. 
  \end{dfn}
  
  \begin{rmrk}
  We will see later that $(1)-(3)$ in Defn. \ref{iso defn} automatically implies $(4)$. 
    \end{rmrk}
 
 The first two basic theorems that is proved before studying $C^{*}$-algebras in detail are stated below. They go under the name `Gelfand-Naimark theorems'.
 
 \begin{thm}[\textbf{Gelfand-Naimark theorem: commutative version}]
 \label{GNC}
 Let $A$ be a commutative $C^{*}$-algebra. Then, $A$ is isomorphic to $C_0(X)$ for some locally compact Hausdorff space $X$. Moreover, the space $X$ is unique up to a homeomorphism. 
   \end{thm}
 In  more fancy language, we can rephrase the above as follows. The map 
 \[
 \Big \{\textrm{~~locally compact Hausdorff spaces~~}\Big \} \ni X \to C_0(X) \in \Big\{~~\textrm{commutative $C^{*}$-algebras}\Big\}\]
 sets up a bijective correspondence between `the category' of locally compact Hausdorff spaces and `the category' of commutative $C^{*}$-algebras. 
 
 \begin{thm}[\textbf{Gelfand-Naimark theorem: non-commutative version}]
 \label{GNN}
 Let $A$ be a $C^{*}$-algebra. Then, $A$ is isomorphic to a $C^{*}$-subalgebra of $B(\clh)$ for some Hilbert space $\clh$. Moreover, if $A$ is unital, then $A$ is isomorphic to a unital $C^{*}$-subalgebra of $B(\clh)$ for some Hilbert space $\clh$. 
 \end{thm}
 
\section{Spectrum and Resolvent}
 
 For the rest of this section, the letter $A$ stands for a unital Banach algebra. We denote the multiplicative identity of $A$ either by $1_A$ or simply by $1$. For $a \in A$, \emph{the spectrum of $a$}, denoted $\sigma(a)$ is defined by 
 \[
 \sigma(a):=\{\lambda \in \bbc: a-\lambda 1_A \textrm{~~is not invertible}\}.\]
 The complement of $\sigma(a)$ is called \emph{the resolvent of $a$} and is denoted $\rho(a)$, i.e.
 \[
 \rho(a):=\{\lambda \in \bbc: a -\lambda \textrm{~~is invertible}\}.\]
The first question that arises is whether $\rho(a)$ and $\sigma(a)$ are non-empty.  It is relatively straightforward to show that $\rho(a)$ is non-empty as the next lemma shows. 

\begin{lmma}
\label{mother of lemmas}
Let $a \in A$ be given. 
If $|\lambda|>||a||$, then $a-\lambda$ is invertible. 
\end{lmma}
\textit{Proof.} Observe that the power series $\displaystyle- \sum_{n=0}^{\infty}\frac{a^n}{\lambda^{n+1}}$ converges and its limit is the inverse of $a-\lambda$. \hfill $\Box$

The series expansion in Lemma \ref{mother of lemmas} will be repeatedly used and it is useful to isolate it as a remark. 
\begin{rmrk}
\label{power series}
Let $a \in A$ be given. If $|\lambda|>||a||$, then 
\[
\frac{1}{\lambda-a}=\sum_{n=0}^{\infty}\frac{a^n}{\lambda^{n+1}}.\]

\end{rmrk}

\begin{xrcs}
\label{continuity of inverse}

Let \[G(A):=\{a \in A: a \textrm{~~is invertible}\}.\] Note that, if $||1-a||<1$, then, by Lemma \ref{mother of lemmas}, $a$ is invertible. Thus, $1=1_A$ is an interior point of $G(A)$. 
\begin{enumerate}
\item[(1)] Prove that every point of $G(A)$ is an interior point, or in other words, $G(A)$ is open. Deduce that for $a \in A$, the resolvent $\rho(a)$ is open and the spectrum $\sigma(a)$ is closed. 
\item[(2)] Show that the map $G(A) \ni a \to a^{-1} \in G(A)$ is continuous. 
\end{enumerate}
\end{xrcs} 

Let $U \subset \bbc$ be a non-empty open set and let $f:U \to A$ be a function. We say $f$ is holomorphic, or more correctly \emph{weakly holomorphic}, if for every $\phi \in A^{*}$, the complex valued function $\phi \circ f$ is holomorphic. Here, $A^{*}$ stands for the dual of $A$, i.e. the space of bounded linear functionals on $A$. 

\begin{ppsn}
\label{spectrum non-empty}
For every $a \in A$, the spectrum $\sigma(a)$ is a non-empty compact subset of $\bbc$. 
\end{ppsn}
\textit{Proof.} We have already seen that $\sigma(a)$ is closed. Moreover, Lemma \ref{mother of lemmas} implies that $\sigma(a)$ is contained in the closed unit ball of radius $||a||$. Thus,  it remains  to 
prove that $\sigma(a)$ is non-empty. 

Suppose not. Then, the resolvent $\rho(a)=\bbc$. Let $f:\mathbb{C} \to A$ be defined by 
\[
f(z):=\frac{1}{z-a}.\]
Let $\phi \in A^{*}$ be given. Note that for $z,w \in \bbc$,
\[
\frac{\phi(f(z))-\phi(f(w))}{z-w}=\phi\Big(\frac{-1}{(z-a)(w-a)}\Big).\]
Use Exercise \ref{continuity of inverse} to conclude that as $z \to w$, \[\frac{\phi(f(z))-\phi(f(w))}{z-w}\longrightarrow \phi\big(\frac{-1}{(w-a)(w-a)}\big).\] Thus, $\phi \circ f$ is entire.  

Thanks to Remark \ref{power series}, for $|z|>||a||$, 
\[
|\phi(f(z))| \leq \sum_{n=0}^{\infty}\frac{||a^n||}{|z|^{n+1}} \leq \sum_{n=0}^{\infty}\frac{||a||^n}{|z|^{n+1}} \leq \frac{1}{|z|-||a||}.\] 
The above forces that $\phi \circ f \to 0$ as $|z| \to \infty$. But, $\phi \circ f$ is entire. Therefore, $\phi(f(z))=0$ for every $z$. 

Since $\phi \circ f=0$ for every $\phi \in A^{*}$, we have $f=0$ which is absurd. Hence, our assumption that $\sigma(a)$ is empty is wrong and the proof is complete. \hfill $\Box$

For, $a \in A$, let 
\[r(a):=r=\sup\{|\lambda|: \lambda \in \sigma(a)\}.\] 
Note that $r(a)$ is a finite non-negative number. The number $r(a)$ is called \emph{the spectral radius of $a$}. The next proposition is a `quantitative version' of Prop. \ref{spectrum non-empty}. 

\begin{ppsn}
\label{spectral radius formula}
For $a \in A$, 
\[
r(a):=\lim_{n \to \infty}||a^n||^{\frac{1}{n}}.\]
\end{ppsn}

\begin{lmma}
\label{poly}
For a polynomial $p$ and $a \in A$, $\sigma(p(a))=\{p(\lambda):\lambda \in \sigma(a)\}$. 
\end{lmma}
\label{Proof.} Let $p$ be a polynomial and $a \in A$ be given.
 It suffices to prove that  $p(a)$ is invertible if and only if $p(\lambda) \neq 0$ for every $\lambda \in \sigma(a)$. 

 Factorise $p$ and write $p$ as 
\[
p(z)=(z-\lambda_1)(z-\lambda_2)\cdots (z-\lambda_n).\]
Then, \[p(a)=(a-\lambda_1)(a-\lambda_2)\cdots (a-\lambda_n).\]
Note the following chain of easily verifiable if and only if statements. Only the first `iff' requires any proof and we leave it to the reader.  
\begin{align*}
\textrm{~$p(a)$ is invertible} & \iff \textrm{~for every $i=1,2,\cdots,n$, $a-\lambda_i$ is invertible}\\
& \iff \textrm{~for every $i=1,2,\cdots,n$, $\lambda_i \in \rho(a)$}\\
& \iff \{\lambda \in \bbc: p(\lambda) =0\} \subset \rho(a).
\end{align*}
This completes the proof. \hfill $\Box$

\textit{Proof of Prop. \ref{spectral radius formula}.} We denote $r(a)$ simply by $r$. Let $f:\rho(a) \to A$ be defined by 
\[
f(z):=(1-az)^{-1}.\]
It can be proved as in Prop. \ref{spectrum non-empty} that $f$ is weakly holomorphic. 

Let $\phi \in A^{*}$ be given. 
Thanks to the power series expansion in Remark \ref{power series}, if $|z|$ is `small', then 
\begin{equation}
\label{pwr}
\phi(f(z))=\sum_{n=0}^{\infty}z^n\phi(a^n).
\end{equation}
But the disk $D_r:=\{z \in \bbc: |z|<\frac{1}{r}\}$ is contained in $\rho(a)$ and $\phi \circ f$ is analytic on $D_r$. Thus, the power series expansion of Eq. \ref{pwr} continues to hold for $z \in D_r$. 

Let $t \in (0,\frac{1}{r})$ be given. Then, the series $\sum_{n=0}^{\infty}t^n\phi(a^n)$ is convergent. In particular, $\phi(t^na^n) \to 0$ as $n \to \infty$. Thus, $t^na^n \to 0$ weakly. By the uniform boundedness principle, there 
exists $M_t >0$ such that $t^n||a^n|| \leq M_t$. Thus, $||a^n||^{\frac{1}{n}} \leq \frac{1}{t}M_t^{\frac{1}{n}}$. Taking $\limsup$, we get
\[
\limsup_{n \to \infty}||a^n||^{\frac{1}{n}} \leq \frac{1}{t}.\]
If we let  $t \to \frac{1}{r}$, we see that 
\[
\limsup_{n \to \infty}||a^n||^{\frac{1}{n}} \leq r.\]

Let $\lambda \in \sigma(a)$. Then, by Lemma \ref{poly}, $\lambda^n \in \sigma(a^n)$. Hence,  $|\lambda|^n \leq ||a^n||$, for every $n \geq 1$. This implies that $|\lambda| \leq ||a^n||^{\frac{1}{n}}$. Taking $\liminf$, we see that 
\[
|\lambda| \leq \liminf_{n \to \infty}||a^n||^{\frac{1}{n}}.\] 
Hence, $r \leq \liminf_{n \to \infty}||a^n||^{\frac{1}{n}}$. The proof is now complete.

\begin{ppsn}
\label{spectral radius for self-adjoint}
Suppose that $A$ is a unital $C^{*}$-algebra. 
Let $a \in A$ be self-adjoint, i.e. $a^{*}=a$. Then, the spectral radius of $a$ is $||a||$. 
\end{ppsn}
\textit{Proof.} The $C^{*}$-identity implies that $||a^2||=||a||^{2}$, $||a^4||=||a^2||^{2}=||a||^{4}$. In general, for any $n$, 
\[
||a^{2^n}||^{\frac{1}{2^n}}=||a||.\]
Therefore, 
\[
\lim_{n \to \infty}||a^n||^{\frac{1}{n}}=\lim_{n \to \infty}||a^{2^n}||^{\frac{1}{2^n}}=||a||.\]
Prop. \ref{spectral radius formula} gives the desired result. \hfill $\Box$

\begin{xrcs}
Let $A:=C(X)$ where $X$ is a compact Hausdorff space. Suppose  $f \in C(X)$. Prove that 
\[
\sigma(f)=\{f(x):x \in X\}.\]
\end{xrcs}

\begin{xrcs}
Let $A:=B(L^{2}(\bbr,d\lambda))$ where $d\lambda$ is the Haar measure. For a bounded measurable function $f:\bbr \to \bbc$, define a bounded operator $M_f$ on $L^{2}(\bbr,d\lambda)$ by 
\[
M_f\xi(t):=f(t)\xi(t).\]
Calculate the spectrum of $M_f$. 
\end{xrcs}

\begin{xrcs}
Let $A:=B(\ell^2(\bbn))$. Suppose $S$ is the `right-shift' operator on $\ell^2(\bbn)$, i.e. 
\[
S(x_0,x_1,x_2,\cdots)=(0,x_0,x_1,\cdots,).\]
Compute the spectrum of $S$ in $A$.
\end{xrcs}

\section{Commutative $C^{*}$-algebras}
In this section, we prove that every unital commutative $C^{*}$-algebra is isomorphic to the algebra of continuous functions on a compact space.

\begin{ppsn}[Gelfand-Mazur]
\label{Gelfand_Mazur}
Let $A$ be a unital Banach algebra. If every non-zero element of $A$ is invertible, then $A$ is $1$-dimensional, i.e. $\{1_A\}$ is a linear basis for $A$. 
\end{ppsn}
\textit{Proof.} Let $a \in A$ be given. Choose $\lambda \in \sigma(a)$. Then, $a-\lambda1_A$ is not invertible. Hence, $a-\lambda1_A=0$, i.e. $a=\lambda 1_A$. The proof is complete. \hfill $\Box$

Assume, until further mention, that $A$ is a unital commutative Banach algebra. 
Let $I$ be an ideal in $A$. Recall that $I$ is said to be maximal if 
\begin{enumerate}
\item[(i)] the ideal $I$ is proper, i.e. $I\neq A$, and
\item[(ii)] if $J$ is an ideal in $A$ containing $I$, then $J=I$ or $J=A$. 
\end{enumerate} 
We denote the set of maximal ideals of $A$ by $\mathcal{M}(A)$, or simply by $\mathcal{M}$. 

Note that maximal ideals are closed. To see this, let $I$ be a maximal ideal of $A$.  Suppose $\overline{I}$ is the closure of $I$. Clearly, $\overline{I}$ is an ideal containing $I$. The maximality of $I$ forces that either $\overline{I}=I$ or $\overline{I}=A$.  Suppose $\overline{I}=A$. Then, $1 \in \overline{I}$. 
This implies that there exists $a \in I$ such that $||1-a||<1$. The inequality $||1-a||<1$ implies that $a$ is invertible. This forces $I=A$ which is a contradiction. Therefore, $\overline{I}=I$, i.e. $I$ is closed. 

The set $\mathcal{M}(A)$ can be identified with a certain subset of $A^{*}$ which we explain next. 

\begin{dfn}
Let $\chi:A \to \bbc$ be a linear map. The map $\chi$ is called a character of $A$ if 
\begin{enumerate}
\item[(1)] the map $\chi$ is multiplicative, i.e. $\chi(ab)=\chi(a)\chi(b)$ for all $a,b \in A$, and
\item[(2)] the map $\chi$ is non-zero.
\end{enumerate}
The set of characters of $A$  is denoted  by $\widehat{A}$. 
\end{dfn}

\begin{rmrk}
 Let $\chi:A \to \bbc$ be a linear map which is multiplicative. Then, $\chi \in \widehat{A}$ if and only if $\chi(1_A)=1$. One direction is obvious. Assume now that $\chi \in \widehat{A}$. 
Choose $a \in A$ such that $\chi(a) \neq 0$. Then, 
\[
\chi(a)=\chi(a.1_A)=\chi(a)\chi(1_A).\]
Cancelling $\chi(a)$, we deduce $\chi(1_A)=1$. 
\end{rmrk}

\begin{rmrk}
\label{maximal ideal space}
 The map 
\[
\widehat{A} \ni \chi \to \ker(\chi) \in \mathcal{M}\] 
is a bijection. We leave it to the reader to check the well-definedness and the injectivity of the prescribed map. Let $I$ be a maximal ideal of $A$. Note that $A/I$ is a Banach algebra, with the quotient norm, and every non-zero element
of $A/I$ is invertible. By Prop. \ref{Gelfand_Mazur}, $\{1_A+I\}$ is a linear basis for $A/I$. 

Let $\pi:A \to A/I$ be the quotient map. Because $\{1_A+I\}$ is a basis for $A/I$, the map $\pi$ is of the form 
\[
\pi(a)=\chi(a)(1_A+I)\]
for some map $\chi:A \to \bbc$. Clearly, the  map $\chi$ is a character and $\ker(\chi)=I$. 

\end{rmrk}

\begin{ppsn}
\label{description of spectrum}
Let $a \in A$ be given. Then, 
\[
\sigma(a)=\{\chi(a): \chi \in \widehat{A}\}.\]
In particular, if $\chi \in \widehat{A}$, then $|\chi(a)| \leq ||a||$ for every $a \in A$. Thus, if $\chi$ is a character, then $||\chi||=1$. 
\end{ppsn}
\textit{Proof.} Let $a \in A$ be given. Suppose $\chi \in \widehat{A}$. Then, $a-\chi(a)1_A \in \ker(\chi)$. 
An invertible element cannot be in a proper ideal. Therefore, $\chi(a) \in \sigma(a)$. Consequently, the set 
$\{\chi(a):\chi \in \widehat{A}\} \subset \sigma(a)$. 

Suppose $\lambda \in \sigma(a)$. Let $I:=A(a-\lambda)$. Then, $I$ is an ideal containing $a-\lambda$. Note that $I$ is proper as $1 \notin I$. By Zorn's lemma, there exists 
a maximal ideal $J$ of $A$ containing $I$. Let $\chi \in \widehat{A}$ be such that $\ker(\chi)=J$. Since $a-\lambda \in J=\ker(\chi)$, we have
$\chi(a)=\lambda \chi(1_A)=\lambda$. This proves that $\sigma(a) \subset \{\chi(a):\chi \in \widehat{A}\}$. Therefore, 
\[
\sigma(a)=\{\chi(a): \chi \in \widehat{A}\}.\]

The remaining assertions follow from Lemma \ref{mother of lemmas} and the fact that $\chi(1_A)=1$. \hfill $\Box$

From the above proposition, we know that $\widehat{A}$ is a subset of the unit ball of $A^{*}$. We endow $\widehat{A}$ with the weak$^*$-topology\footnote{Let $\{\phi_i\}_{i \in A}$ be a net in $A^{*}$ and $\phi \in A^{*}$. Recall that $\phi_i \to \phi$ in the weak$^*$-topology  if and only if for every $a\in A$, $\phi_i(a) \to \phi(a)$.}. Looking at the prescription of $\widehat{A}$, i.e.
\[
\widehat{A}=\Big\{\chi \in \widehat{A}: \chi(1_A)=1~~\textrm{and~}\chi(ab)=\chi(a)\chi(b)~\forall a,b \in A\Big\},\]
it is apparent that $\widehat{A}$ is closed in $\widehat{A}^{*}$, when $A^{*}$ is given the weak$^*$-topology. But $\widehat{A}$ is a subset of the unit ball of $A^{*}$ and the unit ball is weak$^*$-compact. Therefore,  $\widehat{A}$ is a compact Hausdorff space. The space $\widehat{A}$ with the weak$^*$-topology is called \emph{the spectrum of $A$}. 

\textbf{The Gelfand map:} For $a \in A$, let $\widehat{a}:\widehat{A} \to \bbc$ be defined by
\[
\widehat{a}(\chi)=\chi(a).\]
From the definition of the weak$^*$-topology, we see that for any  $a \in A$, $\widehat{a} \in C(\widehat{A})$, i.e. the function $\widehat{a}:\widehat{A} \to \bbc$ is continuous. The map 
\[
A \ni a \to \widehat{a} \in C(\widehat{A})\]
is called \emph{the Gelfand map}. 

\begin{thm}[Gelfand-Naimark]
\label{GN1}
Let $A$ be a unital commutative $C^{*}$-algebra. Then, the Gelfand transformation
\[
A \ni a \to \widehat{a} \in C(\widehat{A})\]
is an isometric $^*$-algebra isomorphism. 
\end{thm}

\begin{lmma}
\label{star structure}
Suppose $A$ is a unital commutative $C^{*}$-algebra. 
Let $a \in A$ and $\chi \in \widehat{A}$ be given. Then, $\chi(a^*)=\overline{\chi(a)}$. 
\end{lmma}
\textit{Proof.} Let $a \in A$ be given. Set $b:=\frac{a+a^*}{2}$ and $c:=\frac{a-a^*}{2i}$. Then, $a=b+ic$. Note that $b$ and $c$ are self-adjoint. It suffices to prove that 
$\chi(b)$ and $\chi(c)$ are real. In other words, it suffices to verify the assertion assuming $a$ is self-adjoint. 

So, suppose $a \in A$ is self-adjoint. Let $\chi \in \widehat{A}$ be given. Write $\chi(a)=\alpha+i\beta$ with $\alpha,\beta \in \bbr$. For $t \in \bbr$, let 
\[
u_t:=e^{ita}:=\sum_{n=0}^{\infty}\frac{(ita)^n}{n!}.\]
Then, $u_su_t=u_{s+t}$, $u_t^{*}=u_{-t}$ and $u_0=1$. In particular, $u_{t}^{*}u_t=1=u_tu_t^{*}$. Consequently, $||u_t||=1$ for every $t \in \bbr$. 

Observe that for every $t \in \bbr$,
\[
e^{-\beta t}=|e^{i(\alpha+i\beta)t}|=|e^{it\chi(a)}|=|\chi(u_t)| \leq 1.\]
Since the above happens for every real number $t$, we have $\beta=0$. Thus, $\chi(a)$ is a real number. This completes the proof. \hfill $\Box$

\begin{lmma}
\label{isometric homo}
Let $A$ and $B$ be $C^{*}$-algebras (not necessarily commutative). Suppose $\pi:A \to B$ is a $^*$-algebra homomorphism.  
The following are equivalent.
\begin{enumerate}
\item[(1)] The homomorphism $\pi$ is an isometry, i.e. $||\pi(a)||=||a||$ for $a \in A$.
\item[(2)] For every self-adjoint $a \in A$, $||\pi(a)||=||a||$.
\end{enumerate}
\end{lmma}
\textit{Proof.} Assume that $(2)$ holds. Let $a \in A$ be given. Observe that $a^{*}a$ is self-adjoint. Hence, 
\[
||\pi(a)||^{2}=||\pi(a)^*\pi(a)||=||\pi(a^*a)||=||a^*a||=||a||^2.\]
This proves that $(2) \implies (1)$. The other implication is a tautology. \hfill $\Box$

\textit{Proof of Thm. \ref{GN1}.} It is clear that the Gelfand transformation is a unital algebra homomorphism. Lemma \ref{star structure} is precisely the assertion 
that the Gelfand map preserves the $^*$-structure. To prove that $A \ni a \to \widehat{a} \in C(\widehat{A})$ is an isometry, we appeal to Lemma \ref{isometric homo}.
Let $a \in A$ be self-adjoint. Thanks to Prop. \ref{spectral radius for self-adjoint} and Prop. \ref{description of spectrum}, we have
\[
||a||=\sup\{|\chi(a)|:\chi \in  \widehat{A}\}=\sup\{|\widehat{a}(\chi)|:\chi \in  \widehat{A}\}=||\widehat{a}||_{\infty}.\]

The only thing that is now left to prove is the surjectivity part. Let $\mathcal{D}:=\{\widehat{a}:a \in A\}$.  Note that $\mathcal{D}$ is a subalgebra,
closed under taking conjugation and contains the constant function $`1'$. Moreover, $\mathcal{D}$ separates points of $\widehat{A}$. By 
Stone-Weierstrass theorem, we conclude that $\mathcal{D}$ is dense in $C(\widehat{A})$. The fact that 
$a \to \widehat{a}$ is an isometry and $A$ is complete implies that $\mathcal{D}$ is norm closed in $C(\widehat{A})$. Hence, $\mathcal{D}=C(\widehat{A})$. 
The proof is now complete. \hfill $\Box$

\begin{xrcs}
Let $X$ be a compact Hausdorff space and set $A:=C(X)$. 
\begin{enumerate}
\item[(1)] For $x \in X$, let $I_x:=\{f \in C(X): f(x)=0\}$. Show that $I_x$ is a maximal ideal and every maximal ideal is of this form. 
               Show that the map 
               \[
               X \ni x \to I_x \in \mathcal{M}(A)\]
               is a bijection.
 \item[(2)] For $x \in X$, let $\epsilon_x:A \to \bbc$ be defined by $\epsilon_x(f)=f(x)$. Prove that $\epsilon_x$ is a character of $A$ and every character arises this way.
               Show that the map, denoted $\epsilon$, 
               \[
               X \ni x \to \epsilon_x \in \widehat{A}\]
               is a homeomorphism. For $f \in A$, what is $\widehat{f} \circ \epsilon$ ? Draw `commutative diagrams'.               
 \item[(3)] Let $F \subset X$ be a closed subset. Define $I_F:=\bigcap_{x \in F}I_x$. Show that $I_F$ is a closed ideal of $A$ and every ideal is of this form for a unique
                 closed set. 
\end{enumerate}
\end{xrcs}

\begin{xrcs}
\label{singly generated}
Let $A$ be a unital commutative $C^{*}$-algebra. Suppose there exists $a \in A$ such that $A=C^{*}\{1,a\}$, i.e. the smallest $C^{*}$-subalgebra generated by $1$ and $a$ is $A$.
Prove that the map 
\[
\widehat{A} \ni \chi \to \chi(a) \in \sigma(a)\]
is a homeomorphism. 
\end{xrcs}

\begin{lmma}
\label{spectral permanence1}
Let $A$ be a unital commutative $C^{*}$-algebra and let $a \in A$ be self-adjoint. Suppose that $a$ is invertible, then there exists a sequence of polynomials $(p_n)$ such that 
\[
a^{-1}=\lim_{n \to \infty}p_n(a).\]
\end{lmma}
\textit{Proof.} By the Gelfand-Naimark theorem, we can assume that $A=C(X)$ for some compact Hausdorff space. Let $a=f$ be self-adjoint, i.e. $f$ is real valued. The fact that 
$f$ is invertible implies that there exists $\epsilon>0$ such that $|f(x)| \geq \epsilon$ for all $x \in X$. 

By the Weierstrass theorem, there exists a sequence of polynomials $p_n$ such that $p_n(t) \to \frac{1}{t}$ uniformly on $[\epsilon,||f||] \cup [-||f||,-\epsilon]$. Then, $p_n(f) \to f^{-1}$ as $n \to \infty$. This completes the 
proof. \hfill $\Box$

\begin{ppsn}
\label{specperm}
Let $B$ be a $C^{*}$-algebra with multiplicative identity $1_B$. Suppose $A$ is a $C^{*}$-subalgebra of $B$ such that $1_B \in A$. 
Let $a \in A$ be such that $a$ is invertible in $B$. Then, $a$ is invertible in $A$. 
\end{ppsn}
\textit{Proof.} Let $a \in A$ be self-adjoint. Denote the inverse of $a$ in $B$ by $a^{-1}$. Define $C:=C^{*}\{1,a\}$ and $D:=C^{*}\{1,a,a^-1\}$. Note that $C$ and $D$ are unital commutative
$C^{*}$-algebras. Also, $C \subset A$ and $D \subset B$. The element $a$ is invertible in $D$. Thanks to Lemma \ref{spectral permanence1}, there exists a sequence of polynomials
$(p_n)$ such that 
\[
a^{-1}=\lim_{n \to \infty}p_n(a).\]
Since $p_n(a) \in C \subset A$ and  $A$ is closed, we have $a^{-1} \in A$. 

Now, let $a \in A$ be an arbitrary element which is invertible in $B$. Then, the self-adjoint elements $a^*a$ and $aa^{*}$ are invertible in $B$ and hence, by what we have established just now, is invertible  in $A$. Set $u:=(aa^*)^{-1}$ and $v=(a^*a)^{-1}$. Then, 
\[
a(a^*u)=aa^*u=1;~~(va^*)a=v(a^*a)=1.\]
The above implies that $a^*u \in A$ is a right inverse of $a$ and $va^* \in A$ is a left inverse of $a$. If an element has a right  inverse and a left inverse in $A$, then it has a two sided inverse in $A$. Hence, 
$a^{-1} \in A$. This completes the proof. \hfill $\Box$

Let $A$ and $B$ be as in  Proposition \ref{specperm}. For $a \in A$, define 
\begin{align*}
\sigma_A(a):&=\{\lambda \in \bbc: a-\lambda \textrm{~~is not invertible in $A$}\}\\
\sigma_B(a):&=\{\lambda \in \bbc: a -\lambda \textrm{~is not invertible in $B$}\}.
\end{align*}
\begin{crlre}
\label{invariance of spectrum}
With the foregoing notation,  for $a \in A$, we have $\sigma_A(a)=\sigma_B(a)$. 
\end{crlre}

\begin{rmrk}
Corollary \ref{invariance of spectrum} means that the spectrum of an element $a$ can be defined in any unital $C^{*}$-algebra containing $a$ and it does not change as long as the multiplicative unit is
kept intact. 

For example, let $T \in B(\clh)$ and let $A=C^{*}\{1,T\}$. We can talk about the spectrum of $T$ irrespective of whether we think $T$ as an element of $A$ or as an element of $B(\clh)$. There is no ambiguity
\end{rmrk}

\textbf{Notation:} Let $K \subset \bbc$ be compact. The function $K \ni z \to z \in \bbc$ will be denoted by $z$ itself. 

\begin{thm}[Continuous functional calculus]
Let $A$ be a unital $C^{*}$-algebra and let $a \in A$ be normal, i.e. $aa^*=a^*a$. Then, there exists a unique $^*$-algebra homomorphism
$\Phi:C(\sigma(a)) \to A$ which is an isometry such that $\Phi(z)=a$. 
\end{thm}
\textit{Proof.} Uniqueness is because the unital $^*$-algebra generated by $z$ is dense in $C(\sigma(a))$. For the existence, let $D:=C^{*}\{1,a\}$. Then, $D$ is commutative and let 
$G:D \to C(\widehat{D})$ be the Gelfand map. By Exercise \ref{singly generated}, we can identify $\widehat{D}$ with $\sigma(a)$. Thus, we can identify $C(\widehat{D})$ with $C(\sigma(a))$. 
Define \[\Phi:= i \circ G^{-1}\] where $i:D \to A$ is the standard inclusion. 
Then, $\Phi$ is the required map. This completes the proof. \hfill $\Box$

\begin{rmrk}
Let $A$ be a unital $C^{*}$-algebra and let $a \in A$ be normal. For $f \in C(\sigma(a))$, the image $\Phi(f)$ is usually denoted $f(a)$. This is a powerful way to construct new elements in a $C^{*}$-algebra. 
For example, for a self-adjoint $a$, the  elements  $|a|$, $|a|^{\frac{1}{2}}$ make sense. 
\end{rmrk}

\begin{xrcs}
Let $A$ be a unital $C^{*}$-algebra and let $a \in A$ be normal. Let $f$ be a continuous complex valued function on $\sigma(a)$. 
Show that 
\begin{enumerate} 
\item[(1)] $||f(a)||=\sup\{|f(\lambda)|:\lambda \in \sigma(a)\}$, 
\item[(2)] $\sigma(f(a))=f(\sigma(a))$. 
\end{enumerate}
\end{xrcs}

\begin{xrcs}[Some special normal elements]
Let $A$ be a unital $C^{*}$-algebra. 
\begin{enumerate}
\item[(1)] An element $u \in A$ is called a unitary if $uu^*=u^*u=1$. Show that $\sigma(u) \subset \bbt$ if $u$ is unitary.
\item[(2)] An element $p \in A$ is called a projection if $p^2=p=p^{*}$. Show that $\sigma(p) \subset \{0,1\}$ if $p$ is a projection.
\item[(3)] For a self-adjoint element $a \in A$, i.e. $a^*=a$, show that $\sigma(a) \subset \bbr$.
\end{enumerate}
\end{xrcs}

\begin{xrcs}
Let $A$ be a unital $C^{*}$-algebra and suppose $u \in A$ is a unitary element. Assume that $||1-u||<1$. Prove that 
$\sigma(u) \in \bbt \backslash \{-1\}$. Show that there exists a self-adjoint element $a \in A$ such that 
$u=e^{ia}$. 
\end{xrcs}
\textit{Hint.} Define a logarithm on $\bbt\backslash \{-1\}$ and use continuous functional calculus. 

\begin{ppsn}
Let $A$ and $B$ be unital $C^{*}$-algebras and let $\pi:A \to B$ be a unital $^*$-homomorphism. Suppose that $\pi$ is injective. Then, $\pi$ is an isometry. 
\end{ppsn}
\textit{Proof.} We apply Lemma \ref{isometric homo}. Let $a \in A$ be self-adjoint. By considering the $C^{*}$-subalgebras generated by $a$ and $\pi(a)$, we 
can assume that $A$ and $B$ are commutative. Then, $A=C(X)$ and $B=C(Y)$ for some compact Hausdorff spaces. 

Let $y \in Y$ be given. The map $\epsilon_y \circ \pi$ is a character of $A$. Here, $\epsilon_y$ is the evaluation map at $y$. Thus, there exists $\phi(y) \in X$ such that 
$\epsilon_y \circ \pi=\epsilon_{\phi(y)}$. This way, we get a map $Y \ni y \to \phi(y) \in X$. A moment's thought with the definitions imply that $\phi$ is continuous. 

Observe that $\pi(f)=f\circ \phi$ for $f \in C(X)$. Thus, it suffices to prove that $\phi$ is onto. Suppose not. Since $Y$ is compact, $\phi(Y)$ is compact. As $\phi$ is not onto, $\phi(Y)$ is  a 
proper closed subset of $X$. By Uryshon's Lemma, there exists a non-zero continuous function $f$ on $X$ that vanishes on $\phi(Y)$. Then, $\pi(f)=0$ which contradicts the hypothesis. 
Hence the proof. \hfill $\Box$

\begin{xrcs}
Let $A$ and $B$ be unital $C^{*}$-algebras and suppose $\pi:A \to B$ is a unital $^*$-algebra homomorphism. Prove that for $a \in A$, $\sigma(\pi(a))\subset \sigma(a)$. 
Use this to prove that $\pi$ is contractive, i.e. for $a \in A$, $||\pi(a)|| \leq ||a||$.
\end{xrcs} 
\textit{Hint.} To prove that $\pi$ is contractive, work with self-adjoint elements first. 

\section{Unitisation}
The Gelfand-Naimark theory for non-unital commutative $C^{*}$-algebras can be obtained from the unital case by adjoining an artificial multiplicative unit. 
The unitisation is usually used when the algebra is non-unital, but it makes perfect sense for unital algebras too. 

Let $A$ be a $C^{*}$-algebra. Define 
\[
A^{+}:=\{(a,\lambda):a \in A,\lambda \in \bbc\}.\]
Then, $A^{+}$ is a $^{*}$-algebra with the multiplication and the $^{*}$-structure defined as follows. 
\begin{align*}
(a,\lambda)(b,\mu)&:=(ab+\lambda b +\mu a, \lambda \mu)\\
(a,\lambda)^{*}&:=(a^{*},\overline{\lambda}).
\end{align*}
Note that $A^{+}$ is a unital $^{*}$-algebra, the unit of $A^{+}$ is $(0,1)$ which we denote by $1_{A^+}$ or simply by $1$. We view $A$ as a subalgebra of $A^{+}$ via 
the embedding $A \ni a \to (a,0) \in A^{+}$. Note that, viewed via this embedding, $A$ is a two sided ideal in $A^{+}$. We abuse notation and write $(a,0)$ as $a$ and write $(0,\lambda)$ as $\lambda$ or $\lambda 1_{A^+}$.

Let $\epsilon:A^{+} \to \bbc$ be the `character' defined by 
\[
\epsilon(a,\lambda)=\lambda.\]
Note that $\ker(\epsilon)=A$ and we have  the following short exact sequence
\[
0 \longrightarrow A \longrightarrow A^{+} \longrightarrow \bbc \longrightarrow \bbc.\]

We now endow $A^{+}$ so that $A^{+}$ becomes a $C^{*}$-algebra. The unital situation is dealt in the following exercise. 

\begin{xrcs}
Let $A$ and $B$ be $C^*$-algebras. Consider the direct sum $A\oplus B$. Then, $A\oplus B$ has an obvious product and a $^*$-structure. Endow $A\oplus B$ with a norm which
makes $A \oplus B$ a $C^{*}$-algebra. 

Let $A$ be a unital $C^*$-algebra. Observe that the map
\[
A^{+} \ni (a,\lambda) \to (a+\lambda1_A,\lambda) \in A \oplus \bbc\]
is a $^*$-isomorphism. The norm on $A^{+}$ is imposed by transporting the norm on $A \oplus \bbc$ via the above isomorphism. 
\end{xrcs}

Let $A$ be a non-unital $C^{*}$-algebra. Since $A$ is an ideal in $A^{+}$, the algebra $A^{+}$ acts on $A$ as follows. For $x=(a,\lambda)$, let $L_x:A \to A$ be the map 
defined by 
\[
L_x(b)=xb=ab+\lambda b.\]
Note that $||L_x|| \leq ||a||+|\lambda|$ and hence $L_x$ is bounded. We define a norm on $A^{+}$ by setting 
\[
||(a,\lambda)||:=||L_{(a,\lambda)}||=\sup\{||ab+\lambda b||: ||b||=1\}.\]

\begin{ppsn}
With the foregoing notation, $(A^{+},||~||)$ is a $C^{*}$-algebra. The norm restricted to $A$ coincides with the original norm on $A$. 
\end{ppsn}
\textit{Proof.} Let $a \in A$ be given. Note that if $||b||=1$, then $||ab||\leq ||a||$. Also, for $b=\frac{a^*}{||a||}$, we have $||ab||=||a||$. Thus, the norm on $A^{+}$ restricted to $A$
agrees with the norm on $A$. 

Next, we verify that $||~||$ satisfies the $C^{*}$-identity. Let $(a,\lambda) \in A^{+}$ be given. Calculate as follows to observe that 
\begin{align*}
||(a,\lambda)||^{2}&=\sup\{||ab+\lambda b||^2: ||b||=1\}\\
&=\sup\{||(ab+\lambda b)^{*}(ab+\lambda b)||: ||b||=1\}\\
&=\sup\{||b^*a^*ab+b^*\lambda a^*b+\overline{\lambda}b^{*}ab+|\lambda|^2b^*b||:||b||=1\}\\
&\leq \sup\{||a^*ab+\lambda a^*b+\overline{\lambda}ab+|\lambda|^2b||:||b||=1\}\\
& \leq ||(a,\lambda)^*(a,\lambda)|| \\
& \leq ||(a,\lambda)^*||||(a,\lambda)||.
\end{align*}
Therefore, $||(a,\lambda)|| \leq ||(a,\lambda)^*||$. By symmetry, we get $||(a,\lambda)^*||\leq ||(a,\lambda)||$. Consequently, $||(a,\lambda)||=||(a,\lambda)^*||$ for $(a,\lambda) \in A^{+}$. 
Therefore, the above chain of inequalities are in fact equalities. Thus, $||(a,\lambda)||^2=||(a,\lambda)^*(a,\lambda)||$ for every $(a,\lambda) \in A^{+}$. 

Let us verify that $||~||$ is indeed a norm. Suppose $x=(a,\lambda)$ is such that $||L_x||=0$. Then, $ab=-\lambda b$ for every $b \in A$. We claim that $\lambda=0$. If not, then $\frac{-a}{\lambda}$ is a left multiplicative identity in $A$. 
Taking adjoints, we see that $A$ has a right multiplicative identity and consequently, $A$ has an identity which is a contradiction to the hypothesis that $A$ is non-unital. 
Therefore, $\lambda=0$. Now the equality $ab=0$ for every $b$ implies that $||a||^{2}=||aa^*||=0$, i.e. $a=0$. Hence $x=0$ and $||~||$ is a norm on $A^{+}$.  

The completeness of $A^{+}$ follows from the completeness of $A$ and the fact that $A$ is  a codimensional one subspace of $A^{+}$. Hence the proof. \hfill $\Box$

Let $A$ be a non-unital commutative $C^{*}$-algebra. We define the spectrum of $A$ exactly like in the unital case. Thus, let $\widehat{A}$ denote the set of non-zero multiplicative maps $\chi:A \to \bbc$. 
For $\chi \in \widehat{A}$, let $\chi^{+}:A^{+} \to \bbc$ be defined by 
\[
\chi^{+}(a,\lambda):=\chi(a)+\lambda.\]
For $\chi \in \widehat{A}$, note that $\chi^{+} \in \widehat{A^+}$. Thus, $\chi \in A^{*}$ for every $\chi \in \widehat{A}$. Endow $\widehat{A}$ with the weak$^*$-topology. The topological space $\widehat{A}$ is called 
\emph{the spectrum of $A$}. 

\begin{xrcs}
Let $X$ be a locally compact Hausdorff space and let $U \subset X$ be an open set. Prove that $C_0(U)$ can be identified with 
$\{f \in C_0(X): f(x)=0 \textrm{~~for all $x \notin U$}\}$. We use this identification without any further mention.

\end{xrcs}

\begin{xrcs}
\begin{enumerate}
\item[(i)] Prove that the map $\widehat{A} \ni \chi \to \chi^{+} \in \widehat{A^+}$ is a topological embedding. This way, we consider $\widehat{A}$ as a subspace of $\widehat{A^+}$. 
\item[(ii)] Show that $\widehat{A^+}=\widehat{A} \cup \{\epsilon\}$. Here, $\epsilon:A^{+} \to \bbc$ is defined by 
\[
\epsilon(a,\lambda)=\lambda.\]
Deduce that $\widehat{A}$ is an open subset of $\widehat{A^+}$. 
\end{enumerate}
\item[(iii)] For $a \in A$, let $\widehat{a}:\widehat{A} \to \bbc$ be defined by 
\[
\widehat{a}(\chi):=\chi(a).\]
Observe that $\widehat{a} \in C_0(\widehat{A})$. 
\end{xrcs}

\begin{thm}[Gelfand-Naimark: the non-unital case]
Let $A$ be a non-unital commutative $C^{*}$-algebra. The Gelfand map 
\[
A \ni a \to \widehat{a} \in C_0(\widehat{A})\]
is an isometric $*$-isomorphism. 
\end{thm}
\textit{ Proof:} Let $G:A^{+} \to C(\widehat{A^+})$ be the Gelfand map. Then, it is routine to argue that $G$ maps $A$ onto $C_0(\widehat{A})$\footnote{The reader should take into account the obvious identifications going on here}. \hfill $\Box$

Let $A$ be a non-unital $C^{*}$-algebra and let $a \in A$. We define the spectrum of $a$, denoted $\sigma(a)$, as
\[
\sigma(a)=\sigma_{A^{+}}((a,0)).\]
Observe that $0 \in \sigma(a)$ for every $a \in A$. Let $C_{0}(\sigma(a))$ be the space of continuous complex valued functions on $\sigma(a)$ that vanish at $0$. 

We have the following continuous functional calculus in the non-unital situation.  
\begin{ppsn}[Continuous functional calculus]
Suppose $A$ is a non-unital $C^{*}$-algebra. Let $a \in A$ be normal, i.e. $aa^*=a^*a$. Then, there exists a unique $^*$-homomorphism $\Phi:C_0(\sigma(a)) \to A$, which is an isometry, such that 
$\Phi(z)=a$. 
\end{ppsn}

\begin{xrcs}
\label{homo iso}
Let $A$ and $B$ be $C^{*}$-algebras and let $\pi:A \to B$ be an injective $^*$-homomorphism. Then, $\pi$ is an isometry, i.e. $||\pi(a)||=||a||$ for every $a \in A$. 
\end{xrcs}
\textit{Hint:} Consider the map $\pi^{+}:A^{+} \to B^{+}$ defined by $\pi^{+}(a,\lambda)=(\pi(a),\lambda)$. 

\begin{xrcs}
Let $A$ and $B$ be $C^{*}$-algebras and let $\pi:A \to B$ be a $^*$-algebra homomorphism. Show that $||\pi(a)|| \leq ||a||$ for every $a \in A$. 
\end{xrcs}

\section{Positivity}
A basic notion in the theory of $C^{*}$-algebras is the notion of positivity and it implements a partial order on the set of self-adjoint elements of a  $C^{*}$-algebra. Let $A$ be a $C^{*}$-algebra and let $a \in A$. Let $\sigma(a)$ be the spectrum of $a$. Recall that, if $A$ is non-unital, $\sigma(a):=\sigma_{A^+}((a,0))$. 

\begin{dfn}
Let $a \in A$ be self-adjoint. Then $a$ is said to be positive, and written $a \geq 0$, if $\sigma(a) \subset [0,\infty)$. For two self-adjoint elements  $x,y \in A$, we say $y \geq x$ if $y-x \geq 0$. 
\end{dfn}
\begin{rmrk}
We have already proved that if $a$ is self-adjoint, then $\sigma(a) \subset \bbr$.
Thus, for a self-adjoint element $a \in A$, $a^2$ is positive. 

\end{rmrk}

We first prove that $\leq x$ is a partial order on the set of self-adjoint elements. 

For the rest of this section, let $A$ be a $C^{*}$-algebra. 

\begin{lmma}
Let $a \in A$ be self-adjoint. Suppose $a \geq 0$ and $a \leq 0$. Then, $a=0$.
\end{lmma}
\textit{Proof.} Note that $\sigma(a)=\{0\}$.  The spectral radius formula gives the desired conclusion. \hfill $\Box$

\begin{ppsn}
\label{sum of positive}
Let $a,b \in A$ be positive. Then, $a+b$ is positive. 
\end{ppsn}
\textit{Proof.} We can assume that $A$ is unital. Clearly, $a+b$ is self-adjoint and hence $\sigma(a+b) \subset \bbr$.  Note that $c$ is positive if and only if $tc$ is positive for some (and every) $t>0$. Thus, by scaling, we can assume that $||a||,||b|| \leq 1$. 
This implies that $\sigma(a)$ and $\sigma(b)$ are contained in $[0,1]$. Consequently, $\sigma(1-a)$ and $\sigma(1-b)$ are contained in $[0,1]$. Therefore, 
$||1-a|| \leq 1$ and $||1-b|| \leq 1$. Calculate as follows to observe that 
\begin{align*}
\Big|\Big|1-\frac{a+b}{2}\Big|\Big|&=\frac{1}{2}||(1-a)+(1-b)||\\
&\leq \frac{1}{2}(||1-a||+||1-b||)\\
& \leq 1.
\end{align*}
Thus, if $\lambda \in \sigma(\frac{a+b}{2})$, then $|1-\lambda| \leq 1$, i.e. $\lambda \in [0,2]$. Hence, $\sigma(\frac{a+b}{2}) \subset [0,\infty)$. This implies $a+b$ is positive. \hfill $\Box$

\begin{crlre}
The order $\leq$ is a partial order on the set of self-adjoint elements of $A$. 
\end{crlre}

Next, we prove the following characterisation of positive elements. 
\begin{thm}
\label{positive char}
Let $a \in A$ be self-adjoint. The following are equivalent.
\begin{enumerate}
\item[(1)] $a$ is positive.
\item[(2)] There exists a self-adjoint element $b \in A$ such that $a=b^2$.
\item[(3)] There exists $b \in A$ such that $a=b^*b$. 
\end{enumerate}
\end{thm}
The non-trivial part is in proving $(3) \implies (1)$ which prove first. 

\begin{lmma}
\label{spectrum of product}
Let $x,y \in A$ be given. Then, $\sigma(xy)\cup\{0\}=\sigma(yx) \cup \{0\}$. 
\end{lmma}
\textit{Proof.} A little bit of thought, playing with scaling $x$ and $y$, will tell us that it suffices to prove that $1-xy$ is invertible if and only if $1-yx$ is invertible. 

Suppose that $1-xy$ is invertible with inverse $v$. Set $u:=1+yvx$. Calculate as follows to observe that 
\begin{align*}
(1-yx)u&=1+yvx-yx-y(xy)vx\\
&=1+yvx-yx+y(1-xy-1)vx \\
&=1+y(1-xy)vx-yx\\
&=1.
\end{align*}
A similar calculation yields $u(1-yx)=1$. Hence, $u$ is the inverse of $1-yx$. We can interchange the roles of $x$ and $y$ to complete the proof. \hfill $\Box$

\begin{rmrk}
The formula for $u$ in the above proof is obtained by playing  with the `formal power series expansion'
\[
\frac{1}{1-z}:=1+z+z^2+\cdots.\]
The reader is encouraged to do the exercise of writing (and simplifying ) the `formal expansion' of $\frac{1}{1-yx}$ and arrive at the formula for $u$. 
\end{rmrk}

\begin{lmma}
Let $a \in A$ be such that $\sigma(a^*a) \subset (-\infty,0]$. Then, $a=0$. 
\end{lmma}
\textit{Proof.} Write $a=b+ic$ with $b,c$ self-adjoint. Then, $a^{*}a+aa^{*}=b^2+c^2$. By Lemma \ref{spectrum of product}, we have $aa^* \leq 0$. Thanks to Prop. \ref{sum of positive}, we have
$a^*a+aa^*=b^2+c^2 \leq 0$. But, $b^2,c^2 \geq 0$. Consequently, $b^2+c^2 \geq 0$ and we deduce $b^2+c^{2}=0$, i.e $b^2=-c^2$. Now, $-c^2 \leq 0$, thus $0 \geq b^2 \geq 0$. Hence,  $b^2=c^2=0$. As $b$ and $c$ are self-adjoint, we have $b=c=0$ which
implies $a=0$. Hence the proof. \hfill $\Box$

\begin{ppsn}
\label{positivity of $a^*a$}
For any $a \in A$, $a^*a \geq 0$. 
\end{ppsn}
\textit{Proof.} Note that $a^*a$ is self-adjoint. Let $f:\sigma(a^*a) \to [0,\infty)$ and $g:\sigma(a^*a) \to [0,\infty)$ be the continuous functions defined by 
\[
f(t)=\max\{t,0\};~ g(t)=\max\{-t,0\}.\]
Note that $f,g \geq 0$ and $f(t)-g(t)=t$ for every $t \in \sigma(a^*a)$. Also, $fg=gf=0$. Set $u:=f(a^*a)$ and $v=g(a*a)$. Then, 
$a^*a=u-v$, $u,v \geq 0$ and $uv=0$. 

Note that $(av)^*(av)=v^*a^*av=v^*(u-v)v=-v^3 \leq 0$. By the previous lemma, we have $av=0$. Hence, $-v^2=a^*av=0$. This implies $v=0$. 
Therefore, $a^*a=u \geq 0$. Hence the proof. \hfill $\Box$

\textit{Proof of Theorem \ref{positive char}.} Let $a \in A$ be self-adjoint. Assume that $(1)$ holds. Then, $\sigma(a) \subset [0,\infty)$. Let $f$ be the function on $\sigma(a)$ defined by 
\[f(t)=\sqrt{t}.\] Set $b=f(a)$. Then, $b^2=a$. This completes the proof of $(1) \implies (2)$. 
The implication $(2) \implies (3)$ is tautological and $(3) \implies (1)$ is exactly Prop. \ref{positivity of $a^*a$}. \hfill $\Box$

\begin{xrcs}
\label{crucial pos}
Let $x,y \in A$ be self-adjoint elements. If $x \geq y$, then $a^*xa \geq a^*ya$. 
\end{xrcs}
\textit{Hint:} Use Theorem \ref{positive char}. 

\begin{ppsn}[Square root]
Let $a \in A$ be a positive element. Then, there exists a unique positive element $b \in A$ such that $a=b^2$. Such an element $b$ is called the square root of $a$ and is  denoted $\sqrt{a}$. 
\end{ppsn}
\textit{Proof.} We provide a proof assuming $A$ is unital. By Gelfand-Naimark, the proposition is true if $A$ is a commutative $C^{*}$-algebra. 
 The proof of Theorem \ref{positive char} provides a square root. Let us prove uniqueness. Suppose $b,c \in A$ are such that $b,c \geq 0$ and $b^2=c^2=a$. 

\textit{Claim:} $b \in C^{*}\{1,b^2\}$. 

Note that $\sigma(b) \subset [0,M]$ for some $M>0$. It follows from Stone-Weierstrass that $span\{t^{2k}: k=0,1,2,\cdots\}$ is dense in $C[0,M]$. Thus, there exists a sequence of 
polynomials $p_n$ such that $p_n(t^2) \to t$ uniformly on $[0,M]$. Then, $p_n(b^2) \to b$. Clearly, $p_n(b^2) \in C^{*}\{1,b^2\}$. Hence, $b \in C^{*}\{1,b^2\}$. This proves the claim. 

Similiarly, $c \in C^{*}\{1,c^2\}$. Thus, $b,c$ are positve elements in the commutative $C^{*}$-algebra $C^{*}\{1,a\}$ whose squares coincide, i.e. $b^2=c^2$. Thus, from the commutative case, $b=c$. 
The proof is complete. \hfill $\Box$

\section{GNS construction}
In this section, we discuss the GNS construction which leads us to the proof of Thm. \ref{GNN}. Let us start with the definition of representations. 
Let $A$ be a $C^{*}$-algebra and suppose  $\clh$ is a Hilbert space. A $^*$-homomorphism $\pi:A \to B(\clh)$ is called \emph{a representation} of $A$ on $\clh$. 
Let $\pi_1:A \to B(\clh_1)$ and $\pi_2:A \to B(\clh_2)$ be two representations. We say that $\pi_1$ and $\pi_2$ are equivalent if there exists a unitary operator $U:\clh_1 \to \clh_2$ such that 
\[
U\pi_1(a)U^*=\pi_2(a)\]
for all $a \in A$.

If $\pi$ is a representation, then $\pi$ is contractive, i.e. $||\pi(a)|| \leq ||a||$ for all $a \in A$. A representation $\pi:A \to B(\clh)$ is said to be non-degenerate if the set \[\pi(A)\clh:=\{\pi(a)\xi:a \in A, \xi \in \clh\}\] is total in $\clh$\footnote{A subset $\mathcal{S} \subset \clh$ is said to be \emph{total} if $span~\mathcal{S}$ is dense in $\clh$}.

\begin{xrcs}
\label{alternate description of non-deg}
Let $\pi:A \to B(\clh)$ be a representation. Prove that the following are equivalent. 
\begin{enumerate}
\item[(1)] The representation $\pi$ is non-degenerate. 
\item[(2)]  The intersection $\displaystyle \bigcap_{a \in A}\ker(\pi(a))=\{0\}$.
\end{enumerate}
\end{xrcs}

\begin{xrcs}
Let $\pi:A \to B(\clh)$ be a representation. A closed subspace $W \subset \clh$ is said to be invariant under $\pi$ if $\pi(a)W \subset W$ for every $a \in A$. 
\begin{enumerate}
\item[(1)] Suppose $W \subset \clh$ is invariant under $\pi$. Then, prove that $W^{\perp}$ is invariant under $\pi$.
\item[(2)] Let $\clh_0:=\overline{span(\pi(A)\clh)}$ and $\clh_1:=\clh_0^{\perp}$. Note that $\clh_0$ and $\clh_1$ are invariant under $\pi$. Show that if we decompose $\clh$ as $\clh=\clh_0\oplus \clh_1$, the operator $\pi(a)$, for $a \in A$, is of the form 
\[\pi(a):=\begin{bmatrix}
\pi_0(a) & 0 \\
0 & 0
\end{bmatrix}.\] 
Prove that $\pi_0$ is a non-degenerate representation of $A$ on $\clh_0$. Thus, every representation can be decomposed as a direct sum of a non-degenerate representation and the zero representation. 
\end{enumerate}
In short, it suffices to  concentrate only on non-degenerate representations if we are concerned with the representation theory of a $C^{*}$-algebra.
\end{xrcs}

\begin{ppsn}
Let $A$ be a unital $C^{*}$-algebra and let $\pi:A \to B(\clh)$ be a representation. Then, the following are equivalent.
\begin{enumerate}
\item[(1)] The representation $\pi$ is unital, i.e. $\pi(1_A)=1_\clh$.
\item[(2)] The representation $\pi$ is non-degenerate. 
\end{enumerate}
\end{ppsn}
\textit{Proof.} Assume that $(1)$ holds. Then, $\xi=\pi(1_A)\xi \in \pi(A)\clh$ for every $\xi \in \clh$. Hence, $\clh=\pi(A)\clh$ in this case. Thus, $(1) \implies (2)$. 
Suppose that $(2)$ holds. Let $\mathcal{S}:=\pi(A)\clh$. Since, $\mathcal{S}$ is total in $\clh$, it suffices to prove that $\pi(1_A)\xi=\xi$ for every $\xi \in \mathcal{S}$. 
Let $\xi=\pi(a)\eta \in \mathcal{S}$ be given. Then, \[\pi(1_A)\xi=\pi(1_A)\pi(a)\eta=\pi(a)\eta=\xi.\] Hence $(1)$ holds. This completes the proof. \hfill $\Box$

\begin{xrcs}[Direct sum of representations]
\begin{enumerate}
\item[(1)] Let $I$ be an indexing set and for each $i \in I$, let $\clh_i$ be a Hilbert space. Denote the algebraic direct sum of $\clh_i$ by $\clh_{alg}:= \bigoplus_{alg}\clh_i$. 
Define an inner product on $\clh$ as follows. For $\xi=\bigoplus \xi_i$ and $\eta=\bigoplus \eta_i$, defined 
\[
\langle \xi|\eta \rangle:=\sum_{i\in I} \langle \xi_i|\eta_i\rangle.\]
Note that the above sum makes sense as $\xi_i$ and $\eta_i$ are zero except for finitely many indices $i$. The completion of $\clh_{alg}$ is called the direct sum of $\clh_i$ and 
is denoted $\displaystyle \bigoplus_{ i \in I}\clh_i$. 

\item[(2)] Keep the notation used in $(1)$. Set $\clh:=\displaystyle \bigoplus_{i \in I}\clh_i$. Suppose $T_i \in B(\clh_i)$ for every $i \in I$ and $\sup\{||T_i||: i \in I\} <\infty$. Then, there exists 
a unique operator denoted $\displaystyle \bigoplus T_i=:T$ such that 
\[
T\Big(\bigoplus \xi_i\Big)=\bigoplus_{i \in I} T_i\xi_i.\]
Prove that $||T||=\sup\{||T_i||:i \in I\}$. 

\item[(3)] Let $A$ be a $C^*$-algebra. Suppose $\pi_i:A \to B(\clh_i)$ is a representation for each $i \in I$. Set $\displaystyle \clh:=\bigoplus_{i \in I}\clh_i$ and for every $a \in A$, let 
\[
\pi(a):=\bigoplus_{i \in I}\pi_i(a).\]
Justify that $\pi(a)$ makes sense. Show that $\pi$ is a representation of $A$ on $\clh$. 
\end{enumerate}
\end{xrcs}

Let $\pi:A \to B(\clh)$ be a representation. A unit vector $\Omega \in \clh$ is said to \emph{cyclic} if $\pi(A)\Omega:=\{\pi(a)\Omega: a \in A\}$ is dense in $\clh$. A representation $\pi$ is said to be 
cyclic if $\pi$ has a cyclic vector. A cyclic representation is obviously non-degenerate.

\begin{ppsn}
Let $\pi:A \to B(\clh)$ be a non-degenerate representation. Then, $\pi$ is unitarily equivalent to a direct sum of cyclic representations. 
\end{ppsn}
\textit{Proof.} Apply Zorn's Lemma to conclude that there is a maximal set of unit vectors $\mathcal{S}$ such that $\pi(A)\Omega$ is orthogonal to $\pi(A)\Omega^{'}$ whenever $\Omega$ and $\Omega^{'}$ are two distinct elements in $\mathcal{S}$. 
For $\Omega \in \mathcal{S}$, let $\clh_{\Omega}:=\overline{\pi(A)\Omega}$. Note that $\clh_\Omega$ is invariant under $\pi$ for every $\Omega \in \mathcal{S}$. Thus, $\pi$ restricts to a representation, denoted $\pi_\Omega$, of $A$ on $\clh_\Omega$ for $\Omega \in \mathcal{S}$.

We claim that $\bigoplus_{alg}\clh_\Omega$ is dense in $\clh$. If not, pick a unit vector $\Omega_0$ orthogonal to $\bigoplus_{alg}\clh_\Omega$. Then, $\pi(A)\Omega_0$ is orthogonal to $\pi(A)\Omega$ for every $\Omega \in \mathcal{S}$. This contradicts the maximality of $\mathcal{S}$. We leave it as an exercise to the reader to fill the proof that $\pi$ is unitarily equivalent $\displaystyle \bigoplus_{\Omega}\pi_\Omega$. \hfill $\Box$

The above proposition means that the study of non-degenerate representations boils down to the study of cyclic representations. How do cyclic representations arise ? They arise from states which is the content of the GNS construction. 
\begin{dfn}
Let $A$ be a $C^{*}$-algebra and suppose $\phi:A \to \bbc$ is a bounded linear functional. Then, $\phi$ is called a state if 
\begin{enumerate}
\item[(1)] $||\phi||=1$, and
\item[(2)] if $x \geq 0$, then $\phi(x) \geq 0$. 
\end{enumerate}
\end{dfn}

\begin{thm}[GNS construction]
\label{GNS state}
Let $A$ be a $C^{*}$-algebra and suppose $\phi$ is a state on $A$. Then, there exists a Hilbert space, a  representation $\pi:A \to B(\clh)$ and a unit vector $\Omega \in \clh$ that is cyclic for $\pi$ such that 
\[
\phi(a)=\langle \pi(a)\Omega|\Omega \rangle.\]
\end{thm}

We will only prove the above theorem in the unital case. The non-unital case requires the notion of ``approximate identities" and we refer the reader to \cite{Arveson_invitation} or \cite{Davidson_Ken} for the non-unital case.

Let $\phi$ be a state on a $C^{*}$-algebra $A$. A triple $(\pi,\clh,\Omega)$ guaranteed by Theorem \ref{GNS state} is usually called a \emph{GNS triple} for the state $\phi$. If we want to stress the dependence of $(\pi,\clh,\Omega)$ on $\phi$, we denote $(\pi,\clh,\Omega)$ by
$(\pi_{\phi},\clh_{\phi},\Omega_{\phi})$. We show next that GNS triples are unique in a certain sense.

\begin{xrcs}
\label{KRP} Suppose $\clh_1$ and $\clh_2$ are Hilbert spaces and $S_1$ and $S_2$ are total subsets of $\clh_1$ and $\clh_2$ respectively. Let $\phi:S_1 \to S_2$ be a map such that $\langle \phi(x)|\phi(y) \rangle=\langle x|y \rangle$ for $x,y \in S_1$. Prove that there exists a unique isometry $V:\clh_1 \to \clh_2$ which extends $\phi$. Moreover, show that if $\phi$ is a bijection, then the isometry $V$ is a unitary.
\end{xrcs}

\begin{ppsn}
Let $\phi$ be a state on a $C^{*}$-algebra $A$. Suppose $(\pi_1,\clh_1,\Omega_1)$ and $(\pi_2,\clh_2,\Omega_2)$ are two GNS triples. Then, there exists a unique unitary $U:\clh_1 \to \clh_2$ such that 
\begin{enumerate}
\item[(1)] for $a \in A$, $U\pi_1(a)=\pi_2(a)U$, and
\item[(2)] $U \Omega_1=\Omega_2$. 
\end{enumerate}
\end{ppsn}
\textit{Proof.} Again we prove assuming $A$ is unital. Let $S_1:=\pi_1(A)\Omega_1$ and $S_2=\pi_2(A)\Omega_2$. Note that for $a,b \in A$, 
\[
\langle \pi_1(a)\Omega_1|\pi_1(b)\Omega_1\rangle=\langle \pi_1(b^*a)\Omega_1|\Omega_1\rangle =\phi(b^*a)=\langle \pi_2(a)\Omega_2|\pi_2(b)\Omega_2 \rangle.\]
By Exercise \ref{KRP}, we have a unitary $U:\clh_1 \to \clh_2$ such that $U\pi_1(a)\Omega_1=\pi_2(a)\Omega_2$. Taking $a=1_A$, we get $U\Omega_1=\Omega_2$. We leave it to the reader 
to verify that $U$ intertwines $\pi_1$ and $\pi_2$, i.e. $U\pi_1(a)=\pi_2(a)U$ for every $a \in A$. \hfill $\Box$

Next, we discuss the existence of `the' GNS triple. Until further mention, the letter $A$ stands for a unital $C^{*}$-algebra. 

\begin{ppsn}
\label{bdd and pos}
Let $\phi:A \to \bbc$ be a linear functional. The following are equivalent. 
\begin{enumerate}
\item[(1)] The functional $\phi$ is positive, i.e. $\phi(x) \geq 0$ whenever $x \geq 0$. 
\item[(2)] The functional $\phi$ is bounded and $||\phi||=\phi(1)$. 
\end{enumerate}
\end{ppsn}
\textit{Proof.} Assume that $(1)$ holds.  Define a sesquilinear form on $A$ by 
\[
\langle x|y \rangle_{\phi}:=\phi(y^*x).\]
Then, by Cauchy-Schwartz inequality, we have for $x,y \in A$, 
\[
|\phi(y^*x)|=|\langle x|y \rangle_\phi| \leq ||x||_{\phi}||y||_{\phi}=\phi(x^*x)^{\frac{1}{2}}\phi(y^*y)^{\frac{1}{2}}.\]
Let $x$ be an arbitrary element of $A$ and set $y=1$ in the above inequality. Then,  $|\phi(x)| \leq \phi(1)^{\frac{1}{2}}\phi(x^*x)^{\frac{1}{2}}$. Note that $x^*x \leq ||x||^2$, i.e. $||x||^2-x^*x=y^*y$ for some element $y$. Since $\phi(y^*y) \geq 0$, we 
get 
\[
|\phi(x)|\leq \phi(1)^{\frac{1}{2}}\phi(||x||^2)^{\frac{1}{2}}=\phi(1)||x||.\]
Therefore, $||\phi|| \leq \phi(1)$. Clearly, $0 \leq \phi(1) \leq ||\phi||~||1||=||\phi||$. The proof of $(1) \implies (2)$ is complete. 

Let us assume that $(2)$ holds. We normalise and assume that $||\phi||=\phi(1)=1$. Let $a \in A$ be self-adjoint. We claim that $\phi(a)$ is real. Write $\phi(a)=\alpha+i\beta$ with $\alpha,\beta \in \bbr$. 
Observe that, for $t \in \bbr$, 
\[
\alpha^2+(\beta+t)^2=|\phi(a+it)|^2 \leq ||a+it||^2=||a||^2+t^2.\]
Since the above happens for every $t \in \bbr$, we can conclude that $\beta=0$. This proves the claim. 

Now, let $a\in A$ be positive. We prove that $\phi(a) \geq 0$. By normalising, we can assume that $||a||=1$. Then, $\sigma(a) \subset [0,1]$. Since $\sigma(1-a)=1-\sigma(a) \subset [0,1]$ and $1-a$ is self-adjoint, we have $||1-a|| \leq 1$. 
Then,
\[
|1-\phi(a)|=|\phi(1-a)|\leq ||1-a|| \leq 1.\]
Therefore, $\phi(a) \subset [0,2]$ and hence the proof. This completes the proof of the implication $(2) \implies (1)$. \hfill $\Box$

Now we prove Thm. \ref{GNS state} when $A$ is unital. To that effect, fix a state $\phi$ on $A$. On $A$, define a sesquilinear form on $A$, as in Prop. \ref{bdd and pos}, by setting
\[
\langle x|y \rangle_{\phi}:=\phi(y^*x).\]
Then, $\langle~|~\rangle_{\phi}$ satisfies all the properties for an inner product except the positive definiteness. Let $N_\phi$ be the set of null vectors, i.e. 
\[
N_{\phi}:=\{x \in A: \phi(x^*x)=0\}.\]
Then, $\langle~|~\rangle_\phi$ descends to a genuine inner product on $A/N_\phi$ and we complete $A/N_\phi$ to get a Hilbert space $\clh_\phi$. 

Let $x,y \in A$ be given. Thanks to Exercise \ref{crucial pos}, we have \[(xy)^{*}(xy)=y^*x^*xy \leq ||x||^{2}y^*y.\] The positivity of $\phi$ implies that 
\begin{equation}
\label{left ideal}
\langle xy|xy \rangle_{\phi}:=\phi((xy)^*(xy)) \leq ||x||^{2}\phi(y^*y)=||x||_{\phi}^{2}||y||^{2}_{\phi}
\end{equation}
Eq. \ref{left ideal} forces 
 that $N_\phi$ is a left-ideal, i.e. if $y \in N_\phi$ and $x \in N_\phi$, then $xy \in N_\phi$. 

Let $x \in A$ be given. Define a linear map $\pi_\phi(x):A/N_\phi \to A/N_\phi$ by setting
\[
\pi_\phi(x)(y+N_\phi)=xy+N_\phi.\]
Thanks to Eq. \ref{left ideal}, we have that $\pi_{\phi}(x)$ is bounded. 
Since $A/N_\phi$ is dense in $\clh_{\phi}$, $\pi_{\phi}(x)$ extends uniquely to a bounded linear operator on $\clh_\phi$, which we again denote by $\phi_\phi(x)$. 

We leave it to the reader to check that $\pi_\phi:A \to B(\clh_\phi)$ is a unital representation. Define $\Omega_{\phi}:=1+N_{\phi}$. That $\Omega_\phi$ is a unit vector is a translation of
the fact that $\phi(1)=1$. Similarly, that $\Omega_\phi$ is cyclic is a consequence of the definition of $\clh_\phi$ as the completion of $A/N_\phi$. By definition, for $x \in A$, 
\[
\phi(x)=\langle \pi_\phi(x)\Omega_\phi|\Omega_\phi\rangle.\]
This completes the proof of Thm. \ref{GNS state} in the unital case. 

We are just a step away from proving Thm. \ref{GNN}. For a unital $C^{*}$-algebra $A$, let $\mathcal{S}(A)$ be the set of states on $A$. 

 \begin{ppsn}
 Let $A$ be a unital $C^{*}$-algebra and let $a \in A$ be self-adjoint. Then, there exists a state $\phi \in \mathcal{S}(A)$ such that 
 $|\phi(a)|=||a||$. 
  \end{ppsn}
  \textit{Proof.} We first give a proof assuming that $A$ is commutative. Thus, let $A=C(X)$, where $X$ is a compact Hausdorff space and let $f \in C(X)$ be a real valued function. 
  Since $X$ is compact and $f$ is continuous, there exists $x \in X$ such that $|f(x)|=||f||$. Let $\epsilon_x:A \to \bbc$ be the character given by evaluation at $x$. Then, $\epsilon_x$ is a 
  state on $A$ and $|\epsilon_x(f)|=||f||$. Thus, the proposition is proved in this case. 
  
  Now, let $A$ be an arbitrary unital $C^{*}$-algebra and let $a \in A$ be self-adjoint. Set $C:=C^{*}\{1,a\}$. Then, $C$ is commutative. Thus, there exists $\omega \in \mathcal{S}(C)$ such that 
  $|\omega(a)|=||a||$. Extend $\omega$ to a bounded linear functional say $\phi$ on $A$ such that $||\phi||=||\omega||$. This is possible, thanks to Hahn-Banach theorem. Since $\phi$ is an extension of $\omega$, 
  we have $\phi(1)=1$. It follows from Prop. \ref{bdd and pos} that $\phi$ is a state. Clearly, $|\phi(a)|=|\omega(a)|=||a||$. Hence the proof. \hfill $\Box$
  
  \begin{thm}
  Let $A$ be a unital $C^{*}$-algebra. Then, there exists a Hilbert space $\clh$ and a faithful unital representation $\pi$ of $A$ on $\clh$. 
    \end{thm}
    \textit{Proof.} Let $A_{sa}$ be the set of self-adjoint elements of $A$. For $a \in A_{sa}$, pick a state $\phi_a$ such that $|\phi_a(a)|=||a||$. For $a \in A_{sa}$, let $\pi_{a}$ be the GNS representation of $\phi_a$ on $\clh_a$. 
    Observe that for $a \in A_{sa}$, 
    \[
    ||a||=|\phi_a(a)|=|\langle \pi_a(a)\Omega_a|\Omega_a| \leq ||\pi_a(a)|| \leq ||a||.\]
    Thus, $||\pi_a(a)|| = ||a||$ for every $a \in A_{sa}$. 
    Define $\displaystyle \clh:=\bigoplus_{a \in A_{sa}}\clh_a$ and for $x \in A$, set $\displaystyle \pi(x)=\bigoplus_{a \in A_{sa}}\pi_a(x)$. 
  
  Let $x \in A_{sa}$ be given.  Note that \[||\pi(x)||=\sup\{||\pi_a(x)||:a \in A_{sa}\}.\] The set $\{||\pi_a(x)||:a \in A_{sa}\}$ is bounded by $||x||$ (since homomorphisms are contractive) and it contains $||x||$. Consequently, 
  $||\pi(x)||=||x||$ for every $x \in A_{sa}$. Applying Lemma \ref{isometric homo}, we obtain that $\pi$ is an isometry. Hence the proof. \hfill $\Box$
    
    \begin{rmrk}
    \label{separable rep}
    Keep the notation of the proof of the above theorem. Suppose $D \subset A_{sa}$ is a norm dense subset. Then, $\displaystyle \bigoplus_{a \in D}\pi_a$ is a faithful unital representation of $A$ on $\displaystyle \bigoplus_{a \in D}\clh_a$. 
    This is because the reasoning employed in the above proof implies that for every $a \in D$, $||\pi(a)||=||a||$. Now, 
    \[
    A_{sa} \ni a \to ||\pi(a)|| \in [0,\infty)~~\textrm{and~~}A_{sa}: \ni a \to ||a||\]
    are continuous functions and they agree on a dense set. Hence, $||\pi(a)||=||a||$ for every $a \in A_{+}$. 
    
    Thus, if $A$ is a unital separable $C^{*}$-algebra, then $A$ can be represented faithfully on a separable Hilbert space.  
        \end{rmrk}
        
        \begin{xrcs}
        Let $A$ be a non-unital $C^{*}$-algebra. Prove that there exists a Hilbert space $\clh$ and a faithful non-degenerate representation $\pi:A \to B(\clh)$.
\end{xrcs}
\textit{Hint:} Consider the unitisation $A^{+}$.

\begin{xrcs}
Let $X$ be a compact metric space and let $A:=C(X)$. Let $\mu$ be a probability measure on $X$. Define a state $\phi:A \to \bbc$ by 
\[
\phi(f)=\int f d\mu.\]
Recall that the Riesz-representation theorem says that every state on $A$ is of this form. 

Set $\clh:=L^{2}(X,\mu)$. For $f \in C(X)$, let $M_f:L^{2}(X,\mu) \to L^{2}(X,\mu)$ be defined by 
\[
M_f\xi(x):=f(x)\xi(x).\]
Show that the map $C(X) \ni f \to M_f \in B(L^{2}(X,\mu))$ is a unital $*$-representation, which we call the multiplication representation associated to $\mu$. Define $\Omega:=1$, i.e the constant function `1'. Prove that $\Omega$ is a cyclic vector for the multiplication representation $M$. Clearly, $\phi(f)=\langle M_f\Omega|\Omega\rangle$ for every $f \in C(X)$. 

Therefore, $(M,\clh,\Omega)$ is the GNS triple for the state $\phi$.

\end{xrcs}

\section{Double commutant theorem}
In this section, we prove von Neumann's double commutant theorem. Let us recall the weak operator topology and the strong operator topology. Let $\clh$ be a Hilbert space. For $\xi \in \clh$, let $||~||_{\xi}$ be 
the seminorm on $B(\clh)$ defined by 
\[
||T||_{\xi}:=||T\xi||.\]
The locally convex topology on $B(\clh)$ induced by the family of seminorms $\{||~||_\xi:\xi \in \clh\}$ is called \emph{the strong operator topology}, abbreviated SOT. Let $(T_i)_{i \in I}$ be a net in $B(\clh)$ and let $T \in B(\clh)$. Then, $T_i \to T$ in SOT if and only if $T_i\xi \to T\xi$ for every $\xi \in \clh$. 

For $\xi,\eta \in \clh$, let $||~||_{\xi,\eta}$ be the seminorm on $B(\clh)$ defined by 
\[
||T||_{\xi,\eta}:=|\langle T\xi|\eta \rangle|.\]
The  topology on $B(\clh)$ induced by the collection of seminorms $\{||~||_{\xi,\eta}:\xi,\eta \in \clh\}$ is called \emph{the weak operator topology}, abbreviated WOT. Let $(T_i)_{i \in I}$ be a net in $B(\clh)$ and let $T \in B(\clh)$. Then, $T_i \to T$ in SOT if and only if $\langle T_i\xi|\eta \rangle  \to \langle T\xi|\eta \rangle$ for every $\xi,\eta \in \clh$. 

Let $\mathcal{S} \subset B(\clh)$. We denote by $\overline{\mathcal{S}}^{sot}$  the closure of $\mathcal{S}$ in the strong operator topology. Similarly, we let $\overline{\mathcal{S}}^{wot}$ be the closure of $\mathcal{S}$ in the weak operator topology. 
Define 
\[
\mathcal{S}^{'}:=\{T \in B(\clh): Ta=aT \textrm{~for all $a \in \mathcal{S}$}\}.\] The set $\mathcal{S}^{'}$ is called \emph{the commutant of $\mathcal{S}$}. 

\begin{xrcs}
Let $\mathcal{S} \subset B(\clh)$. Observe that $\mathcal{S}^{'}$ is a unital subalgebra of $B(\clh)$. Show that $\mathcal{S}^{'}$ is both strongly and weakly closed. Prove that if $\mathcal{S}$ is $^*$-closed, i.e. $a^* \in \mathcal{S}$ whenever $a\in \mathcal{S}$, then $\mathcal{S}^{'}$ is a $^*$-closed. 
\end{xrcs}

\begin{thm}[Double commutant theorem]
\label{double commutant}
Let $A$ be a unital $*$-subalgebra of $B(\clh)$. Then, $\overline{A}^{sot}=A^{''}:=(A^{'})^{'}$. 
\end{thm}

The proof is based on the following two simple exercises.
\begin{xrcs}
\label{calculation of commutant}
Let $\clh$ be a Hilbert space. For $n \geq 1$, let $\clh_n:=\underbrace{\clh \oplus \clh \oplus \cdots \oplus \clh}_{\textrm{n times}}$. For a subset $\mathcal{S} \subset B(\clh)$, let 
\begin{align*}
M_n(\mathcal{S}):&=\{(x_{ij}): x_{ij} \in \mathcal{S}\} \subset B(\clh_n) \\
diag(\mathcal{S}):&=\{(x,x,\cdots,x): x \in \mathcal{S}\}.
\end{align*}
Suppose $A$ is a unital $^*$-algebra of $B(\clh)$. Prove that the commutant of $diag(A)$ is $M_n\big(A^{'}\big)$ and the double commutant $diag(A)^{''}$ is $diag\big(A^{''}\big)$. 
\end{xrcs}

\begin{xrcs}
Let $\mathcal{S} \subset B(\clh)$ be a $^*$-closed subset and suppose $W$ is a closed subspace of $\clh$. Let $P$ be the orthogonal projection onto $W$. Prove that the following are equivalent.
\begin{enumerate}
\item[(1)] The set $\mathcal{S}$ leaves $W$ invariant, i.e. $xW \subset W$ for every $x \in \mathcal{S}$.
\item[(2)] The projection $P$ commutes with $\mathcal{S}$, i.e. $Px=xP$ for every $x \in P$. 
\end{enumerate}
\end{xrcs}

\textit{Proof of Thm. \ref{double commutant}.} Since, $A^{''}$ is strongly closed and contains $A$, $\overline{A}^{sot} \subset A^{''}$. 

Let $\xi \in \clh$ and $T \in A^{''}$  be given. Define $W:=\overline{\{a\xi:a \in A\}}$. Denote the orthogonal projection onto $W$ by $P$. Note that $A$ leaves $W$ invariant. Hence, $P \in A^{'}$. 
Since $P \in A^{'}$ and $T \in A^{''}$, we have $TP=PT$. Consequently, $T$ leaves $W$ invariant. But note that $\xi \in W$ since $A$ contains $1_{\clh}$. Therefore, $T\xi \in W$. Thus, for every $\epsilon>0$, there exists $a \in A$ such that $||T\xi-a\xi|| <\epsilon$. 

Let $\xi_1,\xi_2,\cdots,\xi_n \in \clh$ and $T \in A^{''}$ be given. Set \[\xi:=(\xi_1,\xi_2,\cdots,\xi_n) \in \clh_n:=\underbrace{\clh \oplus \clh \oplus \cdots \oplus \clh}_{\textrm{n times}}.\] 
Then, $diag(T,T,\cdots,T) \in diag\big(A^{''}\big)=(diag(A))^{''}$. By what we have proved, it follows that for any $\epsilon>0$, there exists $a \in A$ such that 
\[
\sum_{i=1}^{n}||T\xi_i-a\xi_i||^{2} \leq \epsilon^{2}.\]
Thus, we have proved that $T \in \overline{A}^{sot}$. Thus, $A^{''} \subset \overline{A}^{sot}$. Hence the proof. \hfill $\Box$

\begin{rmrk}
In fact, it is not necessary to require that $A$ contains the identity operator of $\clh$. It suffices to require that $A$ acts non-degenerately, i.e. $A\clh$ is dense in $\clh$. This could be deduced from the unital case and with the help of `approximate identities', a notion
which we have not discussed and we refer the reader to either \cite{Arveson_invitation} or \cite{Davidson_Ken}. 
\end{rmrk}

\begin{dfn}
Let $\clh$ be a Hilbert space and let $M \subset B(\clh)$ be a unital $^*$-algebra. The algebra $M$ is called a von Neumann algebra if $M$ is strongly closed, or equivalently $M=M^{''}$. 
\end{dfn}
Note that von Neumann algebras are norm closed and hence $C^{*}$-subalgebras of $B(\clh)$. 

\begin{xmpl}
Let $(X,\mathcal{B},\mu)$ be a $\sigma$-finite measure space. For $f \in L^{\infty}(X,\mu)$, let $M_{f}$ be the operator on $\clh:=L^{2}(X,\mu)$ defined by 
\[
M_f\xi(x):=f(x)\xi(x).\]
Set $M:=\{M_f: f \in L^{\infty}(X,\mu)\}$. 
\begin{enumerate}
\item[(1)] Prove that the map $L^{\infty}(X,\mu) \ni f \to M_f \in B(L^{2}(X,\mu))$ is a $^*$-representation which is faithful. 
\item[(2)] Show that $M^{'}=M$. Deduce that $M$ is a von Neumann algebra. 
\item[(3)] Prove directly that $M$ is strongly closed. 
\end{enumerate}

Suppose that $X$ is a compact metric space, $\mathcal{B}$ is the Borel $\sigma$-algebra of $X$ and $\mu$ is a probability measure. Prove that $\{M_f: f \in C(X)\}$ is strongly dense in $M$. 
\end{xmpl}

Does Thm. \ref{double commutant} stay true if we replace SOT by WOT ? The answer is because $A^{''}=\overline{A}^{sot} \subset \overline{A}^{wot}$ and $\overline{A}^{wot} \subset A^{''}$ as $A^{''}$ is weakly closed. However, we don't need double commutant theorem to conclude that $\overline{A}^{sot}=\overline{A}^{wot}$. This is due to the fact that the set of continuous linear functionals on $B(\clh)$ remain the same whether we consider strong operator topology or the weak operator topology.

\begin{xrcs}
Let $V$ be a vector space over $\bbc$ and suppose $f,f_1,f_2,\cdots,f_n$ are linear functionals on $V$. Prove that if $\bigcap_{i=1}^{n}\ker(f_i) \subset \ker(f)$, then $f$ is a linear combination of $f_1,f_2,\cdots,f_n$. 
\end{xrcs}

\begin{lmma}
Let $\omega:B(\clh) \to \bbc$ be a linear functional. Suppose that $\omega$ is continuous w.r.t. to the weak operator topology. Then, there exist $\xi_1,\xi_2,\cdots,x_n,\eta_1,\eta_2,\cdots,\eta_n \in \clh$ such that 
\[
\omega(T)=\sum_{i=1}^{n}\langle T\xi_i|\eta_i\rangle.\]
Moreover, every functional of the above form is continuous w.r.t. the weak operator topology (as well as the strong operator topology). 
\end{lmma}
\textit{Proof.} 
For $\xi,\eta \in \clh$, let $\omega_{\xi,\eta}:B(\clh) \to \bbc$ be defined by 
\[
\omega_{\xi,\eta}(T):=\langle T\xi|\eta \rangle.\]
Let $\omega:B(\clh) \to \bbc$ be a linear functional that is continuous in WOT. Then, there exist $\xi_1,\xi_2,\cdots,\xi_n,\eta_1,\eta_2,\cdots,\eta_n \in \clh$ such that, for every $T \in B(\clh)$,
\[
||\omega(T)|| \leq \sum_{i=1}^{n}|\langle T\xi_i|\eta_i \rangle|.\]
The above inequality implies that $\bigcap_{i=1}^{n}\ker(\omega_{\xi_i,\eta_i}) \subset \ker(\omega)$. Hence, $\omega$ is a linear combination of $\{\omega_{\xi_i,\eta_i}: i=1,2,\cdots,n\}$. 
The proof of the first assertion is over. The remaining assertions are trivial. \hfill $\Box$

\begin{lmma}
Let $\omega:B(\clh) \to \bbc$ be a linear functional. Suppose that $\omega$ is continuous w.r.t. to strong operator topology. Then, there exist $\xi_1,\xi_2,\cdots,\xi_n,\eta_1,\eta_2,\cdots,\eta_n \in \clh$ such that 
\[
\omega(T)=\sum_{i=1}^{n}\langle T\xi_i|\eta_i\rangle.\]
\end{lmma}
\textit{Proof.} Let $\omega$ be as in the statement. Since $\omega$ is strongly continuous, there exist $\xi_1,\xi_2,\cdots,\xi_n \in \clh$ such that 
\[
|\omega(T)| \leq \Big(\sum_{i=1}^{n}||T\xi_i||^2 \Big)^{\frac{1}{2}}.\]
Let $\clh_n:=\clh \oplus \clh \oplus \cdots \clh$ and $\xi:=(\xi_1,\xi_2,\cdots,\xi_n) \in \clh_n$. Consider the subspace 
\[
\mathcal{D}:=\{(T\xi_1,T\xi_2,\cdots,T\xi_n): T \in B(\clh)\}.\]
The map $\Phi:\mathcal{D} \to \bbc$ defined by the rule
\[
\mathcal{D} \ni (T\xi_1,T\xi_2,\cdots,T\xi_n) \to \omega(T) \in \bbc\]
is a well-defined bounded linear function. Thus, it extends to a linear functional on $\clh_n$ which in turn is given by inner product against some vector $\eta=(\eta_1,\eta_2,\cdots,\eta_n)$. 
Thus, we get the equality
\[
\omega(T)=\Phi(T\xi_1,T\xi_2,\cdots,T\xi_n)=\langle diag(T)\xi|\eta\rangle=\sum_{i=1}^{n}\langle T\xi_i|\eta_i \rangle\] 
for every $T \in B(\clh)$. The proof is complete. \hfill $\Box$

A simple application of the Hahn-Banach theorem gives the following corollary. 
\begin{crlre}
Let $\mathcal{S} \subset B(\clh)$ be a convex subset. Then, $\overline{\mathcal{S}}^{sot}=\overline{\mathcal{S}}^{wot}$.

\end{crlre}

\section{A few more facts}
Here, we collect a few more facts that will be used without mention in the chapters that follow. We only state the results and the reader is recommended to consult \cite{Arveson_invitation} for proofs.

\textbf{Approximate identities:} 
Let $A$ be a $C^{*}$-algebra. Suppose $(e_\lambda)_{\lambda \in \Lambda}$ is a net of positive elements in $A$. The net $(e_\lambda)_{\lambda \in \Lambda}$ is called \emph{an approximate identity}
of $A$ if 
\begin{itemize}
\item $e_\lambda \leq e_\mu$ whenever $\lambda \leq \mu$,
\item for every $\lambda \in \Lambda$, $||e_\lambda|| \leq 1$, and
\item for every $x \in A$, $xe_\lambda \to x$ and $e_\lambda x \to x$. 
\end{itemize}
Approximate identities are needed only when $A$ is non-unital. 

\begin{thm}
Every $C^{*}$-algebra has an approximate identity. Moreover, if the $C^{*}$-algebra is separable, we can replace `nets' by `sequences'.
\end{thm}

\begin{xrcs}
Let $A:=C_{0}(\bbr^d)$. Construct an approximate identity of $A$. 
\end{xrcs}

\begin{xrcs}
Let $\clh$ be an infinite dimensional separable Hilbert space. Find an approximate identity for $\clk(\clh)$ where $\clk(\clh)$ is the $C^{*}$-algebra of compact operators on $\clh$. 
\end{xrcs}

\begin{xrcs}
Let $(T_i)$ be a net in $B(\clh)$ and let $T \in B(\clh)$. Let $D \subset \clh$ be a dense subset. Assume further that there exits $M>0$ such that $||T_i|| \leq M$ for every $i \in I$. 
Prove that $(T_i) \to T$ in SOT if and only if $T_i\xi \to T\xi$ for every $\xi \in D$. 

\end{xrcs}

\begin{xrcs}
Let $A$ be a $C^{*}$-algebra and let $(e_\lambda)$ be an approximate identity of $A$. Suppose $\pi:A \to B(\clh)$ is a representation. 
Show that the following are equivalent.
\begin{enumerate}
\item[(1)] The representation $\pi$ is non-degenerate.
\item[(2)] The net $(\pi(e_\lambda))_{\lambda \in \Lambda} \to 1_\clh$ strongly. 
\end{enumerate}
\end{xrcs}

We have earlier seen that for unital $C^{*}$-algebra, the norm of a positive linear functional is attained at $1$. We have the following analog in the non-unital situation.

\begin{ppsn}
Let $A$ be a $C^{*}$-algebra and let $(e_\lambda)$ be an approximate identity of $A$. Let $\phi:A \to \bbc$ be a linear functional. Suppose $\phi$ is positive. Then, $\phi$ is bounded and 
\[
||\phi||=\lim_{\lambda}\phi(e_\lambda).\]

\end{ppsn}

\textbf{Ideals and quotients:} Let $A$ be a $C^{*}$-algebra. For us, by an ideal in $A$, we mean a two sided ideal $I \subset A$ that is norm closed. We need the following facts
concerning ideals and quotients. 
\begin{enumerate}
\item[(1)] Let $I$ be an ideal in $A$. Then, $I$ is $^*$-closed, i.e. $x^{*} \in I$ whenever $x \in I$.
\item[(2)] Let $I$ be an ideal in $A$. Then, the quotient $A/I$ is a $^*$-algebra. Moreover, the quotient norm defined by 
\[
||x+I||:=\inf \{||x+y||: y \in I\}\]
is a $C^{*}$-norm. This way, $A/I$ becomes a $C^{*}$-algebra. 
\end{enumerate}

\begin{crlre}
Let $A$ and $B$ be $C^*$-algebras and suppose $\pi:A \to B$ is a $^*$-homomorphism. Then, $\pi(A)$ is norm closed in $B$. In other words, $\pi(A)$ is a $C^{*}$-subalgebra of $B$. 
\end{crlre}
\textit{Proof.} Let $I:=\ker(\pi)$. Clearly, $I$ is a norm closed two sided ideal in $A$. The map $\widetilde{\pi}:A/I \to B$ defined by
\[
\widetilde{\pi}(x+I):=\pi(x)\]
is an injective $^*$-homomorphism. Thanks to Ex. \ref{homo iso}, $\widetilde{\pi}$ is an isometry. Therefore, $\pi(A)=\widetilde{\pi}(A/I)$ is complete and consequently, $\pi(A)$ is norm closed. 
This completes the proof. \hfill $\Box$

\textbf{Tensor product of Hilbert spaces:} 
We also need the notion of tensor product of Hilbert spaces in the succeeding chapters. Let $\clh_1$ and $\clh_2$ be two Hilbert spaces. Denote the algebraic tensor product $\clh_1 \otimes_{alg} \clh_2$ by 
$\clh_{alg}$. It is not difficult, by making using of `the universal property' of tensor products, to prove that there exists a well defined map $\langle~|~\rangle: \clh_{alg} \times \clh_{alg} \to \bbc$ that is linear in the 
first variable and conjugate linear in the second variable satisfying 
\begin{equation}
\label{inner product tensor}
\langle \xi_1 \otimes \xi_2|\eta_1 \otimes \eta_2 \rangle=\langle \xi_1|\eta_1 \rangle \langle \xi_2|\eta_2 \rangle.
\end{equation}

\begin{ppsn}
The map $\langle~|~\rangle$ defines an inner product on $\clh_{alg}$. 
\end{ppsn}
\textit{Proof.} The only thing that requires verification is positive definiteness. Consider an arbitrary vector $u:=\sum_{i=1}^{n}\xi_i \otimes \eta_i \in \clh_{alg}$.
Let $W_1:=span\{\xi_i:i=1,2,\cdots,n\}$ and $W_2:=span\{\eta_i:i=1,2,\cdots,n\}$. Then, $W_i$ is a finite dimensional subspace of $\clh_i$. Choosing orthonormal 
bases for $W_1$ and $W_2$ and writing $\xi_i,\eta_i$ in terms of the basis elements, we can assume that $u$ is of the form 
\[
u=\sum_{k,\ell}a_{kl}\psi_{k} \otimes \phi_\ell\]
where $a_{kl} \in \bbc$, $\{\psi_k\}$ is an orthonormal basis for $W_1$ and $\{\phi_\ell\}$ is an orthonormal basis for $W_2$. 

Clearly, \[
\langle u|u \rangle=\sum_{k,\ell}|a_{kl}|^2.\]
The above equality implies that $\langle u|u \rangle \geq 0$ and $\langle u|u \rangle=0$ if and only if $u=0$. Hence the proof. \hfill $\Box$

The Hilbert space tensor product of $\clh_1$ and $\clh_2$, denoted $\clh_1 \otimes \clh_2$, is defined as the completion of $\clh_1 \otimes_{alg} \clh_2$ with respect the inner product defined in Eq. \ref{inner product tensor}. 
The reader should do the following exercises to acquaint herself/himself with tensor products. 

\begin{xrcs}
Let $\clh_1$ and $\clh_2$ be Hilbert spaces and let $\clh:=\clh_1 \otimes \clh_2$. 

\begin{enumerate}
\item[(1)] Suppose $\{e_i:i \in I\}$ is an orthonormal basis for $\clh_1$ and $\{f_j:j \in J\}$ is an orthonormal basis for $\clh_2$. Show that $\{e_i \otimes f_j: (i,j) \in I \times J\}$ is an orthonormal basis for $\clh_1 \otimes \clh_2$. 
\item[(2)] Suppose $T \in B(\clh_1)$. Show that there exists a unique operator, denote it by $T \otimes 1$, on $\clh_1 \otimes \clh_2$ such that 
\[
(T \otimes 1)(\xi \otimes \eta)=T\xi \otimes \eta.\]
Prove that the map $B(\clh_1) \ni T \to T \otimes 1 \in B(\clh_1 \otimes \clh_2)$ is an injective $^*$-homomorphism. 
\item[(3)] Suppose $T \in B(\clh_1)$ and $S \in B(\clh_2)$. Show that there exists a unique operator, denote it by $T \otimes S$, on $\clh_1 \otimes \clh_2$ such that 
\[
(T \otimes S)(\xi \otimes \eta)=T\xi \otimes S\eta.\]
\end{enumerate}
\end{xrcs}

\begin{xrcs}
Let $n \geq 1$ be given. Suppose $\clh$ is a Hilbert space. Let \[\clh_n:=\underbrace{\clh \oplus \clh \oplus \cdots\oplus \clh}_{\textrm{n times}}.\]
Show that the map $U:\clh_n \to \clh \otimes \bbc^n$ defined by 
\[
U(\xi_1,\xi_2,\cdots,\xi_n)=\sum_{i=1}^{n}\xi_i \otimes e_i\]
is a unitary. Here, $\{e_i:i=1,2,\cdots,n\}$ is the standard orthonormal basis for $\bbc^n$. 

For $T \in B(\clh)$, show that 
\[
U^{*}(T \otimes 1)U=diag(T,T,\cdots,T).\]

Do the same exercise by setting $n=\infty$. So, $\displaystyle \clh_{\infty}:=\bigoplus_{n \in \bbn}\clh$ and $\bbc^n$ is replaced by $\ell^2(\bbn)$. Define, in a similar fashion, an appropriate unitary $U:\clh_\infty \to \clh \otimes \ell^2(\bbn)$.
Prove that for $T \in B(\clh)$, 
\[
U^{*}(T \otimes 1)U=diag(T,T,T,\cdots).\]
\end{xrcs}

\begin{xrcs}
Let $(X_1,\mathcal{B}_1,\mu_1)$ and $(X_2,\mathcal{B}_2,\mu_2)$ be two $\sigma$-finite measure spaces. Set
$X:=X_1 \otimes X_2$, $\mathcal{B}=\mathcal{B}_1 \otimes \mathcal{B}_2$ and $\mu:=\mu_1 \times \mu_2$. 
For $f \in L^{2}(X_1,\mathcal{B}_1,\mu_1)$ and $g \in L^{2}(X_2,\mathcal{B}_2,\mu_2)$, let $f \otimes g \in L^{2}(X,\mathcal{B},\mu)$ be defined by 
\[
(f \otimes g)(x_1,x_2):=f(x_1)g(x_2).\]
Show that the map 
\[
L^{2}(X_1,\mu_1) \otimes L^{2}(X_2,\mu_2) \ni f \otimes g \to f \otimes g \in L^{2}(X,\mu)\]
is a unitary. 

\end{xrcs}

\chapter{Examples of $C^{*}$-algebras}

\section{Finite dimensional $C^{*}$-algebras}
In this section, we  show that a finite dimensional $C^{*}$-algebra is a direct sum of matrix algebras. The proof is an elegant application
of the double commutant theorem. 

\begin{dfn}
Let $A$ be a $C^{*}$-algebra and let $\pi:A \to B(\clh)$ be a representation. The representation $\pi$ is said to be irreducible if whenever $W$ is a closed subspace of $\clh$ invariant $\pi(A)$, then $W$ is 
either $\clh$ or $\{0\}$. 
\end{dfn}

Let $\pi$ be a representation of $A$ on $\clh$. Suppose $\pi$ is not irreducible. Then, $\clh$ can be written as $\clh_1 \oplus \clh_2$, with $\clh_1$ and $\clh_2$ non-zero, and for $a \in A$, the operator $\pi(a)$ has
the form 
\[
\pi(a)=\begin{bmatrix}
          \pi_1(a) & 0\\
          0 & \pi_2(a)
          \end{bmatrix}.\]
          Clearly, for $i=1,2$, $\pi_i$ is a representation of $A$ on $\clh_i$. 
Note that an irreducible representation is necessarily cyclic. For representations, $\pi_1$ and $\pi_2$ of $A$ on $\clh_1$ and $\clh_2$, let 
\[
L(\pi_1,\pi_2):=\{T \in B(\clh_1,\clh_2): T\pi_1(a)=\pi_2(a)T \textrm{~~for all $a \in A$}\}.\]
The set $L(\pi_1,\pi_2)$ is called \emph{the space of intertwiners} from $\pi_1$ to $\pi_2$. It is easily verifiable that, for representations $\pi_1$ and $\pi_2$ of $A$,  $L(\pi_1,\pi_2)^{*}=L(\pi_2,\pi_1)$. Also, for a representation $\pi$ of $A$, $L(\pi,\pi)$ is the commutant $\pi(A)^{'}$. 
For two representations $\pi_1$ and $\pi_2$ of $A$, we say that $\pi_1$ and $\pi_2$ are \emph{disjoint} if $L(\pi_1,\pi_2)=\{0\}$. 

\begin{lmma}[Schur]
\label{Sch}
Let $A$ be a $C^{*}$-algebra. Suppose $\pi:A \to B(\clh)$ is a representation. Then, the following are equivalent. 
\begin{enumerate}
\item[(1)]  The commutant $\pi(A)^{'}=\mathbb{C}$. 
\item[(2)] The representation $\pi$ is irreducible. 
\end{enumerate}
\end{lmma}
\textit{Proof.} Suppose that $(1)$ holds. Let $W \subset \clh$ be a closed subspace invariant under $\pi$. Denote the orthogonal projection onto $W$ by $P_W$. Then, $P_W \in \pi(A)^{'}$. Hence, $P_W$ is either $0$ or $1$. In other words, $W$ is either $\{0\}$ or $\clh$. Consequently, the representation $\pi$ is irreducible. This completes the proof of $(1) \implies (2)$.

Assume now that $\pi$ is irreducible. Note that $\pi(A)^{'}$ is a $^*$-subalgebra of $B(\clh)$. It suffices to show that every self-adjoint element of $\pi(A)^{'}$ is a scalar, or equivalently it is enough to prove that for a self-adjoint $T \in \pi(A)^{'}$, the spectrum $\sigma(T)$ is singleton.
 Let $T \in \pi(A)^{'}$ be self-adjoint. We claim that   the spectrum $\sigma(T)$ is singleton. 

Suppose not. Let $\lambda$ and $\mu$ be distinct elements in $\sigma(T)$. Choose non-zero functions $f,g \in C(\sigma(T))$ `supported around' $\lambda$ and $\mu$ respectively such that $fg=gf=0$. 
Let $X:=f(T)$ and $Y:=g(T)$. Note that $X,Y \in \pi(A)^{'}$ and $0 \neq Ran(X) \subset Ker(Y)$. The fact that $Y \in \pi(A)^{'}$ implies that $Ker(Y)$ is a closed subspace invariant under $\pi$. Since, $\pi$ is irreducible, this forces $Y=0$, which in turn implies $g=0$ which is a contradiction. Therefore, $\sigma(T)$ is singleton. This proves the claim and hence the proof. \hfill $\Box$

\begin{lmma}[Schur]
\label{Sch1}
Let $A$ be a $C^{*}$-algebra. Suppose $\pi_1$ and $\pi_2$ are irreducible representations of $A$ on $\clh_1$ and $\clh_2$ respectively. Then, the following are equivalent.
\begin{enumerate}
\item[(1)] The representations $\pi_1$ and $\pi_2$ are not unitarily equivalent. 
\item[(2)] The representations $\pi_1$ and $\pi_2$ are disjoint.
\end{enumerate}
\end{lmma}
\textit{Proof.} Clearly, $(2) \implies (1)$. Suppose that $(1)$ holds. Let $T \in L(\pi_1,\pi_2)$ be given. Then, $T^{*}T \in \pi_1(A)^{'}$ and $TT^{*} \in \pi_2(A)^{'}$. Thanks to Lemma \ref{Sch}, we have $T^{*}T=\lambda$ and $TT^{*}=\mu$ for some non-negative scalars $\lambda,\mu$. If either $\lambda$ or $\mu$ is zero, then $T=0$ and we are done. 

Suppose $\lambda,\mu>0$. Set $U:=\frac{1}{\sqrt{\lambda}}T$. Then, $U^{*}U=1$ and $UU^{*}=\frac{\mu}{\lambda}$. This forces that $U$ is an isometry and is onto. Consequently, $U$ is a unitary and clearly $U \in L(\pi_1,\pi_2)$. This means that $\pi_1$ and $\pi_2$ are unitarily equivalent, which is a contradiction to the hypothesis. The proof is complete. \hfill $\Box$

For the rest of this section, the letter $A$ stands for a finite dimensional $C^{*}$-algebra. Until further mention, we assume that $A$ is unital. 

\begin{ppsn}
There exists a Hilbert space $\clh$ which is  finite dimensional and a faithful representation $\pi:A \to B(\clh)$. 
\end{ppsn}
\textit{Proof.} Choose a separable Hilbert space $\clk$ and a faithful unital representation $\widetilde{\pi}:A \to B(\clh)$. This is possible, thanks to Remark \ref{separable rep}. Let $\{\xi_1,\xi_2,\cdots\}$ be an orthonormal basis of $\clk$. Define 
$\omega:A \to \bbc$ by 
\[
\omega(a):=\sum_{i=1}^{\infty}\frac{\langle \widetilde{\pi}(a)\xi_i|\xi_i\rangle}{2^{i-1}}.\]
Then, $\omega$ is a state. Moreover, $\omega$ is faithful, i.e. $\omega(a^*a)=0$ implies $a=0$. 

Let $(\clh_\omega,\pi_\omega,\Omega_\omega)$ be the GNS representation of the state $\omega$. Then, $\clh_\omega$ is finite dimensional and $\pi_\omega$ is faithful. To see the faithfulness of $\pi_\omega$, let $a \in A$ be 
such that $\pi_\omega(a^*a)=0$. Then, $\omega(a^*a)=\langle \pi_\omega(a^*a)\Omega_\omega|\Omega_\omega \rangle=0$. The faithfulness of $\omega$ implies that $a=0$. Hence the proof. \hfill $\Box$

\begin{ppsn}
Let $\pi:A \to B(\clh)$ be a unital representation and assume that $\clh$ is finite dimensional. Then, there exist irreducible representations $\pi_1,\pi_2,\cdots,\pi_r$ such that $\pi$ is equivalent to 
$\pi_1 \oplus \pi_2 \cdots \oplus \pi_r$. 
\end{ppsn}
\textit{Proof.} The proof is by induction on $dim(\clh)$. If $dim(\clh)=1$, then there is nothing to prove. Again, there is nothing to prove if $\pi$ is irreducible. If not, split $\clh$ as $\clh_1 \oplus \clh_2$ as non-zero invariant subspaces and then $\pi(a)$ is of 
the form 
\[
\pi(a)=\begin{bmatrix}
 \pi_1(a) & 0 \\
  0 & \pi_2(a)
  \end{bmatrix}.\]
  Now, $dim(\clh_1),dim(\clh_2)<dim(\clh)$. The proof now follows from induction. \hfill $\Box$

Choose a faithful unital  finite dimensional representation $\pi:A \to B(\clh)$. Decompose $\pi$ as $\pi=\pi_1 \oplus \pi_2 \oplus \cdots \pi_r$ into irreducible representations. Define an equivalence relation on $I:=\{1,2,\cdots,r\}$ as
$i \sim j$ if $\pi_i$ and $\pi_j$ are equivalent. Note that if $i \sim j$, then $\ker(\pi_i)=\ker(\pi_j)$. 

For each equivalence class $[i]$, pick a representative $s([i])$. Define $\displaystyle \widetilde{\pi}:=\bigoplus_{[i] \in I/\sim}\pi_{s([i])}$. 
Then,
\[
\ker(\widetilde{\pi})=\bigcap_{[i] \in I/\sim}\ker(\pi_{s([i])})=\bigcap_{i \in I}\ker(\pi_i)=\ker(\pi)=\{0\}.\]

We have, in particular, proved the following. 

\begin{ppsn}
\label{faithful irrep}
There exist irreducible representations $\pi_1,\pi_2,\cdots,\pi_s$ of $A$ such that 
\begin{enumerate}
\item[(1)] for $i\neq j$, $\pi_i$ and $\pi_j$ are disjoint, and
\item[(2)] the direct sum $\displaystyle \bigoplus_{i=1}^{s}\pi_i$ is faithful.
\end{enumerate}
\end{ppsn}

\begin{thm}
\label{structure of finite dimensional}
Let $A$ be a unital finite dimensional $C^{*}$-algebra. Then, $A$ is isomorphic to a direct sum of matrix algebras, i.e. there exist positive integers $n_1,n_2,\cdots,n_s$ such that 
$A$ is isomorphic to $M_{n_1}(\bbc) \oplus M_{n_2}(\bbc) \oplus \cdots \oplus M_{n_s}(\bbc)$. 
\end{thm}
\textit{Proof.} Let $\pi_1,\pi_2,\cdots,\pi_s$ be irreducible representations as in Prop. \ref{faithful irrep}. Let $\clh_i$ be the Hilbert space on which $\pi_i$ acts and
set $\displaystyle \clh:=\bigoplus_{i=1}^{s}\clh_i$. Define, 
\[
\pi:=\bigoplus_{i=1}^{s}\pi_i.\]
Since, $\pi$ is faithful, $A$ is isomorphic to $\pi(A)$. 

The reader should be able to do the following `matrix calculations'. As $L(\pi_i,\pi_j)=\{0\}$ if $i \neq j$ and $L(\pi_i,\pi_i)=\bbc$, we have 
\[
\pi(A)^{'}=\{(x_{ij}: x_{ij} \in L(\pi_j,\pi_i)\}=\{diag(\lambda_1,\lambda_2,\cdots,\lambda_s): \lambda_i \in \bbc\}.\]
Then, 
\[
\pi(A)^{''}=\{diag(T_1,T_2,\cdots,T_s): T_i \in B(\clh_i)\} \cong M_{n_1}(\bbc) \oplus M_{n_2}(\bbc) \oplus \cdots \oplus M_{n_s}(\bbc).\]
Here, $n_i:=dim(\clh_i)$. Thanks to the double commutant theorem, we have 
\[
A \cong \pi(A)=\pi(A)^{''} \cong M_{n_1}(\bbc) \oplus M_{n_2}(\bbc) \oplus \cdots \oplus M_{n_s}(\bbc).\]
Hence the proof. \hfill $\Box$

\begin{xrcs}
Let $I$ be a two sided ideal of $M_n(\bbc)$. Prove that $I$ is either $\{0\}$ or $M_n(\bbc)$. 
\end{xrcs}

\begin{xrcs}
Let $A_1,A_2,\cdots,A_n$ be $C^{*}$-algebras and set $\displaystyle A:=\bigoplus_{i=1}^{n}A_i$. Prove that if $I$ is an ideal of $A$, then $I$ is of the form
$\displaystyle I=\bigoplus_{i=1}^{n}I_i$ where $I_i$ is an ideal in $A_i$. 

\end{xrcs}

\begin{crlre}
A finite dimensional $C^{*}$-algebra is unital.
\end{crlre}
\textit{Proof.} Let $A$ be a finite dimensional $C^{*}$-algebra. Its unitisation $A^{+}$ is a unital finite dimensional $C^{*}$-algebra. 
By Thm. \ref{structure of finite dimensional}, it follows that $A^{+}$ is a direct sum of matrix algebras. Use the last two exercises to conclude 
that $A$ is a direct sum of matrix algebras which is clearly unital. \hfill $\Box$

\begin{xrcs}
Let $A$ be a $C^{*}$-algebra. Recall the partial order on the set of self-adjoint elements $A_{sa}$ given by $x \leq y$ if $y-x$ is positive. 
Suppose $p$ and $q$ are projections in $A$. Prove that $p \leq q$ if and only $pq=qp=p$. 
\end{xrcs}

For a $C^{*}$-algebra $A$, let 
\[
Z(A):=\{z \in A: za=az \textrm{~~for all $a \in A$}\}.\]
Note that $Z(A)$ is a commutative $C^{*}$-algebra. The $C^{*}$-algebra $Z(A)$ is called the center of $A$. A projection in $Z(A)$ is usually called a \emph{central projection}.
Denote the set of central projections of $A$ by $\mathcal{P}(Z(A))$. A projection $p \in \mathcal{P}(Z(A))$ is said to be a \emph{minimal central projection} if it satisfies the following. If $q \leq p$ with $q \in \mathcal{P}(Z(A))$, then $q$ is either $0$ or $p$. 

For a central projection $z \in Z(A)$, define 
\[
zAz:=\{a \in A: za=az=a\}.\]
Clearly, $zAz$ is a unital $C^{*}$-algebra where the multiplicative unit is $z$. 

\begin{xrcs}
Let $A=M_{m_1}(\bbc)\oplus M_{m_2}(\bbc) \cdots \oplus M_{m_r}(\bbc)$. For $i=1,2,\cdots,r$, let 
\[
z_i:=(0,0,\cdots,1_{m_i},0,\cdots,0).\]
\begin{enumerate}
\item[(1)] Show that $dim(Z(A))=r$. 
\item[(2)] Prove that the set of non-zero minimal central projections of $A$ is $\{z_1,z_2,\cdots,z_r\}$. 
\item[(3)] Observe that for $i=1,2,\cdots,r$, $z_iAz_i \cong M_{m_i}(\bbc)$. 
\item[(4)] Suppose $A \cong M_{n_1} \oplus M_{n_2}(\bbc) \oplus \cdots \oplus M_{n_s}(\bbc)$. Use $(1)$-$(3)$ to conclude that $r=s$ and there exists a permutation $\sigma$ of $\{1,2,\cdots,r\}$ such that $n_i=m_{\sigma(i)}$. 
\end{enumerate}
Thus, the decomposition of a finite dimensional $C^{*}$-algebra as a direct sum of matrix algebras, as in Theorem \ref{structure of finite dimensional}, is `unique'. 

\end{xrcs}

\section{The $C^{*}$-algebra of compact operators}

In this section, we discuss  $C^{*}$-algebra of  compact operators on a separable Hilbert space. We assume throughout that all the Hilbert spaces that we consider are separable.
Our convention is that the inner product is linear in the first variable and conjugate linear in the second variable. Let us recall the following facts about compact operators. 
Let $\clh$ be a separable Hilbert space. 

\begin{enumerate}
\item[(1)] A bounded linear operator $T:\clh  \to \clh $ is said to be compact if the following condition is satisfied. If$(x_n)$ is a bounded sequence in $\clh$, then $(Tx_n)$ has a convergent
subsequence in $\clh$. 
\item[(2)] Denote the set of compact operators on $\clh$ by $\clk(\clh)$. Then $\clk(\clh)$ is a norm closed two sided ideal in $B(\clh)$. Moreover, $\clk(\clh)$ is closed under taking adjoints. In short, $\clk(\clh)$ is a $C^{*}$-subalgebra of $B(\clh)$. 
\item[(3)] An operator $T:\clh \to \clh$ is said to be of finite rank if $Ran(T)$ is finite dimensional. Denote the set of finite rank operators on $\clh$ by $\mathcal{F}(\clh)$. Then $\mathcal{F}(\clh)$ is norm dense in $\clk(\clh)$. 
\end{enumerate}

For $\xi,\eta \in \clh$, let $\theta_{\xi,\eta} \in B(\clh)$ be defined by $\theta_{\xi,\eta}(\gamma)=\xi \langle \gamma|\eta \rangle$. Clearly, $\theta_{\xi,\eta}$ is of rank one and hence compact. Note the following relations. 
\begin{align*}
\theta_{\xi,\eta}^{*}&=\theta_{\eta,\xi} \\
\theta_{\xi_1,\eta_1}\theta_{\xi_2,\eta_2}&=\langle \xi_2|\eta_1 \rangle \theta_{\xi_1,\eta_2}\\
T\theta_{\xi,\eta}&=\theta_{T\xi,\eta} \\
\theta_{\xi,\eta}T&=\theta_{\xi,T^{*}\eta}
\end{align*}
for $\xi,\eta,\xi_1,\xi_2,\eta_1,\eta_2 \in \clh$ and $T \in B(\clh)$. 
The above relations imply that the linear span of $\{\theta_{\xi,\eta}: \xi,\eta \in \clh\}$ is a $^*$-closed subalgebra of $\clk(\clh)$. 
\begin{xrcs}
Prove that the linear span of $\{\theta_{\xi,\eta}: \xi,\eta \in \clh\}$ is $\mathcal{F}(\clh)$. 
\end{xrcs}

Observe that the map $\clh \times \clh \ni (\xi,\eta) \to \theta_{\xi,\eta} \in \mathcal{K}(\clh)$ is linear in the first variable and  conjugate linear in the second variable. Note that for $\xi,\eta \in \clh$, $||\theta_{\xi,\eta}||=||\xi||||\eta||$. These two facts and the fact that $\mathcal{F}(\clh)$ is dense in $\clk(\clh)$ together imply that 
if $D$ is a countable dense subset of $\clh$, then $\{\theta_{\xi,\eta}: \xi,\eta \in D\}$ is total in $\clk(\clh)$. Thus, $\clk(\clh)$ is a separable $C^{*}$-subalgebra of $B(\clh)$.

The two most important results regarding the $C^{*}$-algebra of compacts is that $\clk(\clh)$ is simple, i.e. it has no nontrivial two sided ideals\footnote{For us, ideals mean closed two sided ideals.} and that $\clk(\clh)$ has only one irreducible representation up to unitary equivalence. 

\begin{thm}
\label{simplicity of compacts}
Let $\clh$ be a separable Hilbert space. Then, $\clk(\clh)$ is simple. 
\end{thm}
First we prove the result assuming $\clh$ is finite dimensional. Suppose $\dim(\clh)=n$. Then, every linear operator on $\clh$ is compact and consequently, $\clk(\clh)$ is isomorphic to $M_n(\mathbb{C})$. 

\begin{lmma}
For $n \geq 1$, $M_n(\mathbb{C})$ is simple. 
\end{lmma}
\textit{Proof.} Let $n \geq 1$ be given. For $i,j \in \{1,2,\cdots,n\}$, let $e_{ij}$ be the matrix with $1$ at the $(i,j)^{th}$-entry and zero elsewhere. Then $\{e_{ij}\}_{i,j}$ forms a basis for $M_n(\mathbb{C})$. Note the following relations.
\begin{align*}
e_{ij}e_{kl}&=\delta_{jk}e_{il} \\
e_{ij}^{*}&=e_{ji}.
\end{align*}
for $i,j,k,l \in \{1,2,\cdots,n\}$. Let $I$ be a non-zero two sided ideal in $M_n(\mathbb{C})$. Pick a non-zero element $X \in I$. Write $X=\sum_{i,j}x_{ij}e_{ij}$. There exists $k,l$ such that $x_{kl} \neq 0$. Note that $e_{kk}Xe_{ll}=x_{kl}e_{kl}$. Hence $e_{kl} \in I$ as $I$ is a two sided ideal.   Let $i,j \in \{1,2,\cdots,n\}$ be given. Note that $e_{ij}=e_{ik}e_{kl}e_{lj}$. Hence $e_{ij} \in I$ for every $i,j$. But $\{e_{ij}\}_{i,j}$ is a basis for $M_{n}(\mathbb{C})$. As a consequence, it follows that $I=M_{n}(\mathbb{C})$. This completes the proof. \hfill $\Box$

\begin{lmma}
Let $A$ be a $C^{*}$-algebra. The following are equivalent.
\begin{enumerate}
\item[(1)] For every non-zero representation $\pi$, $||\pi(a)||=||a||$.
\item[(2)] The $C^{*}$-algebra $A$ is simple. 
\end{enumerate}
\end{lmma}
\textit{Proof.} Suppose $(1)$ holds. Let $I$ be a non-zero ideal in $A$. Let $\pi:A/I \to B(\clh_\pi)$ be a faithful representation of $A/I$. Denote the quotient map $A \to A/I$ by $q$. Then for every $a \in A$, $||\pi \circ q(a)||=||a||$. In other words, $||a||=||a+I||$ for every $a \in A$. Pick a non-zero element $x \in I$. Then, the previous equality implies that $||x||=||x+I||=0$ which is a contradiction. This proves $(1) \implies (2)$. 

Suppose $(2)$ holds. Let $\pi$ be a non-zero representation of $A$. Then, $ker(\pi)=\{0\}$. Hence $\pi$ is injective. But any injective $^*$-homomorphism is isometric. Thus $||\pi(a)||=||a||$ for every $a \in A$. Thus $(2) \implies (1)$ is proved. Hence the proof. \hfill $\Box$

Fix an orthonormal basis $\{\xi_1,\xi_2,\cdots \}$ for $\clh$. For $i,j$, let $E_{ij}=\theta_{\xi_i,\xi_j}$. Observe the following. 
\begin{align}
\label{matrix units}
E_{ij}E_{kl}&=\delta_{jk}E_{il}\\
E_{ij}^{*}&=E_{ji}
\end{align}
for $i,j,k,l\in \mathbb{N}$. Let $\cla_n$ be the linear span of $\{E_{ij}: i,j \in \{1,2,\cdots,n\}\}$. The above relations imply that $\cla_n$ is a $^*$-subalgebra of $\clk(\clh)$. Since $\cla_n$ is finite dimensional, it follows that $\cla_n$ is norm closed. 
Moreover the map $e_{ij} \to E_{ij}$ from $M_{n}(\mathbb{C}) \to \cla_n$ is an isometric $^*$-isomorphism (Why?). Thus, $\cla_n$ is simple. Observe that $\cla_n \subset \cla_{n+1}$ and $\cla:= \bigcup_{n \geq 1}\cla_n$ is norm dense in $\clk(\clh)$ (Why?).

\textit{Proof of Theorem \ref{simplicity of compacts}.} Let $\pi:\clk(\clh) \to B(\clh_0)$ be a non-zero representation. Since $\cla$ is dense in $\clk(\clh)$, it follows that there exists $n_0$ such that $\pi$ is non-zero on $\cla_{n_0}$. Since $\cla_{n_0} \subset \cla_{n}$ for $n \geq n_0$, it follows that $\pi$ restricted to $\cla_n$ is non-zero. But $\cla_n$ is simple. Consequently $\pi$ is isometric on $\cla_n$ for $n \geq n_0$. Since $\cla=\bigcup_{n \geq n_0}\cla_n$, it follows that $||\pi(a)||=||a||$ for every $a \in \cla$. Since $\cla$
is dense in $\clk(\clh)$, it follows that $||\pi(a)||=||a||$ for every $a \in \clk(\clh)$. Hence $\clk(\clh)$ is simple. This completes the proof. \hfill $\Box$

Next we derive a ``universal picture" of $\clk(\clh)$. Keep the foregoing notation. 
\begin{ppsn}
\label{universal compact}
Let $A$ be a $C^{*}$-algebra. Suppose there exists a system of matrix units $\{e_{ij}:i,j \in \mathbb{N}\}$ in $A$, i.e. the set $\{e_{ij}:i,j \in \mathbb{N}\}$ satisfies the following relations.
\begin{align*}
e_{ij}e_{kl}&=\delta_{jk}e_{il}\\
e_{ij}^{*}&=e_{ji}
\end{align*}
for $i,j,k,l \in \mathbb{N}$. Then there exists a unique $^*$-homomorphim $\pi:\clk(\clh) \to A$ such that for $i,j \in \mathbb{N}$, $\pi(E_{ij})=e_{ij}$. 
\end{ppsn}
\textit{Proof.} Note that $\{E_{ij}:i,j \in \{1,2,\cdots,n\}\}$ is a basis for $\cla_n$ for every $n$. Thus, there exists a linear map $\pi_n:\cla_n \to A$ such that $\pi_n(E_{ij})=e_{ij}$ for $i,j \in \{1,2,\cdots,n\}$. Clearly, $\pi_n$ is a $^*$-homomorphism. Since $\cla_{n}$ is simple, it follows that $\pi_n$ is isometric. The maps $(\pi_n)$'s are consistent, i.e. $\pi_{n+1}|_{\cla_n}=\pi_n$. Thus, there exists a $^*$-homomorphism $\pi:\cla \to A$ such that $\pi|_{\cla_n}=\pi_n$. As each $\pi_n$ is isometric, it follows that $\pi$ is isometric. 
Consequently, $\pi$ extends to a $^*$-homomorphism to the closure of $\cla$ which is $\clk(\clh)$. We denote the extension again by $\pi$. It is clear that $\pi$ is the required map. Uniqueness of $\pi$ is obvious. \hfill $\Box$

Derive the following ``coordinate free" description of the universal picture of $\clk(\clh)$. 
\begin{xrcs}
\label{co-ordinate free universal picture}
Let $D$ be a dense subspace of $\clh$ and $A$ be a $C^{*}$-algebra. Suppose that for $\xi,\eta \in D$, there exists $e_{\xi,\eta} \in A$ such that 
\begin{align*}
e_{\xi,\eta}^{*}&=e_{\eta,\xi} \\
e_{\xi_1,\eta_1}e_{\xi_2,\eta_2}&=\langle \xi_2|\eta_1 \rangle e_{\xi_1,\eta_2}
\end{align*}
for $\xi,\xi_1,\xi_2,\eta,\eta_1,\eta_2 \in D$. 
Prove that there exists a unique $^*$-homomorphism $\pi:\clk(\clh) \to A$ such that $\pi(\theta_{\xi,\eta})=e_{\xi,\eta}$ for $\xi,\eta \in D$. 
\end{xrcs}

Next, we study the representation theory of the algebra of  compact operators. The crucial facts regarding the representation theory of compacts are the following.
\begin{enumerate}
\item[(1)] Any non-degenerate representation of $\clk(\clh)$ is a direct sum of irreducible representations. 
\item[(2)] The only irreducible representation, up to unitary equivalence, of $\clk(\clh)$ is the identity representation. 
\end{enumerate}
This is the content of the next theorem.

\begin{xrcs}
Keep the foregoing notation. Let $E_n=\sum_{i=1}^{n}E_{ii}$. Note that $E_n$ is the projection onto the subspace spanned by $\{\xi_1,\xi_2,\cdots,\xi_n\}$. Clearly, $E_n \leq E_{n+1}$ for $n \geq 1$. Show the following.
\begin{enumerate}
\item[(1)] The sequence $(E_n) \to 1$ strongly, i.e. $E_n \xi \to \xi$ for every $\xi \in \clh$.
\item[(2)] For every finite rank operator $T$ on $\clh$, $TE_n \to T$ and $E_nT \to T$ in norm. 
\item[(3)] For every compact operator $T$ on $\clh$, $TE_n \to T$ and $E_n T \to T$ in norm. In other words, $(E_n)$ is an approximate identity of $\clk(\clh)$. 
\end{enumerate}
\end{xrcs}

\begin{lmma}
The identity representation of $\clk(\clh)$ on $\clh$ is irreducible. 
\end{lmma}
 \textit{Proof.} Let $W$ be a non-zero closed subspace of $\clh$ which is invariant under $\clk(\clh)$. Pick a unit vector $\eta \in W$. Note that $\theta_{\xi,\eta}(\eta)=\xi$. Thus $\xi \in W$ for every $\xi \in \clh$. This implies that $W=\clh$. 
 Hence the proof. \hfill $\Box$
\begin{thm}
\label{representation of compacts}
Let $\pi:\clk(\clh) \to B(\widetilde{\clh})$ be a non-degenerate representation. Then, there exists a Hilbert space $\clh_0$ and a unitary $U:\clh \otimes \clh_0 \to B(\widetilde{\clh})$ such that 
\[
\pi(A)=U(A \otimes 1)U^{*}\]
for $A \in \clk(\clh)$.
\end{thm}
\textit{Proof.} Since $\clk(\clh)$ is simple, $\pi$ is an isometry and therefore injective. Thus, for every $i,j$, $\pi(E_{ij})\neq 0$. Let $\clh_0$ be the range space of $\pi(E_{11})$. Denote the dimension of $\clh_0$ by $d$ and let $\{\eta_i\}_{i=1}^{d}$ be an orthonormal basis for $\clh_0$. We claim that $\{\pi(E_{i1})\eta_{j}\}_{i,j}$ is total in $\widetilde{\clh}$. Denote the closed linear span of $\{\pi(E_{i1})\eta_{j}\}_{i,j}$ by $\clh_1$. It is clear that $\pi(E_{rs})$ leaves $\clh_1$ invariant for every $r,s$. Since the linear span of $\{E_{rs}\}$ is norm dense in $\clk(\clh)$, it follows that $\clh_1$ is invariant under $\pi$ and so is $\clh_1^{\perp}$.

 Suppose $\clh_1^{\perp} \neq \{0\}$. Clearly, $\clh_0 \subset \clh_1$. Hence $\clh_1^{\perp} \subset \clh_0^{\perp}=Ker(\pi(E_{11}))$. Thus $\pi(E_{11})=0$ on $\clh_1^{\perp}$. On $\clh_1^{\perp}$, $\pi(E_{i1})^{*}\pi(E_{i1})=\pi(E_{11})=0$. This implies that $\pi(E_{i1})=0$, on $\clh_1^{\perp}$, for every $i$. This implies that $\pi(E_{ii})=\pi(E_{i1})\pi(E_{i1})^{*}=0$ on $\clh_1^{\perp}$. 
 For $n \in \bbn$, set $E_n:=\sum_{i=1}^{n}E_{ii}$. Then, $\pi(E_n)=0$ on $\clh_1^{\perp}$. However, $(E_n)$ is an approximate identity of $\clk(\clh)$ and therefore, $\pi(E_n) \to 1$ strongly. This forces that $\clh_1^{\perp}=0$. This proves the claim that 
 $\{\pi(E_{i1})\eta_{j}\}_{i,j}$ is total in $\widetilde{\clh}$.
 

Let $r,s \in \mathbb{N}$ and $j,k \in \{1,2,\cdots,d\}$ be given. Calculate as follows to observe that
\begin{align*}
\langle \pi(E_{r1})\eta_{j}|\pi(E_{s1})\eta_k\rangle & = \langle \pi(E_{1s})\pi(E_{r1})\eta_j|\eta_k \rangle \\
 &= \delta_{rs}\langle \pi(E_{11})\eta_{j}|\eta_{k} \rangle \\
 &=\delta_{rs} \langle \eta_j|\eta_k \rangle \\
 &=\delta_{rs}\delta_{jk}.
\end{align*}
The above calculation together with Exercise \ref{KRP} and the fact that $\{\pi(E_{i1})\eta_{j}\}_{i,j}$ is total in $\widetilde{\clh}$ ensures that there exists a unitary $U:\clh \otimes \clh_0 \to \widetilde{\clh}$ such that $U(\xi_i \otimes \eta_j)=\pi(E_{i1})\eta_j$. 
A direct calculation reveals that $U(E_{rs} \otimes 1)U^{*}=\pi(E_{rs})$. The proof is now completed by appealing to the fact that linear span of $\{E_{ij}: i,j\}$ is dense in $\clk(\clh)$. \hfill $\Box$

\begin{xrcs}
Let $\pi:\clk(\clh) \to B(\widetilde{\clh})$ be a non-degenerate representation. Suppose there exists a Hilbert space $\clh_0$ and a unitary $U:\clh \otimes \clh_0 \to \widetilde{\clh}$ such that 
\[
\pi(A)=U(A \otimes 1)U^{*}\]
for $A \in \clk(\clh)$. Show that $dim(\clh_0)$ is the dimension of the range space of $\pi(p)$ where $p$ is any rank one projection. The dimension of $\clh_0$ is called  \textbf{the multiplicity} of the identity representation in $\pi$. 
\end{xrcs}

We usually identify all infinite dimensional separable Hilbert space and reserve the letter $\clk$ to indicate the $C^{*}$-algebra of compact operators on a separable infinite dimensional Hilbert space. 

\section{Universal $C^{*}$-algebras}
Often, a $C^{*}$-algebra is prescribed in terms of generators and relations. We have already seen one such example, where in Prop. \ref{universal compact}, we have shown that   the $C^{*}$-algebra of compact operators $\clk$ is the universal $C^{*}$-algebra generated by a system of ``matrix units"
$\{e_{ij}:i,j \in \mathbb{N}\}$. We make this idea precise here. This section is based on the lectures given by Cuntz in 2014 at Oberwolfach. 

Let $\cla$ be a $^*$-algebra.  Let $p:\cla \to [0,\infty)$ be a map. We say that $p$ is a $C^{*}$-seminorm if 
\begin{enumerate}
\item[(1)] $p$ is seminorm on $\cla$, 
\item[(2)] $p$ satisfies the $C^*$-identity, i.e  for $x \in \cla$, $p(x^{*}x)=p(x)^{2}$, and
\item[(3)] for $x,y \in \cla$, $p(xy) \leq p(x)p(y)$.
\end{enumerate}
For $x \in \cla$, define \[||x||:=\sup\{p(x): ~\textrm{$p$ is a $C^{*}$-seminorm on $\cla$}\}.\] It is quite possible that $||x||$ is infinite for some $x \in \cla$.  Assume that $||x||<\infty$ for every $x \in \cla$. Let 
\[
I:=\{x \in \cla: ||x||=0\}.\]
Condition $(3)$ implies that $I$ is an ideal in $\cla$. Consider the quotient $\cla/I$. The seminorm $||~||$ descends to a $C^{*}$-norm on $\cla/I$. The completion of $\cla/I$ with respect to this norm is called
the universal $C^{*}$-algebra of $\cla$, or the enveloping $C^{*}$-algebra of $\cla$ and is usually denoted $C^{*}(\cla)$. 
\begin{xrcs}
Keep the foregoing notation. Show that for every $x \in \cla$, 
\begin{align*}
||x||&=\sup\{||\pi(x)||: \textrm{~$\pi$ is a $^*$-homomorphism from $\cla$ to a $C^{*}$-algebra}\}\\
&=\sup\{||\pi(x)||: \textrm{~$\pi$ is a non-degenerate representation of $\cla$} \}
\end{align*}
(Recall that a representation $\pi:\cla \to B(\clh)$ is said to be non-degenerate if $\pi(\cla)\clh$ is dense in $\clh$. If $\cla$ is unital, non-degenerate representations
are precisely unital representations).
\end{xrcs}
 
\begin{rmrk}
Note that $C^{*}(\cla)$ exists if and only if $||x||<\infty$ for every $x \in \cla$. 
\end{rmrk}
Consider the natural map $\cla \to C^{*}(\cla)$. We abuse notation and write the image of an element $x \in \cla$ under this map by $x$ itself. The $C^{*}$-algebra $C^{*}(\cla)$ is called the ``universal $C^{*}$-algebra of $\cla$" because it satisfies the following universal property. Keep the foregoing notation.

\begin{ppsn}
\label{universal property}
Let $B$ be a $C^{*}$-algebra and let  $\pi:\cla \to B$ be a $^*$-homomorphism. Then, there exists a unique $^*$-homomorphism $\widetilde{\pi}:C^{*}(\cla) \to B$ such that $\widetilde{\pi}(x)=\pi(x)$ for every $x \in \cla$. 
\end{ppsn}
\textit{Proof.} Uniqueness is obvious. For existence, let $p:\cla \to [0,\infty)$ be defined by $p(x)=||\pi(x)||$. Note that $p$ is a $C^{*}$-seminorm on $\cla$. Thus $||\pi(x)|| \leq ||x||$ for every $x \in \cla$. This implies that $\pi$ descends to a $^*$-homomorphim, say $\widetilde{\pi}:\cla/I \to B$. It is clear that $\widetilde{\pi}$ is bounded. Denote the extension to $C^{*}(\cla)$ again by $\widetilde{\pi}$. Then, $\widetilde{\pi}$ is the required map. \hfill $\Box$

Often, the algebra $\cla$ itself is given by generators and relations.  For example, consider the following statements

\begin{enumerate}
\item[(1)] Let $\cla$ be the universal unital $^*$-algebra generated by a single element $u$ such that $u^{*}u=1$ and $uu^{*}=1$. 
\item[(2)] Let $\cla$ be the universal unital $^*$-algebra generated by $v$ such that $v^{*}v=1$.
\item[(3)] Let $\cla$ be the universal unital $^*$-algebra generated by $P,Q$ such that the Heisenberg commutation relation $PQ-QP=1$ is satisfied. 
\item[(4)] Let $\cla$ be the universal  $^*$-algebra generated by $\{p_i\}_{i=1}^{n}$ such that $p_i^{2}=p_i=p_{i}^{*}$ and $p_{i}p_{j}=\delta_{ij}p_i$. 
\end{enumerate}
What is the meaning of each statement? For example in $(1)$, we mean that  there exists a $^*$-algebra $\cla$ which is generated by a single element $u$ and that $\cla$ has the following universal property : 
Suppose $\mathcal{B}$ is a unital $^*$-algebra and $w \in \mathcal{B}$ is such that $w^{*}w=ww^{*}=1$. Then, there exists a unique $^*$-homomorphism $\pi:\cla \to \mathcal{B}$ such that $\pi(u)=w$. It is not difficult to prove that the universal algebra $\cla$ is unique up to a unique isomorphism. We do the same for $(2)$, $(3)$ and $(4)$. The justification for the existence of such an $\cla$ is always by `abstract nonsense' which we do not discuss here. 

Consider now the following statement. 
Let $\mathcal{T}$ be the universal unital $C^{*}$-algebra generated by a single element $v$ such that $v^{*}v=1$. What do we mean by this ? First, we take the universal unital $^*$-algebra generated by $v$ such that $v^{*}v=1$. Denote it by $\cla$. Then, $\mathcal{T}=C^{*}(\cla)$. The $C^{*}$-algebra $\mathcal{T}$ is called the \emph{Toeplitz algebra} in the literature. 

But, does $\mathcal{T}$ exist ? Yes, it exists. For suppose $p$ is a $C^{*}$-seminorm on $\cla$. Define $I_{p}:=\{a \in \cla: p(a)=0\}$. Then $p$ descends to a $C^{*}$-norm on $\cla/I_{p}$. Let $A_{p}$ be the completion of $\cla/I_{p}$. Since $v+I_{p}$ is an isometry in $A_p$, we have $||v+I_{p}|| \leq 1$. Thus any word in $v$ and $v^{*}$ has $p$-norm at most $1$. Now, let $x$ be an element in $\cla$. Write $x=\sum_{\alpha}x_\alpha w_{\alpha}$ where $w_{\alpha}$ is a word in $v$ and $v^{*}$. Then $p(x) \leq \sum_{\alpha}|x_{\alpha}|$ and the latter bound is independent of $p$. Consequently $||x||<\infty$ for every $x$. 

\begin{rmrk}
\label{existence of universal $C^{*}$-algebra}
The argument outlined above works in the following situtation. Suppose $\cla$ is a $^*$-algebra generated by $\{x_{i}\}$ and each $x_{i}$ has $p$-norm atmost $1$ for every $C^{*}$-seminorm $p$ on $\cla$. Then $C^{*}(\cla)$ exists. However, $C^{*}(\cla)$ might be zero. 
\end{rmrk}

 The Toeplitz algebra has the following ``universal property"
\begin{xrcs}
Suppose $B$ is a unital $C^{*}$-algebra and $w \in B$ is such that $w^{*}w=1$. Then, there exists a unique $^*$-homomorphims $\pi:\mathcal{T} \to B$ such that $\pi(v)=w$. 
\end{xrcs}

Let us now show that the Toeplitz algebra is non-zero. We show this by producing a non-zero representation of $\mathcal{T}$. Consider the Hilbert space $\ell^{2}(\mathbb{N})$. Let $\{\delta_n\}_{n \geq 1}$ be the standard orthonormal basis for $\ell^{2}(\mathbb{N})$. Let $S:\ell^{2}(\mathbb{N}) \to \ell^{2}(\mathbb{N})$ be the unique operator such that $S(\delta_n)=\delta_{n+1}$. Then $S$ is an isometry, i.e. $S^{*}S=1$. 
The universal property of $\mathcal{T}$ guarantees that there exists a unique $^*$-homomorphism $\pi:\mathcal{T} \to C^{*}(S)$ such that $\pi(v)=S$. This guarantees that $v \neq 0$ and hence $\mathcal{T} \neq 0$. Later, we will show that $\pi$ is an isomorphism. 

\textit{Warning:} When one talks of the universal $C^{*}$-algebra given in terms of generators and relations, one should be cautious and decide first whether it exists or not and whether it is zero or  is non-zero. 
 We usually, but not always, apply Remark \ref{existence of universal $C^{*}$-algebra} to justify the existence of the universal $C^{*}$-algebra. To show, it is non-zero, we need to find a non-zero representation of the universal $^*$-algebra $\cla$ on a Hilbert space, or equivalently a non-zero $^*$-homomorphism from $\cla$ to a $C^{*}$-algebra. Use this to do the following exercises.

\begin{xrcs}
Show that the universal unital $C^{*}$-algebra generated by $u$ such that $u^{*}u=uu^{*}=1$ exists. 
\end{xrcs}

\begin{xrcs}
Show that the universal unital $C^{*}$-algebra generated by $\{p_{i}\}_{i=1}^{n}$ satisfying the relations $p_{i}^{2}=p_{i}=p_{i}^{*}$ and $p_{i}p_j=\delta_{ij}p_i$ exists. 

\end{xrcs}
Let us identify the universal $C^{*}$-algebras considered in the above two exercises concretely. 

\begin{ppsn}
\label{continuous function on the circle}
The algebra of continuous functions on the circle $\mathbb{T}$ denoted $C(\mathbb{T})$ is the universal $C^{*}$-algebra generated by $u$ such that $u^{*}u=uu^{*}=1$.
\end{ppsn}
\textit{Proof.} We denote the function $\mathbb{T} \ni z \to z \in \mathbb{C}$ by $z$ itself. Let $A$ be the universal $C^{*}$-algebra generated by $u$ such that $u^{*}u=uu^{*}=1$. 
Note that $u$ is a unitary in $A$. The continuous functional calculus gives a $^*$-homomorphism $C(\mathbb{T}) \to A$ which maps $z \to u$. Call it $\rho$. The universal property of $A$ gives a map $\pi:A \to C(\mathbb{T})$ such that $\pi(u)=z$. It is clear that $\pi \circ \rho(z)=u$ and $\rho \circ \pi(u)=z$. Since $u$ and $z$ generates $A$ and $C(\mathbb{T})$ respectively, it follows that $\pi$ and $\rho$ are inverses of each other. This completes the proof. \hfill $\Box$

\begin{ppsn}
\label{orthogonal projections universal}
Let $A$ be the universal $C^{*}$-algebra generated by $\{p_{i}:i \in \mathbb{N}\}$ such that $p_{i}^{2}=p_i=p_{i}^{*}$ and $p_ip_j=\delta_{ij}p_i$. 
Then $A \simeq C_{0}(\mathbb{N})$. 
\end{ppsn}
\textit{Proof.} We leave the proof that $A$ exists to the reader. Let $e_{i} \in C_{0}(\mathbb{N})$ be such that the $i$th coordinate of $e_i$ is $1$ and the rest of the coordinates are zero. It is clear that $e_{i}^{2}=e_i=e_i^{*}$ and $e_{i}e_{j}=\delta_{ij}e_i$. 
By the universal property, there exists a $^*$-homomorphism $\pi:A \to C_{0}(\mathbb{N})$ such that $\pi(p_i)=e_i$. 

\textit{Claim:} Let $B$ be a $C^{*}$-algebra and $q_1,q_2,\cdots,q_n$ be a finite sequence of orthogonal projections. Then for every $\lambda_1,\lambda_2,\cdots,\lambda_n \in \mathbb{C}$, 
\[
||\sum_{i=1}^{n}\lambda_iq_i|| \leq \sup_{1 \leq i \leq n}|\lambda_i|.\]
By representing $B$ faithfully on a Hilbert space say $\clh$, we can assume that $q_1,q_2,\cdots,q_n$ are operators on $\clh$. Then for a unit vector $\xi \in \clh$, we have
\begin{align*}
||\sum_{i=1}^{n}\lambda_i q_i \xi||^{2}&=\sum_{i=1}^{n}|\lambda_i|^{2}\langle q_{i}\xi|\xi \rangle \\
&\leq (\sup_{1 \leq i \leq n}|\lambda_i|)^{2}(\langle \sum_{i=1}^{n}q_i\xi|\xi \rangle)\\
& \leq (\sup_{1 \leq i \leq n}|\lambda_i|)^{2} ~(\textrm{ since $\sum_{i=1}^{n}q_i$ is a projection}).
\end{align*}
Hence, \begin{equation}
\label{estimate}
||\sum_{i=1}^{n}\lambda_iq_i|| \leq \sup_{1 \leq i \leq n}|\lambda_i|.
\end{equation}

Consider the dense $^*$-subalgebra $C_{c}(\mathbb{N})$ of $C_{0}(\mathbb{N})$. Note that $\{e_i: i \in \mathbb{N}\}$ is a basis for $C_{c}(\mathbb{N})$. Let $\rho:C_{c}(\mathbb{N}) \to A$ be the linear map such that $\rho(e_i)=p_i$.
Clearly $\rho$ is a $^*$-homomorphism. The estimate \ref{estimate} implies that $\rho$ is bounded. Denote the extension of $\rho$ to $C_0(\mathbb{N})$ by $\rho$ itself. Then, $\rho \circ \pi$ agrees with the identity map on the generators. 
Consequently, $\rho \circ \pi$ is identity. Similarly $\pi \circ \rho$ is identity. This shows that $\rho$ and $\pi$ are inverses of each other. Hence $\pi$ is an isomorphism. This completes the proof. \hfill $\Box$

Let us give a non-example. The universal $C^{*}$-algebra generated by two elements $P,Q$ such that $PQ-QP=1$ is zero. It suffices to show the following.
\begin{ppsn}
Let $\clh$ be a non-zero Hilbert space. Then, there does not exist bounded operators $P$ and $Q$ on $\clh$ such that $PQ-QP=1$. 
\end{ppsn}
\textit{Proof.} Suppose, on the contrary, assume that there exist $P,Q \in B(\clh)$ such that the commutator $[P,Q]= PQ-QP=1$. For a bounded operator $T$, let $\sigma(T)$ be the spectrum of $T$. Recall that for bounded operators $T,S$, $\sigma(TS) \cup \{0\}=\sigma(ST) \cup \{0\}$. 

Choose $\lambda \in \sigma(QP)$. Note that \[\lambda+1 \in \sigma(QP+1)=\sigma(PQ) \subset \sigma(PQ) \cup \{0\} \subset \sigma(QP) \cup \{0\}.\] The compactness of $\sigma(QP)$ implies that $\lambda+k=0$ for some positive integer $k$. 
 The fact that $\lambda \in \sigma(QP) \implies \lambda+1 \in \sigma(QP)\cup \{0\}$ implies that $-1 \in \sigma(QP)$. 

Then, $-1 \in \sigma(PQ)$. The relation $PQ-QP=1$ implies that  $-2 \in \sigma(QP)$ which in turn implies $-2 \in \sigma(PQ)$. By induction, we obtain $-k \in \sigma(PQ)$ for every positive integer $k$ which 
contradicts the fact that $\sigma(PQ)$ is bounded. This completes the proof. \hfill $\Box$

The notion of universal $C^{*}$-algebra is very handy and allows us to quickly define  group $C^{*}$-algebras and crossed products of discrete groups which provide important examples of $C^{*}$-algebras. 

\textbf{Group $C^{*}$-algebras:} Let $G$ be a discrete group.\footnote{We always assume some sort of separability hypothesis. For instance, we mostly assume topological spaces are   second countable, discrete groups are  countable, Hilbert spaces are separable etc... We do things as if this hypothesis is always there and make no explicit mention of this.} \emph{The full group $C^{*}$-algebra of $G$}, denoted $C^{*}(G)$,  is defined to be the universal $C^{*}$-algebra generated by $\{u_{s}:s \in G\}$ that satisfy the following relations.
\begin{align*}
u_su_t&=u_{st} \\
u_{s}^{*}&=u_{s^{-1}}
\end{align*}
for $s,t \in G$. 
Note that the above relations imply that $u_{e}$ is the multiplicative identity of $C^{*}(G)$, where $e$ is the identity element of $G$. Moreover $\{u_{s}:s \in G\}$ is a family of ``unitaries" and consequently $||u_s|| \leq 1$ for every $s \in G$. Thus, by Remark \ref{existence of universal $C^{*}$-algebra}, the $C^{*}$-algebra $C^{*}(G)$ exists. The next thing to show is that $C^{*}(G)$ is non-zero. 

Consider the Hilbert space $\ell^{2}(G)$ and let $\{\epsilon_t: t \in G\}$ be the standard orthonormal basis for $\ell^{2}(G)$. For $s \in G$, let $\lambda_s$ be the unitary operator on $\ell^{2}(G)$ such that 
\[
\lambda_{s}(\epsilon_t)=\epsilon_{st}\]
for $t \in G$. The map $G \ni s \to \lambda_{s} \in B(\ell^{2}(G))$ is called \emph{the left regular representation of $G$}.  Clearly, $\lambda_{s}\lambda_t=\lambda_{st}$ and $\lambda_{s}^{*}=\lambda_{s^{-1}}$. Thus, there exists a unique unital $^*$-homomorphism $\widetilde{\lambda}:C^{*}(G) \to B(\ell^{2}(G))$ such that $\widetilde{\lambda}(u_s)=\lambda_s$. This shows that $C^{*}(G)$ is non-zero. 

The image of $\widetilde{\lambda}$ is a $C^{*}$-subalgebra of $B(\ell^{2}(G))$ and is called \emph{the reduced $C^{*}$-algebra of $G$} and is denoted $C_{red}^{*}(G)$. Note that $C_{red}^{*}(G)$ is the $C^{*}$-algebra generated by $\{\lambda_s:s \in G\}$. Sometimes, we abuse notation and write $\widetilde{\lambda}$ simply by $\lambda$. It is natural to ask whether $\widetilde{\lambda}$ is an isomorphism. It turns out that $\widetilde{\lambda}$ is an isomorphism if and only if the group $G$ is \emph{amenable}. Abelian groups are amenable. An example of a non-amenable group is the free group on $2$ generators $\mathbb{F}_{2}$. 

Let us take a closer look at $C^{*}(G)$. Let $\cla$ be the universal $^*$-algebra generated by $\{u_{s}\}_{s \in G}$ such that $u_{s}u_t=u_{st}$ and $u_{s}^{*}=u_{s^{-1}}$. We first obtain a concrete description of $\cla$. Let $C_{c}(G)$ denote the space of finitely supported complex valued functions on $G$.  Define a $^*$-algebraic structure on $C_{c}(G)$ as follows. 
For $f,g \in C_{c}(G)$, let $f*g:G \to \mathbb{C}$  be defined by 
\[
f*g(s)=\sum_{t\in G}f(st)g(t^{-1}).\]
Note that the above sum is a finite sum. Note that $f*g \in C_{c}(G)$. The multiplication operation defined above is called \emph{the convolution}. Define a $^*$-operation on $C_{c}(G)$ by 
$f^{*}(s)=\overline{f(s^{-1})}$. 
\begin{xrcs}
Show that $C_{c}(G)$ with the convolution and the $^*$-operation defined above is a $^*$-algebra. 
\end{xrcs}

The algebra $C_{c}(G)$ is usually called \emph{the group algebra of $G$} and the usual notation is $C[G]$.
For $s \in G$, let $\delta_{s} \in C_{c}(G)$ be given by 
\begin{equation} 
\delta_s(t):=\begin{cases}
 1  & \mbox{ if
} t=s,\cr
   &\cr
    0 &  \mbox{ if } t \neq s.
         \end{cases}
\end{equation}
Observe that $\delta_{s}*\delta_{t}=\delta_{st}$ and $\delta_{s}^{*}=\delta_{s^{-1}}$. Note that $\delta_{e}$ is the multiplicative identity of $C_{c}(G)$. Thus, there exists a $^*$-homomorphism $\pi:\cla \to C_{c}(G)$ such that $\pi(u_s)=\delta_s$ for $s \in G$. 
\begin{lmma}
The map $\pi$ is an isomorphism. 
\end{lmma}
\textit{Proof.} We define the inverse map directly by setting $\rho(\delta_s)=u_s$. This is possible provided we can show that $\{\delta_s:s \in G\}$ is a basis for $C_{c}(G)$. We claim that $\{\delta_s: s \in G\}$ is a basis for $C_{c}(G)$. Let $f \in C_{c}(G)$ be given.
Then $f=\sum_{s \in G}f(s)\delta_s$. Moreover if $f=\sum_{s \in G}a_{s}\delta_s$, then applying the equality at an arbitrary point $t$, we get $f(t)=a_t$. This proves our claim. Let $\rho:C_{c}(G) \to \cla$ be the linear map such that $\rho(\delta_s)=u_s$. Then, clearly $\rho \circ \pi$ and $\pi \circ \rho$ agree with identity maps on the generators and hence agree with the identity maps everywhere. This shows that $\rho$ and $\pi$ are inverses of each other. Hence the proof. \hfill $\Box$

Now $C^{*}(G)$ is the completion of $C_{c}(G)$ where the norm on $C_{c}(G)$ is given by 
\[
||f||:=\sup\{\pi(f): \textrm{~$\pi$ is a unital representation of $C_{c}(G)$ on a Hilbert space}\}.\]

How do representations of $C_{c}(G)$ arise ? Let $\clh$ be a Hilbert space and $U:G \to B(\clh)$ be a map. We say that $U$ is a unitary representation if 
\begin{enumerate}
\item[(1)] for $s,t \in G$, $U_{s}U_{t}=U_{st}$, and
\item[(2)] for $s \in G$, $U_{s}$ is a unitary. 
\end{enumerate}
The set of unitaries $\mathcal{U}(\clh)$ is a group and a unitary representation of $G$ on $\clh$ is simply a group homomorphism from $G$ to $\mathcal{U}(\clh)$. 
Let $U:G \to \mathcal{U}(\clh)$ be a unitary representation. Then $U_{s}^{*}=U_{s^{-1}}$ for $s \in G$. Thus, there exists a unique unital $^*$-homomorphism, denoted $\pi_{U}:C_{c}(G) \to B(\clh)$ such that 
\[
\pi_{U}(\delta_s)=U_s
\]
for $s \in G$. Conversely, suppose $\pi$ is a unital representation of $C_{c}(G)$ on a Hilbert space $\clh$. Set $U_{s}=\pi(\delta_s)$. Then $\{U_{s}\}_{s \in G}$ is a unitary representation of $G$. Clearly $\pi_{U}$ and $\pi$ agree on $\{\delta_{s}:s \in G\}$. Since $\{\delta_{s}:s \in G\}$ is a basis for $C_{c}(G)$, it follows that $\pi=\pi_{U}$. Thus, representations of $C_{c}(G)$ are the ``same" as the unitary representations of the group $G$. 
Thus for $f \in C_{c}(G)$, 
\[
||f||_{C^{*}(G)}:=\sup\{||\pi_{U}(f)||:  U \textrm{~is a unitary representation of $G$}\}.\]  

In particular, we have proved that representations of $C^{*}(G)$ are in one-one correspondence with representations of $C_{c}(G)$. To summarise, \textit{$C^{*}(G)$ is the $C^{*}$-algebra that  captures the representation theory of the group $G$}.

\begin{rmrk}
The map $U \to \pi_{U}$ respects unitary equivalence, irreducibility, direct sum, etc..... Thus, the study of the representation theory of groups is equivalent to the study of the representation theory of the associated full group $C^{*}$-algebra. 
This has advantages, for then we can use (operator) algebraic techniques. The following theorem is a typical application.
\end{rmrk}

\begin{xrcs}
\label{finite irrep}
Let $A_1,A_2,\cdots,A_r$ be $C^{*}$-algebras such that for every $i=1,2,\cdots,r$, $A_i$ has exactly one irreducible representation. Show that $A:=\bigoplus_{i=1}^{r}A_i$ has exactly $r$ irreducible representations. 
\end{xrcs}

\begin{thm}
Let $G$ be a finite group. Then, up to unitary equivalence, $G$ has only finitely many irreducible representations. 
\end{thm}
\textit{Proof.} It suffices to prove that $C^{*}(G)$ has only finitely many irreducible representations. Note that $C^{*}(G)$ is finite dimensional. Therefore, there exists positive integers $n_1,n_2,\cdots,n_r$ such that $C^{*}(G)$ is isomorphic to $M_{n_1}(\bbc) \oplus M_{n_2}(\bbc) \oplus \cdots \oplus M_{n_r}(\bbc)$. The proof is completed by appealing to Ex. \ref{finite irrep} and the fact that $M_n(\bbc)$ has only one irreducible representation. This completes the proof. \hfill $\Box$

Let us identify the full $C^{*}$-algebra of a discrete abelian group. Let $G$ be a discrete abelian group. Denote the `unit circle' by $\bbt$, i.e. 
\[
\bbt:=\{z \in \bbc: |z|=1\}.\]
Note that $\bbt$ is a multiplicative group. A homomorphism $\chi:G \to \bbt$ is called \emph{a character of $G$}. We denote the set of characters of $G$ by $\widehat{G}$. 

The set $\widehat{G}$ has a group structure, where the group multiplication is pointwise multiplication. The map $G \ni s \to 1 \in \mathbb{T}$ is the identity element of $G$. For $\chi \in \widehat{G}$, the inverse of $\chi$ is $\overline{\chi}$. We endow $\widehat{G}$ with the topology of pointwise convergence, i.e. the product topology. The convergence of nets is as follows. Suppose $(\chi_{\alpha})$ is a net in $\widehat{G}$ and $\chi \in \widehat{G}$. Then $\chi_\alpha \to \chi$ if and only if $\chi_{\alpha}(s) \to \chi(s)$ for every $s \in G$. By Tychonoff theorem,  $\widehat{G}$ is compact. It is routine to check that $\widehat{G}$ is a topological group. 

\begin{ppsn}
Let $G$ be a discrete abelian group. Then $C^{*}(G)$ is isomorphic to $C(\widehat{G})$. 
\end{ppsn}
\textit{Proof.} Since the group $G$ is abelian, it follows that $A:=C^{*}(G)$ is commutative. Moreover $A$ is unital. It suffices to show that $\widehat{A}$ is homeomorphic to $\widehat{G}$. Let $\chi \in \widehat{G}$ be given. Then, by the universal property, there exists a homomorphism $\omega_\chi:A \to \mathbb{C}$ such that $\omega_\chi(u_s)=\chi(s)$. Since $\{u_{s}:s \in G\}$ generates $C^{*}(G)$, the map 
\[
\widehat{G} \ni \chi \to \omega_\chi \in \widehat{A}\]
is $1$-$1$. 
Let $\omega:A \to \mathbb{C}$ be a character. Since $\{u_{s}:s \in G\}$ is a set of unitaries, it follows that for every $s \in G$, $\omega(u_s) \in \mathbb{T}$. Define $\chi:G \to \mathbb{T}$ by $\chi(s)=\omega(u_s)$. It is clear that $\chi$ is a character of $G$ and $\omega=\omega_\chi$. This proves that the map $\widehat{G} \ni \chi \to \omega_\chi \in \widehat{A}$ is onto. 

Let $(\chi_{\alpha})$ be a net in $\widehat{G}$ such that $(\chi_{\alpha}) \to \chi \in \widehat{G}$. Note that $\{\omega_{\chi_\alpha}\}$ is uniformly bounded. Thus to show $\omega_{\chi_\alpha} \to \omega_{\chi}$, it suffices to check $\omega_{\chi_\alpha}(x) \to \omega_{\chi}(x)$ for $x$ in a total set $F$ of $A$. Set $F:=\{u_{s}: s \in G\}$ and observe that $F$ is total in $A$. Clearly $\omega_{\chi_\alpha}(x) \to \omega_{\chi}(x)$ for every $x \in F$. Hence, the map $\widehat{G} \ni \chi \to \omega_{\chi} \in \widehat{A}$ is continuous. Since $\widehat{G}$ and $\widehat{A}$ are both compact and Hausdorff, it follows that the map $\chi \to \omega_{\chi}$ is a homeomorphism. This completes the proof. \hfill $\Box$

\textbf{Crossed products:} Let $G$ be a discrete group and let $A$ be a $C^{*}$-algebra. By an \emph{action of $G$ on $A$}, we mean a family $\alpha:=\{\alpha_{s}\}_{s \in G}$ of automorphisms of $A$ such that $
\alpha_{s} \circ \alpha_{t}=\alpha_{st}$ for $s,t \in G$. Such a triple $(A,G,\alpha)$ is called a \emph{ $C^{*}$-dynamical system}. Here is an example of a dynamical system. 
\begin{xmpl}
Let $X$ be a locally compact Hausdorff space and suppose $G$ is a discrete group that acts on $X$ via homeomorphisms on the left. For $s \in G$ and $f \in C_{0}(X)$, define
\[
\alpha_{s}(f)(x)=f(s^{-1}x).\]
Then $\alpha:=\{\alpha_s\}_{s \in G}$ is an action of $G$ on $C_{0}(X)$. Note that $G$ leaves $C_{c}(X)$ invariant, where $C_{c}(X)$ is the dense subalgebra of compactly supported continuous functions on $X$. 

\end{xmpl}
Let $(A,G,\alpha)$ be a $C^{*}$-dynamical system. Assume that $A$ is unital. The \emph{full crossed product}, denoted $A \rtimes_{\alpha} G$, is defined to be the universal unital $C^{*}$-algebra generated by a copy of $A$ and unitaries $\{u_{s}\}_{s \in G}$ such that $u_{s}u_{t}=u_{st}$ and $u_{s}au_{s}^{*}=\alpha_{s}(a)$ for $s \in G$ and $a \in A$. Note that, if $A=\mathbb{C}$ and $\alpha$ is the trivial action, then $A \rtimes_{\alpha} G \simeq C^{*}(G)$. 

Let us take a closer look at the $C^{*}$-algebra $A \rtimes_{\alpha} G$. First, let $\mathcal{B}$ be the universal $^*$-algebra generated by a copy of $A$ and unitaries $\{u_{s}\}_{s \in G}$ that satisfy the relations mentioned above. The relations imply that the linear span of $\{a_{s}u_{s}:a_{s} \in A, s \in G\}$ is $\mathcal{B}$. Our experience with group $C^{*}$-algebras suggest that we should treat $\mathcal{B}$ as the algebra of functions defined on $G$ but now taking values in $A$. Thus consider $C_{c}(G,A)$, i.e. the set of functions $f:G \to A$ such that $f$ is finitely supported. 

We make $C_{c}(G,A)$ into a $^*$-algebra by defining the multiplication and the $^*$-operation as follows: for $f,g \in C_{c}(G,A)$, 
\begin{align*}
f*g(s)&=\sum_{t \in G}f(t)\alpha_{t}(g(t^{-1}s)) \\
f^{*}(s)&=\alpha_{s}(f(s^{-1})^{*})
\end{align*}
for $f,g \in C_{c}(G,A)$. It is  a routine exercise to verify that the multiplication and the $^*$-operation defined above makes $C_{c}(G,A)$ into a $^*$-algebra. Note that for $f,g \in G$, 
$f*g(s)=\displaystyle \sum_{t \in G}f(st^{-1})\alpha_{st^{-1}}(g(t))$. 
For $a \in A$ and $s \in G$, denote the element of $C_{c}(G,A)$ which vanishes at points other than $s$ and whose value at $s$ is $a$ by $a \otimes \delta_s$. Note that for $f \in C_{c}(G,A)$, 
$f=\sum_{s \in G}f(s) \otimes \delta_s$. 

\begin{xrcs}
Keep the foregoing notation. Prove that for $a,b \in A$ and $s,t \in G$,
\begin{align*}
(a \otimes \delta_s)*(b \otimes \delta_t)&=a\alpha_{s}(b)\otimes \delta_{st} \\
(a \otimes \delta_s)^{*}&=\alpha_{s^{-1}}(a^{*}) \otimes \delta_{s^{-1}}
\end{align*}
\end{xrcs}
The above relations and the universal property of $\mathcal{B}$ together imply that there exists a $^*$-homomorphism $\lambda:\mathcal{B} \to C_{c}(G,A)$ such that $\lambda(a)=a \otimes \delta_{e}$ and $\lambda(u_s)=1 \otimes \delta_s$. 
Let $\mu:C_{c}(G,A) \to \mathcal{B}$ be defined by $\mu(f)=\sum_{s \in G}f(s)u_{s}$. Note that $\mu$ is a $^*$-homomorphism. (The multiplication and the $^*$-operation are defined in such a way on $C_{c}(G,A)$ precisely to make this map a homomorphism). Clearly, $\lambda \circ \mu=Id$. Note that  $\mu \circ \lambda$ agrees with the identity map on $A$ and $\{u_s:s \in G\}$ which generates $\mathcal{B}$ as an algebra. Thus $\mu \circ \lambda=Id$. Hence $\lambda$ and $\mu$ are inverses of each other. 

Then, the full crossed product $A \rtimes G$ is the enveloping $C^{*}$-algebra of $C_{c}(G,A)$. Note that the $^*$-algebra $C_{c}(G,A)$ makes sense, even if $A$ is not unital. 

\begin{dfn}
Let $(A,G,\alpha)$ be a $C^*$-dynamical system.  The full crossed product, denoted $A \rtimes_{\alpha} G$, is defined to be the enveloping $C^{*}$-algebra of $C_{c}(G,A)$. 
\end{dfn}

We need to show that $A \rtimes_{\alpha} G$ exists and is non-zero. This requires us to prove that the universal norm is finite and we are forced to understand non-degenerate representations of $C_{c}(G,A)$ in more concrete terms. We will make use of the following remark in the sequel. 
\begin{rmrk}
\label{KRP}
 Suppose $\clh_1$ and $\clh_2$ are Hilbert spaces and $S_1$ and $S_2$ are total subsets of $\clh_1$ and $\clh_2$ respectively. Let $\phi:S_1 \to S_2$ be a map such that $\langle \phi(x)|\phi(y) \rangle=\langle x|y \rangle$ for $x,y \in S_1$. Then there exists a unique isometry $V:\clh_1 \to \clh_2$ which extends $\phi$. Moreover, if $\phi$ is a bijection, the isometry $V$ is a unitary.
\end{rmrk}

Let $\lambda:C_{c}(G,A) \to B(\clh)$ be a non-degenerate representation. Since the linear span of  $\{a \otimes \delta_{s}: a \in A, s \in G\}$ is $C_{c}(G,A)$, it follows that $\{\lambda(a \otimes \delta_s)\xi:a \in A, s \in G, \xi \in \clh\}$ is total in $\clh$. Fix $r \in G$. For $a,b \in A$, $s,t \in G$ and $\xi,\eta \in \clh$, calculate as follows to observe that 
\begin{align*}
&\langle \lambda(\alpha_{r}(a) \otimes \delta_{rs})\xi|\lambda(\alpha_{r}(b) \otimes \delta_{rt})\eta \rangle \\
&= \langle \xi|\lambda(\alpha_{r}(a) \otimes \delta_{rs})^{*}\lambda(\alpha_{r}(b) \otimes \delta_{rt})\eta \rangle \\
&=\langle \xi|\lambda((\alpha_{s^{-1}r^{-1}}(\alpha_r(a^{*})) \otimes \delta_{s^{-1}r^{-1}})*(\alpha_r(b) \otimes \delta_{rt}))\eta \rangle \\
&=\langle \xi|\lambda(\alpha_{s^{-1}}(a^{*}b) \otimes \delta_{s^{-1}t})\eta \rangle \\
&=\langle \xi|\lambda(a \otimes \delta_{s})^{*}\lambda(b \otimes \delta_t)\eta \rangle \\
&=\langle \lambda(a \otimes \delta_s)\xi|\lambda(b \otimes \delta_t)\eta \rangle. 
\end{align*}
Appealing to Remark \ref{KRP}, we conclude that there exists a unique unitary, denoted $U_{r}$, such that $U_{r}(\lambda(a \otimes \delta_s)\xi)=\lambda(\alpha_{r}(a) \otimes \delta_{rs})\xi$ for $a \in A$, $s \in G$ and $\xi \in \clh$. By evaluating on the total set 
$\{\lambda(a \otimes \delta_s)\xi: a \in A, s \in G, \xi \in \clh\}$, we conclude that $U_{r}U_{s}=U_{rs}$ for every $r,s \in G$. Thus, $U:=\{U_s\}_{s \in G}$ is a unitary representation of $G$ on $\clh$. 

Define for $a \in A$, $\pi(a)=\lambda(a \otimes \delta_{e})$. Then $\pi$ is a $^*$-representation of $A$ on $\clh$. 

\begin{xrcs}
 By evaluating on the total set $\{\lambda(a \otimes \delta_s)\xi:a \in A, s \in G,\xi \in \clh\}$, prove that 
\begin{enumerate}
\item[(1)] the representation $\pi$ is a $^*$-representation,
\item[(2)] the family $U:=\{U_{s}\}_{s \in G}$ is a unitary representation of $G$,
\item[(3)] for $a \in A$ and $s \in G$, $U_{s}\pi(a)U_{s}^{*}=\pi(\alpha_s(a))$, or equivalently $U_{s}\pi(a)=\pi(\alpha_s(a))U_{s}$. 
\end{enumerate}
Such a pair $(\pi,U)$ is called a covariant representation of the dynamical system. 
\end{xrcs}
Keep the foregoing notation. The representation $\lambda$ can be recovered from the pair $(\pi,U)$. Note, again by evaluating on the total set $\{\lambda(a \otimes \delta_s)\xi:a \in A, s \in G,\xi \in \clh\}$, that $\lambda(a \otimes \delta_s)=\pi(a)U_{s}$. Since the set $\{\lambda(a \otimes \delta_s)\xi:a \in A, s \in G, \xi \in \clh\}$ is total, it follows that $\{\pi(a)U_{s}\xi: a \in A, s \in G, \xi \in \clh\}$ is total in $\clh$. This implies that the representation $\pi$ is non-degenerate. 

We can reverse the above process. First a definition. 
\begin{dfn}
\label{covariant}
Consider a $C^{*}$-dynamical system $(A,G,\alpha)$. Let $\pi:A \to B(\clh)$ be a representation and $U:G \to \mathcal{U}(\clh)$ be a unitary representation. We say that the pair $(\pi,U)$ is a covariant representation of the dynamical system $(A,G,\alpha)$ if for $a \in A$, $s \in G$, \[U_{s}\pi(a)U_{s}^{*}=\pi(\alpha_s(a)).\] 
We always assume that $\pi$ is non-degenerate. 
\end{dfn}
Let $(\pi,U)$ be a covariant representation of the dynamical system $(A,G,\alpha)$. Define a map $\lambda:C_{c}(G,A) \to B(\clh)$ by \[\lambda(f)=\sum_{s \in G}\pi(f(s))U_{s}.\] It is clear that $\lambda(f*g)=\lambda(f)\lambda(g)$ and $\lambda(f)^{*}=\lambda(f^{*})$ if $f$ and $g$ are of the form $a \otimes \delta_s$. But since $\{a\otimes \delta_s\}$ spans $C_{c}(G,A)$, it follows that $\lambda$ is a $^*$-homomorphism. Note that $\lambda(a \otimes \delta_{e})=\pi(a)$. Hence $\lambda$ is non-degenerate. We denote this map $\lambda$ by $\pi \rtimes U$. 
\begin{xrcs}
Let $(\pi,U)$ be a covariant representation of $(A,G,\alpha)$. Prove that the covariant representation that we obtain  if we apply the process described before Definition \ref{covariant} to the non-degenerate representation $\pi \rtimes U$  is $(\pi,U)$.
\end{xrcs} 

Thus, the  non-degenerate representations of the $^*$-algebra $C_{c}(G,A)$ are in $1$-$1$ correspondence with the covariant representations of the dynamical system $(A,G,\alpha)$. Therefore, the universal norm on $C_{c}(G,A)$ is given by 
\[
||f||=\sup\{||(\pi \rtimes U)(f)||: (\pi,U) \textrm{~is a covariant representation of $(A,G,\alpha)$}\}. \]
For $f \in C_{c}(G,A)$, let \[||f||_{1}:=\sum_{s \in G}||f(s)||.\]

 Let $(\pi,U)$ be a covariant representation of $(A,G,\alpha)$. Note that for $f \in C_{c}(G,A)$, 
\[||(\pi \rtimes U)(f)||=||\sum_{s \in G}\pi(f(s))U_{s}|| \leq \sum_{s \in G}||\pi(f(s))||||U_s|| \leq \sum_{s \in G}||f(s)||= ||f||_{1}.\] Hence $||f|| \leq ||f||_{1}$ for every $f \in C_{c}(G,A)$. This proves that $||~||$ is a genuine $C^{*}$-seminorm on $C_{c}(G,A)$. Next we show that $||~||$ is indeed a norm on $C_{c}(G,A)$ by exhibiting a covariant representation. 

Let $\pi:A \to B(\clh)$ be a faithful representation. Set $\widetilde{\clh}:=\clh \otimes \ell^{2}(G)$. Let $\{\epsilon_{t}: t \in G\}$ be the standard orthonormal basis for $\ell^{2}(G)$. For $a \in A$, let $\widetilde{\pi}(a)$ be the bounded operator on $\widetilde{\clh}$ given by the equation \[\widetilde{\pi}(a)(\xi \otimes \epsilon_{t})=\pi(\alpha_{t}^{-1}(a)) \otimes \epsilon_{t}.\] Let $\{\lambda_{s}:s \in G\}$ be the left regular representation of $\ell^{2}(G)$. For $s \in G$, set $\widetilde{\lambda_s}=1 \otimes \lambda_s$. 

\begin{xrcs}
Verify that $(\widetilde{\pi},\widetilde{\lambda})$ is a covariant representation of $(A,G,\alpha)$. 
\end{xrcs}

\begin{ppsn}
The map $C_{c}(G,A) \ni f \to (\widetilde{\pi} \rtimes \widetilde{\lambda})(f) \in B(\widetilde{\clh})$ is injective. 
\end{ppsn}
\textit{Proof.} Suppose $(\widetilde{\pi} \rtimes \widetilde{\lambda})(f)=0$. For $s \in G$, set $a_{s}=f(s)$. Then for every $\xi,\eta \in \clh$ and $r,t \in G$, we have 
$\displaystyle \sum_{s \in G} \langle \widetilde{\pi}(a_s)\widetilde{\lambda_s}(\xi \otimes \epsilon_r)|\eta \otimes \epsilon_t \rangle =0$. This implies that for $\xi,\eta \in \clh$ and $r,t \in G$, 
\[
\sum_{s \in G} \langle \pi(\alpha_{sr}^{-1}(a_s))\xi \otimes \epsilon_{sr}|\eta \otimes \epsilon_t \rangle=0.\]  
Fix $s \in G$. In the previous expansion, substitute $r=s^{-1}$ and $t=e$ to obtain that 
$\langle \pi(a_s)\xi|\eta \rangle=0$ for every $\xi,\eta \in \clh$. But $\pi$ is faithful. This implies that $a_s=0$. Hence $f=0$. This completes the proof. \hfill $\Box$

Keep the foregoing notation. For $f \in C_{c}(G,A)$, define \[||f||_{red}=||(\widetilde{\pi} \rtimes \widetilde{\lambda})(f)||.\] By what we have shown, it follows that $||~||_{red}$ is a $C^{*}$-norm on $C_{c}(G,A)$. By definition, $||f||_{red} \leq ||f||$ for $f \in C_{c}(G,A)$. Hence, the universal norm $||~||$ is a $C^{*}$-norm. 
\begin{dfn}
The completion of $C_{c}(G,A)$ with respect to the universal norm $||~||$ is called the full crossed product and is denoted $A \rtimes_{\alpha} G$.

\end{dfn}
\begin{rmrk}
It is a remarkable fact that $||~||_{red}$ is independent of the chosen faithful representation $\pi$. We will prove this in  the next chapter. The norm $||~||_{red}$ is called the reduced norm on $C_{c}(G,A)$. The completion of $C_{c}(G,A)$ with respect to the reduced norm is called the reduced crossed product and is denoted $A \rtimes_{r,\alpha} G$. 

Clearly, there is a natural surjective homomorphism from $A \rtimes_{\alpha} G \to A \rtimes_{r,\alpha}G$. Unless there is some `amenability' hypothesis, we cannot expect the above map to be an isomorphism. 
\end{rmrk}

\begin{xrcs}
\label{a dense subalgebra of a crossed product}
Let $(A,G,\alpha)$ be a $C^{*}$-dynamical system. Suppose $\cla$ is a dense $^*$-algebra of $A$. Assume that $\alpha_s(\cla) \subset \cla$ $\forall s \in G$. Prove that \[\cla \rtimes_{\alpha} G:=span\{a \otimes \delta_s: a \in \cla, s \in G\}\] is a dense $^*$-subalgebra of $A \rtimes_{\alpha} G$. 
\end{xrcs}
Let us identity one example of a crossed product explicity. Let $G$ be a discrete group and let $G$ acts on the topological space $G$ by left translations. Consider the induced action $\alpha$ of $G$ on $C_{0}(G)$. For $s \in G$, let $\chi_{s}$ be the characteristic function at $s$. Then $\chi_{s} \in C_{c}(G)$ and $\alpha_{s}(\chi_t)=\chi_{st}$. 
\begin{ppsn}
\label{discrete stone-von Neumann}
The crossed product $C_{0}(G) \rtimes_{\alpha} G$ is isomorphic to $\mathcal{K}(\ell^{2}(G))$. 
\end{ppsn}
\textit{Proof.} The algebra $\mathcal{K}(\ell^{2}(G))$ has a universal picture. Thus, it suffices to exhibit appropriate matrix units in $C_{0}(G) \rtimes_{\alpha}G$. Let $\{E_{s,t}:s ,t \in G\}$ be the natural system of matrix units in $\mathcal{K}(\ell^{2}(G))$ which correspond to the standard orthonormal basis $\{\epsilon_{t}: t \in G\}$ of $\ell^{2}(G)$.

For $s,t \in G$, let 
$e_{s,t}:=\chi_{s} \otimes \delta_{st^{-1}}$.  For $q,r,s,t \in G$, calculate as follows to observe that 
\begin{align*}
e_{q,r}*e_{s,t}&=(\chi_{q} \otimes \delta_{qr^{-1}})*(\chi_{s} \otimes \delta_{st^{-1}}) \\
&=\chi_{q}\chi_{qr^{-1}s}\otimes \delta_{qr^{-1}st^{-1}} \\
&=\delta_{q,qr^{-1}s}\chi_{q} \otimes \delta_{qt^{-1}}\\
&=\delta_{r,s}\chi_{q} \otimes \delta_{qt^{-1}}\\
&=\delta_{r,s}e_{q,t}.
\end{align*}
For $s,t \in G$,  observe that 
\[e_{s,t}^*=(\chi_s \otimes \delta_{st^{-1}})^{*}=\alpha_{ts^{-1}}(\chi_s) \otimes ts^{-1}=\chi_{t} \otimes \delta_{ts^{-1}}=e_{t,s}.\]
Thus $\{e_{s,t}:s ,t \in G\}$ forms a system of matrix units. Thus, by the universal property of $\mathcal{K}(\ell^{2}(G))$, there exists a $^*$-homomorphism $\lambda:\mathcal{K}(\ell^{2}(G)) \to C_{0}(G) \rtimes_{\alpha} G$ such that $\lambda(E_{s,t})=e_{s,t}$. By Exercise \ref{a dense subalgebra of a crossed product}, it follows that $\lambda$ is onto. Since $\mathcal{K}(\ell^{2}(G))$ is simple, it follows that $\lambda$ is one one. Hence $\lambda$ is an isomorphism. This completes the proof. \hfill $\Box$

\begin{xrcs}
Consider the Hilbert space $\ell^{2}(G)$. For $f \in C_{0}(G)$, let $M(f)$ be the bounded operator on $\ell^{2}(G)$ defined by the equation
\[
M(f)\xi(s):=f(s)\xi(s)\]
for $\xi \in \ell^{2}(G)$. Show that $M:C_{0}(G) \to B(\ell^{2}(G))$ is a non-degenerate $^*$-representation. Let $\lambda:=\{\lambda_{s}\}_{s \in G}$ be the left regular representation 
of $G$ on $\ell^{2}(G)$. Prove that $(M,\lambda)$ is a covariant pair. 

Use the fact that $C_{0}(G) \rtimes_{\alpha} G$ is simple to show that $M \rtimes \lambda$ implements an isomorphism between $C_{0}(G) \rtimes_{\alpha} G$ and $\mathcal{K}(\ell^{2}(G))$.

\end{xrcs}

Let us end this section by listing out a few examples of universal $C^{*}$-algebras which have played crucial roles in the development of the subject. 
\newline

\textbf{Cuntz algebra $O_n$ :} Let $n \geq 2$. The Cuntz algebra $O_n$ is defined to be the universal unital $C^{*}$-algebra generated by isometries $s_1,s_2,\cdots,s_n$ such that 
\[
\sum_{i=1}^{n}s_is_{i}^{*}=1.\]
Note that $\{s_{i}s_{i}^{*}:i=1,2,\cdots,n\}$ is a family of projections which add up to $1$ which is again a projection. Thus, the projections $\{s_{i}s_{i}^{*}\}_{i=1}^{n}$ form a family of 
mutually orthogonal projections which is equivalent to saying $s_{i}^{*}s_{j}=0$ if $i \neq j$. It is an important fact that the Cuntz algebra $O_n$ is simple. 
\newline

\textbf{The non-commutative torus $A_{\theta}$ :} Let $\theta \in \mathbb{R}$. The non-commutative torus $A_{\theta}$ is defined to be the universal $C^{*}$-algebra generated by 
two unitaries $u,v$ such that \[uv=e^{2\pi i \theta}vu.\] Define $R_{\theta}:\mathbb{T} \to \mathbb{T}$ by $R_{\theta}(z)=e^{-2\pi i \theta}z$. Note that $R_{\theta}$ is a homeomorphism of $\bbt$.
Consequently, this gives rise to an action of the cyclic group $\mathbb{Z}$ on $C(\mathbb{T})$. Show that  the crossed product $C(\mathbb{T}) \rtimes \mathbb{Z}$ is isomorphic to $A_{\theta}$. If $\theta$ is irrational, 
then $A_{\theta}$ is simple. 

The computation of the $K$-theoretic invariants for the two $C^{*}$-algebras listed above were significant breakthroughs in operator $K$-theory. The non-commutative torus still remains one of 
the widely studied example in noncommutative geometry. 
\newline

\textbf{The odd dimensional quantum sphere :} Let $0<q<1$ be given and $\ell \geq 0$. 
The $C^*$-algebra $C(S_q^{2\ell+1})$ of the quantum
sphere $S_q^{2\ell+1}$
is the universal $C^*$-algebra generated by
elements
$z_1, z_2,\ldots, z_{\ell+1}$
satisfying the following relations:
\bean
z_i z_j & =& qz_j z_i,\qquad 1\leq j<i\leq \ell+1,\\
z_i^* z_j & =& q z_j z_i^* ,\qquad 1\leq i\neq j\leq \ell+1,\\
z_i z_i^* - z_i^* z_i +
(1-q^{2})\sum_{k>i} z_k z_k^* &=& 0,\qquad \hspace{2em}1\leq i\leq \ell+1,\\
\sum_{i=1}^{\ell+1} z_i z_i^* &=& 1.
\eean
 Note that for $\ell=0$, the $C^*$-algebra
$C(S_q^{2\ell+1})$ is the algebra of continuous functions
$C(\bbt)$ on the torus and for $\ell=1$, it is denoted $C(SU_q(2))$. The $C^{*}$-algebra $C(SU_q(2))$ is one of the first examples in Woronowicz theory of compact quantum groups and is one of the first examples whose 
representation theory was explicitly worked out.


\section{The Toeplitz algebra and the unilateral shift}

In this section, we discuss  the $C^{*}$-algebra generated by the unilateral shift on $\ell^{2}(\mathbb{N})$. We prove Coburn's theorem which asserts that it is the universal $C^{*}$-algebra generated by a 
single isometry. Coburn's theorem is a fundamental theorem and we will see its importance when we discuss Cuntz' proof of Bott periodicity in $K$-theory. 

\begin{dfn}
Let $\mathcal{T}$ be the universal unital $C^{*}$-algebra generated by $v $ such that $v^{*}v=1$. The $C^{*}$-algebra $\mathcal{T}$ is called the Toeplitz algebra. 
\end{dfn}
By Remark \ref{existence of universal $C^{*}$-algebra}, the $C^{*}$-algebra $\mathcal{T}$ exists. Consider the Hilbert space $\ell^{2}(\mathbb{N})$. Let  $\{\delta_n: n \geq 0\}$ be the standard orthonormal basis 
for $\ell^{2}(\mathbb{N})$. Let $S$ be the bounded operator on $\ell^{2}(\mathbb{N})$ such that $S(\delta_n)=\delta_{n+1}$. The operator $S$ is called \emph{the unilateral shift} on $\ell^{2}(\mathbb{N})$. Clearly $S^{*}S=1$. Thus, by the universal property of $\mathcal{T}$, there exists a unique $^*$-homomorphism $\mathcal{T} \to C^{*}(S)$ which maps $v \to S$. Here, $C^{*}(S)$ denotes the $C^{*}$-algebra generated by $S$. Coburn's theorem asserts that this map is indeed an isomorphism. 

Let us take a closer look at the $C^{*}$-algebra $C^{*}(S)$. For $m,n \geq 0$, let $E_{m,n}=\theta_{\delta_m,\delta_n}$. Set $P:=1-SS^{*}$. Note that $P=E_{0,0}$ and $E_{m,n}=S^{m}PS^{*n}$. Since the linear span of $\{E_{m,n}: m,n \in \mathbb{N}\}$ is dense in $\mathcal{K}(\ell^{2}(\mathbb{N}))$, it follows that $\mathcal{K}(\ell^{2}(\mathbb{N}))$ is contained in $C^{*}(S)$.  Hence, $\clk:=\mathcal{K}(\ell^{2}(\mathbb{N}))$ is an ideal in $C^{*}(S)$. Note that $\dot{S}$, the image of $S$ under the canonical surjection, in the quotient $C^{*}(S)/\mathcal{K}$ is a unitary. Thus, the quotient is generated by a single unitary $\dot{S}$

\begin{lmma}
\label{spectrum}
The spectrum of $\dot{S}$ in $C^{*}(S)/\mathcal{K}$ is $\mathbb{T}$.
\end{lmma}
\textit{Proof.} For $z \in \mathbb{T}$, let $U_{z}$ be the unitary on $\ell^2(\bbn)$ defined $U_{z}(\delta_n)=z^{n}\delta_n$. Note that $U_{z}SU_{z}^{*}=zS$ for every $z \in \mathbb{T}$. Fix $z \in \mathbb{T}$. The map $T \to U_{z}TU_{z}^{*}$ defines an automorphism of $C^{*}(S)$ which leaves the ideal $\mathcal{K}$ invariant. Thus, it descends to an automorphism, let us denote it by $\alpha_{z}$, on the quotient $C^{*}(S)/\clk$. Note that $\alpha_{z}(\dot{S})=z\dot{S}$. 

Denote the spectrum of $\dot{S}$ by $\sigma(\dot{S})$. Fix $z \in \mathbb{T}$. Since $\alpha_z$ is an automorphism, it follows that $\sigma(\alpha_{z}(\dot{S}))=\sigma(\dot{S})$. But $\alpha_{z}(\dot{S})=z\dot{S}$. Hence $\sigma(\alpha_{z}(\dot{S}))=z \sigma(\dot{S})$. This implies that $\sigma(\dot{S})$ is invariant under multiplication by every element of $\mathbb{T}$. Since $\dot{S}$ is a unitary, it follows that $\sigma(\dot{S})$ is contained in $\mathbb{T}$. Hence $\sigma(\dot{S})=\mathbb{T}$. This completes the proof. 
\hfill $\Box$

Denote the function $\mathbb{T} \ni z \to z \in \mathbb{C}$ by $z$ itself. Now continuous functional calculus and Lemma \ref{spectrum} implies that there exists a $^*$-homomorphism $C^{*}(S) \to C(\mathbb{T})$ which maps $S \to z$. Summarising our discussion so far, we have the following exact sequence
\[
0 \longrightarrow \mathcal{K}(\ell^{2}(\mathbb{N})) \longrightarrow C^{*}(S) \longrightarrow C(\mathbb{T}) \longrightarrow 0
\]
 where the map $C^{*}(S) \to C(\mathbb{T})$ sends $S\to z$ and the map $\mathcal{K} \to C^{*}(S)$ is the natural inclusion.  

The first step towards the proof of Coburn's theorem is to derive a similar exact sequence for the Toeplitz algebra $\mathcal{T}$. We imitate what we did for $C^{*}(S)$. Set $p:=1-vv^{*}$. Then $p \neq 0$. For $m,n \geq 0$, set $e_{m,n}=v^{m}pv^{*n}$. Note that $v^{*}p=0$. Hence $v^{*n}p=0$ for every $n \geq 1$. Taking adjoints, we get $pv^{m}=0$ if $m \geq 1$. Note that $v^{*n}v^{m}=v^{m-n}$ if $m\geq n$ and if $m<n$ then $v^{*n}v^{m}=v^{*(n-m)}$.

Let $m_1,n_1,m_2,n_2 \geq 0$ be given. Suppose $n_1>m_2$. Then \[
e_{m_1,n_1}e_{m_2,n_2}=v^{m_1}pv^{*n_1}v^{m_2}pv^{*n_2}=v^{m_1}pv^{*(n_1-m_2)}pv^{*n_2}=0.\]
A similar calculation reveals that if $m_2>n_1$ then $e_{m_1,n_1}e_{m_2,n_2}=0$. Clearly if $m_2=n_1$, then $e_{m_1,n_1}e_{m_2,n_2}=e_{m_1,n_2}$. Hence 
\[
e_{m_1,n_1}e_{m_2,n_2}=\delta_{m_2,n_1}e_{m_1,n_2}.\]
Clearly, $e_{m,n}^{*}=e_{n,m}$ for $m,n \geq 0$. Thus, $\{e_{m,n}:m,n \geq 0\}$ is a system of matrix units. Let $I$ be the closed linear span of $\{e_{m,n}:m,n \geq 0\}$. Note that $I$ is an ideal in $\mathcal{T}$. Indeed $I$ is the ideal generated by $p$. 
By the universal property, there exists a homomorphism $\lambda: \mathcal{T} \to C(\mathbb{T})$ such that $\lambda(v)=z$. Since $C(\mathbb{T})$ is the universal $C^{*}$-algebra generated by a single $v$ such that $v^{*}v=1$ and $1-vv^{*}=0$, it follows that 
the kernel of $\lambda$ is $I$. 

Thus, we obtain a short exact sequence 
\[
0 \longrightarrow I \longrightarrow \mathcal{T} \longrightarrow C(\mathbb{T}) \longrightarrow 0.\]

\begin{thm}[Coburn]
\label{Coburn's theorem}
The natural map $\mathcal{T} \to C^{*}(S)$ which sends $v \to S$ is an isomorphism.
\end{thm}
\textit{Proof.} 
The two short exact sequences that we have derived fit into the following commutative diagram. 
\begin{equation}
\label{commutative}
\def\labelstyle{\scriptstyle}
\xymatrix@C=25pt@R=20pt{
0\ar[r] & I \ar[d]\ar[r]& \mathcal{T} \ar[d] \ar[r]& C(\mathbb{T})\ar[d]\ar[r] &0  \\
0\ar[r] & \mathcal{K} \ar[r]& C^*(S)\ar[r]& C(\mathbb{T})\ar[r] &0 } 
\end{equation}
The maps in the vertical directions are as follows. The map $\mathcal{T} \to C^{*}(S)$ is the map that sends $v$ to $S$ and the map $C(\mathbb{T}) \to C(\mathbb{T})$ is the identity map. Now an application of the five lemma yields the proof. \hfill $\Box$

 Next, we prove \emph{the Wold decomposition} for a single isometry, a result which describes how a generic isometry looks like. Consider the isometry $S$ with multiplicity, i.e. consider a Hilbert space $\clk$ and look at $S \otimes 1$ on $\ell^{2}(\mathbb{N}) \otimes \clk$. Suppose $U$ is a unitary on a different Hilbert space say $\clh_1$. Then
$ \begin{bmatrix}
 S \otimes 1 & 0 \\
 0 & U \\
 \end{bmatrix}$ is an isometry on $(\ell^{2}(\mathbb{N})\otimes \clk)\oplus \clh_1$. The Wold decomposition asserts that every isometry, up to a unitary equivalence,  is of this form. 

\begin{thm}[Wold decomposition]
\label{Wold}
Let $\clh$ be a separable Hilbert space and $V$ be an isometry on $\clh$. Then there exists Hilbert spaces $\clk$ and $\clh_1$, a unitary $U$ on $\clh_1$ and a unitary $W:\clh \to (\ell^{2}(\mathbb{N})\otimes \clk) \oplus \clh_1$ such that 
\[
WVW^{*}=\begin{bmatrix}
S \otimes 1 & 0 \\
0 & U \\
\end{bmatrix}. \]
\end{thm}

\begin{lmma}
Let $A$ be a $C^{*}$-algebra and let $I \subset A$ be an ideal. Suppose  $\pi$ is   a non-degenerate representation of $I$ on a Hilbert space $\clh$. Then, there exists a unique representation $\widetilde{\pi}:A \to B(\clh)$ such that $\widetilde{\pi}(x)=\pi(x)$ for $x \in I$. 
\end{lmma}
\textit{Proof.} Any non-degenerate representation can be written as a direct sum of cyclic representations. A moment's thought reveals that it suffices to prove the lemma when $\pi$ is cyclic. Thus, let $\pi$ be cyclic and $\xi \in \clh$ be a cyclic vector, i.e. $\{\pi(x)\xi:x \in I\}$ is dense in $\clh$. Fix $a \in A$. Calculate as follows to observe that  for $x \in I$,
\begin{align*}
\langle \pi(ax)\xi|\pi(ax)\xi \rangle &=\langle \pi(x^{*}a^{*}ax)\xi|\xi \rangle \\
& \leq ||a||^{2} \langle \pi(x^{*}x)\xi|\xi \rangle ~(\textrm{since $x^{*}a^{*}ax \leq ||a||^{2}x^{*}x$}) \\
& \leq ||a||^{2} \langle \pi(x)\xi|\pi(x)\xi \rangle.
\end{align*}
The above calculation implies that there exists a unique bounded operator, denoted $\widetilde{\pi}(a)$, such that $\widetilde{\pi}(a)\pi(x)\xi=\pi(ax)\xi$ for every $x \in I$. If $a \in I$, then \[\widetilde{\pi}(a)\pi(x)\xi=\pi(ax)\xi=\pi(a)\pi(x)\xi\] for every $x \in I$.
Since $\{\pi(x)\xi:x \in I\}$ is dense in $\clh$, it follows that $\widetilde{\pi}(a)=\pi(a)$. Evaluating on the dense set $\{\pi(x)\xi:x \in I\}$, it is routine to see that $\widetilde{\pi}$ is a $*$-representation. 

Uniqueness follows from the fact that  $\{\pi(x)\eta: x \in I, \eta \in \clh\}$ is total in $\clh$. This completes the proof. \hfill $\Box$

\textit{Proof of Theorem \ref{Wold}}. Let $V$ be an isometry on $\clh$. By Coburn's theorem, there exists a representation $\pi:C^{*}(S) \to B(\clh)$ such that $\pi(S)=V$. Let $I:=\mathcal{K}(\ell^{2}(\mathbb{N}))$. 
Set $\clh_0=\pi(I)\clh$ and $\clh_1=\clh_0^{\perp}$. Since $I$ is an ideal of $C^{*}(S)$, it follows that $\clh_0$ and $\clh_1$ are invariant under $\pi$. Note that $\pi(x)$ vanishes on $\clh_1$ if $x \in I$. Thus, $\pi(V)|_{\clh_1}$ is a unitary. Set $U:=\pi(V)|_{\clh_1}$.

Restrict the representation $\pi$ to $I$ on $\clh_0$. Then $\pi|_{I}$ is non-degenerate. Hence there exists a Hilbert space $\clk$ and a unitary $W_0: \clh_0 \to \ell^{2}(\mathbb{N}) \otimes \clk$ such that $W_0\pi(x)W_{0}^{*}=x \otimes 1$ for $x \in \mathcal{K}(\ell^{2}(\mathbb{N}))$.  The representation $W_{0}\pi(.)W_{0}^{*}$ and $x \to x \otimes 1$ are both extensions to $A$ of the representation $W_0\pi(.)W_{0}$ defined on $I$. Hence $W_{0}\pi(x)W_{0}^{*}=x\otimes 1$ for every $x \in A$. Define \[W: \clh_0 \otimes \clh_1 \to (\ell^{2}(\mathbb{N}) \otimes \clk) \oplus \clh_1\] by $W=W_0 \oplus Id$.
Then $WVW^{*}=W\pi(S)W^{*}=\begin{bmatrix}
                            S \otimes 1 & 0 \\
                             0 & U 
                             \end{bmatrix}$. This completes the proof. \hfill $\Box$
 
 \begin{rmrk}                            
Here, we have derived Wold decomposition from Coburn's theorem. We could for instance first prove Wold decomposition and derive Coburn's theorem as a corolllary. The derivation undertaken here is more operator algebraic in nature.     
In the next chapter, we prove Wold decomposition, a result due to Cooper (\cite{Cooper}), in the continuous case.                          
\end{rmrk}

\section{Measure theoretic preliminaries}
We end this chapter by discussing group $C^*$-algebras in the locally compact case. We first collect the measure theoretic preliminaries needed in this section. 
Let $X$ be a second countable, locally compact,  Hausdorff topological space. Denote the Borel $\sigma$-algebra of $X$, i.e. the $\sigma$-algebra generated
by open subsets of $X$ by $\mathcal{B}_{X}$. On a locally compact space, we always consider this Borel $\sigma$-algebra. Let $\mu$ be a measure on $(X,\mathcal{B}_X)$. 
We say that $\mu$ is finite on compact sets if $\mu(K)<\infty$ for every compact $K \subset X$. Measures which are finite on compact sets are also called \emph{Radon measures}. 
Let $\mu$ be a Radon measure on $(X,\mathcal{B}_X)$. Denote the algebra of continuous complex valued functions on $X$ with compact support by $C_{c}(X)$.  We have the following.
\begin{enumerate}
\item[(1)] The measure $\mu$ is regular, i.e. for every $E \in \mathcal{B}_X$, 
\begin{align*}
\mu(E)&=\sup\{\mu(K): K \subset E, \textrm{$K$ is compact}\} \\
           &=\inf \{\mu(V): E \subset V, \textrm{$V$ is open}\}.
\end{align*}
\item[(2)] The fact that $\mu$ is finite on compact sets implies that for $f \in C_{c}(X)$, $f$ is integrable with respect to $\mu$. Moreover, $C_{c}(X) \subset L^{p}(X)$ for every $1 \leq p \leq \infty$.
\item[(3)] The fact that $\mu$ is regular has the consequence that $C_{c}(X)$ is dense in $L^{p}(X)$ for every $1 \leq p< \infty$. 
\end{enumerate}
A linear functional $\phi:C_{c}(X) \to \mathbb{C}$ is said to be positive if $f \geq 0$, then $\phi(f) \geq 0$. Denote the set of positive linear functionals on $C_{c}(X)$ by $C_{c}(X)_{+}^{*}$. Let $\mu$ be a Radon measure
on $(X,\mathcal{B}_X)$. Define $\phi_{\mu}:C_{c}(X) \to \mathbb{C}$ by
\[
\phi_{\mu}(f)=\int f(x)d\mu(x)
\]
for $f \in C_{c}(X)$. Then, $\phi_{\mu}$ is a positive linear functional. Denote the set of Radon measures on $X$ by $\mathcal{M}(X)$.

\begin{thm}[Riesz representation theorem]
The map \[\mathcal{M}(X) \ni \mu \to \phi_{\mu} \in C_{c}(X)_{+}^{*}\] is a bijection. 
\end{thm}

\textbf{Push forward measure:} Let $(X,\mathcal{B}_X)$ and $(Y,\mathcal{B}_{Y})$ be measurable spaces. Suppose $\mu$ is a measure on $(X,\mathcal{B}_X)$ and $T:X \to Y$ is a measurable map.  For $E \in \mathcal{B}_{Y}$, define
\[
T_{*}\mu(E):=\mu(T^{-1}(E)).\]
Then, $T_{*}\mu$ is a measure on $(Y,\mathcal{B}_{Y})$. The measure $T_{*}\mu$ is called \emph{the push forward of $\mu$ by $T$}.
Keep the foregoing notation. Suppose $f:Y \to [0,\infty]$ is measurable. It is not difficult to prove that  
\[
\int f d(T_{*}\mu)=\int (f \circ T) d\mu.\]
Thus, a measurable function $f:Y \to \mathbb{C}$ is integrable if and only if $f \circ T$ is integrable, and in that case
\[
\int fd(T_{*}\mu)=\int (f \circ T) d\mu.\]

\textbf{Inductive limit topology:} Let $X$ be a locally compact Hausdorff second countable topological space. Suppose $(f_n)$ is a sequence in $C_{c}(X)$ and $f \in C_{c}(X)$. We say that $f_n \to f$ in \emph{the inductive limit topology} if 
there exists a compact set $K \subset X$ such that 
\begin{enumerate}
\item[(1)] for every $n \geq 1$, $supp(f_n) \subset K$, and
\item[(2)] the sequence $f_n \to f$ uniformly on $X$. 
\end{enumerate}
Suppose $V$ is a topological vector space and $T:C_{c}(X) \to V$ is a linear map. We say that $T$ is continuous with respect to the inductive limit topology, if whenever $(f_n)$ is a sequence in $C_{c}(X)$ which converges to $f$ in the inductive limit topology, then $T(f_n) \to T(f)$ in $V$. 

\begin{xrcs}
Let $X$ be a second countable, locally compact, Hausdorff topological space and let $\mu$ be a Radon measure on $X$. Prove that the functional \[
C_{c}(X) \ni f \to \int f(x)d\mu(x) \in \mathbb{C}
\]
is continuous with respect to the inductive limit topology. 
\end{xrcs}

\begin{xrcs}
Let $X$ be a second countable, locally compact, Hausdorff topological space and let $\mu$ be a Radon measure on $X$. Prove that the ``natural" map
$C_{c}(X) \to L^{p}(X)$ is continuous with respect to the inductive limit topology for every $1 \leq p < \infty$.

\end{xrcs}

\textbf{Haar measure: }
Let us now discuss the basics of Haar measure on a locally compact group. The letter $G$ stands for a locally compact, second countable, Hausdorff topological group. The Borel $\sigma$-algebra of $G$ is denoted $\mathcal{B}_{G}$. 
For $s \in G$, let $\sigma_{s}:G \to G$ be defined by $\sigma_{s}(t)=st$ and let $\rho_{s}:G \to G$ be defined by $\rho_s(t)=ts$. For $s \in G$, let $L_{s}, R_{s}:C_{c}(G) \to C_{c}(G)$ be defined by
\begin{align*}
L_{s}f(t)&:=f(s^{-1}t)\\
R_{s}f(t)&:=f(ts).
\end{align*}

A measure $\mu$ on $(G,\mathcal{B}_{G})$ is said to be left invariant if $(\sigma_{s})_{*}(\mu)=\mu$ for every $s \in G$. By Riesz representation theorem, a Radon measure $\mu$ is left invariant if and only if 
\[
\int f(s^{-1}t)d\mu(t)=\int f(t)d\mu(t)
\]
for every $f \in C_{c}(G)$ and $s \in G$. We accept the following theorem without proof. We refer the reader to \cite{Folland} for a proof. 

\begin{thm}[Haar measure]
Let $G$ be a second countable, locally compact, Hausdorff topological group. Then, there exists a non-zero Radon measure $\mu$ which is left invariant. Moreover, if $\mu$ and $\nu$ are two non-zero left invariant
Radon measures, then there exists $c>0$ such that $\nu=c\mu$.

\end{thm}

\begin{dfn}
A left invariant non-zero Radon measure on $G$ is called a Haar measure on $G$. 
\end{dfn}
Note that any two Haar measures differ by a scalar. Of course, we could talk about a right Haar measure. 
Fix a left Haar measure $\mu$ on $G$. 
\begin{ppsn}
If $U$ is a non-empty open subset of $G$, then $\mu(U)>0$. 
\end{ppsn}
\textit{Proof.} Let $U$ be a non-empty open subset of $G$. Suppose that $\mu(U)=0$. By left invariance of $\mu$, it follows that $\mu(xU)=0$ for every $x \in G$.
Note that $\{xU: x \in G\}$ is an open cover of $G$. Since $G$ is second countable, there exists a countable subset $\{x_1,x_2,\cdots,\} \subset G$ such that 
$G=\bigcup_{i=1}^{\infty}x_iU$. Then, $\mu(G)=0$ which is a contradiction. Hence the proof. \hfill $\Box$

\begin{xrcs}
 Suppose $f \in C_{c}(G)$ is non-negative and $\int f d\mu=0$. Prove that $f$ is identically zero.
\end{xrcs}

For $s \in G$, let $\mu_s=(\rho_{s^{-1}})_{*}\mu$. Then, $\mu_{s}$ is a left invariant Radon measure on $G$. Thus, there exists a positive scalar $\Delta(s)$ such that 
$\mu_{s}=\Delta(s)\mu$. A moment's thought reveals that the map $G \ni s \to \Delta(s) \in (0,\infty)$ does not depend on the chosen Haar measure $\mu$, and 
it is in fact a homomorphism. The function $\Delta$ is called \emph{the modular function} of the group $G$. 

\begin{xrcs}
Let $\mu$ be a Haar measure on $G$. Prove that for $f \in G$, 
\[
\int f(ts^{-1})d\mu(t)=\Delta(s)\int f(t)d\mu(t)\]
for $f \in C_{c}(G)$ and $s \in G$. 
\end{xrcs}

\begin{ppsn}
\label{continuity of left action}
Fix $f \in C_{c}(G)$. The map $G \ni s \to L_{s}(f) \in C_{c}(G)$ and the map $G \ni s \to R_{s}(f) \in C_{c}(G)$ are continuous, when $C_{c}(G)$ is given the inductive limit topology, i.e. 
if $s_n \to s$, then $L_{s_n}(f) \to L_{s}(f)$ and $R_{s_n}(f) \to R_{s}(f)$ in the inductive limit topology. 
\end{ppsn}
\textit{Proof.}  Let $s_n$ be a sequence in $G$ such that $s_n \to s$. Choose a compact set $L$ such that $L$ contains $\{s_n: n \geq 1\}\cup \{s\}$. Denote the support of $f$ by $K$. 
Note that for $t \in G$, $supp(L_t(f)) \subset tK$. Hence, $supp(L_{s_n}(f))$ and $supp(L_{s}(f))$ are contained in $LK$ and $LK$ is compact.

 Suppose $L_{s_n}(f)$ does not converge to $L_{s}(f)$ uniformly.
Then there exists $\epsilon>0$ and a subsequence $(t_{n_k}) \in G$ such that 
\[|L_{s_{n_k}}(f)(t_{n_k})-L_{s}(f)(t_{n_k})| \geq \epsilon.\]
The above inequality implies that $t_{n_k} \in LK$. But $LK$ is compact. By passing to a subsequence, if necessary, we can assume that $t_{n_k}$ converges, say, to $t$. Then, \[L_{s_{n_k}}(f)(t_{n_k})=f(s_{n_k}^{-1}t_{n_k}) \to f(s^{-1}t).\] 
Similarly, $L_{s}(f)(t_{n_k}) \to f(s^{-1}t)$ which contradicts the fact that for every $k$
\[|L_{s_{n_k}}(f)(t_{n_k})-L_{s}(f)(t_{n_k})| \geq \epsilon.\]
Hence the proof. \hfill $\Box$

\begin{xrcs}
Suppose $G$ acts on a locally compact space $X$ on the left. For $s \in G$, let $L_{s}:C_{c}(X) \to C_{c}(X)$ be defined by $L_{s}(f)=f(s^{-1}x)$. Show that for $f \in C_{c}(X)$, the map $G \ni s \to L_{s}(f) \in C_{c}(X)$ is continuous, where $C_{c}(X)$ is given the inductive limit topology. 

Show that for $f \in C_{0}(X)$, the map $G \ni s \to L_{s}(f) \in C_{0}(X)$ is continuous, where $C_{0}(X)$ is given the norm topology. Here, the map $L_{s}:C_{0}(X) \to C_{0}(X)$ is defined in the same fashion. 
\end{xrcs}
\textit{Hint: } The inclusion $C_{c}(X) \to C_{0}(X)$ is continuous and has dense range. 
\begin{lmma}
The modular function $\Delta$ is continuous. 
\end{lmma}
\textit{Proof.} Choose $f \in C_{c}(G)$ such that $\int f(t)d\mu(t)=1$. Note that for $ s \in G$, \[\int R_{s^{-1}}f(t)d\mu(t)=\int f(ts^{-1})d\mu(s)=\Delta(s)\int f(t)d\mu(t)=\Delta(s).\] 
The above equality together with Proposition \ref{continuity of left action} imply that $\Delta$ is continuous. \hfill $\Box$

\begin{rmrk}
A locally compact group $G$ is called unimodular if $\Delta=1$. Abelian groups are clearly unimodular. Compact groups are unimodular. For, if $G$ is compact, then the image $\Delta(G)$ is a compact subgroup of $(0,\infty)$ and
the only compact subgroup of $(0,\infty)$ is $\{1\}$. 
\end{rmrk}

For $E \in \mathcal{B}_{G}$, define $\overline{\mu}(E)=\mu(E^{-1})$. It is clear that $\overline{\mu}$ is a right invariant Haar measure. Note that $\overline{\mu}$ is the push forward measure of $\mu$ under the map $s \to s^{-1}$. 
The proof of the next proposition is taken from \cite{Faraut_Lie}.
\begin{ppsn}
The measure $\overline{\mu}$ and $\mu$ are absolutely continuous with respect to each other. Moreover, $\frac{d\overline{\mu}}{d\mu}(s)=\Delta(s)^{-1}$. 

Equivalently, for $f \in C_{c}(G)$, \[\int f(s^{-1})d\mu(s)=\int f(s)\Delta(s)^{-1}d\mu(s).\]
\end{ppsn}
\textit{Proof.} Let $I:C_{c}(G) \to \mathbb{C}$ be the positive linear functional defined by the equation
\[
I(f)=\int f(s^{-1})\Delta(s^{-1})d\mu(s).\]
Let $f \in C_{c}(G)$ and $t \in G$ be given. Define $g(s)=f(s^{-1})\Delta(s^{-1})$. Calculate as follows to observe that 
\begin{align*}
I(L_{t}(f))&=\int L_{t}(f)(s^{-1})\Delta(s^{-1})d\mu(s) \\
&= \int f(t^{-1}s^{-1})\Delta(s^{-1})d\mu(s) \\
&=\Delta(t) \int g(st)d\mu(s)\\
&=\int g(s)d\mu(s) \\
&=I(f).
\end{align*}
Hence the measure $\nu$ associated to the linear functional $I$, via Riesz representation theorem, is left invariant. Hence there exists $C>0$ such that for $f \in C_{c}(G)$,
\[
\int f(s^{-1})\Delta(s)^{-1}d\mu(s)=C \int f(s)d\mu(s).\]
Let $f \in C_{c}(G)$ be given. Apply the above equality to the function $s \to f(s^{-1})\Delta(s^{-1})$ to see that 
\[
\int f(s)d\mu(s)=C \int f(s^{-1})\Delta(s^{-1})d\mu(s)=C^{2}\int f(s)d\mu(s).\]
Hence $C=1$. This completes the proof. \hfill $\Box$

Let us end this section  by describing the Haar measure for a few examples of topological groups. 
\begin{xmpl}
 The Lebesgue measure on $\mathbb{R}^{d}$ is a Haar measure. 
\end{xmpl}

\begin{xmpl}
Let $\mathbb{T}$ be the unit circle. Note that $\mathbb{T}$ is a compact group with respect to multiplication. We can identify $C(\mathbb{T})$ with the space of periodic functions on $\bbr$, i.e. \[
C(\mathbb{T})=\{f:\mathbb{R} \to \mathbb{C}: ~~f~\text{is continuous and ~}f(x+1)=f(x),~ \forall x \in \mathbb{R}\}.\]
Define $\phi:C(\mathbb{T}) \to \mathbb{C}$ by the formula
\[
\phi(f):=\int_{0}^{1}f(x)dx.\]
Then, $\phi$ is a positive linear functional on $C(\mathbb{T})$. Thus, there exists a measure $\mu$ on $\mathbb{T}$ such that $\int f d\mu=\phi(f)$ for every $f \in C(\mathbb{T})$. Show that $\mu$ is a Haar measure on $\bbt$. 
\end{xmpl}

\begin{xmpl}
Let $G$ be a countable discrete group. Then, the counting measure on $G$ is a Haar measure on $G$.

\end{xmpl}

\begin{xmpl}
The \textbf{ax+b-group:}
Let \[G:=\Big\{\begin{bmatrix}
                  a & b \\
                  0 & 1 
                 \end{bmatrix}: a \neq 0, b \in \mathbb{R}\Big\}.\] Show that $G$ is a closed subgroup of $GL_{2}(\mathbb{R})$. As a set $G=\mathbb{R}\backslash \{0\} \times \mathbb{R}$. Consider the two measures $\frac{da db}{|a|}$ and $\frac{da db}{a^{2}}$ on $G$. One of them is right invariant and the other is left invariant. Determine which one is right invariant and which one is left invariant. 
                    Compute the modular function and show that the group $G$, also called the \textbf{ax+b}-group, is not unimodular. 
\end{xmpl}

\section{Group $C^{*}$-algebras}
Let $G$ be an arbitrary, locally compact, second countable, Hausdorff topological group fixed for the rest of this section. Fix a Haar measure $\mu$. We write $\int f(s)d\mu(s)$ simply as $\int f(s)ds$. Recall the following formulas. For $f \in C_{c}(G)$ and $t \in G$, 
\begin{align*}
\int f(t^{-1}s)ds&=\int f(s)ds \\
\int f(st)ds&=\Delta(t)^{-1}\int f(s)ds \\
\int f(s^{-1})ds&=\int f(s)\Delta(s^{-1})ds.
\end{align*}

For $f,g \in C_{c}(G)$, let $f*g:G \to \mathbb{C}$ be defined by 
\[
f*g(s)=\int f(st)g(t^{-1})dt=\int f(t)g(t^{-1}s)ds.\]
The function $f*g$ is called \emph{the convolution of $f$ and $g$}\footnote{One could equally convolve $L^{1}$ functions. But most of the time, it suffices to work with the dense subspace $C_{c}(G)$}.  Define an involution $^*$ on $C_{c}(G)$ as follows.  
For $f \in C_{c}(G)$, let $f^{*} \in C_{c}(G)$ be defined by
$
f^{*}(s)=\Delta(s)^{-1}\overline{f(s^{-1})}.$

\begin{xrcs}
Show that for $f,g \in C_{c}(G)$, $f*g \in C_{c}(G)$. If $K$ denotes the support of $f$ and $L$ denotes the support of $g$, prove that the support of $f*g$ is contained in $KL$. 
\end{xrcs}

\begin{ppsn}
The space $C_c(G)$ with convolution as multiplication and $^*$ as involution is a $^*$-algebra. 
\end{ppsn}
\textit{Proof.} The proof is really a straightforward application of Fubini's theorem and the left invariance of the Haar measure. For the reader's benefit, let us verify that the convolution is associative and the involution $^*$ is anti-multiplicative. 
Let $f,g,h \in C_{c}(G)$ be given. For $s \in G$, calculate as follows to observe that 
\begin{align*}
((f*g)*h)(s)&= \int (f*g)(st)h(t^{-1})dt \\
 &= \int \Big(\int f(str)g(r^{-1})dr\Big)h(t^{-1})dt \\
 &= \int \Big( \int f(sr)g(r^{-1}t)dr\Big)h(t^{-1})dt ~(\textrm{ left invariance of the Haar measure})\\
 &=\int f(sr)\Big(g(r^{-1}t)h(t^{-1})dt \Big)dr ~(\textrm{Fubini's theorem}) \\
 &= \int f(sr)(g*h)(r^{-1})dr\\
 &=(f*(g*h))(s).
\end{align*}
This proves that the convolution is associative. Let $f,g \in C_{c}(G)$ be given. For $f,g \in C_{c}(G)$ and $s \in G$, calculate as follows to observe that 
\begin{align*}
(f*g)^{*}(s)&=\Delta(s)^{-1}\overline{(f*g)(s^{-1})} \\
&= \Delta(s)^{-1}\int \overline{f(s^{-1}t)}\overline{g(t^{-1})}dt\\
&=\int g^{*}(t)f^{*}(t^{-1}s)dt\\
&=(g^{*}*f^{*})(s).
\end{align*}
This completes the proof. \hfill $\Box$

If $G$ is discrete, then the characteristic functions $\delta_{s} \in C_{c}(G)$ and $\{\delta_{s}:s \in G\}$ spans $C_{c}(G)$. Moreover in the discrete case, the multiplication and the involution of basis elements are as follows.
\begin{align*}
\delta_s*\delta_t&=\delta_{st} \\
\delta_{s}^{*}&=\delta_{s^{-1}}.
\end{align*}
If $G$ is not discrete, then the characteristic functions are no longer  elements of $C_{c}(G)$. The trick to overcome this  is to use  \emph{approximate identities} instead. 

Let $\{U_n\}_{n=1}^{\infty}$ be a decreasing sequence of  open sets containing the identity element $e$. We assume that $U_n$ is symmetric around $e$, i.e. $U_n^{-1}=U_n$. Suppose that $\{U_n\}_{n=1}^{\infty}$ is a basis at $e$, i.e. given an 
open set $U$ containing $e$, there exists $N$ such that $U_{N} \subset U$. Note that such a sequence of open sets can always be constructed. For, $G$ is metrisable and we can
let $V_n$ be the open ball (with respect to a metric inducing the topology of $G$) of radius $\frac{1}{n}$ centered at $e$. Then, set $U_n=V_n \cap V_n^{-1}$. 
For each $n$, choose $\phi_n \in C_{c}(G)$ such that $ supp(\phi_n) \subset U_n$, $\phi_n^{*}=\phi_n$, $\phi_n \geq 0$ and $\int \phi_n(s)ds=1$. Such a sequence $\{\phi_n\}_{n=1}^{\infty}$ is called \emph{an approximate identity}
of $C_{c}(G)$. The justification for the name \emph{approximate identity} is due to the following proposition. 

\begin{ppsn}
\label{approximate}
Keep the foregoing notation. For $f \in C_{c}(G)$, the sequence $\{\phi_n*f\}_{n=1}^{\infty}$ and the sequence $\{f*\phi_n\}_{n=1}^{\infty}$ converge to $f$ in the inductive limit topology. 
\end{ppsn}
\textit{Proof.} Since $\phi_n$ is self-adjoint, it suffices to prove that $\phi_n*f \to f$ in the inductive limit topology. Let $K$ be a compact neighbourhood at $e$. For large $n$, $U_n \subset K$ and consequently $supp(\phi_n*f)\subset supp(\phi_n)supp(f) \subset Ksupp(f)$. Thus $\{\phi_n *f \}_{n=1}^{\infty}$ is supported inside a common compact set. 

Let $\epsilon>0$ be given. 
 Since the map $G\ni t \to L_{t}f \in C_{c}(G)$ is continuous, when $C_{c}(G)$ is given the inductive limit topology,  there exists $N$ large such that for $t \in U_{N}$, $||L_{t}f-f||_{\infty} \leq \epsilon$. For $s \in G$ and $n \geq N$, calculate as follows to observe that 
 \begin{align*}
 |\phi_n*f(s)-f(s)|&=\Big|\int \phi_{n}(t)f(t^{-1}s)dt - f(s)\Big| \\
 &= \Big|\int_{U_n} \phi_{n}(t)f(t^{-1}s)dt - \int \phi_n(t)f(s)dt \Big|\\
 &=\Big|\int_{U_n}\phi_{n}(t)(f(t^{-1}s)-f(s))dt\Big|\\
 & \leq \int_{U_n}\phi_{n}(t)||L_{t}f-f||_{\infty}dt\\
 & \leq \epsilon.
 \end{align*}
 Hence, the sequence $\{\phi_n*f\}_{n=1}^{\infty}$ converges to $f$ in the inductive limit topology. \hfill $\Box$
 
 \begin{xrcs}
 Suppose $G$ is not discrete. Prove that $C_{c}(G)$ has no multiplicative identity.   \end{xrcs}
 \textit{Hint:}Use the fact that $C_{c}(G)$ has an approximate identity.

 For $f \in C_{c}(G)$, define \[||f||_{1}:=\int |f(s)|ds.\] It is clear that $||f^{*}||_{1}=||f||_{1}$ for $f \in C_{c}(G)$. Let $f,g \in C_{c}(G)$ be given. 
 Calculate as follows to observe that 
 \begin{align*}
 \int |f*g(s)|ds &\leq \int \int |f(t)||g(t^{-1}s)|dtds \\
 &= \int |f(t)|\big(\int |g(t^{-1}s)|ds\big)dt \\
 & =\int |f(t)|||g||_{1}dt  (\textrm{~Haar measure is left invariant})\\
 &=||f||_{1}||g||_{1} .
  \end{align*}
 Hence $||f*g|| \leq ||f||_{1}||g||_{1}$ for $f,g \in C_{c}(G)$. In other words, $(C_{c}(G), ||~||_{1})$ is a normed $^*$-algebra. 
 
 \begin{dfn}
 Let $\cla$ be a $^*$-algebra and $||~||$ be a norm on $\cla$. We say that the pair $(\cla,||~||)$ is a normed $^*$-algebra if for $a,b \in \cla$, 
 \begin{align*}
 ||ab||& \leq ||a||||b|| \\
 ||a^{*}||&=||a||.
 \end{align*}
 
 \end{dfn}
 Let $(\cla,||~||_1)$ be a normed $^*$-algebra. The enveloping $C^{*}$-algebra of $\cla$ is defined in the same fashion as in Section 2, the only difference here is that to define the  universal $C^{*}$-seminorm, we consider only representations
 which are bounded w.r.t. $||~||_1$. Suppose $\pi:\cla \to B(\clh)$ be a representation. We say that $\pi$ is  bounded w.r.t. $||~||_1$ if for every $x \in \cla$, $||\pi(x)|| \leq ||x||_1$.
 
 Define a $C^{*}$-seminorm $||~||$ on $\cla$ as 
 \[
 ||x||:=\sup\{||\pi(x)||: \pi \textrm{~is a bounded $*$-representation}\}\]
 for $x \in \cla$. Suppose that $||x||<\infty$ for every $x \in \cla$. Let $I:=\{x \in \cla: ||x||=0\}$. Then $||~||$ descends to a genuine $C^{*}$-norm on $\cla/I$ and the completion of $\cla/I$ is called the enveloping $C^{*}$-algebra of $\cla$ and is denoted $C^{*}(\cla)$. 
 Note that $*$-representations of $C^{*}(\cla)$ are in one-one correspondence with bounded representations of $\cla$. 
 
 \begin{dfn}
 The full group $C^{*}$-algebra, denoted $C^{*}(G)$,  is defined as the enveloping $C^{*}$-algebra of $C_{c}(G)$. 
  \end{dfn}
 Of course, we need to show that $C^{*}(G)$ exists and is non-zero. This requires us to study bounded $^*$-representations of $C_{c}(G)$ in more detail. Just like in the discrete setting, first we show that non-degenerate bounded $^*$-representations of $C_{c}(G)$ are in
 $1$-$1$ correspondence with strongly continuous unitary representations of $G$.
 
 \begin{dfn}
 \label{strongly continuous}
 Let $\clh$ be a Hilbert space and let $U:G \to B(\clh)$ be a map. We denote the image of an element $s$ under $U$ by $U_{s}$. We say that $U$ is a strongly continuous unitary representation of $G$ on $\clh$
 if 
 \begin{enumerate}
 \item[(1)] for $s,t \in G$, $U_sU_t=U_{st}$, 
 \item[(2)] for $s \in G$, $U_{s}$ is a unitary, and
 \item[(3)] for $\xi \in \clh$, the map $G \ni s \to U_{s}\xi \in \clh$ is continuous, where $\clh$ is given the norm topology. 
  \end{enumerate}
   \end{dfn}
   Since we will not consider unitary representations that are not strongly continuous, we drop the modifying term "strongly continuous".  Note that to check $(3)$, it suffices to check the continuity requirement for vectors $\xi$ in a total set. 
 
 \begin{xrcs}
 Show that Condtion $(3)$ in Defn. \ref{strongly continuous} can be replaced by the following condition.
 \begin{enumerate}
 
 \item[$(3)^{'}$] For $\xi,\eta \in \clh$, the map $G \ni s \to \langle U_s\xi|\eta \rangle \in \mathbb{C}$ is continuous. 
 \end{enumerate}
 Condition $(3)^{'}$ is equivalent to saying that $G \ni s \to U_s \in B(\clh)$ is continuous, when $B(\clh)$ is given the weak topology. 
 \end{xrcs}
 \begin{xmpl}
 Consider the Hilbert space $L^{2}(G,\mu)$ where $\mu$ is the Haar measure on $G$. For $s \in G$, let $\lambda_s$ be the unitary on $L^{2}(G)$ defined by the equation
 \[
 \lambda_s(f)(t)=f(s^{-1}t).\]
 Then, $\lambda:=\{\lambda_s\}_{s \in G}$ is a unitary representation of $G$ and is called \emph{the left regular representation} of $G$. The only non-trivial thing to verify is the continuity condition. 
 
 Fix $f \in C_{c}(G)$. Note that the map $G \ni s \to \lambda_s(f) \in L^{2}(G)$ is the composite of the map $G \ni s \to \lambda_s(f) \in C_{c}(G)$ and the inclusion $C_{c}(G) \to L^{2}(G)$. 
 But both these maps are continuous. Hence $G \ni s \to \lambda_s(f) \in L^{2}(G)$ is continuous. 
 
 For $s \in G$, let $\rho_s$ be the unitary on $L^{2}(G)$ defined by the equation
 \[
 \rho_s(f)(t)=f(ts).\]
 Then $\rho:=\{\rho_s\}_{s \in G}$ is a unitary representation of $G$ and is called the right regular representation of $G$. Note that $\lambda(G)$ and $\rho(G)$ commute with each other. It is a  fact that the commutant of $\lambda(G)$ is the von Neumann algebra generated by $\rho(G)$. Similarly, $\rho(G)^{'}$ is the von Neumann algebra generated by $\lambda(G)$.
  \end{xmpl}
To explain the bijective correspondence between the set of unitary representations of $G$ and the set of non-degenerate bounded $^*$-representations of $C_{c}(G)$, we need to recall how to integrate operator valued functions. 

\textbf{Operator valued integration:} Let $(X,\mathcal{B})$ be a measurable space and $\clh$ be a separable Hilbert space. A function $f:X \to B(\clh)$ is said to be \emph{weakly measurable} if for $\xi,\eta \in \clh$, the map \[
X \ni x \to \langle f(x)\xi|\eta \rangle \in \mathbb{C}\]
is measurable. 
 
\begin{lmma}
Keep the foregoing notation. Let $f:X \to B(\clh)$ be weakly measurable. Then, the map 
$
X \ni x \to ||f(x)|| \in \mathbb{C}$
is measurable.
\end{lmma}  
\textit{Proof.} Let $D$ be a countable dense subset of the unit ball of $\clh$. Note that for $x \in X$,
\[
||f(x)||=\sup_{\xi,\eta \in D}|\langle f(x)\xi|\eta \rangle|.\] 
 The proof is clear now. \hfill $\Box$
 
 Let $\mu$ be a measure on $(X,\mathcal{B})$. A weakly measurable map $f:X \to B(\clh)$ is said to be integrable w.r.t. $\mu$ if $\int ||f(x)||d\mu(x)<\infty$. 
 \begin{ppsn}
 Let $\mu$ be a measure on $(X,\mathcal{B})$ and $f:X \to B(\clh)$ be integrable w.r.t. $\mu$. Then, there exists a unique bounded linear operator on $\clh$, denoted $\int f(x)d\mu(x)$, such that for $\xi,\eta \in \clh$,
 \[
 \Big\langle \big(\int f(x)d\mu(x)\big)\xi\big|\eta \Big\rangle=\int \langle f(x)\xi|\eta \rangle d\mu(x)\]
 for $\xi,\eta \in \clh$. Also, \[
 \Big|\Big|\int f(x)d\mu(x)\Big|\Big| \leq \int ||f(x)||d\mu(x).\]
  \end{ppsn}
 \textit{Proof.} Let $B:\clh \times \clh \to \mathbb{C}$ be defined by \[B(\xi,\eta)=\int \langle f(x)\xi|\eta \rangle d\mu(x).\] Then, $B$ is a bounded sesquilinear form on $\clh$. Consequently, there exists a unique bounded linear operator, denote it by $\int f(x)d\mu(x)$, such that \[\Big\langle \big(\int f(x)d\mu(x)\big)\xi\big|\eta \Big\rangle=\int \langle f(x)\xi|\eta \rangle d\mu(x)\] for $\xi,\eta \in \clh$. The estimate about the norm is obvious. Hence the proof. \hfill $\Box$
 
 Note the following properties about operator valued integration. 
 \begin{enumerate}
 \item[(1)] Suppose $f,g:X \to B(\clh)$ are integrable. Then $\alpha f+g$ is integrable for every $\alpha \in \bbc$ and in that case 
 \[
 \int (\alpha f(x)+g(x))d\mu(x)=\alpha \int f(x)d\mu(x)+\int g(x)d\mu(x).\]
 In short, the integral is linear. 
 \item[(2)] Suppose $f:X \to B(\clh)$ is integrable and $T \in B(\clh)$. Then, the maps \[x \to Tf(x); ~~\textrm{~and~}x \to f(x)T\] are integrable. Also, we have the equality
 \[
 \int Tf(x)d\mu(x)=T\Big(\int f(x)d\mu(x)\Big)\] and  \[\Big(\int f(x)d\mu(x)\Big)T=\int f(x)Td\mu(x).\] 
 \item[(3)] Suppose $f:X \to B(\clh)$ is integrable. Then, $X \ni x \to f(x)^{*} \in B(\clh)$ is integrable and
 \[
 \int f(x)^{*}d\mu(x)=\Big(\int f(x)d\mu(x)\Big)^{*}.\]
 \end{enumerate}
 
 \begin{xrcs}
 Formulate a version of the dominated convergence theorem and prove it. 
  \end{xrcs}
   Let us now return to the study of non-degenerate bounded $^*$-representations of $C_{c}(G)$.  Let $U:G \to B(\clh)$ be a unitary representation of $G$ on $\clh$. For $f \in C_{c}(G)$, let $\pi_{U}(f)$ be the bounded operator given by the equation
  \[
  \pi_{U}(f)=\int f(s)U_s ds.\]
  Note that $\pi_{U}(f)$ exists, since $s \to f(s)U_s$ is weakly continuous and $s \to ||f(s)U_{s}||=|f(s)|$ is integrable. Clearly, for $f \in C_{c}(G)$, $||\pi_{U}(f)|| \leq ||f||_{1}$. 
 
 \begin{ppsn}
 \label{unitary to repn}
 Keep the foregoing notation. The map $\pi_{U}:C_c(G) \to B(\clh)$ is a bounded non-degenerate $^*$-representation of $C_{c}(G)$. The representation $\pi_{U}$ is continuous w.r.t. the inductive limit topology.

 \end{ppsn}
 \textit{Proof.} First let us check that $\pi_{U}$ preserves the adjoints. Let $f \in C_{c}(G)$ and $\xi, \eta \in \clh$ be given.
Then 
\begin{align*}
\langle \pi_{U}(f^{*})\xi|\eta \rangle&=\int f^{*}(s)\langle U_{s}\xi|\eta \rangle ds \\
                                                      &= \int \Delta(s^{-1})\overline{f(s^{-1})}\langle U_{s}\xi|\eta \rangle ds \\
                                                      &=\int \Delta(s^{-1})\langle \xi|f(s^{-1})U_{s^{-1}}\eta \rangle ds \\
                                                      &=\int \langle \xi|f(s)U_{s}\eta \rangle ds \\
                                                      &=\int \overline{f(s) \langle U_{s}\eta|\xi \rangle}ds \\
                                                      &=\overline{\langle \pi_{U}(f)\eta|\xi\rangle} \\
                                                      &= \langle \pi_{U}(f)^{*}\xi|\eta \rangle.
\end{align*}
Hence $\pi_{U}(f^{*})=\pi_{U}(f)^{*}$.

Let $f,g \in C_{c}(G)$ and $\xi, \eta \in \clh$ be given. Calculate as follows to observe that 
\begin{align*}
\langle \pi_{U}(f*g)\xi|\eta \rangle&= \int f*g(s)\langle U_{s}\xi|\eta \rangle ds \\
                                                     & = \int \Big (\int f(t)g(t^{-1}s)dt \Big) \langle U_{s}\xi|\eta \rangle ds \\
                                                     &= \int f(t)\Big (\int g(t^{-1}s)\langle U_{s}\xi|\eta \rangle ds\Big)dt \\
                                                     &= \int f(t) \Big( \int g(s) \langle U_{ts}\xi|\eta \rangle ds \Big)dt \\
                                                     &= \int f(t) \Big(\int g(s) \langle U_{s}\xi|U_{t}^{*}\eta \rangle ds \Big)dt \\
                                                     & = \int f(t) \langle \pi_{U}(g)\xi|U_{t}^{*}\eta \rangle dt \\
                                                     &= \int f(t)\langle U_{t}\pi_{U}(g)\xi|\eta \rangle dt \\
                                                     &= \langle \pi_{U}(f)\pi_{U}(g)\xi|\eta \rangle.
\end{align*}
Consequently, we have $\pi_{U}(f*g)=\pi_{U}(f)\pi_{U}(g)$. This proves that $\pi_{U}$ is a $^*$-representation. We have already noted that $\pi_{U}$ is bounded. The continuity w.r.t. the inductive limit topology follows as a consequence. 
 
 Let $\{\phi_n\}_{n=1}^{\infty}$ be the approximate identity constructed in Prop. \ref{approximate}. Keep the notation used in Prop. \ref{approximate}. We claim that $\pi_{U}(\phi_n)$ converges strongly to $Id$. 
  Let $\xi \in \clh$ be given. Suppose $\epsilon>0$ is given. Choose $N$ large such that for $s \in U_{N}$, $||U_{s}\xi-\xi|| \leq \epsilon$. For $\eta \in \clh$ and $n \geq N$, calculate as follows to observe that 
  \begin{align*}
  |\langle \pi_{U}(\phi_n)\xi-\xi|\eta\rangle|&= \Big|\int \phi_n(s)\langle U_{s}\xi|\eta \rangle ds-\int \phi_n(s)\langle \xi|\eta \rangle ds\Big|\\
  &\leq \Big|\int \phi_n(s)\langle U_s\xi-\xi|\eta \rangle ds\Big| \\
  & \leq \int \phi_n(s)||U_s\xi-\xi||||\eta||ds \\
  & \leq \epsilon ||\eta||. 
    \end{align*}
 Hence, for $n \geq N$, $||\pi_{U}(\phi_n)\xi-\xi|| \leq \epsilon$. Consequently, $\pi_{U}(\phi_n)\xi \to \xi$ for every $\xi \in \clh$. This proves our claim. Thus, $\pi_{U}$ is non-degenerate. This completes the proof. \hfill $\Box$
 
 The representation $\pi_{U}$ constructed in the previous proposition is called \emph{the integrated form of $U$}. 
 \begin{xrcs}
 \label{uniqueness of unitary rep}
 Keep the notation of the previous proposition. Prove that 
 \begin{enumerate}
 \item[(1)] for $s \in G$ and $f \in C_{c}(G)$, $U_s\pi_{U}(f)=\pi_{U}(L_sf)$, and
 \item[(2)]  for $s \in G$, $\pi_{U}(L_s\phi_n) \to U_{s}$ in the strong operator topology.
  \end{enumerate}
 Thus, if $U$ and $V$ are two unitary representations of $G$ on the same Hilbert space, then $\pi_{U}=\pi_{V}$ if and only if $U=V$. 
  \end{xrcs}
 
 Next we show that every non-degenerate bounded $*$-representation of $C_{c}(G)$ is of the form $\pi_{U}$ for a unique unitary representation $U$ of $G$. We need to invoke vector valued integration at a crucial point and it is worthwhile to digress a bit into vector valued integration. We start by recalling the Krein-Smulian theorem whose proof can be found for instance in \cite{Conway}.
 
 \begin{rmrk}[Krein-Smulian]
Let $E$ be a separable Banach space and $\phi:E^{*} \to \mathbb{C}$ be a linear functional. Then $\phi$ is weak  $^*$-continuous if and only if $\phi$ is weak $^*$-sequentially continuous.
We refer the reader to Corollary 12.8 of \cite{Conway}.
\end{rmrk}

 \textbf{Vector valued integration:} Suppose $E$ is a separable Banach space and let $(X,\mathcal{B})$ be a measurable space. 
\begin{enumerate}
\item[(1)] A map $f:X \to E$ is said to be \emph{weakly measurable} if $\phi \circ f$ is measurable for every $\phi \in E^{*}$. 
\item[(2)] Suppose $f:X \to E$ is weakly measurable. Then, the map $X \ni x \to ||f(x)|| \in \mathbb{C}$ is measurable. This is because, since $E$ is separable, the unit ball of $E^{*}$ w.r.t. to the weak $^*$-topology is a compact metrisable space. 
\item[(3)] Let $\mu$ be a measure on $(X,\mathcal{B})$ and  let $f:X \to E$ be a weakly measurable map. We say that $f$ is integrable w.r.t $\mu$ if $x \to ||f(x)||$ is integrable. Suppose $f$ is integrable. Define $F:E^{*} \to \mathbb{C}$ by
\[
F(\phi)=\int \phi(f(x))d\mu(x).\]
An application of the Krein-Smulian theorem implies that $F$ is weak $^*$-continuous. Thus, there exists a unique element, denoted $\int f(x)d\mu(x) \in E$, such that 
\[
\phi\Big(\int f(x)d\mu(x)\Big)=\int \phi(f(x))d\mu(x)\]
for every $\phi \in E^{*}$. 
We call $\int f(x) d\mu(x)$, the integral of $f$ w.r.t the measure $\mu$. The $\int$ satisfies the usual linearity properties and the dominated convergence theorem.
\end{enumerate}

Now let $\mu$ be a Haar measure on $G$. Let $f,g \in C_{c}(G)$ be given. Consider $g$ as an element of $L^{1}(G)$. Note that the map $G \ni s \to L_{s}(g) \in L^{1}(G)$ is continuous, when $L^{1}(G)$ is given the norm topology. Consequently, the vector valued integral, in the sense explained above, $\int f(s)L_{s}(g)ds$ exists. 

\begin{lmma}
\label{integral expression}
With the foregoing notation, we have $f*g=\int f(s)L_s(g)ds$ in $L^{1}(G)$. 
\end{lmma}
\textit{Proof.} For $\phi \in L^{\infty}(G)$, let $\omega_{\phi}:L^{1}(G) \to \mathbb{C}$ be defined by $\omega_{\phi}(f)=\int f(s)\phi(s)ds$. The map $\phi \to \omega_{\phi}$ identifies $L^{\infty}(G)$ with the dual of $L^{1}(G)$.  
 It suffices to show that for every $\phi \in L^{\infty}(G)$, $\omega_{\phi}(f*g)=\omega_{\phi}\Big(\int f(s)L_{s}(g)ds \Big)$. Fix $\phi \in L^{\infty}(G)$. Calculate as follows to observe that 
 \begin{align*}
 \omega_{\phi}(f*g)&=\int f*g(t)\phi(t)dt \\
 &=\int \Big(\int f(s)g(s^{-1}t)ds)\phi(t)dt \\
 &=\int f(s)\Big(\int \phi(t)g(s^{-1}t)dt\Big)ds\\
 &=\int f(s)\Big(\int \phi(t)L_{s}(g)(t)dt\Big)ds\\
 &=\int f(s)\omega_{\phi}(L_{s}(g))ds\\
 &=\omega_{\phi}\Big(\int f(s)L_{s}(g)ds\Big).
  \end{align*}
  Hence the proof. \hfill $\Box$
  
  \begin{ppsn}
  Let $\pi:C_{c}(G) \to B(\clh)$ be a  bounded $^*$-representation that is non-degenerate. Then, there exists a unique unitary representation $U:G \to B(\clh)$ such that $\pi=\pi_{U}$. 
    \end{ppsn}
  \textit{Proof.} Uniqueness follows from Exercise \ref{uniqueness of unitary rep}. Note that $\{\pi(f)\xi: f \in C_{c}(G), \xi \in \clh\}$ is total in $\clh$. Let $s \in G$ be given. Note that for $f,g \in C_{c}(G)$, $(L_{s}f)^{*}*L_{s}g=f^{*}*g$. This has the consequence that for $\xi,\eta \in \clh$ and $f,g \in C_{c}(G)$, 
  \[
  \langle \pi(L_sf)\xi|\pi(L_{s}g)\eta \rangle=\langle \pi(f)\xi|\pi(g)\eta \rangle.\]
  Hence, there exists a unique unitary operator, denote it by $U_s$, such that \[U_{s}\pi(f)\xi=\pi(L_sf)\xi\] for $f \in C_{c}(G)$ and $\xi \in \clh$. Evaluating on the total set $\{\pi(f)\xi: f \in C_{c}(G), \xi \in \clh\}$, it is straightforward to verify that $U_sU_t=U_{st}$ for $s,t \in G$. To check that $\{U_{s}\}_{s \in G}$ is strongly continuous, it is sufficient to verify that for $f \in C_{c}(G)$ and $\xi \in \clh$, the map $G \ni s \to \pi(L_sf)\xi \in \clh$ is continuous. But the last assertion follows, as $\pi$ is continuous with respect to the inductive limit topology.
  
  We claim that $\pi=\pi_{U}$. Since $\{\pi(g)\xi: g \in C_{c}(G), \xi \in \clh\}$ is total in $\clh$, It suffices to show that for $f,g \in C_{c}(G)$ and $\xi,\eta \in \clh$,
  \[
  \langle \pi_{U}(f)\pi(g)\xi|\eta \rangle=\langle \pi(f)\pi(g)\xi|\eta \rangle.\]
  Let $f,g \in C_{c}(G)$ and $\xi,\eta \in \clh$ be given. Denote the linear extension of $\pi$ to $L^{1}(G)$ by $\widetilde{\pi}$. Define $\omega:L^{1}(G) \to \mathbb{C}$ by $
  \omega(h)=\langle \widetilde{\pi}(h)\xi|\eta \rangle.$
  Calculate as follows to observe that 
  \begin{align*}
  \langle \pi(f)\pi(g)\xi|\eta \rangle&=\langle \pi(f*g)\xi|\eta \rangle \\
    &=\omega(f*g) \\
   &=\int f(s)\omega(L_sg)ds ~~(\textrm{by Lemma \ref{integral expression}})\\
   &=\int f(s)\langle \pi(L_{s}g)\xi|\eta \rangle ds\\
   &= \int f(s)\langle U_{s}\pi(g)\xi|\eta \rangle ds\\
   &=\langle \Big(\int f(s)U_{s}ds\Big)\pi(g)\xi|\eta \rangle\\
   &=\langle \pi_{U}(f)\pi(g)\xi|\eta \rangle. 
    \end{align*}
  This completes the proof. \hfill $\Box$

  Next, we show that the universal $C^{*}$-norm is indeed a norm by exhibiting a faithful representation of $C_{c}(G)$. Let $\lambda:=\{\lambda_s\}_{s \in G}$ be the left regular representation of $G$ on $L^{2}(G)$. We denote the integrated form of the left regular representation by $\lambda$ itself. Thus, for $f \in C_{c}(G)$, $\lambda(f)=\int f(s)\lambda_s ds$. 
  Let $f,g,h \in C_{c}(G)$ be given. Calculate as follows to observe that 
  \begin{align*}
  \langle \lambda(f)g|h \rangle&=\int f(s)\langle \lambda_s(g)|h \rangle ds\\
  &= \int f(s) \Big(\int g(s^{-1}t)\overline{h(t)}dt\Big)ds\\
  &=\int \overline{h(t)}\Big(\int f(s)g(s^{-1}t)ds\Big)dt \\
  &=\int (f*g)(t)\overline{h(t)}dt\\
  &=\langle f*g|h \rangle.
   \end{align*}
  Hence for $f,g \in C_{c}(G)$, $\lambda(f)g=f*g$. 
  
  \begin{lmma}
 The representation $\lambda$ is faithful.
  \end{lmma}
 \textit{Proof.} Let $f \in C_{c}(G)$ be such that $\lambda(f)=0$. Consider an approximate identity $\{\phi_n\}_{n=1}^{\infty}$ of $C_{c}(G)$. Considered as an element of $L^{2}(G)$, we have $f*\phi_n=\lambda(f)\phi_n=0$. Hence ,$f*\phi_n=0$ in $C_{c}(G)$.
 But $f*\phi_n \to f$ in the inductive limit topology. Consequently, it follows that $f=0$. Hence $\lambda$ is faithful. \hfill $\Box$

  An immediate consequence of the previous lemma is that the universal $C^{*}$-seminorm on $C_{c}(G)$ given by 
  \[
  ||f||:=\sup\{||\pi_U(f)||: \textrm{~$U$ is a unitary representation of $G$}\}
  \]
  for $f \in C_{c}(G)$ is indeed a norm. 
  
  \begin{dfn}
  The completion of $C_{c}(G)$ with respect to the universal $C^{*}$-norm is called the full group $C^{*}$-algebra and is denoted $C^{*}(G)$. 
  For $f \in C_{c}(G)$, let \[
  ||f||_{red}:=||\lambda(f)||.\]
  Here, $\lambda$ is the left regular representation.   Then, $||~||_{red}$ is a $C^{*}$-norm on $C_{c}(G)$ and its completion is called the reduced $C^{*}$-algebra of $G$ and is denoted $C_{red}^{*}(G)$. Note that $C_{red}^{*}(G)$ is the $C^{*}$-subalgebra of $B(L^{2}(G))$ generated by
  $\{\lambda(f):f \in C_{c}(G)\}$. 
   The map $C_{c}(G) \ni f \to \lambda(f) \in C_{red}^{*}(G)$ extends to a surjective homomorphism from $C^{*}(G)$ onto $C_{red}^{*}(G)$. 
    \end{dfn}

 Note that non-degenerate $*$-representations of $C^{*}(G)$ are in one-one correspondence with bounded non-degenerate $*$-representations of $C_{c}(G)$, which in turn is in one-one correspondence with unitary representations of $G$. 
 In other words, the map \[U \to \pi_{U}\] identifies strongly continuous unitary representations of $G$ and non-degenerate representations of $C^{*}(G)$. Moreover, this map preserves unitary equivalence, direct sum, irreducibility, etc.....
 Thus, in principle, studying the representation theory of a locally compact group is equivalent to studying the representation theory of $C^{*}(G)$. We prove Raikov's theorem, which asserts that every locally compact group has sufficiently many irreducible representations, as an application of this principle. 
 
 We need the following which is the corollary to Theorem 1.7.2 of \cite{Arveson_invitation}. 
 \begin{ppsn}
 Let $A$ be a $C^{*}$-algebra and $a \in A$ be non-zero. Then, there exists an irreducible representation $\pi$ such that $\pi(a) \neq 0$. In other words, irreducible representations of $A$ separate points of $A$. 
 
 \end{ppsn}
 
 \begin{thm}[Raikov's theorem]
 Let $G$ be a locally compact second countable Hausdorff topological group. For every $s \in G$, with $s\neq e$, there exists an irreducible unitary representation $U$ of $G$ such that $U_s \neq Id$. 
  \end{thm}
  \textit{Proof.} Let $s \in G$ be such that $s \neq e$. Suppose, on the contrary, assume that for every irreducible unitary representation $U$ of $G$, $U_s = Id$. Choose $f \in C_{c}(G)$ such that $L_sf \neq f$. Then, there exists an irreducible representation $\pi$ of $C^{*}(G)$, say on the Hilbert space $\clh$, such that $\pi(L_sf)\neq \pi(f)$. Let $U$ be the unitary representation of $G$ such that $\pi=\pi_{U}$. Then, $U$ is irreducible. 
  Note that $U_{s}\pi_{U}(f)=\pi_{U}(L_sf) \neq \pi_{U}(f)$. Hence $U_s \neq I$.  \hfill $\Box$
  
 We end this section by identifying the $C^{*}$-algebra of an abelian group. For the rest of this section, let $G$ be a locally compact, second countable, topological group which we assume is abelian. 
  Note that $C^{*}(G)$ is commutative. First, we identify the spectrum of $C^{*}(G)$. Let $\chi:G \to \mathbb{T}$ be a continuous map. We say that $\chi$ is a \emph{character of $G$} if 
  $\chi(st)=\chi(s)\chi(t)$ for every $s,t \in G$. Denote the set of characters of $G$ by $\widehat{G}$. For $\chi_1,\chi_2 \in \widehat{G}$, define $\chi_1.\chi_2:G \to \mathbb{T}$ 
  by the formula \[
  \chi_1.\chi_2(s)=\chi_1(s)\chi_2(s)\]
  for $s \in G$. Then $\chi_1.\chi_2 \in \widehat{G}$. With this multiplication, $\widehat{G}$ is an abelian group. 
  
  We endow $\widehat{G}$ with the topology of uniform convergence on compact sets. 
  The convergence of nets in $\widehat{G}$ is as follows. Let $(\chi_i)$ be a net in $\widehat{G}$ and let $\chi \in \widehat{G}$ be given. Then $\chi_{i} \to \chi$ if and only if for every
  compact set $K \subset G$, the net $(\chi_i)$ converges uniformly to $\chi$ on $K$. Endowed with the topology of convergence on compact sets, $\widehat{G}$ is a topological
  group. 
  
  Set $A:=C^{*}(G)$. Let $\chi \in \widehat{G}$ be given. Then, $\chi$ is a $1$-dimensional unitary representation of $G$. Thus, there exists a $^*$-homomorphism denoted $\omega_{\chi}:A \to \mathbb{C}$ such that 
  \[\omega_{\chi}(f)=\int f(s)\chi(s)ds\] for $f \in C_{c}(G)$. 
  
  \begin{xrcs}
  \label{convergence on compact sets}
  Let $E$ be a Banach space and $\{\phi_{i}\}$ be a bounded net in $E^{*}$. Suppose $\phi \in E^{*}$ and $\phi_{i} \to \phi$ in the weak $^*$-topology. Let $K$ be a compact set in $E$. Then $\phi_{i} \to \phi$ uniformly on $K$.

  \end{xrcs}

  \begin{thm}
  \label{abelian group algebras}
  With the foregoing notation, the map $\widehat{G} \ni \chi \to \omega_{\chi} \in \widehat{A}$ is a homeomorphism. As a consequence, it follows that $C^{*}(G) \simeq C_{0}(\widehat{G})$. 
    \end{thm}
   \textit{Proof.} First, we prove that $\chi \to \omega_{\chi}$ is continuous.  Suppose $\chi_i$ is a net in $\widehat{G}$ and $\chi_i \to \chi \in \widehat{G}$. Since $\{\omega_{\chi_i}\}$ is uniformly bounded, it suffices to 
  prove that for $f \in C_{c}(G)$, $\omega_{\chi_i}(f) \to \omega_{\chi}(f)$. Let $f \in C_{c}(G)$ be given. Denote the support of $f$ by $K$. Let $\epsilon>0$ be given. Choose $i_0$ such that for $i \geq i_0$, $|\chi_i(x)-\chi(x)| \leq \epsilon$ for $x \in K$. 
  Calculate as follows to observe that for $i \geq i_0$,
  \begin{align*}
  |\omega_{\chi_i}(f)-\omega_{\chi}(f)|&=|\int (f(s)\chi_i(s)-f(s)\chi(s))ds \\
  & \leq \int_{K} |f(s)||\chi_i(s)-\chi(s)|ds\\
  & \leq \epsilon ||f||_1.
  \end{align*}
  This proves that $\omega_{\chi_i} \to \omega_\chi$. Hence, the map $\chi \to \omega_{\chi}$ is continuous. 
  
  Let $\omega$ be a character of $A$. View $\omega$ as a $1$-dimensional representation on the Hilbert space $\mathbb{C}$. Then there exists a unitary representation $\chi:G \to \mathcal{U}(\mathbb{C})\simeq \mathbb{T}$ such that for $f \in C_{c}(G)$, 
  $\omega(f)=\int f(s)\chi(s)ds=\omega_{\chi}(f)$. Since $C_{c}(G)$ is dense in $C^{*}(G)$, it follows that $\omega=\omega_\chi$. This proves that the map $\chi \to \omega_\chi$ is onto. The injectivity of the map follows from Exercise \ref{uniqueness of unitary rep}.
  
  Consider a net $(\omega_{\chi_i}) \to \omega_{\chi}$. Then $\omega_{\chi_i\chi^{-1}} \to \omega_1$. It suffices to prove that $\chi_i\chi^{-1} \to 1$. Thus, with no loss of generality, we can assume that $\chi$ is the trivial character. Denote $\omega_\chi$ by $\omega_0$. Let $K$ be a compact subset of $G$ and let $\epsilon>0$ be given. Choose $f \in C_{c}(G)$ such that $f \geq 0$ and $\int f(s)ds=1$. Note that $f*\chi(s)=1$ for every $s \in G$.
  
  Note that the inclusion $C_{c}(G) \to A$ is continuous. Hence the set $\{L_{s^{-1}}f: s \in K\}$ is a compact subset of $A$. By Exercise \ref{convergence on compact sets}, there exists $i_0$ such that for $i \geq i_0$ and $s \in K$, \[|f*\chi_i(s)-\chi(s)|= |f*\chi_i(s)-f*\chi(s)|=|\omega_{\chi_i}(L_{s^{-1}}f)-\omega_0(L_{s^{-1}}f)| \leq \epsilon.\] 
  Choose $i_1$ such that for $i \geq i_1$, $|\omega_{\chi_i}(f)-\omega_0(f)| \leq \epsilon$. Choose $i_2 \geq i_0, i_1$. Let $i \geq i_2$ and $s \in G$ be given. Calculate as follows to observe that 
  \begin{align*}
  |f*\chi_i(s)-\chi_i(s)|&=|\int f(t)\chi_{i}(t^{-1}s)ds-\chi_{i}(s)\int f(t)dt| \\
  & \leq |\int \chi_i(s)\overline{(\chi_i(t)-1)}f(t)dt|\\
  & \leq |\chi_i(s)| \Big|\int f(t)(\chi_i(t)-1)dt\Big| \\
  & \leq |\omega_{\chi_i}(f)-\omega_0(f)|\\
  & \leq \epsilon.
    \end{align*}
  Combining the above two inequalities, we see that for $i \geq i_2$ and $s \in K$, \[|\chi_i(s)-\chi(s)| \leq 2\epsilon.\] This proves that $\chi_i \to \chi$. Hence, the map $\widehat{G} \ni \chi \to \omega_\chi \in \widehat{A}$ is a homeomorphism. This completes the proof. \hfill $\Box$
  
  We end this section by   discussing Plancherel's theorem for abelian groups. Let $G$ be a locally compact, second countable, Hausdorff topological group, which we assume is abelian. For $f \in L^{1}(G)$, let $\widehat{f}:\widehat{G} \to \mathbb{C}$ be defined by 
  \[
  \widehat{f}(\chi)=\int f(s)\overline{\chi(s)}ds\] 
  for $\chi \in \widehat{G}$. The function $\widehat{f}$ is called \emph{the Fourier transform of $f$}. Note that for $f \in C_{c}(G)$, $\widehat{f}(\chi)=\omega_{\overline{\chi}}(f)$. Hence, $\widehat{f} \in C_{0}(\widehat{G})$ for $f \in C_{c}(G)$. Using the fact that $C_{c}(G)$ is dense in $L^{1}(G)$, it follows at once that $\widehat{f} \in C_0(\widehat{G})$ for $f \in L^{1}(G)$. 
  
\begin{thm}[Plancherel's theorem]
\begin{enumerate}
\item[(1)] For $f \in L^{1}(G) \cap L^{2}(G)$, $\widehat{f} \in L^{2}(\widehat{G})$.
\item[(2)] There exists a unique Haar measure $\overline{\mu}$ on $\widehat{G}$ such that the map \[L^{1}(G)\cap L^{2}(G) \ni f \to \widehat{f} \in L^{2}(\widehat{G})\] extends to a unitary map from $L^{2}(G)$ onto $L^{2}(\widehat{G})$. The unitary $f \to \widehat{f}$ is usually denoted $\mathcal{F}$. 
\end{enumerate}

\end{thm}
Let $\sigma:C^{*}(G) \to C_{0}(\widehat{G})$ be the $^*$-homomorphism such that for $f \in C_{c}(G)$ and $\chi \in \widehat{G}$,
\[
\sigma(f)(\chi)=\omega_{\overline{\chi}}(f)=\int f(s)\chi(s)ds.\]
Theorem \ref{abelian group algebras} and the Gelfand-Naimark theorem asserts that $\sigma$ is well-defined and is a $^*$-isomorphism. Let $\overline{\mu}$ be a Haar measure on $\widehat{G}$ as in Plancherel's theorem.   Let $M:C_{0}(\widehat{G}) \to B(L^{2}(\widehat{G}))$ be the multiplication representation, i.e. for $f \in C_{0}(\widehat{G})$ the operator $M(f)$ is given by
\[
M(f)\xi(\chi)=f(\chi)\xi(\chi)
\]
for $\xi \in L^{2}(\widehat{G})$. Then, $M$ is a faithful representation. Let $\lambda:C^{*}(G) \to B(L^{2}(G))$ be the left regular representation. For $f \in C_{c}(G)$ and $\xi \in C_{c}(G)$, calculate as follows to observe that 
\begin{align*}
\mathcal{F}\lambda(f)\xi(\chi)&=\widehat{f*\xi}(\chi) \\
 &=\int (f*\xi)(s)\chi(s)ds \\
 &= \int (\int f(t)\xi(t^{-1}s)dt)\chi(s)ds \\
 &=\int f(t)\chi(t)\Big(\int \xi(t^{-1}s)\chi(t^{-1}s)ds\Big)dt \\
 &= \int f(t)\chi(t)\Big(\int \xi(s)\chi(s)ds\Big)dt\\
 &=\int f(t)\chi(t)\widehat{\xi}(\chi)dt \\
 &=\omega_{\overline{\chi}}(f)\mathcal{F}\xi(\chi)\\
 &=\sigma(f)(\chi)\mathcal{F}\xi(\chi).
 \end{align*}
  Hence $\mathcal{F}\lambda(f)=M(\sigma(f))\mathcal{F}$ for $f \in C_{c}(G)$. Hence $\mathcal{F}\lambda(.)\mathcal{F}^{*}=M \circ \sigma$. But $M \circ \sigma$ is a faithful representation of $C^{*}(G)$. Thus we obtain the following corollary.
  
  \begin{crlre}
  The left regular representation $\lambda:C^{*}(G) \to C_{red}^{*}(G)$ is an isomorphism. 
    \end{crlre}
    
    We finish this section by stating the Pontraygin duality theorem. We refer the reader to Chapter 4 of \cite{Folland} for a proof. Let $G$ be a locally compact abelian group. Denote the dual group by $\widehat{G}$. For $s \in G$, let $\widehat{s}:\widehat{G} \to \mathbb{T}$ be defined by
    \[
    \widehat{s}(\chi)=\chi(s)\]
    for $\chi \in \widehat{G}$. Then $\widehat{s} \in \widehat{\widehat{G}}$. Moreover, the map $G \ni s \to \widehat{s} \in \widehat{\widehat{G}}$ is continuous. Raikov's theorem implies that the map $s \to \widehat{s}$ is indeed one-one. 
  
  \begin{thm}[Pontraygin duality]
  The map $\widehat{G} \ni s \to \widehat{s} \in \widehat{\widehat{G}}$ is a homeomorphism. 
  \end{thm}

   \begin{xrcs}
  In this  exercise, we identify the duals of a few concrete abelian groups. 
   
  \begin{enumerate}
  
  \item[(1)] For $\xi \in \mathbb{R}$, let $\chi_{\xi}:\mathbb{R} \to \mathbb{T}$ be defined by $\chi_\xi(x)=e^{2\pi i x\xi}$. Prove that  the map $\mathbb{R} \ni \xi \to \chi_{\xi} \in \widehat{\mathbb{R}}$ is a homeomorphism of topological groups.
  \item[(2)]  For $z \in \mathbb{T}$, let $\chi_z: \mathbb{Z} \to \mathbb{T}$ be defined by $\chi_z(n)=z^{n}$. Show that $\mathbb{T} \ni z \to \chi_z \in \widehat{\mathbb{Z}}$ is a homeomorphism of topological groups.
  \item[(3)] For $ m \in \mathbb{Z}$, let $\chi_m:\mathbb{T} \to \mathbb{T}$ be defined by $\chi_m(z)=z^{m}$. Show that $\mathbb{Z} \ni m \to \chi_m \in \widehat{\mathbb{T}}$ is a homeomorphism of topological groups. 
  \item[(4)] Identify the duals of $\mathbb{R}^{d}$, $\mathbb{Z}^{d}$ and $\mathbb{T}^{d}$ for $d \geq 2$. 
  \end{enumerate}
  \end{xrcs}

  \chapter{Crossed products and Hilbert $C^{*}$-modules}
  \section{Crossed products}
  In this section, we discuss the notion of crossed products of $C^{*}$-algebras associated with actions of topological groups. We will omit  proofs as we have the discussed the case of group $C^{*}$-algebras and discrete crossed products in complete detail. 
  
  Let $A$ be a $C^{*}$-algebra and let $G$ be a locally compact, second countable, Hausdorff topological group. By an action of $G$ on $A$, we mean a map $\alpha:G \to Aut(A)$, the image of an element $s$ under $\alpha$ is denoted by $\alpha_s$, such that 
  \begin{enumerate}
  \item[(1)] for $s \in G$, $\alpha_s:A \to A$ is an automorphism, 
    \item[(2)] for $s,t \in G$, $\alpha_s \circ \alpha_t=\alpha_{st}$, and
  \item[(3)] for $a \in A$, the map $G \ni s \to \alpha_s(a) \in A$ is continuous, when $A$ is given the norm topology.
   \end{enumerate}
  The triple $(A,G,\alpha)$ is called a $C^{*}$-dynamical system. 
  
  \begin{xrcs}
  Let $(A,G,\alpha)$ be a $C^{*}$-dynamical system. Prove that the map
  \[
  G \times A \ni (s,a) \to \alpha_s(a) \in A\]
  is continuous. 
   \end{xrcs}
  
  \begin{xmpl}
  Let $X$ be a left $G$-space, where $X$ is a locally compact, second countable, Hausdorff topological space. Define for $s \in G$, $\alpha_s:C_{0}(X) \to C_{0}(X)$ by 
  \[
  \alpha_s(f)(x):=f(s^{-1}x)\]
  for $f \in C_{0}(X)$. Then, $(C_{0}(X),G,\alpha)$ is a $C^{*}$-dynamical system. 
    \end{xmpl}
  Let $(A,G,\alpha)$ be a $C^{*}$-dynamical system. Let $ds$ be a left Haar measure on $G$ and $\Delta$ be the modular function of $G$. 
  Consider the vector space $C_{c}(G,A)$, i.e.
  \[
  C_{c}(G,A):=\{f:G \to A: \textrm{$f$ is continuous and compactly supported}\}.\]
  Define the convolution and the involution on $C_{c}(G,A)$ as follows : for $f,g \in C_{c}(G,A)$,
  \begin{align*}
  f*g(s)&=\int f(t)\alpha_{t}(g(t^{-1}s))dt \\
  f^*(s)&=\Delta(s^{-1})\alpha_{s}(f(s^{-1})^{*}).\
    \end{align*}
    Then $C_{c}(G,A)$ becomes a $^*$-algebra. For $f \in C_{c}(G,A)$, let 
    \[
    ||f||_{1}:=\int ||f(s)||ds.
    \]
    With the norm defined above $C_{c}(G,A)$ is normed $^*$-algebra. The crossed product $A \rtimes_{\alpha} G$ is defined as the enveloping $C^{*}$-algebra of $C_{c}(G,A)$. 
    
    \begin{dfn}
    Let $(A,G,\alpha)$ be a $C^{*}$-dynamical system. By a covariant representation of $(A,G,\alpha)$ on a Hilbert space $\clh$, we mean a pair $(\pi,U)$ such that 
    \begin{enumerate}
    \item[(1)] $\pi$ is a $^*$-representation of $A$, 
    \item[(2)] $U:=\{U_s\}_{s \in G}$ is a strongly continuous representation of $G$, and
    \item[(3)] for $s \in G$ and $a \in A$, the covariance condition $U_s\pi(a)U_{s}^{*}=\pi(\alpha_s(a))$ is satisfied. 
   \end{enumerate}
   We say that $(\pi,U)$ is non-degenerate if $\pi$ is non-degenerate. 
  \end{dfn}
  
  The following theorem characterises bounded $^*$-representations of $C_{c}(G,A)$. 
  
  \begin{thm}
  \label{correspondence for crossed products}
  Let $(A,G,\alpha)$ be a $C^{*}$-dynamical system. Let $(\pi,U)$ be a non-degenerate covariant representation of $(A,G,\alpha)$ on a Hilbert space $\clh$. For $f \in C_{c}(G,A)$, let 
  \[
  (\pi \rtimes U)(f):=\int \pi(f(s))U_s ds.\]
  Then, $\pi \rtimes U$ is a non-degenerate bounded $*$-representation of $C_{c}(G,A)$. 
  
   Suppose $\widetilde{\pi}$ is a bounded non-degenerate $^*$-representation of $C_{c}(G,A)$ on a Hilbert space $\clh$. Then, there exists a unique covariant representation $(\pi,U)$ of $(A,G,\alpha)$ on $\clh$ which is non-degenerate such that 
  for $f \in C_{c}(G,A)$, 
  \[
  \widetilde{\pi}(f)=(\pi \rtimes U)(f)=\int \pi(f(s))U_s ds.\]

  \end{thm}
  \begin{rmrk}
   The map \[(\pi,U) \to \pi \rtimes U\] establishes a bijective correspondence between non-degenerate covariant representations of $(A,G,\alpha)$ and non-degenerate bounded $^*$-representations of $C_{c}(G,A)$. 
  Moreover, the above correspondence preserves direct sum, irreducibility etc.....
  \end{rmrk}
  For $f \in C_{c}(G,A)$, let 
  \begin{align*}
  ||f|| :=  \sup\{& ||(\pi \rtimes U)(f)||: (\pi,U) \\
  & \textrm{~is a non-degenerate covariant representation of $(A,G,\alpha)$}\}.
  \end{align*}
    We will show that $||~||$ is indeed a genuine $C^{*}$-norm on $C_{c}(G,A)$. The completion of $C_{c}(G,A)$ with respect to this universal norm is called the \textit{full crossed product} and is denoted $A \rtimes G$. 
  
  Next, we exhibit a concrete covariant representation of $(A,G,\alpha)$.  Let $\pi$ be  a faithful non-degenerate $*$-representation of $A$ on a Hilbert space $\clh$. 
 Set $\widetilde{\clh}:=L^{2}(G,\clh)$. Recall that $\widetilde{\clh}$ consists of weakly measurable square integrable $\clh$-valued functions. The inner product on $\widetilde{\clh}$ is given by
 \[
 \langle \xi|\eta \rangle:=\int \langle \xi(s)|\eta(s) \rangle ds\]
 for $\xi,\eta \in \widetilde{\clh}$. The proof that $\widetilde{\clh}$ is a Hilbert space is similar to the case when $\clh=\mathbb{C}$. 
 
 For $s \in G$ and $a \in A$, let $\lambda_s$ and $\widetilde{\pi}(a)$ be the bounded operators on $\widetilde{\clh}$ defined by
 \begin{align*}
 \lambda_s\xi(t)&=\xi(s^{-1}t) \\
 \widetilde{\pi}(a)\xi(t)&=\pi(\alpha_{t}^{-1}(a))\xi(t)
  \end{align*}
  Then, $(\widetilde{\pi},\lambda)$ is a covariant representation of $(A,G,\alpha)$. 
  
  \begin{xrcs}
  Show that $(\widetilde{\pi},\lambda)$ is non-degenerate. 
    \end{xrcs}
  
  \begin{ppsn}
  The representation $\widetilde{\pi} \rtimes \lambda$ is a faithful representation of $C_{c}(G,A)$. 
    \end{ppsn}
  \textit{Proof.} Let $f \in C_{c}(G,A)$ be such that $(\widetilde{\pi} \rtimes \lambda)(f)=0$. For $\xi, \eta \in C_{c}(G)$ and $u,v \in \clh$, let $\xi_0(s)=\xi(s)u$ and $\eta_0(s)=\eta(s)v$.  Note that $\xi_0,\eta_0 \in \widetilde{\clh}$. 
 Set \[K(s,t)=\langle \pi(\alpha_{t}^{-1}(f(s)))u|v \rangle.\] Calculate as follows to observe that 
 \begin{align*}
 0 & =\langle (\widetilde{\pi} \rtimes \lambda)(f)\xi_0|\eta_0 \rangle\\
 &=\int \langle \widetilde{\pi}(f(s))\lambda_s\xi_0|\eta_0 \rangle ds \\
 &=\int \Big ( \int \langle (\widetilde{\pi}(f(s))\lambda_s\xi_0)(t)|\eta_0(t) \rangle dt\Big) ds \\
 &=\int \Big (\int \langle \pi(\alpha_t^{-1}(f(s)))\xi(s^{-1}t)u|\eta(t)v \rangle dt \Big) ds \\
 &= \int \overline{\eta(t)} \Big(\int \xi(s^{-1}t) K(s,t)ds\Big)dt.
  \end{align*}
  Since $\eta$ is arbitrary and $t \to \int \xi(s^{-1}t)K(s,t)ds$ is continuous, it follows that for every $t$, 
  \[
  \int \xi(s^{-1}t)K(s,t)ds=0.\]
  The arbitrariness of $\xi$ implies $K(s,t)=0$ for every $s,t$. Hence $K(s,e)=0$. This implies that $\langle \pi(f(s))u|v \rangle=0$ for every $s$, $u,v \in \clh$. But $\pi$ is faithful and hence $f(s)=0$ for every $s \in G$. 
  This implies that $f=0$ and the proof is complete. \hfill $\Box$
  
  The above proposition allows us to define the reduced $C^{*}$-norm on $C_{c}(G,A)$. For $f \in C_{c}(G,A)$, let 
  \[
  ||f||_{red}=||(\widetilde{\pi} \rtimes \lambda)(f)||.\]
  The faithfulness of $\widetilde{\pi} \rtimes \lambda$ implies that $||~||_{red}$ is a $C^{*}$-norm on $C_{c}(G,A)$. Moreover, for $f \in C_{c}(G,A)$, the reduced $C^{*}$-norm of $f$ is atmost the full $C^{*}$-norm of $f$, i.e. 
  \[
  ||f||_{red} \leq ||f||.\]
  The completion of $C_{c}(G,A)$ with respect to the norm $||~||_{red}$ is called the \textit{reduced crossed product} and is denoted $A \rtimes_{red} G$. There is a natural surjection from $A \rtimes G$ onto $A \rtimes_{red} G$ which need not be one-one, if we don't assume `amenability' hypothesis. 
  
  A priori, it looks as if the reduced $C^{*}$-norm depends on the chosen faithful representation $\pi$. But it is in fact independent of the chosen representation. The proof of this requires us to digress into the theory of Hilbert C*-modules which we undertake next. 
  
  \section{Hilbert $C^{*}$-modules}
  Hilbert $C^{*}$-modules are analogues of Hilbert spaces, where the inner product takes values in a $C^{*}$-algebra. Rieffel in his seminar paper \cite{Rieffel}  succesfully demonstrated the use of Hilbert $C^{*}$-modules to understand the imprimitivity theorems of Mackey. Kasparov's development of KK-theory utilises Hilbert C*-modules in an essential way and is now an indispensable tool in several areas of operator algebras. For more on K or KK-theory, see \cite{Bla} and \cite{Jensen-Thomsen}.

  Let $B$ be a $C^{*}$-algebra. Suppose $E$ is a vector space. We say that $E$ is a right $B$-module, if $E$ has a right $B$-action that satisfy the usual consistency conditions. 
  
  \begin{dfn}
  Let $E$ be a right $B$-module. By a $B$-valued inner product on $E$, we mean a map $\langle ~|~ \rangle: E \times E \to B$ such that 
  \begin{enumerate}
  \item[(1)] $\langle~|~\rangle$ is linear in the second variable and conjugate linear in the first variable, 
  \item[(2)] for $b \in B$ and $x,y \in E$, $\langle x|yb \rangle=\langle x|y \rangle b$,
  \item[(3)] for $x,y \in E$, $\langle x|y \rangle ^{*}=\langle y|x \rangle$, 
  \item[(4)] for $x \in E$, $\langle x|x \rangle $ is a positive element of $B$, and
  \item[(5)] if $\langle x|x \rangle=0$, then $x=0$. 
   \end{enumerate}
   \end{dfn}
  
  Let $E$ be a right $B$-module with a $B$-valued inner product. For $x \in E$, set \[||x||:=||\langle x| x\rangle||^{\frac{1}{2}}.\] 
  
  \begin{ppsn}[Cauchy-Schwarz inequality]
  For $x,y \in E$, \[||\langle x|y \rangle|| \leq ||x|| ||y||.\] 
    \end{ppsn}
    \textit{Proof.} Let $\rho$ be a state on $B$. The map \[E \times E \ni (e,f) \to \rho(\langle e|f \rangle) \in \mathbb{C}\] defines a semi-definite inner product on $E$. Applying the usual Cauchy-Schwarz inequality to this semi-definite inner product and  by taking $e=x\langle x|y\rangle$ and $f=y$, we see
    \begin{align*}
    \rho(\langle x|y \rangle ^{*} \langle x|y \rangle) & = \rho(\langle e|f \rangle) \\
    & \leq \rho(\langle e|e \rangle)^{\frac{1}{2}}\rho(\langle f|f \rangle)^{\frac{1}{2}} \\
    & \leq \rho(\langle x|y \rangle^{*} \langle x|x \rangle \langle x|y \rangle)^{\frac{1}{2}}\rho(\langle y|y\rangle)^{\frac{1}{2}}\\
    & \leq ||x|| ||y|| \rho(\langle x|y \rangle^{*}\langle x|y \rangle)^{\frac{1}{2}}.
        \end{align*}
  On simplification, we get $\rho(\langle x|y \rangle^{*} \langle x|y \rangle) \leq ||x||^{2}||y||^{2}$ for every state $\rho$. But, for a positive element $a \in B$, \[||a||=\sup\{\rho(a): \textrm{$\rho$ is a state of $B$}\}.\]
  Therefore $||\langle x|y \rangle||^{2}=||\langle x|y \rangle^{*} \langle x|y \rangle|| \leq ||x||^{2}||y||^{2}$. Taking square roots, we have \[||\langle x| y \rangle|| \leq ||x||||y||.\] The proof is complete. \hfill $\Box$      
  
  \begin{xrcs}
  \label{continuity of right action}
  Show that for $x \in E$ and $b \in B$, $||xb|| \leq ||x||||b||$. 
    \end{xrcs}
  Once we have the Cauchy-Schwarz inequality, it is proved as in the Hilbert space setting that $||~||$ defines a norm on $E$. 
  
  \begin{dfn}
  Let $E$ be a right $B$-module with a $B$-valued inner product. We say that $E$ is a Hilbert $B$-module if $E$ is complete with respect to the norm $||~||$, where for $x \in E$, 
  \[
  ||x||=||\langle x| x \rangle_{B}||^{\frac{1}{2}}.\]
   \end{dfn}
  
  \begin{xmpl}
  Hilbert $\mathbb{C}$-modules are just Hilbert spaces. The only difference is that now the inner product is linear in the second variable as opposed to our usual convention. 
    \end{xmpl}
  
  \begin{xmpl}
  Let $B$ be a $C^{*}$-algebra and let $E:=B$. The right multiplication by $B$ makes $E$ into a right $B$-module. For $x,y \in E$, define $\langle x|y \rangle=x^{*}y$. Then $E$ is a Hilbert $B$-module. The norm on $E$ induced by the inner product coincides with the $C^*$-norm on $B$. 
  \end{xmpl}
  
  \begin{xmpl}
  Let $B$ be a $C^{*}$-algebra. Set 
  \[
  H_B:=\{(b_1,b_2,\cdots,): \sum_{n=1}^{\infty}b_n^{*}b_n \textrm{~converges in $B$} \}.\]
  The $C^{*}$-algebra $B$ acts on the right by coordinatewise multiplication. 
  
  For $b:=(b_1,b_2,b_3,\cdots,)$ and $c:=(c_1,c_2,c_3,\cdots,)$, set 
  \[
  \langle b|c \rangle:=\sum_{n=1}^{\infty}b_n^{*}c_n.\]
  Then $H_B$ is a Hilbert $B$-module. 
  \end{xmpl}
  
  \begin{dfn}
  Let $E_1$ and $E_2$ be Hilbert $B$-modules. Suppose $T:E_1 \to E_2$ is a map. We say that $T$ is adjointable with adjoint $T^{*}$ if there exists a map (which is necessarily unique) $T^{*}:E_2 \to E_1$ such that 
  \[
  \langle Tx|y \rangle = \langle x|T^{*}y\rangle\]
  for $x \in E_1$ and $y \in E_2$. 
  The set of adjointable operators from $E_1$ to $E_2$ is denoted by $\mathcal{L}_{B}(E_1,E_2)$. When $E_1=E_2=E$, we write $\mathcal{L}_{B}(E,E)$ as $\mathcal{L}_{B}(E)$. 
    \end{dfn}
  
  \begin{xrcs}
  Let $T:E_1 \to E_2$ be an adjointable operator. Show that 
  \begin{enumerate}
  \item[(1)] $T$ is $\mathbb{C}$-linear, 
  \item[(2)] the map $T$ is $B$-linear, and
  \item[(3)] the map $T:E_1 \to E_2$ is bounded. 
    \end{enumerate}
  Show that $\mathcal{L}_B(E_1,E_2)$ is a norm closed subspace of $B(E_1,E_2)$. Here, $B(E_1,E_2)$ denotes the set of bounded operators from $E_1$ to $E_2$.  
  \end{xrcs}
  \textit{Hint:} To prove $(3)$, use the adjoint equation and the closed graph theorem. For other assertions, use the adjoint equation. 
  
  \begin{ppsn}
  Let $E$ be a Hilbert $B$-module. Then $\mathcal{L}_{B}(E)$ is a $C^{*}$-algebra. 
    \end{ppsn}
  \textit{Proof.} The only thing that requires proof is that the operator norm satisfies the $C^{*}$-identity. First, note that by the Cauchy-Schwarz inequality, we have for $x \in E$, 
  \[
  ||x||=\sup\{||\langle x|y \rangle||: ||y||=1\}.\]
  
  For $T \in \mathcal{L}_{B}(E)$, $||T||=\sup\{||\langle Tx|y \rangle||: ||x||=1=||y||\}$. Thus, it is clear that $||T^{*}||=||T||$. Since the operator norm is submultiplicative, it follows that for $T \in \mathcal{L}_{B}(E)$, 
 $ ||T^{*}T|| \leq ||T||^{2}$. 
 
 Let $T \in \mathcal{L}_{B}(E)$ be given. Then 
 \begin{align*}
 ||T||^{2}&=\sup\{||Tx||^{2}: ||x||=1\} \\
             &=\sup\{||\langle Tx|Tx \rangle||: ||x||=1\} \\
             &=\sup\{||\langle T^{*}Tx|x \rangle ||: ||x||=1\} \\
             & \leq ||T^{*}T|| ~(\textrm{ by Cauchy-Schwarz inequality}).
  \end{align*}
  The proof is now complete. \hfill $\Box$
  
  Unlike in the case of Hilbert spaces, it is not true that bounded operators between Hilbert modules are adjointable. Here is an example. 
  \begin{xmpl}
  Let $B:=C[0,1]$ and $J:=\{f \in B: f(0)=0\}$. Let $E_1=J$ and $E_2=B$. Both $E_1$ and $E_2$ are Hilbert $B$-modules. Consider the inclusion $T: E_1 \to E_2$, i.e. $T(x)=x$.
  Then ,$T$ is not adjointable. Suppose not and let $S$ be the adjoint of $T$. Let $h:=S(1)$. Then $h \in J$. Calculate as follows to observe that for $f \in J$,
  \begin{align*}
  \overline{f}&=\langle T(f)|1 \rangle \\
  &=\langle f|S(1) \rangle \\
  &=\overline{f}h.
    \end{align*}
    In other words, $h$ is a multiplicative identity of the non-unital $C^{*}$-algebra $J$ which is absurd. 
       \end{xmpl}
       
  \begin{rmrk}
  One needs to exercise caution while dealing with Hilbert modules. For example, it is not true that if $F$ is a submodule of $E$, then $E=F \oplus F^{\perp}$. Can you construct a counterexample ? 
      \end{rmrk}
      \textit{Hint:} What about the last example ?

      In practice, it is essential to complete right $B$-modules to get genuine Hilbert modules. The following proposition helps in achieving this. 
       \begin{ppsn}
      Let $B_0$ be a dense $^*$-subalgebra of a $C^{*}$-algebra $B$. Suppose $E_0$ is a right $B_0$-module with a $B_0$-valued inner product. 
      $E_0$ is usually called a pre-Hilbert $B_0$-module. 
      Denote the completion of $E_0$ by $E$. Then, the $B_0$-module structure on $E_0$ lifts uniquely to make $E$ into a Hilbert $B$-module. 
            \end{ppsn}
            \textit{Proof.} The proof is routine and makes essential use of Exercise \ref{continuity of right action}. \hfill $\Box$

  Let us construct the Hilbert module of interest associated to a $C^{*}$-dynamical system. Suppose $(A,G,\alpha)$ is a $C^{*}$-dynamical system. 
  Let $E_0:=C_{c}(G,A)$.  We make $E_0$ into a right $A$-module as follows. For $f \in E_0$ and $a \in A$, define
  \[
  (f.a)(s):=f(s)\alpha_s(a).\]
  The $A$-valued inner product on $E_0$ is given by 
  \[
  \langle f|g \rangle_{A}:=(f^{*}*g)(e)=\int \Delta(s^{-1})\alpha_s(f(s^{-1})^*)\alpha_s(g(s^{-1}))ds.\]
  Verify that $E_0$ is a pre-Hilbert $A$-module. We obtain a genuine Hilbert $A$-module upon completion which we denote by $E$. 
  
  For $a \in A$, let $i_{A}(a):E_0 \to E_0$ be defined by \[i_{A}(a)f(s)=af(s).\] For $s \in G$, let $i_{G}(s):E_0 \to E_0$ be defined by 
  \[
  i_{G}(s)f(t)=\alpha_{s}(f(s^{-1}t)).\]
  
  \begin{xrcs}
  \label{reduced covariant}
  For $a \in A$ and $s \in G$, show that $i_{A}(a)$ and $i_{G}(s)$ extend to bounded operator on $E$. We denote the extensions by the same symbols. 
  Verify that $i_{A}(a)^{*}=i_{A}(a)$ and $i_{G}(s)^{*}=i_{G}(s^{-1})$.  Prove that 
  \begin{enumerate}
  \item[(1)] the map $i_{A}: A \to \mathcal{L}_{A}(E)$ is a non-degenerate $^*$-representation, 
  \item[(2)] the map $s \to i_{G}(s)$ is a strongly continuous unitary representation of $G$ on $E$, and
  \item[(3)] for $s \in G$ and $a \in A$, the covariance condition $i_{G}(s)i_{A}(a)i_{G}(s)^{*}=i_{A}(\alpha_s(a))$ is satisfied. 
    \end{enumerate}
  In short, the pair $(i_A,i_G)$ is a covariant representation of $(A,G,\alpha)$ on the Hilbert $A$-module $E$. 
    \end{xrcs}
  
  Just like we can integrate covariant representations on a Hilbert space to obtain a representation of the crossed product on the same Hilbert space, we could do the same with Hilbert modules. 
  To do so, we need to discuss vector valued integration one more time. 
  
  Suppose $E$ is a Hilbert $B$-module. For $x \in E$, define a seminorm $||~||_{x}$ on $\mathcal{L}_{B}(E)$ by \[||T||_{x}=||Tx||+||T^*x||.\] The topology on $\mathcal{L}_{B}(E)$ induced by the family of seminorms $\{||~||_{x}: x \in E\}$ is called the topology of $^*$-strong convergence on $\mathcal{L}_{B}(E)$. Let $(T_i)$ be a net in $\mathcal{L}_{B}(E)$ and $T \in \mathcal{L}_{B}(E)$ be given. Then $T_i \to T$ in the topology of $^*$-strong convergence if and only if for every $x \in X$, $T_{i}x \to Tx$ and $T_{i}^{*}x \to T^{*}x$.

  \begin{ppsn}
  Let $X$ be a locally compact second countable Hausdorff topological space and $\mu$ be a Radon measure on $X$. Let $f:X \to \mathcal{L}_{B}(E)$ be continuous, when $\mathcal{L}_{B}(E)$ is given the topology of $^*$-strong convergence. 
  Suppose that the map $X \ni x \to ||f(x)|| \in [0,\infty)$ is integrable. Then, there exists a unique adjointable operator on $E$ denoted $\int f(x) d\mu(x)$ such that for $u,v \in E$, 
  \[
  \langle  \Big(\int f(x)d\mu(x)\Big)u|v \rangle = \int \langle f(x)u|v \rangle d\mu(x).\]
      \end{ppsn}
  \textit{Proof.} Fix $u \in E$. The map $X \ni x \to f(x)u \in E$ is continuous and integrable. Define 
  \[
  \Big(\int f(x)d\mu(x)\Big)u:=\int f(x)u d\mu(x).\]
  The assertion follows from the defining properties of vector valued integration. \hfill $\Box$
  
  \begin{rmrk}
  Theorem \ref{correspondence for crossed products} stays true with Hilbert spaces replaced by Hilbert modules. 
    \end{rmrk}
  
  \textbf{Interior tensor product:} One of the most important construction in Hilbert modules is the notion of interior tensor product which we discuss next. Let $E$ be a Hilbert $B$-module and let $F$ be a Hilbert $C$-module. 
  Suppose $\pi:B \to \mathcal{L}_{C}(F)$ is a $^*$-homomorphism. Think of $F$ as a left $B$-module and $E$ as a right $B$-module. Consider the algebraic tensor product $X:=E \otimes_B F$. Note that $C$ acts naturally on the right  
  on $X$. Moreover, in $X$, we have the equality $eb \otimes f= e \otimes \pi(b)f$ for $e \in E$, $f \in F$ and $b \in B$. Define a $C$-valued semi-definite inner product on $X$ by the formula
  \begin{equation}
  \label{interior inner product}
  \langle e_1 \otimes f_1|e_2 \otimes f_2 \rangle = \langle f_1|\pi(\langle e_1|e_2))f_2 \rangle.\end{equation}
  The fact that the inner product is positive requires a bit of work. 
  
  \begin{lmma}
  \label{positivity in Hilbert modules}
  Let $E$ be a Hilbert $B$-module. Suppose $T \in \mathcal{L}_{B}(E)$. The following are equivalent.
  \begin{enumerate}
  \item[(1)] $T$ is a positive element of $\mathcal{L}_{B}(E)$.
  \item[(2)] For every $x \in E$, $\langle Tx|x \rangle \geq 0$. 
    \end{enumerate}
   \end{lmma}
  \textit{Proof.} Suppose $(1)$ holds. Write $T=S^{*}S$ with $S \in \mathcal{L}_{B}(E)$. Then for $x \in E$, \[\langle Tx|x \rangle = \langle S^{*}Sx|x \rangle = \langle Sx|Sx \rangle \geq 0.\]
  This proves that $(1)$ implies $(2)$.  
  
  Now suppose that $(2)$ holds. Note that for $x \in E$, $\langle Tx|x \rangle=\langle Tx|x \rangle^{*}=\langle x|Tx \rangle$. 
  For $x,y \in E$, let $[x,y]=\langle Tx|y \rangle$ and $[x,y]'=\langle x|Ty \rangle$. Both $[~,~]$ and $[~,~]^{'}$ are sesquilinear $B$-valued forms and $[x,x]=[x,x]^{'}$ for $x \in E$. By the polarisation
  identity, it follows that the forms $[~,~]$ and $[~,~]'$ agree. Consequently $\langle Tx|y \rangle = \langle x|Ty \rangle$. In other words, $T=T^{*}$. 
  
  Write $T=R-S$ with $R,S \geq 0$ and $SR=RS=0$. For $x \in E$, calculate as follows to observe that 
  \begin{align*}
  0 & \leq \langle TSx| Sx \rangle \\
  & \leq -\langle S^{3}x|x \rangle\\
  & \leq 0 (~\textrm{as $S^{3}$ is positive }).
    \end{align*}
    Hence $\langle S^{3}x|x \rangle=0$ for every $x \in E$. By the polarisation identity, it follows that for $x,y \in E$, $\langle S^{3}x|y \rangle=0$. Hence, $S^{3}=0$ which forces $S=0$.
    Thus $T=R \geq 0$. This completes the proof. \hfill $\Box$
    
    \begin{xrcs}
    \label{norm of positive operator}
    Let $T \in \mathcal{L}_{B}(E)$ be such that $T \geq 0$. Prove that 
    \[
    ||T||=\sup\{ ||\langle x|Tx \rangle||: ||x||=1\}.\]
    
    \end{xrcs}
  
  \begin{lmma}
  Let $E$ be a Hilbert $B$-module and $e_1,e_2,\cdots,e_n \in E$ be given. Then the matrix $(\langle e_i|e_j \rangle)$ is a positive element of $M_n(B)$. 
    \end{lmma}
    \textit{Proof.} Consider the Hilbert $B$-module $B^{n}:=B \oplus B \oplus \cdots \oplus B$. Think of elements of $B^{n}$ as column vectors. For $A \in M_{n}(B)$, let 
    $L_{A}: B^{n} \to B^{n}$ be defined by 
    \[
    L_{A}(x)=Ax.\]
    Then $L_{A}$ is an adjointable operator on $B^{n}$ and the map $M_{n}(B) \ni A \to L_{A} \in \mathcal{L}_{B}(B^{n})$ is an injective $^*$-homomorphism (Justify!). 
    
    Let $A:=(\langle e_i|e_j \rangle)$. It suffices to show that $L_{A}$ is positive. In view of Lemma \ref{positivity in Hilbert modules}, it suffices to show that 
    $\langle x|L_{A}(x) \rangle \geq 0$. Let $x:=(b_1,b_2,\cdots,b_n)^{t}$ be given. 
    Note that 
    \[
    \langle x|L_{A}(x) \rangle=\sum_{i=1}^{n}b_{i}^{*}(\sum_{j=1}^{n}\langle e_i|e_j \rangle b_j)=\sum_{i,j}\langle e_{i}b_i|e_jb_j\rangle=\langle \sum_ie_ib_i|\sum_i e_ib_i \rangle \geq 0.\]
    This completes the proof. \hfill $\Box$
    
    We now prove that Eq. \ref{interior inner product} defines a  semi-definite inner product. Keep the notation used in the paragraph preceding Eq. \ref{interior inner product}.
    Let $x:=\sum_{i=1}^{n}e_i \otimes f_i$ be an arbitrary element in $X$. The representation $\pi$ ``amplifies naturally" to a representation of $M_{n}(B)$ on the Hilbert $C$-module $F^{n}:=F \oplus F \oplus \cdots \oplus F$. 
  Since $(\langle e_i|e_j \rangle)$ is a positive element in $M_n(B)$, the operator $T:=(\pi(\langle e_i|e_j \rangle))$ is a positive operator on $F^{n}$. 
  Set $f:=(f_1,f_2,\cdots,f_n)^{t}$. Then,
  \[
  \langle x| x \rangle= \langle f|Tf \rangle \geq 0.\]
  Thus, $\langle ~|~\rangle$ defines a positive semi-definite $C$-valued inner product on $X$. We mod out the null vectors and complete it to obtain a genuine Hilbert $C$-module. We denote the resulting $C$-module by
  $E \otimes_{\pi} F$. The module $E \otimes_{\pi} F$ is called the interior tensor product, or the internal tensor product of $E$ and $F$.

  \begin{ppsn}
  \label{towards Rieffel induction}
  Keep the foregoing notation. Suppose $T \in \mathcal{L}_{B}(E)$. Then, there exists a unique adjointable operator, denoted $T \otimes 1$, on $E \otimes_\pi F$ such that 
  \[
  (T \otimes 1)(e \otimes f)=Te \otimes f\]
  for $e \in E$ and $f \in F$. 
    \end{ppsn}
  \textit{Proof.} Let $x:=\sum_{i=1}^{n}e_i \otimes f_i$ be given. We claim that 
  \begin{equation}
  \label{defining inequality}
  \sum_{i,j} \langle f_i|\pi(\langle Te_i|Te_j \rangle)f_j \rangle \leq ||T||^{2} \sum_{i,j}\langle f_i|\pi(\langle e_i|e_j \rangle)f_j\rangle.
    \end{equation}
  We leave it to the reader to convince herself that once Eq. \ref{defining inequality} is established, the conclusion follows. Argue as in the previous lemma, with the aid of the next exercise, that the matrix $(\langle Te_i|Te_j \rangle) \leq ||T||^{2}(\langle e_i|e_j \rangle)$.
  
   Let $A:=(\pi(\langle Te_i|Te_j \rangle))$ and $B:=||T||^{2}(\pi(\langle e_i |e_j\rangle))$. Then $A$ and $B$ are adjointable operators on $F^{n}$ and $A \leq B$. Set $f:=(f_1,f_2,\cdots,f_n)^{t}$.  Inequality \ref{defining inequality} follows from the  fact that $\langle f|Af \rangle \leq \langle f|Bf \rangle$. This completes the proof. \hfill $\Box$
  
  \begin{xrcs}
  Let $E$ be a Hilbert $B$-module and $T \in \mathcal{L}_{B}(E)$ be given. Prove that for $x \in E$, $\langle Tx|Tx \rangle \leq ||T||^{2} \langle x|x \rangle$. 
    \end{xrcs}
    \textit{Hint:} Write $||T||^{2}-T^{*}T$ as $S^{*}S$ for some $S \in \mathcal{L}_{B}(E)$.

  \begin{rmrk}
  We have the following. 
  \begin{enumerate}
  \item[(1)] The map $\mathcal{L}_{B}(E) \ni T \to T \otimes 1 \in \mathcal{L}_{B}(E \otimes_\pi F)$ is a $^*$-homomorphism. If $\pi$ is injective, then the map $T \to T \otimes 1$ is injective. 
  \item[(2)] Suppose $(T_i)$ is a bounded net which converges to $T$ in the $^*$-strong topology. Then $T_i \otimes 1 \to T\otimes 1$ in the $^*$-strong topology. 
      \end{enumerate}
    \end{rmrk}
  
  Proposition \ref{towards Rieffel induction} leads us to a very important notion of induced representations due to Rieffel (\cite{Rieffel}). The data given is as follows. Suppose $E$ is a Hilbert $B$-module and  let $\phi:A \to \mathcal{L}_{B}(E)$ be a representation. 
  $E$ is usually called a Hilbert $A$-$B$ bimodule. Suppose $\pi$ is a representation of $B$ on a Hilbert space $\clh$. Consider the Hilbert space $\clh_{\pi}:=E \otimes_\pi \clh$. For $a \in A$, define $Ind(\pi)(a)$ on $\clh_\pi$ by 
  \[
  Ind(\pi)(a)=\phi(a) \otimes 1.\]
  Then $Ind(\pi)$ is a representation of $A$ on the Hilbert space $\clh_\pi$. The representation $Ind(\pi)$ is called the representation induced by $\pi$ via the bimodule $E$. 
  
  Suppose $\pi_1$ and $\pi_2$ are representations of $B$ on the Hilbert spaces $\clh_1$ and $\clh_2$. Suppose $T:\clh_1 \to \clh_2$ is an intertwiner, i.e. $T\pi_1(b)=\pi_2(b)T$. Because $T$ commutes with the action of $B$, the map $1 \otimes T$ is well-defined first on the algebraic level and then extends to give a genuine adjointable operator from $E \otimes_{\pi_1} \clh_1 \to E \otimes_{\pi_2} \clh_2$. It is clear that $1 \otimes T$ commutes with the left action of $A$. Or in other words, $1\otimes T$ intertwines $Ind(\pi_1)$ and $Ind(\pi_2)$. 
    
  Summarising the above discussion, we observe that 
    \[
  \pi \to Ind(\pi)\]
  is a functor from the category of representations of $B$ to the category of representations of $A$. 
  
  \begin{rmrk}
  \label{primitive Morita equivalence}
  We can try to seek an inverse for the above functor which leads us to the notion of Morita equivalence which we will discuss later in the course. The idea is that we can define an inverse if there is a $B$-$A$ Hilbert bimodule $F$ such that  $E \otimes_B F \cong A$ and $F \otimes_A E \cong B$ as bimodules.   If such bimodules exist, then we say that $A$ and $B$ are Morita equivalent in the sense of Rieffel (see the  Section 10). Rieffel's induction then ensures that the representation theory of $A$ and that of $B$ are the same if $A$ and $B$ are Morita-equivalent.     As a warm-up exercise, the reader can try to prove that for a $C^{*}$-algebra $B$, $B$ and $M_n(B)$ are Morita-equivalent. An obvious $M_n(B)$-$B$-bimodule should stare at our face.     
     
     The first remarkable result due to Rieffel is that   $C_{0}(G/H) \rtimes G$ is Morita equivalent to $C^{*}(H)$, where $G$ is a locally compact group and $H$ is a closed subgroup of $G$. 
     Note that $G$ acts by left translations on the space of cosets $G/H$. 
    \end{rmrk}
    
    Let us return back to our original motivation for considering Hilbert $C^{*}$-modules which is to establish that the reduced $C^{*}$-norm does not depend on the choice of the faithful representation. 
    Let $(A,G,\alpha)$ be a $C^{*}$-dynamical system. Let $E_0:=C_{c}(G,A)$ be the pre-Hilbert $A$-module constructed earlier. The inner product and the right action of $A$ are given by
    \begin{align*}
    \langle f|g \rangle_{A}:&=\int \Delta(t)^{-1}\alpha_{t}(f(t^{-1})^{*})\alpha_t(g(t^{-1}))dt \\
    (f.a)(s)&:=f(s)\alpha_s(a)
       \end{align*}
    for $f,g \in E_0$ and $a \in A$. Denote the completion of $E_0$ by $E$. 
    
    Let $(i_A,i_G)$ be the covariant representation of $(A,G,\alpha)$ on $E$ as in Exercise \ref{reduced covariant}. Recall that the formulas for $i_{A}$ and $i_{G}$ are given by 
    \begin{align*}
    i_{A}(a)f(t)&=af(t) \\
    i_{G}(s)f(t)&:=\alpha_{s}(f(s^{-1}t))
        \end{align*}
    for $f \in E_0$, $a \in A$ and $s \in G$. 
    
    Fix a faithful non-degenerate representation $\pi$ of $A$ on $\clh$. Set $\widetilde{\clh}:=L^{2}(G,\clh)$. For $a \in A$ and $s \in G$, let $\widetilde{\pi}$ and $\lambda_s$ be the bounded operators on $\widetilde{\clh}$ defined by
    \begin{align*}
 \lambda_s\xi(t)&=\xi(s^{-1}t) \\
 \widetilde{\pi}(a)\xi(t)&=\pi(\alpha_{t}^{-1}(a))\xi(t).
  \end{align*}
  For $f \in C_{c}(G,A)$, the reduced $C^{*}$-norm is $||(\widetilde{\pi} \rtimes \lambda)(f)||$. We prove that \[||f||_{red}=||(i_A \rtimes i_{G})(f)||\] and the right side does not depend on the representation $\pi$. 

 The trick is to use Rieffel's induction. Note that $\pi$ is a representation of $A$ on $\clh$. Thus, we can form the interior tensor product $E \otimes_{\pi} \clh$ which is a Hilbert space and we show that the latter 
 Hilbert space is identified with $\widetilde{\clh}$ by a specific unitary. Once this identification is made, then we show $i_{A}(a) \otimes 1=\widetilde{\pi}(a)$ and $i_{G}(s) \otimes 1=\lambda_s$.  On integration, we 
 get for $f \in C_{c}(G,A)$, $(i_A \rtimes i_{G})(f) \otimes 1=(\widetilde{\pi} \rtimes \lambda)(f)$.  Since $\pi$ is faithful, the map $\mathcal{L}_{B}(E) \ni T \to T \otimes 1 \in B(\widetilde{\clh})$ is $1$-$1$, and hence preserves the
 norm. Consequently, $||f||_{red}=||(i_A \rtimes i_G)(f)||$. The verification is carried out in the next exercise. 
  
  \begin{xrcs}
  Assume that $G$ is discrete. Show that there exists a unique unitary operator $U:E \otimes_{\pi} \clh \to \ell^{2}(G) \otimes \clh$ such that 
  \[
  U((a \otimes \delta_t) \otimes \xi)=\delta_t \otimes \pi(\alpha_{t}^{-1}(a))\xi.\]
Prove that $U(i_{A}(a) \otimes 1)U^{*}=\widetilde{\pi}(a)$ and $U(i_{G}(s) \otimes 1)U^{*}=\lambda_s$ for $a \in A$ and $s \in G$.  Conclude that for $f \in C_{c}(G,A)$, 
\[||f||_{red}=||(i_A \rtimes i_G)(f)||.\] 
Treat the topological case similarly. 
    \end{xrcs}
  
  \section{Stone-von Neumann theorem}
  As an application of the material developed so far, we determine in this section, the irreducible representations of the Heisenberg group. The determination of the irreducible representations of the Heisenberg group goes under the name ``Stone-von Neumann theorem". For $n \geq 1$,  let \[H_{2n+1}:=\mathbb{R}^{n} \times \mathbb{R}^{n} \times \mathbb{R}.\] 
 The group law on $H_{2n+1}$ is defined by
 \[
 (x_1,y_1,z_1)(x_2,y_2,z_2):=(x_1+x_2,y_1+y_2,z_1+z_2+\langle x_1|y_2 \rangle).\] 
  Verify that $H_{2n+1}$ together with the group law defined above is indeed a topological group. 
  
  \begin{xrcs}
  Let $Z$ be the center of $H_{2n+1}$. Show that $Z=\{(0,0,z): z \in \mathbb{R}\}$. 
   \end{xrcs}
  
  The reader who has already seen Lemma \ref{Sch} should not have difficulty with the following proposition. 
  \begin{ppsn}[Schur's lemma]
  Let $G$ be a locally compact, Hausdorff, second countable topological group. Suppose $\pi: G \to \mathcal{U}(\clh)$ is a strongly continuous unitary representation. Then, $\pi$ is irreducible if and only if  the commutant $\pi(G)^{'}=\mathbb{C}$. 
    \end{ppsn}

 \begin{ppsn}
 \label{irrep of Heisenberg}
 Suppose $\pi:H_{2n+1} \to \mathcal{U}(\clh)$ be a strongly continuous unitary representation. Assume that $\pi$ is irreducible. Then, there exists a unique $\lambda \in \mathbb{R}$ such that for $z \in \mathbb{R}$,
 \[\pi(0,0,z)=e^{i\lambda z}.\] 
 For $x \in \mathbb{R}^{n}$ and $y \in \mathbb{R}^{n}$, set $U_x:=\pi(x,0,0)$ and $V_y=\pi(0,y,0)$.
 Then, $\{U_x\}_{x \in \mathbb{R}^{n}}$ and $\{V_y\}_{y \in \mathbb{R}^{n}}$ are strongly continuous group of unitaries such that 
  \begin{equation}
 \label{Weyl}
 U_xV_y=e^{i\lambda \langle x|y \rangle}V_yU_x
  \end{equation}
 for $x \in \mathbb{R}^{n}$ and $y \in \mathbb{R}^{n}$. Moreover the commutant  $\{U_x,V_y: x ,y \in \mathbb{R}^{n}\}^{'}=\mathbb{C}$.  
 
 Conversely, suppose $\{U_x\}_{x \in \mathbb{R}^{n}}$ and $\{V_y\}_{y \in \mathbb{R}^{n}}$ are strongly continuous unitary representations on a Hilbert space $\clh$ which satisfy Eq. \ref{Weyl}. Assume in addition that $\{U_x,V_y: x,y \in \mathbb{R}^{n}\}'=\mathbb{C}$. For $(x,y,z) \in H_{2n+1}$, set 
 \[
 \pi(x,y,z)=e^{i\lambda z}V_yU_x.\]
Then, $\pi$ defines an irreducible representation of $H_{2n+1}$. 
 
 \end{ppsn}
    \textit{Proof.} Note that $\pi(0,0,z)$ commutes with $\pi(H_{2n+1})$ for every $z \in \mathbb{R}^{n}$. By Schur's lemma, it follows that $\pi(0,0,z) \in \mathbb{T}$. Since  $\pi$ is a strongly continuous representation, it follows that the map  $\mathbb{R} \ni z \to \pi(0,0,z) \in \mathbb{T}$ is a continuous group homomorphism. Hence, there exists a unique $\lambda \in \mathbb{R}$ such that $\pi(0,0,z)=e^{i \lambda z}$. 
    
    Note that in $H_{2n+1}$, \[(x,0,0)(0,y,0)=(0,y,0)(x,0,0)(0,0,\langle x|y \rangle).\] Hence Eq. \ref{Weyl} is satisfied. Note that $\pi(H_{2n+1})^{'}=\{U_x,V_y:x,y \in \mathbb{R}^{n}\}'$. Since $\pi$ is irreducible, it follows that $\{U_x,V_y: x,y \in \mathbb{R}^{n}\}'=\mathbb{C}$. The proof of the converse part is routine. \hfill $\Box$
    
\begin{rmrk}
The relation \ref{Weyl} when $\lambda=1$ is usually called the ``Weyl commutation" relation.

\end{rmrk}    
\begin{dfn}
Let $\clh$ be a Hilbert space and $U:=\{U_x\}_{x \in \mathbb{R}^{n}}$ and $V:=\{V_y\}_{y \in \mathbb{R}^{n}}$ be strongly continuous group of unitaries. We say that the pair $(U,V)$ is a Weyl family of unitaries with phase factor $\lambda$ if Eq. \ref{Weyl} is satisfied. We call the pair $(U,V)$ irreducible, if the commutant $\{U_x,V_y: x,y \in \mathbb{R}^{n}\}'=\mathbb{C}$. 

\end{dfn}    
In view of Prop. \ref{irrep of Heisenberg}, the  problem of determining the irreducible representations of the Heisenberg group reduces to determine  Weyl families of unitaries which are irreducible. 

\textbf{Case 1: $\lambda=0$} In this case, a Weyl family $(U,V)$ corresponds to two unitary representations of $\mathbb{R}^{n}$ which commute. Equivalently, in this case, a Weyl family corresponds to a unitary 
representation of the cartesian product $\mathbb{R}^{n} \times \mathbb{R}^{n}$. Thus, irreducible Weyl families are precisely the irreducible representations of the abelian group $\mathbb{R}^{n} \times \mathbb{R}^{n}$, or in other words the characters of $\mathbb{R}^{2n}$. We know that the dual of $\mathbb{R}^{2n}$ is $\mathbb{R}^{2n}$. 

Let $\mu \in \mathbb{R}^{2n}$ be given. Write $\mu:=(\mu_1,\mu_2)$ with $\mu_1,\mu_2 \in \mathbb{R}^{n}$. Set $U_x:=e^{i \langle \mu_1|x \rangle}$ and $V_{y}:=e^{i \langle \mu_2|y \rangle}$. 
  Then $(U,V)$ is an irreducible Weyl family of unitaries on the one dimensional Hilbert space $\mathbb{C}$ with phase factor $\lambda=0$.   Up to unitary equivalence, every such Weyl family arises this way. 
    Moreover, for distinct values of $\mu$, the corresponding Weyl families are inequivalent. 
    
    The non-trivial case is when $\lambda \neq 0$. For the rest of  our discussion, $\lambda$ will be a fixed non-zero real number. The first order of business is to exhibit an irreducible Weyl family with phase factor $\lambda$. Let $\clh:=L^{2}(\mathbb{R}^{n})$. 
 For $x,y \in \mathbb{R}^{n}$, let $U_x$ and $V_y$ be the unitary operators on $\clh$ defined by the following equation
 \begin{align*}
 U_xf(t)&:=f(t-x) \\
 V_yf(t)&:=e^{-i \lambda \langle y|t \rangle}f(t)
 \end{align*}
 
 \begin{ppsn}
 \label{irreducible Weyl family}
 Keep the foregoing notation. The pair $(U,V)$ is an irreducible Weyl family with phase factor $\lambda$. 
  \end{ppsn}
 \textit{Proof.} It is routine to verify that $U$ and $V$ are strongly continuous unitary representations of $\mathbb{R}^{n}$ and they satisfy Eq. \ref{Weyl}. Consider $L^{\infty}(\mathbb{R}^{n})$ and let $M:L^{\infty}(\mathbb{R}^{n}) \to B(\clh)$ be the multiplication representation, i.e. for $\phi \in L^{\infty}(\mathbb{R}^{n})$ and $f \in \clh$, 
 \[
 M(\phi)f(t)=\phi(t)f(t).\]
 Note that $L^{\infty}(\mathbb{R}^{n}) \ni \phi \to M(\phi) \in B(\clh)$ is continuous, when $L^{\infty}(\mathbb{R}^{n})$ is given the weak $*$-topology (after identifying $L^{\infty}(\mathbb{R}^{n})$ with the dual of $L^{1}(\mathbb{R}^{n})$) and when $B(\clh)$ is given the weak operator topology. 
 
 \textit{Claim:} The linear span of $\{e^{-i \lambda \langle y|t \rangle}: y \in \mathbb{R}^{n}\}$ is weak $^*$-dense in $L^{\infty}(\mathbb{R}^{n})$. Suppose not. Then, there exists a non-zero $f \in L^{1}(\mathbb{R}^{n})$ such that for $y \in \mathbb{R}^{n}$, 
 \[
 \int f(t)e^{-i\lambda \langle y|t \rangle}dt=0.\]
 In other words, the Fourier transform of $f$ is zero which in turn implies $f=0$. This is a contradiction. This proves our claim. 
 
 Suppose $T \in B(\clh)$ is such that $TU_x=U_xT$ and $TV_y=V_yT$ for $x,y \in \mathbb{R}^{n}$. The density of $\{e^{-i\lambda \langle y|t \rangle}: y \in \mathbb{R}^{n}\}$ in $L^{\infty}(\mathbb{R}^{n})$ implies that $TM(\phi)=M(\phi)T$ for every $\phi \in L^{\infty}(\mathbb{R}^{n})$. It is a good exercise to prove that the commutant of $L^{\infty}(\mathbb{R}^{n})$ is $L^{\infty}(\mathbb{R}^{n})$ (see, for instance, Theorem 2.2.1 of \cite{Arveson_invitation}).
 Hence, there exists $\phi \in L^{\infty}(\mathbb{R}^{n})$ such that $T=M(\phi)$. 
 
 Now the equation $U_xTU_x^{*}=T$ implies that for every $x \in \mathbb{R}^{n}$, $\phi(t+x)=\phi(t)$ for almost all $t \in \mathbb{R}^{n}$. Let $\omega_{\phi}:C_{c}(\mathbb{R}^{n}) \to \mathbb{R}$ be defined by
 \[
 \omega_{\phi}(f)=\int f(t)\phi(t)dt.\]
 To show that $\phi$ is a scalar, it suffices to show that  $\omega_{\phi}$ is a scalar multiple of the linear functional $I:C_{c}(\mathbb{R}^{n}) \to \mathbb{C}$ defined by the equation
 \[
 I(f)=\int f(t)dt.\]
 Let $g \in C_{c}(\mathbb{R}^{n})$ be such that $I(g)=0$. Choose $f \in C_{c}(\mathbb{R}^{n})$ such that $\int f(t)dt=1$. 
 Calculate as follows to observe that 
 \begin{align*}
 \omega_{\phi}(g)&=\int f(s) (\int g(t)\phi(t)dt)ds \\
 &=\int f(s)(\int g(t)\phi(t-s)dt)ds \\
 &= \int g(t) (\int f(s)\phi(t-s)ds)dt \\
 &=\int g(t)(\int f(s)\phi(-s)ds)dt \\
 &=(\int g(t)dt) (\int f(s)\phi(-s)ds)\\
 &=0.
  \end{align*} 
 Hence, $Ker(I) \subset Ker(\omega_\phi)$. This shows that $\omega_\phi$ is a scalar multiple of $I$. Hence $\phi$ is a scalar and consequently $T$ is a scalar. This proves that the commutant $\{U_x,V_y: x,y \in \mathbb{R}^{n}\}'=\mathbb{C}$. The proof is now complete. \hfill $\Box$
 
 \begin{thm}[Stone-von Neumann]
 Let $(\widetilde{U},\widetilde{V})$ be an irreducible Weyl family of unitaries with phase factor $\lambda$ on a Hilbert space $\widetilde{\clh}$. Denote the Weyl family constructed in Proposition \ref{irreducible Weyl family} by $(U,V)$. Then $(\widetilde{U},\widetilde{V})$ is unitarily equivalent to $(U,V)$. This means that there exists a unitary $T:\widetilde{\clh} \to \clh$ such that $T\widetilde{U}_xT^{*}=U_x$ and $T\widetilde{V}_{y}T^{*}=V_{y}$. 
  \end{thm}
  The proof of Stone-von Neumann's theorem relies on the following steps. 
  \begin{enumerate}
  \item[(1)] First we   show that Weyl family of unitaries are in $1$-$1$ correspondence with covariant representations of the dynamical system $(C_{0}(\mathbb{R}^{n}),\mathbb{R}^{n},\alpha)$ where the action $\alpha$ is by translations. Moreover the correspondence respects irreducibility. 
  
  \item[(2)] Thus the problem reduces to the determination of irreducible covariant representations of $(C_0(\mathbb{R}^{n}),\mathbb{R}^{n},\alpha)$ or in other words determining the irreducible representations of the crossed product $C_{0}(\mathbb{R}^{n}) \rtimes \mathbb{R}^{n}$. 
  \item[(3)] Next, we  show that $C_{0}(\mathbb{R}^{n}) \rtimes \mathbb{R}^{n}$ is isomorphic to $\mathcal{K}(L^{2}(\mathbb{R}^{n}))$.\footnote{For the discrete version, see Prop. \ref{discrete stone-von Neumann}.} Since the algebra of compact operators has only one irreducible representation, up to unitary equivalence, the theorem follows. 
  \end{enumerate}
  
  For $y \in \mathbb{R}^{n}$, let $\alpha_{y}:C^{*}(\mathbb{R}^{n}) \to C^{*}(\mathbb{R}^{n})$ be defined by
  \[
  \alpha_{y}f(x)=e^{-i\lambda \langle x|y \rangle}f(x).\]
  It is routine to verify that $\alpha:=\{\alpha_y\}_{ y \in \mathbb{R}^{n}}$ is an action of $\mathbb{R}^{n}$ on $C^{*}(\mathbb{R}^{n})$. Let $(U,V)$ be a Weyl family of unitaries with phase factor $\lambda$. Denote the integrated form of $U$ by $\pi_{U}$. 
  Then $(\pi_{U},V)$ is a covariant representation of the dynamical system $(C^{*}(\mathbb{R}^{n}),\mathbb{R}^{n},\alpha)$. 
  
  Conversely, suppose $(\pi,V)$ is a non-degenerate covariant representation of $(C^{*}(\mathbb{R}^{n}), \mathbb{R}^{n},\alpha)$. Let $U:=\{U_{x}\}_{x \in \mathbb{R}^{n}}$ be the strongly continuous unitary representation of $\mathbb{R}^{n}$ whose integrated form is $\pi$. Note that $(\pi,V)$ is covariant implies that for $y \in \mathbb{R}^{n}$, $f \in C_{c}(\mathbb{R}^{n})$ and vectors $\xi,\eta$, we have 
  \[
  \int f(x)\langle V_yU_xV_y^{*}\xi|\eta \rangle dx= \int e^{-i\lambda \langle x|y \rangle}f(x)\langle U_x\xi|\eta \rangle dx.\]
  Since the above equality is true for every $f \in C_{c}(\mathbb{R}^{n})$, it follows that $V_yU_xV_y^{*}=e^{-i\lambda \langle x|y\rangle}U_x$ for $x,y \in \mathbb{R}^{n}$. In other words, it follows that $(U,V)$ is a Weyl family of unitaries with phase factor $\lambda$. 
  
  \begin{xrcs}
  Prove that the correspondence $(U,V) \to (\pi_U,V)$ preserves irreducibility. 
   \end{xrcs}
  Since the dual of $\mathbb{R}^{n}$ is $\mathbb{R}^{n}$, it follows from Gelfand-Naimark theorem (see Theorem \ref{abelian group algebras}) that the map $\sigma:C^{*}(\mathbb{R}^{n}) \to C_{0}(\mathbb{R}^{n})$ defined by the equation
  \[
  \sigma(f)(s)=\int e^{i \lambda \langle s|x \rangle}f(x)dx\]
  is an isomorphism. Note that for $y \in \mathbb{R}^{n}$ and $f \in C_{c}(\mathbb{R}^{n})$, $\sigma(\alpha_y(f))=L_{y}(\sigma(f))$ where for $g \in C_{0}(\mathbb{R}^{n})$, $L_{y}(g)(s)=g(s-y)$. 
  Thus the dynamical system $(C^{*}(\mathbb{R}^{n}),\mathbb{R}^{n},\alpha)$ is isomorphic to $(C_{0}(\mathbb{R}^{n}),\mathbb{R}^{n},L)$. Thus, the final step in the proof of Stone-von Neumann theorem is the fact that $C_{0}(\mathbb{R}^{n}) \rtimes \mathbb{R}^{n}$ is isomorphic to $\mathcal{K}(L^{2}(\mathbb{R}^{n}))$. We prove the latter assertion for a general locally compact group.\footnote{We give a proof only in the unimodular case and leave the intricacies with the modular function to the interested reader.}
 
 Let $A$ be a $C^{*}$-algebra and $\cla$ be a dense $^*$-algebra of $A$. Suppose $X$ is a second countable locally compact Hausdorff topological space. For $a \in \cla$ and $f \in C_{c}(G)$, let $f\otimes a \in C_{c}(X,A)$ be defined by
 \[
 (f\otimes a)(x):=f(x)a.\]
 
 \begin{ppsn}
 \label{density in inductive limit}
 The linear span of $\{f \otimes a: f \in C_{c}(X), a \in \cla\}$ is dense in $C_{c}(X,A)$ with respect to the inductive limit topology. 
  \end{ppsn}
 \textit{Proof.} Let $F \in C_{c}(X,A)$ be given. Denote the support of $F$ by $K$. Choose an open set $U$ such that $K \subset U$ and $\overline{U}$ is compact. Fix $n \geq 1$. For $x \in K$, choose an open set $U_x \subset U$ such that for $y \in U_x$, 
 $||F(y)-F(x)|| \leq \frac{1}{2n}$. Choose $a_{x} \in \cla$ such that $||F(x)-a_x|| \leq \frac{1}{2n}$. Note that for $y \in U_x$, $||F(y)-a_x|| \leq \frac{1}{n}$. The family $\{U_x: x \in K\}$ covers $K$. Choose a finite subcover $\{U_{x_i}:i=1,2,\cdots,N\}$. For $i=1,2,\cdots,N$, let $a_i=a_{x_i}$.  
 Let $\{\phi_1,\phi_2,\cdots,\phi_n\}$ be a family in $C_{c}(X)$ such that 
 \begin{enumerate}
 \item[(a)] $supp(\phi_i) \subset U_{x_i}$ and $0 \leq \phi_i \leq 1$, and
 \item[(b)] for $x \in K$, $\sum_{i=1}^{N}\phi_{i}(x)=1$.
  \end{enumerate}
  Since $\sum_{i=1}^{N}\phi>0$ on $K$ and $K$ is a compact set, it follows that there exists an open subset $V \subset U$ such that $K \subset V$, $\sum_{i=1}^{N}\phi>0$ on $V$ and $\overline{V}$ is compact. Let $\chi \in C_{c}(X)$ be such that $0 \leq \chi \leq 1$, $\chi=1$ on $K$
  and $supp(\chi) \subset V$. Set $\psi:=\frac{\chi}{\sum_{i=1}^{N}\phi_i}$. 
  
 Set $F_{n}:=\sum_{i=1}^{N}\psi \phi_i \otimes a_i$. Then $supp(F_n) \subset \overline{U}$. Let $x \in X$ be given. Calculate as follows to observe that 
 \begin{align*}
 ||F(x)-F_n(x)||&=||\sum_{i=1}^{N}\psi(x)\phi_i(x)(F(x)-a_i)||\\
 & \leq \psi(x) \sum_{i=1}^{N} \phi_{i}(x)||F(x)-a_i||\\
 & \leq \frac{1}{n} \psi(x)\sum_{i=1}^{N}\phi_i(x) ~~(\textrm{since $||F(x)-a_i|| \leq \frac{1}{n}$ if $\phi_{i}(x)>0$}) \\
 & \leq \frac{1}{n}.
  \end{align*}
 This shows that $F_n \to F$ in the inductive limit topology and the proof is complete. \hfill $\Box$
 
 Let $X$ be a locally compact second countable Hausdorff space on which $G$ acts on the left. For $s \in G$ and $f \in C_{0}(X)$, define
 \[
 \alpha_s(f)(x)=f(s^{-1}x).\]
 Then $\alpha:=\{\alpha_s\}_{s \in G}$ is an action of $G$ on $C_{0}(X)$. Consider the vector space $C_{c}(X \times G)$. The vector space $C_{c}(X \times G)$ has a $^*$-algebra structure where the multiplication and involution are defined by
 \begin{align*}
 F*G(x,s)&=\int F(x,t)G(t^{-1}x,t^{-1}s)dt \\
 F^{*}(x,s)&=\Delta(s)^{-1}\overline{F(s^{-1}x,s^{-1})}.
  \end{align*}
 for $F,G \in C_{c}(X \times G)$. 
 
 For $F \in C_{c}(X \times G)$, let $\widetilde{F} \in C_{c}(G,C_0(X))$ be defined by
 \[
 \widetilde{F}(s)(x):=F(x,s).\]
 The map $C_{c}(X \times G) \ni F \to \widetilde{F} \in C_{c}(G,C_0(X))$ is an embedding and preserves the $^*$-algebra structure. By Prop. \ref{density in inductive limit}, it follows that $C_{c}(X \times G)$ is dense in $C_{c}(G,C_0(X))$. Consequently, $C_{c}(X \times G)$ is a dense $*$-subalgebra of the crossed product $C_{0}(X) \rtimes G$. 
 
 In the following, we let $G$ act on itself by left multiplication. For $s \in G$ and $f \in C_0(G)$, we let $\alpha_s(f) \in C_0(G)$ be defined by 
 \[
 \alpha_s(f)(x)=f(s^{-1}x).\]
 Then, $\alpha:=\{\alpha_s\}_{s \in G}$ defines an action of $G$ on $C_0(G)$. 
 
 \begin{thm}[Stone-von Neumann theorem: the abstract version]
 \label{Stone abstract}
 Let $G$ be a unimodular group. The crossed product $C_{0}(G) \rtimes_{\alpha} G$ is isomorphic to $\mathcal{K}(L^{2}(G))$. 
  \end{thm}
 \textit{Proof.} In view of Exercise \ref{co-ordinate free universal picture}, it suffices to exhibit a family $\{\widetilde{\theta}_{f,g}: f,g \in C_{c}(G)\}$ in $C_{0}(G) \rtimes G$ such that 
 \begin{enumerate}
 \item[(1)] for $f_1,f_2,g_1,g_2 \in C_{c}(G)$, $\widetilde{\theta}_{f_1,g_1}\widetilde{\theta}_{f_2,g_2}=\langle f_2|g_1 \rangle \widetilde{\theta}_{f_1,g_2}$, 
 \item[(2)] for $f,g \in C_{c}(G)$, $\widetilde{\theta}_{f,g}^{*}=\widetilde{\theta}_{g,f}$, and
 \item[(3)] the linear span of $\{\widetilde{\theta}_{f,g}: f,g \in C_{c}(G)\}$ is dense in $C_0(G) \rtimes G$. 
  \end{enumerate}
 For $f,g \in C_{c}(G)$, let $\widetilde{\theta}_{f,g} \in C_{c}(X \times G)$ be defined by
 \[
 \widetilde{\theta}_{f,g}(x,s)=f(x)\overline{g(s^{-1}x)}.\]
 Here $X=G$. It is routine to check $(1)$ and $(2)$.  By Prop. \ref{density in inductive limit} and the fact that the map $X \times G \ni (x,s) \to (x,s^{-1}x) \to X \times X$ is a homeomorphism, it follows that the linear span of $\{\widetilde{\theta}_{f,g}: f,g \in C_{c}(G)\}$ is dense in $C_{0}(G) \rtimes G$.  
 The proof is now complete. \hfill $\Box$
 
 Let $M:C_0(G) \to B(L^{2}(G))$ be the multiplication representation, i.e. for $f \in C_0(G)$, let $M(f) \in B(L^{2}(G))$ be defined by 
 \[
 M(f)\xi(s):=f(s)\xi(s).\]
 Let $\lambda:=\{\lambda_s\}_{s \in G}$ be the left regular representation of $G$ on $L^{2}(G)$. For $f,g \in C_{c}(G)$, let $\theta_{f,g}$ be the rank one operator on $L^{2}(G)$ defined by 
 \[
 \theta_{f,g}(\xi)=f\langle \xi|g \rangle.\]

 \begin{xrcs}[Stone-von Neumann: another version]
 \label{Stone another}
 Keep the foregoing notation and the notation of Thm. \ref{Stone abstract}.
 \begin{enumerate}
 \item[(1)] Prove that $(M,\lambda)$ is a covariant representation of the dynamical $(C_0(G),G,\alpha)$. 
 \item[(2)] Show that for $f,g \in C_{c}(G)$, 
 \[
 (M \rtimes \lambda)(\widetilde{\theta}_{f,g})=\theta_{f,g}.\]
 Conclude that the range of $M \rtimes \lambda$ is $\mathcal{K}(L^{2}(G))$. It is now clear that $M \rtimes \lambda$ is an irreducible representation of $C_0(G) \rtimes G$. 
 \item[(3)] Suppose $(\pi,U)$ is a covariant representation of $(C_0(G),G,\alpha)$ on $\clh$. Suppose that $(\pi,U)$ is irreducible i.e. $\{\pi(f),U_s:f \in C_0(G), s \in G\}^{'}=\bbc$ (or equivalently $\pi \rtimes U$ is irreducible). Then, there exists a unitary $X:\clh \to L^{2}(G)$ such that for $f \in C_0(G)$ and $s \in G$, 
 \[
 X\pi(f)X^{*}=M(f);~~\textrm{and~}XU_s X^{*}=\lambda_s.\]
 \item[(4)] Suppose $(\pi,U)$ is a non-degenerate covariant representation of $(C_0(G),G,\alpha)$ on a Hilbert space $\clh$.  Then, there exists a Hilbert space $\mathcal{L}$ and a unitary $X:\clh \to L^{2}(G) \otimes \mathcal{L}$ such that 
 \[
 X\pi(f)X^{*}=M(f) \otimes 1;~~\textrm{and~}XU_sX^{*}=\lambda_s \otimes 1\]
 for $f \in C_0(G)$ and $s \in G$. 
  \end{enumerate}
  \end{xrcs}
 To prove $(3)$ and $(4)$, use $(2)$, apply Thm. \ref{Stone abstract}, Thm. \ref{representation of compacts} and the fact that covariant representations of a dynamical system are in bijective correspondence with representations of the associated crossed product. 
 
 \begin{rmrk}
 The reader interested to know more about the history of Stone-von Neumann theorem and its role in subsequent developments in the $C^{*}$-algebra theory should consult  the excellent article \cite{Rosenberg} by Jonathan Rosenberg.
  \end{rmrk}
 
 \section{Cooper's theorem}
 In this section, we derive the Wold decomposition picture of a  $1$-parameter semigroup of isometries in the continuous case. Recall that, in the discrete case, we have already derived Wold decomposition in Thm. \ref{Wold}. Let $\bbr_{+}:=[0,\infty)$. Note that $\bbr_{+}$ is a semigroup. 
 Let us start with a few definitions. 
  \begin{dfn}
 Let $\clh$ be a Hilbert space and let $V:[0,\infty) \to B(\clh)$ be a map. We denote the image of an element $t \in [0,\infty)$ under $V$ by $V_t$. We say that $V:=\{V_t\}_{t \in \bbr_{+}}$ is a strongly continuous semigroup of isometries if 
 \begin{enumerate}
 \item[(1)] for every $t \in \bbr_{+}$, $V_t$ is an isometry,
 \item[(2)] for $s,t \in \bbr_{+}$, $V_{s+t}=V_s V_t$, and
 \item[(3)] for $\xi \in \clh$, the map $\bbr_{+} \ni t \to V_t \xi \in \clh$ is continuous. 
 \end{enumerate}
  \end{dfn}
  Sometime, we call a strongly continuous $1$-parameter semigroup of isometries as an isometric representation of $\bbr_{+}$. Let $V:=\{V_t\}_{t \geq 0}$ be an isometric representation of $\bbr_{+}$ on  a Hilbert space $\clh$. Then, $V$ is said to be \emph{pure} if $V_{t}^{*} \to 0$ strongly, i.e. for $\xi \in \clh$, $V_{t}^{*}\xi \to 0$ as $t \to \infty$. If each $V_t$ is  unitary, then we say that $V$ is a unitary representation of $\bbr_{+}$. 
  
  \begin{xrcs}
  Let $U:=\{U_t\}_{t \geq 0}$ be a unitary representation  of $\bbr_{+}$ on a Hilbert space $\clh$. Then, $U$ can be extended to a unitary representation of the group $\bbr$. 
    \end{xrcs}
  \textit{Hint:} Simply set $U_{-t}:=U_{t}^{*}$ for $t>0$.

 \begin{ppsn}
Let $V:=\{V_t\}_{t \geq 0}$ be an isometric representation of $\bbr_{+}$ on a Hilbert space $\clh$. For $t \geq 0$, let $E_t$ be the range projection of $V_t$, i.e. $E_t:=V_{t}V_{t}^{*}$. The following are equivalent. 
\begin{enumerate}
\item[(1)] The isometric representation $V$ is pure. 
\item[(2)] For $\xi \in \clh$, $E_{t}\xi \to 0$ as $t \to \infty$. 
\item[(3)] The intersection $\displaystyle \bigcap_{t \geq 0}V_t\clh=\{0\}$. 
\end{enumerate}
\end{ppsn}  
\textit{Proof.} Let $t \in \bbr_{+}$ be given. Observe that for $\xi \in \clh$, $||E_t\xi|| =||V_tV_{t}^{*}\xi|| \leq ||V_t^{*}\xi||$. Also, for $\xi \in \clh$, 
$
||V_{t}^{*}\xi||=||V_{t}^{*}V_{t}V_{t}^{*}\xi|| \leq ||E_t\xi||$. The equivalence between $(1)$ and $(2)$ is now clear. 

Let $\clh_{\infty}:=\bigcap_{t \geq 0}V_t\clh$. Suppose that $(2)$ holds. Let $\xi \in \clh_{\infty}$. Then, $E_{s}\xi=\xi$ for every $s \in \bbr_{+}$. Since $E_{s}\xi \to 0$ as $s \to \infty$, it follows that $\xi=0$. 
Hence $(2) \implies (3)$. 

 Suppose that $\clh_{\infty}=\{0\}$. Taking orthocomplements, it follows that $\bigcup_{t \geq 0}\ker(V_t)^{*}$ is dense in $\clh$. Since $\{E_t\}_{t \geq 0}$ is norm bounded, it suffices to verify $(2)$ when $\xi \in \bigcup_{s \geq 0}\ker(V_s)^{*}$. 
 Let $s>0$ and let $\xi \in Ker(V_s^*)$ be given. Note that for $t>s$, $Ran(V_t) \subset Ran(V_s)$. Thus, for $t>s$, $E_t\xi=E_tE_s\xi=0$. Consequently, $E_t\xi \to 0$ as $t \to \infty$. This completes the proof of $(3) \implies (2)$. 
 \hfill $\Box$

 Let $V:=\{V_t\}_{t \geq 0}$ be an isometric representation of $\bbr_{+}$ on a Hilbert space $\clh$. Set $\clh_{\infty}:=\displaystyle \bigcup_{t \geq 0}V_t\clh$. Note that for any $s>0$, $\displaystyle \clh_{\infty}=\bigcap_{t >s}V_t\clh$. More generally, if $F$ is a cofinal subset of $[0,\infty)$, then $\clh_{\infty}=\displaystyle \bigcap_{t \in F}V_t\clh$. It is clear that for $t \geq 0$, $V_t\clh_{\infty}=\clh_{\infty}$ and $V_{t}^{*}\clh_\infty=\clh_\infty$. Thus, $\clh_{\infty}$ is invariant under $\{V_t,V_{t}^{*}: t \geq 0\}$. Also, $\clh_{\infty}^{\perp}$ is invariant under $\{V_t,V_{t}^{*}:t \geq 0\}$.

 For $t \geq 0$, let $U_t$ be the restriction of $V_t$ to $\clh_\infty$ and let $W_t$ be the restriction of $V_t$ to $\clh_{\infty}^{\perp}$. If we decompose $\clh$ as $\clh=\clh_{\infty}^{\perp}\oplus \clh_{\infty}$, then $V_t$ is of the form 
\[V_t= \begin{bmatrix}
           W_t & 0 \\
           0 & U_t
           \end{bmatrix}.\] By construction, $U:=\{U_t\}_{t \geq 0}$ is a unitary representation of $\bbr_{+}$ on $\clh_{\infty}$. 
\begin{xrcs}
With the foregoing notation, prove that  $W:=\{W_t\}_{t \geq 0}$ is a pure isometric representation of $\bbr_{+}$ on $\clh_{\infty}^{\perp}$. 
 \end{xrcs} 
 
 Thus, we have proved the following `baby version' of Wold decomposition.
\begin{ppsn}
\label{Baby Wold}
Any isometric representation of $\bbr_{+}$ is unitarily equivalent to a direct sum of a pure isometric representation  of $\bbr_{+}$ and a unitary representation of $\bbr_{+}$. 
\end{ppsn}

Thus, to understand $1$-parameter semigroups of isometries, we need to understand pure isometric representations and unitary representations. We have already talked about  
 unitary representations of $\bbr$. They are in bijective correspondence with non-degenerate representations of $C^{*}(\bbr) \cong C_0(\bbr)$. What can we say about pure semigroups of isometries ? 
 The pleasant feature is that there are only countably many of them and we can list them. This is the content of Cooper's theorem. 

For $t \geq 0$, let $S_t$ be the isometry on $L^{2}((0,\infty))$ defined by \begin{equation}
\label{isometries}
S_{t}(\xi)(s):=\begin{cases}
 \xi(s-t)  & \mbox{ if
} s  \geq t,\cr
   &\cr
    0 &  \mbox{ if } s<t
         \end{cases}
\end{equation}
for $\xi \in L^{2}((0,\infty))$. Then, $S:=\{S_t\}_{t \geq 0}$ is a strongly continuous, pure semigroup of isometries, usually called the shift semigroup on $L^{2}((0,\infty))$. The main result of this section 
is the following theorem. 
\begin{thm}[Cooper]
\label{Cooper}
Let $V:=\{V_t\}_{t \geq 0}$ be a strongly continuous semigroup of isometries on a Hilbert space $\clh$. Suppose that $V$ is pure, i.e. $V_t^{*} \to 0$ strongly as $t \to \infty$. Then, 
there exists a Hilbert space $\mathcal{L}$ and a unitary $U:\clh \to L^{2}((0,\infty)) \otimes \mathcal{L}$ such that 
 for $t \geq 0$, 
$
UV_tU^{*}=S_t \otimes 1$.
\end{thm}

The usual proofs available in the literature make essential use of unbounded operators. The proof that is provided here is more algebraic.  The idea behind the proof is to look for  an algebra whose non-degenerate
representations are in $1$-$1$ correspondence with pure semigroups of isometries. It turns out that the algebra that we seek is $C_0(\bbr) \rtimes \bbr$, where the action of $\bbr$ is by translation.
Thanks to the abstract version of the Stone-von Neumann theorem (Thm. \ref{Stone abstract}), $C_0(\bbr) \rtimes \bbr \cong \mathcal{K}(L^{2}(\bbr))$ and the latter algebra has only countably
many non-degenerate representations up to unitary equivalence. We make this precise in what follows.

\textbf{Unitary dilations:} Let $V:=\{V_t\}_{t \geq 0}$ be an isometric representation of $\bbr_+$ on a Hilbert space $\clh$. Suppose $U:=\{U_t\}_{t \in \bbr}$ is a unitary representation on a Hilbert space $\mathcal{K}$ and suppose $j:\clh \to \mathcal{K}$ is an isometry. We say that the triple $(U,\mathcal{K},j)$ is a unitary dilation of $V$ if for $t \geq 0$ and $\xi \in \clh$, 
\[
U_t(j(\xi))=j(V_t\xi).\]
We usually suppress the notation $j$ and regard $\clh$ as a closed subspace of $\mathcal{K}$ and by abusing notation, we write the above equality as $V_t\xi=U_t\xi$ whenever $t \geq 0$ and $\xi \in \clh$. Let $(U,\mathcal{K})$ be a unitary dilation of $V$. We say that $(U,\mathcal{K})$ is \emph{minimal} if $\bigcup_{t \geq 0}U_{-t}\clh$ is dense in $\mathcal{K}$.  

Let $(U^{(1)},\clk_1)$ and $(U^{(2)},\clk_2)$ be two unitary dilations of $V$. We say that they are unitarily equivalent if there exists a unitary $X:\clk_1 \to \clk_2$ such that $X\clh=\clh$ and for $t \in \bbr$, $XU^{(1)}_tX^{*}=U^{(2)}_{t}$. 

\begin{thm}
Let $V:=\{V_t\}_{t \geq 0}$ be an isometric representation of $\bbr_{+}$ on a Hilbert space $\clh$. Then, $V$ admits a minimal unitary dilation. Moreover, any two minimal unitary dilations of $V$ are unitarily equivalent. 
\end{thm}
\textit{Proof.} First, we prove the uniqueness. Let $(U^{(1)},\clk_1)$ and $(U^{(2)},\clk_2)$ be two minimal unitary dilations of $V$. Let $\mathcal{D}_{i}:=\bigcup_{t \geq 0}U^{(i)}_{-t}\clh$. 
Then, $\mathcal{D}_i$ is a dense subspace of $\clk_i$. If we make use of Remark \ref{KRP}, we see that there exists a unique unitary $X:\clk_1 \to \clk_2$ such that 
\[
XU^{(1)}_{-t}\xi=U^{(2)}_{-t}\xi\]
for $t \geq 0$ and $\xi \in \clh$. It is routine to verify that $X$ implements the required equivalence between $(U^{(1)},\clk_1)$ and $(U^{(2)},\clk_2)$. 

Next, we construct a dilation that is minimal. Define an equivalence relation on $\clh \times [0,\infty)$ as follows. We say $(\xi,s) \sim (\eta,t)$ if $V_t\xi=V_s\eta$. Denote the set of equivalence classes by $\mathcal{K}_{\infty}$. Then, $\mathcal{K}_{\infty}$ has the structure of an inner product space where addition, scalar multiplication and inner product are given by 
\begin{align*}
[(\xi,s)]+[(\eta,t)]&=[(V_t\xi+V_s\eta,s+t)]\\
\lambda[(\xi,s)]&=[(\lambda\xi,s)]\\
\langle [(\xi,s)]|[(\eta,t)]\rangle&=\langle V_t\xi|V_s\eta \rangle.
\end{align*}
Denote the completion of $\clk_{\infty}$ by $\clk$. Let $j:\clh \to \clk$ be defined by $j(\xi):=[(\xi,0)]$. Then, $j$ is an isometry. Through $j$, we view $\clh$ as a closed subspace of $\clk$. 

For $t \geq 0$, let $U_t:\clk_\infty \to \clk_\infty$ be defined by 
\[
U_t([(\xi,s)]=[(V_t\xi,s)].\]
It is clear that $U_t$ is well defined, linear and preserves the inner product. Clearly, for $\xi \in \clh$, $U_t\xi=V_t\xi$. More precisely, $U_t(j(\xi))=j(V_t\xi)$ for $t \geq 0$ and $\xi \in \clh$.

We claim that for $t \geq 0$,  $U_t$ is onto. Let $[(\eta,s)] \in \clk_\infty$ be given. Then, by definition, 
\begin{equation}
\label{formula for inverse}
[(\eta,s)]=[(V_t\eta,s+t)]=U_t([(\eta,s+t)]).\end{equation}
The above equality implies that $U_t$ is onto. 
Thus, $U_t$ extends to a unitary operator on $\clk$ which we denote again by $U_t$. Clearly, $U_sU_t=U_{s+t}$ for $s,t \geq 0$. Thanks to Eq. \ref{formula for inverse}, we have $\bigcup_{t \geq 0}U_{t}^{-1}\clh=\bigcup_{t \geq 0}\{[(\xi,t):\xi \in \clh\}=\clk_\infty$ and the latter set is dense in $\clk$ by definition. 

Now, we check the strong continuity of the map $[0,\infty) \ni t \to U_t \in B(\clh)$. It suffices to prove that for every $\xi \in \clh$ and $s \geq 0$, the map $[0,\infty) \ni t \to U_t([(\xi,s)]) \in \clk$ is continuous. This follows from the strong continuity of $V$ and the inequality 
\[
|| U_{t_1}([(\xi,s)])-U_{t_2}([(\xi,s)])||=||V_{t_1}\xi-V_{t_2}\xi||.\]
Thus, $U:=\{U_t\}_{t \geq 0}$ is a strongly continuous group of unitaries. Extend $U$ to the whole of $\bbr$ and denote the extension again by $U$. Then, $U$ is clearly `the minimal unitary dilation' of $V$. 
\hfill $\Box$

\begin{xrcs}
Let $S:=\{S_t\}_{t \geq 0}$ be the shift semigroup on $L^{2}((0,\infty))$. Suppose $\lambda:=\{\lambda_t\}_{t \in \bbr}$ is the regular representation of $\bbr$ on $L^{2}(\bbr)$. Prove that $\lambda$ is the minimal unitary dilation of $S$. 
\end{xrcs}
\textit{Hint:} View $L^{2}((0,\infty))$ as a subspace of $L^{2}(\bbr)$ in the obvious way. It suffices to check that $\lambda:=\{\lambda_t\}_{t \in \bbr}$ has `the defining properties' for it to be the minimal unitary dilation. 

We will proceed towards the proof of Cooper's theorem. Let $V:=\{V_t\}_{t \geq 0}$ be a strongly continuous semigroup of isometries on a Hilbert space $\clh$. Assume that $V$ is pure. The isometric representation $V$ will be fixed for the rest of this section. 
Let $(U,\clk)$ be the minimal unitary dilation of $V$. For $t \in \bbr$, let $E_t$ be the orthogonal projection onto $U_t\clh$. 

Observe the following properties regarding the family $\{E_t:t \in \bbr\}$ of projections. 
\begin{enumerate}
\item[(i)] Since $U_t\clh \subset \clh$ for every $ t\geq 0$, it follows that if $s\leq t$, then $U_t\clh \subset U_s\clh$. This means that if $s \leq t$, then $E_s \geq E_t$. In particular, the family $\{E_t: t \in \bbr\}$ is a commuting family of projections. 
\item[(ii)] For $s,t \in \bbr$, $U_{s+t}\clh=U_s(U_t\clh)$. Thus, we have the following covariance relation. For $s \in \bbr$ and $t \in \bbr$, 
\[
U_sE_tU_{s}^{*}=E_{s+t}.\]
\item[(iii)] The map $\bbr \ni t \to E_{t} \in B(\clk)$ is strongly continuous. This follows from the strong continuity of $\{U_t\}_{t \in \bbr}$ and from the fact that $E_t=U_tE_0U_{t}^{*}$. 
\end{enumerate}

The key ingredient in the proof of Cooper's theorem is to understand the algebra generated by $\{E_t:t \in \bbr\}$. However, if we take the $C^{*}$-algebra generated by the family $\{E_t: t \in \bbr\}$, then it is not separable, and in particular its spectrum
will be `too big'. Taking inspiration from the study of group $C^{*}$-algebra, we look at the $C^{*}$-algebra generated by $\{\int f(t)E_t dt: f \in C_{c}(\bbr)\}$. 

Let us introduce a few notation. Adjoint $\infty$ at the right end of $\bbr$ and endow $(-\infty,\infty]$ with the usual topology. For $\phi \in C_0((-\infty,\infty])$ and $s \in \bbr$, let $R_s(\phi) \in C_{0}((-\infty,\infty])$ be defined by 
\[
R_s\phi(x)=\phi(x-s).\]
As usual, the convention is that $\infty\pm t=\infty$ for all $t \in \bbr$. Then, $R:=\{R_s\}_{s \in \bbr}$ defines an action of $\bbr$ on $C_0((-\infty,\infty])$.

For $f \in C_{c}(\bbr)$, let $\widetilde{f}:(-\infty,\infty] \to \bbc$ be defined by 
\[
\widetilde{f}(x):=\int_{-\infty}^{x}f(t)dt.\]
Then, $\widetilde{f} \in C_0((-\infty,\infty])$ for every $f \in C_{c}(\bbr)$. 

Let $\mathcal{D}$ be the $C^{*}$-algebra generated by $\Big \{\int f(t)E_t dt: f \in C_{c}(\bbr)\Big\}$. Note that $\mathcal{D}$ is commutative. 

\begin{lmma}
\label{description of characters}
Let $\chi$ be a character of the commutative $C^{*}$-algebra $\mathcal{D}$. Then, there exists a unique $x=:x_{\chi} \in (-\infty,\infty]$ such that 
\[
\chi\Big(\int f(t)E_tdt\Big)=\int_{-\infty}^{x}f(t)dt\]
for all $f \in C_{c}(\bbr)$. 
\end{lmma}
\textit{Proof.} Let $\chi$ be a character of $\mathcal{D}$. Note that the map 
\[
C_{c}(\bbr) \ni f \to \chi\Big(\int f(t)E_tdt\Big) \in \bbc\]
admits an extension to $L^{1}(\bbr)$. Therefore, there exists $\phi \in L^{\infty}(\bbr)$ such that 
\[
\chi\Big(\int f(t)E_tdt\Big)=\int f(t)\phi(t)dt\]
for $f \in C_{c}(\bbr)$. 

Since $\chi$ is non-zero, there exists real numbers $a$ and $b$ with $a<b$ and a function $g \in C_{c}(\bbr)$ such that $supp(g) \subset (a,b)$ and $\int g(s)\phi(s)dt=1$. 
We claim that $\phi(t)=1$ for almost all $t \in (-\infty,a)$. 

Let $f \in C_{c}(\bbr)$ be any function such that $supp(f) \subset (-\infty,a)$. 
Note that $E_tE_s=E_s$ whenever $t<a<s$. Therefore, 
\begin{align*}
\Big(\int f(t)E_tdt\Big)\Big(\int g(s)E_sds\Big)&=\Big(\int_{-\infty}^{a}f(t)E_tdt\Big)\Big(\int_{a}^{b}g(s)E_sds\Big)\\
&=\Big(\int_{(t,s) \in (-\infty,a) \times (a,b)} f(t)g(s)E_tE_sdtds\Big)\\
&=\Big(\int_{-\infty}^{a}f(t)dt\Big)\int_{-\infty}^{\infty}g(s)E_sds.
\end{align*}
Applying $\chi$ to the above equation and appealing to the fact that $\chi$ is a character, we deduce that 
\begin{align*}
\int_{-\infty}^{a}f(t)dt& = \Big(\int_{-\infty}^{a}f(t)dt\Big)\Big(\int_{a}^{b} g(s)\phi(s)ds\Big)\\
&=\Big(\int_{-\infty}^{a}f(t)dt\Big)\chi\Big(\int g(s)E_s ds\Big)\\
&=\chi\Big(\Big(\int_{-\infty}^{a}f(t)dt\Big)\Big(\int g(s)E_sds\Big)\Big) \\
&=\chi\Big(\Big(\int f(t)E_tdt \Big)\Big(\int g(s)E_sds\Big)\Big)\\
&=\chi\Big(\int f(t)E_tdt \Big)\chi\Big(\int g(s)E_sds\Big)\\
&=\int_{-\infty}^{a}f(t)\phi(t)dt \int_{a}^{b} g(s)\phi(s)ds \\
&=\int_{-\infty}^{a}f(t)\phi(t)dt. 
\end{align*}
Since $f$ is arbitrary, we can conclude that $\phi(t)=1$ for almost all $t \in (-\infty,a)$.

Define \[A:=\{a \in \bbr: \phi(t)=1 \textrm{~~for almost all $t \in (-\infty,a)$}\}.\] We have shown that $0 \in A$. Set $x:=\sup(A)$. Clearly, $\phi(t)=1$ for almost all $t \in (-\infty,x]$.  If $x=\infty$, then $\phi=1$ a.e. and we are done. 
Suppose that $x<\infty$. 

We claim that $\phi(t)=0$ for almost all $t \in (x,\infty)$. Suppose not. Then, there exist real numbers $a,b$ with $b>a>x$ and a function $g \in C_{c}(\bbr)$ with $supp(g) \subset (a,b)$ and $\int_{a}^{b} g(s)\phi(s)ds=1$. Arguing as before, we 
can conclude that $\phi(t)=1$ for almost all $t \in (-\infty,a)$. Then, $A \ni a >x$ which is a contradiction. This proves the claim. In other words, $\phi(t)=1_{(-\infty,x]}(t)$ a.e. Hence the proof. \hfill $\Box$

\begin{ppsn}
\label{covariant Cooper}
Keep the foregoing notation. 
\begin{enumerate}
\item[(1)] There exists a unique $*$-homomorphism $\pi:C_0((-\infty,\infty]) \to B(\clk)$ such that for $f \in C_{c}(\bbr)$, 
\[
\pi(\widetilde{f})=\int f(t)E_t dt.\]
\item[(2)] The pair $(\pi,U)$ is covariant, i.e. for $\phi \in C_0((-\infty,\infty])$ and $s \in \bbr$, 
\[
U_s\pi(\phi)U_s^{*}=\pi(R_s\phi).\]
\end{enumerate}
\end{ppsn}
\textit{Proof.} Uniqueness of the homomorphism $\pi$ follows from the fact that $\{\widetilde{f}:f \in C_{c}(\bbr)\}$ generates $C_0((-\infty,\infty])$. This can be proved using Stone-Weierstrass theorem. For the existence, let us denote the map 
\[
\widehat{D} \ni \chi \to x_\chi \in (-\infty,\infty]\]
by $T$. Note that the map $T$ is continuous. For $\phi \in C_0((-\infty,\infty])$, set \[\pi(\phi):=G(\phi \circ T).\] Here, $G: C(\widehat{\mathcal{D}}) \to \mathcal{D} \subset B(\clk)$ is the Gelfand transformation. It follows from definitions that for $f \in C_{c}(\bbr)$, 
\[
\pi(\widetilde{f})=\int f(t)E_t dt.\]

Let $s \in \bbr$ and $f \in C_{c}(\bbr)$ be given. Note that $\widetilde{R_sf}=R_s(\widetilde{f})$. Calculate as follows to observe that 
\begin{align*}
U_s\pi(\widetilde{f})U_s^{*}&=U_{s} \Big (\int f(t)E_tdt\Big)U_{s}^{*}\\
&=\int f(t)U_sE_tU_s^{*}dt\\
&=\int f(t)E_{s+t}dt\\
&=\int f(t-s)E_tdt\\
&=\int R_s(f)(t)E_tdt\\
&=\pi(\widetilde{R_s(f)})=\pi(R_s(\widetilde{f})).
\end{align*}
Since $\{\widetilde{f}:f \in C_{c}(\bbr)\}$ generates $C_0((-\infty,\infty])$, it follows that $U_s\pi(\phi)U_{s}^{*}=\pi(R_s\phi)$ for $s \in \bbr$ and $\phi \in C_0((-\infty,\infty])$. This completes the proof. \hfill $\Box$

There are a couple of crucial observations that we need to make. Recall that we have assumed $V$ is pure. 
\begin{lmma}
\label{Cooper vanishing at infinity}
For every $\xi \in \clk$, $E_t\xi \to 0$ as $t \to +\infty$.
\end{lmma}
\textit{Proof.} Note that for $s \geq 0$, $V_s=E_0U_sE_0$ and $V_{s}^{*}=E_0U_s^*E_0$. 
As $\{E_t:t \in \bbr\}$ is norm bounded,  It suffices to verify the statement when $\xi \in \bigcup_{s \geq 0}U_{s}^{*}\clh$. Let $s \geq 0$ and $\eta \in \clh$ be given. Note that for large $t$.
\[
E_{t}U_{s}^{-1}\xi=U_{t}E_0U_{t-s}^{*}E_0\xi=U_{t}V_{t-s}^{*}\xi.\]
Hence, for large $t$, 
\[
||E_{t}U_{s}^{-1}\xi|| \leq ||V_{t-s}^{*}\xi||.\]
The purity of $V$ implies that $V_{t-s}^{*}\xi \to 0$ as $t \to \infty$. Hence the proof. \hfill $\Box$

\begin{lmma}
\label{non-degenerate Cooper}
Let $\pi$ be the representation of $C_0((-\infty,\infty])$ described in Thm. \ref{covariant Cooper}. Then, $\pi$ restricted to $C_{0}(\bbr)$ is non-degenerate. 
\end{lmma}
\textit{Proof.} Suppose not. Then, there exists a non-zero vector $\eta$ such that 
\[
\langle \pi(\phi)\xi|\eta\rangle=0\]
for every $\phi \in C_{c}(\bbr)$ and $\xi \in \clk$. 

Let $\xi \in \clk$ be arbitrary. Let $\omega_\xi:C_{c}(\bbr) \to \bbc$ be defined by 
\[
\omega_\xi(f)=\int f(t)\langle E_t\xi|\eta \rangle dt.\]
Define $I:C_{c}(\bbr) \to \bbc$ by \[I(f)=\int f(t)dt.\] Suppose $f \in \ker(I)$. Then, $\widetilde{f} \in C_{c}(\bbr)$. The fact that $\langle \pi(\widetilde{f})\xi|\eta \rangle=0$ implies that 
$\int f(t)\langle \xi|E_t\eta\rangle=0$, i.e $\omega_\xi(f)=0$. Consequently, $\ker(I)=\ker(\omega_\xi)$. Thus, there exists $c_{\xi} \in \bbc$ such that for all $f \in C_{c}(\bbr)$, 
\[
\int f(t)\langle E_t\xi|\eta \rangle dt=c_{\xi}\int f(t)dt.\]
As the above equality happens for every $f \in C_{c}(\bbr)$ and the map $\bbr \to t \to \langle E_t\xi|\eta \rangle \in \bbc$ is continuous, we conclude that $\langle E_t\xi|\eta\rangle =c_{\xi}$ for every $t \in \bbr$. 
Lemma \ref{Cooper vanishing at infinity} implies that $\langle E_t\xi|\eta \rangle=0$ for every $t \in \bbr$. 
 
 Thus, $\eta$ is orthogonal to $U_t\clh$ for every $t \in \bbr$. However, $\bigcup_{t \geq 0}U_{-t}\clh$ is dense in $\clk$. Thus, $\eta$ is orthogonal to every vector of $\clk$ forcing $\eta=0$. This contradiction
 completes the proof. \hfill $\Box$
 
 The final step that is required is the following lemma. Define
 \[
 \pi(C_{c}((0,\infty)))\clk:=\overline{span\{\pi(\phi)\xi:\phi \in C_{c}((0,\infty)), \xi \in \clk\}}.\]
 
 \begin{lmma}
 \label{projection onto the given space}
 The orthogonal projection onto $\pi(C_{c}((0,\infty)))\clk$ is $E_0$, or in other words, $\clh=\pi(C_{c}((0,\infty)))\clk$.  
  \end{lmma}
  \textit{Proof.} Let $C_{c}^{\infty}((0,\infty))$ be the algebra of smooth functions $f$ on $\bbr$ that vanishes outside a compact subset contained in $(0,\infty)$. Clearly, $\pi(C_{c}^{\infty}((0,\infty)))\clk=\pi(C_{c}((0,\infty)))\clk$. 
 Let $\phi \in C_{c}^{\infty}((0,\infty))$ be given. Observe that 
 \[
 E_0\pi(\phi)=E_0\Big(\int_{0}^{\infty}\phi^{'}(t)E_tdt\Big)=\int_{0}^{\infty} \phi^{'}(t)E_0E_tdt=\int \phi^{'}(t)E_tdt=\pi(\phi).\]
 The above equality implies that $\pi(C_{c}^{\infty}((0,\infty)))\clk \subset \clh$. 
 
 Suppose that the inclusion $\pi(C_c^{\infty}((0,\infty)))\clk \subset \clh$ is proper. Then, there exists a non-zero vector $\eta \in \clh$ such that for every $\phi \in C_{c}^{\infty}((0,\infty))$ and $\xi \in \clk$,
 \[
 \int_{0}^{\infty} \phi^{'}(t)\langle E_t\xi|\eta \rangle dt=0.\]
 Let $f \in C_{c}^{\infty}((0,\infty))$ be such that $\int f(t)dt=0$. Define $\phi(x):=\int_{-\infty}^{x}f(t)dt$. Then, $\phi^{'}(x)=f(x)$. The hypothesis implies that 
 \[
 \int f(t)\langle E_t\xi|\eta \rangle dt=0.\]
 Thus, arguing as in Lemma \ref{non-degenerate Cooper}, we can conclude that $\langle E_t\xi|\eta\rangle =0$ for every $t \in (0,\infty)$. Therefore, $E_t\eta=0$ for every $t>0$. Since $E_t\eta \to E_0\eta=\eta$ as $t \to 0$, we conclude that $\eta=0$, which is a contradiction. Hence the proof. \hfill $\Box$

 \textit{Proof of Thm. \ref{Cooper}:} Let $(\pi,U)$ be the covariant representation of the dynamical system $(C_0((-\infty,\infty]),\bbr,R)$ constructed in Prop. \ref{covariant Cooper}. Thanks to Lemma \ref{non-degenerate Cooper}, $\pi$ restricted to $C_0(\bbr)$ is non-degenerate. Thus, we obtain a non-degenerate covariant representation $(\pi,U)$ of the dynamical system $(C_0(\bbr),\bbr,R)$ on $\clk$. By the Stone-von Neumann theorem (Ex. \ref{Stone another}), there exists a Hilbert space $\mathcal{L}$ and a unitary $X:\clk \to L^{2}(\bbr) \otimes \mathcal{L}$ such that for $\phi \in C_0(\bbr)$ and $s \in \bbr$, 
 \begin{align*}
 X\pi(\phi)X^{*}&=M(\phi)\otimes 1\\
 XU_sX^{*}&=\lambda_s \otimes 1.
 \end{align*}
 Here, $M:C_0(\bbr) \to B(L^{2}(\bbr))$ is the multiplication representation and $\lambda:=\{\lambda_s\}_{s \in \bbr}$ is the regular representation of $\bbr$ on $L^{2}(\bbr)$. It is clear from the above equation and Lemma \ref{projection onto the given space} that $X$ takes $\clh$ onto $L^{2}((0,\infty))\otimes \mathcal{L}$. Let $Y$ be the restriction of $X$ to $\clh$. Then, $Y$ is a unitary between $\clh$ and $L^{2}((0,\infty))\otimes \mathcal{L}$. 
 
 Let $S:=\{S_t\}_{t \geq 0}$ be the shift semigroup on $L^{2}((0,\infty))$. 
 For $t \geq 0$, $S_t \otimes 1$ is the restriction of $\lambda_t \otimes 1$ to $L^{2}((0,\infty))\otimes \mathcal{L}$. As $U:=\{U_t\}_{t \geq 0}$ is a dilation of $V:=\{V_t\}_{t \geq 0}$, for $t \geq 0$, $V_t$ is the restriction of $U_t$ to $\clh$. As $X$ intertwines $\{U_t\}_{t \geq 0}$ and $\{\lambda_t\otimes 1\}_{t \geq 0}$, it follows that $Y$ intertwines $\{V_t\}_{t \geq 0}$ and $\{S_t \otimes 1\}_{t \geq 0}$. The proof is complete. \hfill $\Box$
 
 With Cooper's theorem in hand, the Wold decomposition can now be restated as follows. 
 
 \begin{thm}[Wold-decomposition]
 Let $V:=\{V_{t}\}_{t \geq 0}$ be a strongly continuous semigroup of isometries on a Hilbert space $\clh$. Then, there exist Hilbert spaces $\mathcal{L}$ and $\clk$, a  unitary representation $U:=\{U_t\}_{t \geq 0}$ on $\clk$ and a unitary $X:\clh \to \big(L^{2}((0,\infty))\otimes \mathcal{L} \big)\oplus \mathcal{K}$ such that for every $t \geq 0$,
 \[XV_tX^{*}=\begin{bmatrix}
 S_t \otimes 1 & 0 \\
 0 & U_t
 \end{bmatrix}.\]

 \end{thm}

  \section{The non-commutative torus $A_\theta$}
  In this section, we discuss the simplicity of the $C^{*}$-algebra, called the non-commutative torus, associated to irrational rotations on the circle. The non-commutative torus is probably the widely studied example in the field of non-commutative geometry. 
 
 \begin{dfn}
  Let $\theta \in \mathbb{R}$ be given. The non-commutative torus $A_{\theta}$ is the universal $C^{*}$-algebra generated by two unitaries $u$ and $v$ such that 
  \[
  uv=e^{2\pi i \theta}vu.\]
   \end{dfn}
   The reason $A_{\theta}$ is called the non-commutative torus is because when $\theta=0$,  $A_{\theta}$ is isomorphic to $C(\mathbb{T}^{2})$. First, we realise $A_\theta$ as a crossed product. Let $\alpha:C(\mathbb{T}) \to C(\mathbb{T})$ be defined by 
   \[
   \alpha(f)(z)=f(e^{-2\pi i\theta}z).\]
   Let $\widetilde{u}:=z \in C(\mathbb{T})$. Then $\alpha(\widetilde{u})=e^{-2\pi i \theta}\widetilde{u}$. 
   Clearly, $\alpha$ is an automorphism of $C(\mathbb{T})$. The automorphism $\alpha$ induces an action of the cyclic group $\mathbb{Z}$ on $C(\mathbb{T})$. Consider the crossed product $C(\mathbb{T}) \rtimes_{\alpha} \mathbb{Z}$. By the definition of the crossed product, $C(\mathbb{T}) \rtimes_{\alpha} \mathbb{Z}$ is the universal $C^{*}$-algebra generated by a copy of $C(\mathbb{T})$ and a unitary $\widetilde{v}$ such that 
   \[
   \widetilde{v}f\widetilde{v}^{*}=\alpha(f)\]
   for $f \in C(\mathbb{T})$. However, $C(\mathbb{T})$ is generated by $\widetilde{u}$ and the equation $\widetilde{u}\widetilde{v}=e^{2 \pi i \theta}\widetilde{v}\widetilde{u}$ is satisfied in $C(\mathbb{T}) \rtimes \mathbb{Z}$. Consequently, there exists a surjective $^*$-homomorphism $\phi:A_{\theta} \to C(\mathbb{T}) \rtimes \mathbb{Z}$ such that $\phi(u)=\widetilde{u}$ and $\phi(v)=\widetilde{v}$. 
   
   Note that $C(\mathbb{T})$ is the universal $C^{*}$-algebra generated by the unitary $\widetilde{u}$. Consequently, there exists a homomorphism $\pi:C(\mathbb{T}) \to A_{\theta}$ such that $\pi(\widetilde{u})=u$. The relation \[uv=e^{2\pi i \theta}vu\] implies that the pair $(\pi,v)$ is a covariant representation of $(C(\mathbb{T}),\mathbb{Z},\alpha)$. Denote the homomorphism  $\pi \rtimes v$ from $ C(\mathbb{T}) \rtimes \mathbb{Z} \to A_{\theta}$ by $\psi$. It is clear that $\psi(\widetilde{u})=u$ and $\psi(\widetilde{v})=v$. 
   Hence $\psi$ and $\phi$ are inverses of each other. This proves that $A_{\theta}$ is isomorphic to $C(\mathbb{T}) \rtimes \mathbb{Z}$. 
   
   The main theorem of this section is that if $\theta$ is irrational, then $A_{\theta}$ is simple, i.e. it has no non-trivial closed two sided ideals. The proof makes use of a very useful notion called conditional expectation.
   
   \begin{dfn}
   Let $A$ be a $C^{*}$-algebra and $B \subset A$ be a $C^{*}$-subalgebra. A linear map $E: A \to B$ is called a conditional expectation of $A$ onto $B$ if 
   \begin{enumerate}
   \item[(1)] for $b \in B$, $E(b)=b$, 
   \item[(2)] for $b_1,b_2 \in B$ and $a \in A$, $E(b_1ab_2)=b_1E(a)b_2$, and
   \item[(3)] for $a \in A$, $E(a^{*}a) \geq 0$. 
      \end{enumerate}
   The conditional expectation $E$ is said to be faithful if whenever $E(a^{*}a)=0$, $a=0$. 
      \end{dfn}
   \begin{lmma}
   Let $E:A \to B$ be a conditional expectation. Then $E$ is contractive, i.e. for $a \in A$,  $||E(a)|| \leq ||a||$. 
     \end{lmma}
   \textit{Proof.} Consider the Hilbert $B$-module $B$. Represent $B$ as adjointable operators on $B$ by left multiplication. That is, for $b \in B$, let $L_{b}:B \to B$ be defined by $L_{b}(c)=bc$. Then, $L_{b}$ is adjointable for every $b \in B$. Moreover, $B \ni b \to L_{b} \in \mathcal{L}_{B}(B)$ is an injective $^*$-homomorphism. Hence for $b \in B$, $||b||=||L_{b}||$. Suppose $a \in A$ is a positive element. Calculate, using Exercise \ref{norm of positive operator}, as follows to observe that
   \begin{align*}
   ||E(a)||&=||L_{E(a)}|| \\
              &= \sup \{||\langle b|E(a)B \rangle||: b \in B, ||b||=1\}\\
              &=\sup\{||b^{*}E(a)b||: b \in B, ||b||=1\}\\
              &=\sup \{||E(b^{*}ab)||: b \in B, ||b||=1\} \\
              & \leq \sup \{||E(b^{*}||a||b)||: b \in B, ||b||=1\} \\
              & \leq ||a||.
     \end{align*} 
     
     For $a_1,a_2 \in A$, let $\langle a_1|a_2 \rangle=E(a_1^{*}a_2)$. Then $\langle~|~\rangle$ is a $B$-valued semi-definite inner product. Hence by the Cauchy-Schwarz inequality, we have for $b \in B$ and $a \in A$, 
     \[
     ||b^{*}E(a)||=||E(b^{*}a)|| \leq ||E(b^{*}b)||^{\frac{1}{2}}||E(a^{*}a)||^{\frac{1}{2}} \leq ||b||||a||.\]
   Let $(e_\lambda)$ be an approximate identity of $B$. The above equation implies that \[||e_\lambda E(a)|| \leq ||a||.\] But $e_{\lambda}E(a) \to E(a)$ in $B$. Consequently, $||E(a)|| \leq ||a||$ for every $a \in A$. This completes the proof. \hfill $\Box$
   
   \begin{rmrk}
   A theorem due to Tomiyama asserts that if $E:A \to B$ is a contractive linear map such that $E(b)=b$ for every $b \in B$, then $E$ is a conditional expectation of $A$ onto $B$. 
     \end{rmrk}
  Here is an example of a conditional expectation. Suppose $A$ is a a $C^{*}$-algebra and $\alpha:=\{\alpha_s\}$ is an action of a compact group $G$ on $A$. We choose the Haar measure on $G$ to be a probability measure. 
   The fixed point algebra of $\alpha$, denoted $A^{\alpha}$, is defined by
   \[
   A^{\alpha}:=\{a \in A: \alpha_s(a)=a ~~\forall s \in G\}.\]
   Note that $A^{\alpha}$ is a $C^{*}$-subalgebra of $A$. Define $E:A \to A$ by \[
   E(a)=\int \alpha_s(a)ds.\]
   \begin{xrcs}
   Verify that $E$ is a faithful conditional expectation of $A$ onto $A^{\alpha}$. 
      \end{xrcs}
   Let $(A,G,\alpha)$ be a $C^{*}$-dynamical system. Assume that $G$ is discrete and abelian. Note that $\widehat{G}$ is a compact group. For the sake of simplicity, assume that $A$ is unital. Recall that the crossed product $A \rtimes_\alpha G$ is the universal $C^{*}$-algebra generated by a copy of $A$ and unitaries $\{u_s: s \in G\}$ such that $u_su_t=u_{st}$ and $u_s a u_s^{*}=\alpha_s(a)$ for $a \in A$. This universal picture reveals that for $\chi \in \widehat{G}$, there exists a unique $^*$-homomorphism $\beta_{\chi}: A \rtimes G \to A \rtimes G$ such that 
   \begin{align*}
   \beta_{\chi}(a)&=a \\
   \beta_{\chi}(u_s)&=\chi(s)u_s.
     \end{align*}
   Then $\beta:=\{\beta_\chi\}_{\chi \in \widehat{G}}$ is an action of $\widehat{G}$ on the crossed product $A \rtimes G$. The action $\beta$ is called the dual action on the crossed product $A \rtimes G$. We claim that the fixed point algebra of $\beta$ is $A$. 
   Set $B=A \rtimes G$. It is clear that $A \subset B^{\beta}$. Let $E:B \to B^{\beta}$ be the conditional expectation given by
   \[
   E(b):=\int \beta_{\chi}(b) d\chi.\]
   
    It suffices to show that $E(b) \in A$ for every $b \in B$. Since $A$ is closed in $A \rtimes G$ and $E$ is continuous, it suffices to show that $E(b) \in A$ whenever $b$ is the form $b=\sum_{s \in G}a_su_s$. It is clear that it suffices to show that 
    $E(au_s)=0$ if $s \neq e$. Let $s \neq e$ and $a \in A$ be given. Note that \[E(au_s)=\int a\chi(s)u_s d\chi= (\int \chi(s)d\chi)au_s.\] Since $s \neq e$, by Raikov's theorem,  there exists a character $\chi_0$ of $G$ such that $\chi_0(s) \neq 1$. Calculate as follows to observe that 
    \begin{align*}
    \int \chi(s)d\chi&=\int (\chi_0 \chi)(s) d\chi ~~(\textrm{by the left invariance of the Haar measure})\\
    &= \int \chi_0(s) \chi(s) d\chi \\
    &= \chi_0(s) ( \int \chi(s) d\chi).
        \end{align*}
    Since $\chi_0(s) \neq 1$, it follows that $\int \chi(s) d\chi=0$. This proves that $E(au_s)=0$ if $s \neq e$. Hence $A=B^{\beta}$. 
    
    \begin{rmrk}
    If $b=\sum_{s \in G}a_su_s$ then $E(b)=a_{e}$. 
        \end{rmrk}
    Let us turn our attention back to $A_{\theta}$. For the rest of this section, assume that $\theta$ is irrational. Write $A_{\theta}=C(\mathbb{T}) \rtimes \mathbb{Z}$. Denote the generating unitary of $C(\mathbb{T})$ by $u$, and the unitary corresponding to the generator $1$ of the group $\mathbb{Z}$ by $v$. Then, $uv=e^{2\pi i \theta}vu$. Let $E:A_\theta \to C(\mathbb{T})$ be the conditional expectation constructed out of the dual action of $\widehat{\mathbb{Z}}=\mathbb{T}$ on $A_{\theta}$.     
        For $n \geq 1$, let $E_{n}: A_\theta \to A_\theta$ be defined by
        \[
        E_{n}(x):=\frac{1}{n+1}\sum_{k=0}^{n}u^{k}xu^{*k}.\]
        The crucial fact that we need to conclude the simplicity of $A_\theta$ is the following. 
  
  \begin{lmma}
  For $x \in A_{\theta}$, $E_n(x) \to E(x)$. 
    \end{lmma}
    \textit{Proof.} It suffices to check that for $x=u^{r}v^{s}$ with $r,s \in \mathbb{Z}$, $E_n(x) \to E(x)$. This is because the sequence $\{E_n\}_{n \geq 1}$ is norm bounded. Let $r,s \in \mathbb{Z}$ be given and let $x=u^{r}v^{s}$. It is clear that if $s=0$, $E_n(x)=u^{r}=E(x)$. It suffices to consider the case when $s \neq 0$ which we henceforth assume.  Set $z=e^{2 \pi i s \theta}$. Since $\theta$ is irrational, it follows that $z \neq 1$.  
    Now a simple calculation using the relation $u^{k}v^{s}=e^{2 \pi i ks \theta}v^{s}u^{k}$ implies that 
    \[
    E_{n}(x)=\frac{1}{n+1}(\sum_{k=0}^{n}e^{2 \pi i ks \theta})u^{r}v^{s}=\frac{1}{n+1}\Big(\frac{1-z^{n+1}}{1-z}\Big)u^{r}v^{s}.\] 
    Thus, as $n \to \infty$, $E_n(x) \to 0 = E(x)$. This completes the proof. \hfill $\Box$
    
    \begin{thm}
    Let $\theta$ be an irrational number. The $C^{*}$-algebra $A_{\theta}$ is simple. 
        \end{thm}
        \textit{Proof.} Let $I \subset A_{\theta}$ be a closed two sided non-zero ideal of $A_{\theta}$. 
        Denote $E(I)$ be $J$.  A consequence of the previous lemma is that $J \subset I$.  Note that $J=I \cap C(\mathbb{T})$. 
        Since $E$ is faithful, it follows that $J$ is non-zero. Moreover, the fact that $E$ is a conditional expectation implies that $J$ is a two sided ideal of $C(\mathbb{T})$. 
        
        For $x \in J$, $\alpha^{k}(x)=v^{k}xv^{*k} \in I$. Clearly $\alpha^{k}(x) \in C(\mathbb{T})$. Thus, $J$ is an $\alpha$-invariant ideal of $C(\mathbb{T})$. In other words, $\alpha(J)=J$. 
    Let $F \subset \mathbb{T}$ be a  closed subset such that \[J=\{f \in C(\mathbb{T}): f \textrm{~vanishes on $F$}\}.\] The fact that $\alpha(J)=J$ implies that $e^{2\pi i k\theta}F=F$ for every $k \in \mathbb{Z}$. 
    It is well known that for every $x_0$, $\{e^{2 \pi i k\theta}x_0: k \in \mathbb{Z}\}$ is dense in $\mathbb{T}$. We claim that $F$ is empty. Suppose not. Since $F$ is closed and $e^{2 \pi i k \theta}F=F$, we have $F=\mathbb{T}$. 
    This however  forces $J=0$ which is a contradiction. Hence $F=\emptyset$. Consequently, $J=C(\mathbb{T}) \subset I$. But then the ideal generated by $C(\mathbb{T})$ is $A_{\theta}$. Hence $I=A_{\theta}$. This completes the proof. \hfill $\Box$
  
\section{Mackey's imprimitivity theorem : the discrete case}
  This section is devoted to a discussion on Mackey's imprimitivity theorem cast in Rieffel's language of Hilbert $C^{*}$-modules. We only give a proof in the discrete setting and refer the reader to the monographs \cite{Raeburn_Williams} and \cite{Williams_Dana} for the topological case. In Rieffel's language, Mackey's imprimitivity theorem reads as follows. 
  \begin{thm}
  \label{Rieffel-Mackey}
  Let $G$ be a second countable locally compact topological group and $H$ be a closed subgroup of $G$. The crossed product $C_{0}(G/H) \rtimes G$ is  Morita equivalent to  $C^{*}(H)$.
    \end{thm}
  
   First we proceed towards defining  the notion of strong Morita equivalence due to Rieffel. 
   \begin{dfn}
   \label{Imprimitivity}
   Let $A$ and $B$ be $C^{*}$-algebras. An $A$-$B$ imprimitivity bimodule is a vector space $E$ which has a left $A$-action and a right $B$-action together with $A$-valued and $B$-valued inner products satisfying the following
   \begin{enumerate}
   \item[(1)] the $A$-valued inner product is linear in the first variable and conjugate linear in the second variable,
   \item[(2)] the $B$-valued inner product is linear in the second variable and conjugate linear in the first variable,
   \item[(3)] for $x,y \in E$ and $a \in A$, $\langle ax|y \rangle_{B}=\langle x|a^{*}y\rangle_{B}$,
   \item[(4)] for $x,y \in E$ and $b \in B$, $\langle xb|y \rangle_{A}=\langle x|yb^{*} \rangle_{A}$, 
   \item[(5)] for $x,y,z \in E$, $\langle x|y \rangle_{A} z=x \langle y|z \rangle_{B}$, 
   \item[(6)] the linear span of $\{\langle x|y \rangle_{A}:x,y \in E\}$ is dense in $A$, 
      \item[(7)] the linear span of $\{\langle x|y\rangle_{B}: x,y \in E\}$ is dense in $B$, and
      \item[(8)] $E$ is complete with respect to the norm induced by both the $A$-valued and the $B$-valued inner products. 
             \end{enumerate}
    \end{dfn}
    Let $E$ be an $A$-$B$ imprimitivity bimodule. For $x \in E$, define \begin{align*}
        ||x||_{A}:&=||\langle x|x \rangle_{A}||^{\frac{1}{2}}\\
    ||x||_{B}:&=||\langle x|x \rangle_{B}||^{\frac{1}{2}}. 
  \end{align*}
  
  \begin{ppsn}
  \label{equality of norms}
  With the foregoing notation, we have $||x||_{A}=||x||_{B}$ for every $x \in E$. 
    \end{ppsn}
    \textit{Proof.} View $E$ as a Hilbert $B$-module. For $a \in A$, let $L_{a}:E \to E$ be defined by $L_{a}(x)=a.x$. Then for every $a \in A$, $L_{a}$ is adjointable and the map $A \ni a \to L_{a} \in \mathcal{L}_{B}(E)$ is a $^*$-homomorphism. 
    We claim that $a \to L_{a}$ is injective. Suppose $L_{a}=0$. Then $\langle ax|y \rangle_{A}=0$ for every $x,y \in E$. Consequently $a \langle x|y \rangle_{A}=0$. But $\{\langle x|y \rangle_{A}: x,y \in E\}$ is dense in $A$. Hence $ab=0$ for every $b \in A$. This shows that $a=0$. This proves our claim. 
    
    Let $x \in E$ be given. Define $\theta_{x,x}:E \to E$ by $\theta_{x,x}(y)=x \langle x|y \rangle_{B}$. Note that $\theta_{x,x}=L_{\langle x| x \rangle_A}$. Thus to complete the proof, it suffices to show that for $x \in E$, 
    \[
    ||\theta_{x,x}||=||x||_{B}^{2}.\]
    It follows from Cauchy-Schwarz inequality that $||\theta_{x,x}|| \leq ||x||_{B}^{2}$. Set $y:=\frac{x}{||x||_{B}}$ and calculate as follows to observe that 
    \begin{align*}
    ||\theta_{x,x}(y)||^{2}&=||\langle x\langle x|y \rangle| x \langle x|y \rangle \rangle||\\
    &= \frac{1}{||x||^{2}}||\langle x|x \rangle \langle x|x \rangle \langle x|x \rangle|| \\
    &=\frac{1}{||x||^{2}}||x||^{6} \\
    &=||x||^{4}.
       \end{align*}
    Hence $||x||_{B} ^{2} \leq ||\theta_{x,x}||$. Consequently, $||\theta_{x,x}||=||x||_{B}^{2}$. This completes the proof. \hfill $\Box$
   
   \begin{dfn}
   Let $A$ and $B$ be $C^{*}$-algebras. We say that $A$ and $B$ are  Morita equivalent if there exists an $A$-$B$ imprimitivity bimodule. 
    \end{dfn}
    
    \begin{xmpl}
    Let $A$ be a $C^{*}$-algebra. Set $E:=A$ and $B:=A$. Then $E$ is a Hilbert $B$-module. The $C^{*}$-algebra $A$ acts on $E$ by left multiplication. Define an $A$-valued (left) inner product on $E$ by
    \[
    \langle x|x \rangle_{A}=xy^{*}.\]
    Then $E$ is an $A$-$A$ imprimitivity bimodule. 
    \end{xmpl}
    
    \begin{xmpl}
    Let $A$ be a $C^{*}$-algebra and $p \in A$ be a projection. Suppose the ideal generated by $p$ is $A$. Let $B:=pAp$. Set $E:=pA$. Define a $B$-valued inner product on $E$ by 
    \[
    \langle x|y \rangle=xy^{*}.\]
    Then $E$ is a $B$-$A$ imprimitivity bimodule. 
        \end{xmpl}
    
    \textbf{The algebra of compact operators of a Hilbert module:} Let $E$ be a Hilbert $B$-module. For $x,y \in E$, let $\theta_{x,y}:E \to E$ be defined by 
    \[
    \theta_{x,y}(z)=x \langle y|z \rangle.\]
    Note that for $x,y \in E$, $\theta_{x,y}$ is an adjointable operator and $\theta_{x,y}^{*}=\theta_{y,x}$. Moreover for $T \in \mathcal{L}_{B}(E)$, $T\theta_{x,y}=\theta_{Tx,y}$ and $\theta_{x,y}T=\theta_{x,T^{*}y}$. The $C^{*}$-algebra of compact operators on $E$, denoted $\mathcal{K}_{B}(E)$, is defined to be the closed linear span of $\{\theta_{x,y}:x,y \in E\}$. 
    
    \begin{rmrk}
    Let $A$ be a $C^{*}$-algebra and consider the Hilbert $A$-module $E:=A$. For $a \in A$, let $L_{a}:E \to E$ be defined by $L_{a}(x)=ax$.     
    Note that for $x,y \in E$, $\theta_{x,y}=L_{xy^{*}}$. This implies that $A \ni a \to L_{a} \in \mathcal{K}_{B}(E)$ is an isomorphism. 
    
    Suppose $A$ is unital. Then $L_{1}$ is the identity operator which is not compact in the sense of Banach space theory unless the algebra $A$ is finite dimensional. Thus compact operators in the sense of Hilbert $C^{*}$-modules need not be compact in the usual Banach space theory sense. 
        \end{rmrk}
    
    \begin{xrcs}
    Let $E:=A^{n}$ be n copies of $A$. Show that $M_{n}(A)$ is isomorphic to $\mathcal{K}_{A}(A^{n})$. 
       \end{xrcs}
   
   Let $E$ be a Hilbert $B$-module and set $A:=\mathcal{K}_{B}(E)$. The $C^{*}$-algebra $A$ acts on $E$ on the left by the formula: $T.x=Tx$ for $T \in \mathcal{K}_{B}(E)$ and $x \in E$. Define an $A$-valued inner product on $E$ by
   \[
   \langle x|y \rangle_{A}=\theta_{x,y}.\]
   The proof of Proposition \ref{equality of norms} imply that for $x \in E$, $||\theta_{x,x}||=||x||^{2}$. It is routine to verify that $E$ satisfies all the axioms, except $(7)$, of Definition \ref{Imprimitivity}. Note that for a Hilbert $B$-module $E$, the linear span of $\{\langle x|y \rangle: x,y \in E\}$ is always a two sided ideal. 
   
   \begin{dfn}
   Suppose $E$ is a Hilbert $B$-module. The Hilbert module $E$ is said to be full if the closed linear span of $\{\langle x|y \rangle: x,y \in E\}$ is $B$. 
    \end{dfn}

   \begin{ppsn}
   Let $A$ and $B$ be $C^{*}$-algebras. The following are equivalent.
   \begin{enumerate}
   \item[(1)] The $C^{*}$-algebras $A$ and $B$ are Morita equivalent. 
   \item[(2)] There exists a full Hilbert $B$-module and a faithful representation $\phi:A \to \mathcal{L}_{B}(E)$ such that $\phi(A)=\mathcal{K}_{B}(E)$. 
      \end{enumerate}
      \end{ppsn}
      \textit{Proof.} Suppose $(1)$ holds.  Let $E$ be an $A$-$B$ imprimitivity bimodule. For $a \in A$, let $\phi(a):E \to E$ be defined by $\phi(a)(x)=ax$. In the proof of Prop. \ref{equality of norms}, we observed that $A \ni a \to \phi(a) \in \mathcal{L}_{B}(E)$ is injective. 
      Note that for $x,y \in E$, $\phi(\langle x|y \rangle_{A})=\theta_{x,y}$. Since $\{\langle x|y \rangle_{A}:x,y \in E\}$ is dense in $A$, it follows that the image of $\phi$ is $\mathcal{K}_{B}(E)$. Axiom $(7)$ implies that $E$ is a full Hilbert $B$-module. 
      
      We have already observed that if $E$ is a full Hilbert $B$-module then $E$ is a $\mathcal{K}_{B}(E)$-$B$ imprimitivity bimodule. Thus $(2) \implies (1)$ is clear. \hfill $\Box$
   
   Next we show that Morita equivalence is indeed an equivalence relation on $C^{*}$-algebras.
   \begin{ppsn}
   Morita equivalence is an equivalence relation.
     \end{ppsn}
     \textit{Proof.} We have already observed that $A$ is an $A$-$A$ imprimitivity bimodule. Suppose $E$ is an $A$-$B$ imprimitivity bimodule. Denote the conjugate vector space by $\overline{E}$. Then as a set $\overline{E}$ is just $E$. For an element $x \in E$, when we regard $x$ as an element of $\overline{E}$, we write $j(x)$ for $x$. The addition and scalar multiplication are defined by 
     \begin{align*}
     j(x)+j(y)&=j(x+y)\\
     \lambda.j(x)&=j(\overline{\lambda}x).
          \end{align*}
          Then $\overline{E}$ is a $B$-$A$ imprimitivity bimodule where the right action of $A$, the left action of $B$ and the inner products are given by 
          \begin{align*}
          j(x).a&=j(a^{*}x) \\
          b.j(x)&=j(xb^{*})\\
          \langle j(x)|j(y) \rangle_A&=\langle x|y \rangle_{A} \\
          \langle j(x)|j(y) \rangle_{B}&=\langle x|y \rangle_B
            \end{align*}
   for $x,y \in E$ and $a \in A$, $b \in B$. 
   Suppose $E$ is an $A$-$B$ imprimitivity bimodule and $F$ is a $B$-$C$ imprimitivity bimodule then the interior tensor product $E \otimes_{B} F$ is an $A$-$C$ imprimitivity bimodule. \hfill $\Box$       
   
   \begin{rmrk}
   Let $E$ be a $A$-$B$ imprimitivity bimodule and $\overline{E}$ be a conjugate $B$-$A$ imprimitivity bimodule constructed in the previous proposition. Then the maps 
   \[E \otimes_{B} \overline{E} \ni x \otimes j(y) \to \langle x|y \rangle_{A} \in A \] and \[\overline{E} \otimes_{A} E \ni j(x) \otimes y \to \langle x|y \rangle_{B} \in B\] are isomorphisms of Hilbert modules.
   Thus $E \otimes_{B} \overline{E} \cong A$ and $\overline{E} \otimes_{A} E \cong B$. 
      Thus Morita equivalent $C^{*}$-algebras have the same representation theory (See Remark \ref{primitive Morita equivalence}).
    \end{rmrk}
   
   In practice, imprimitivity bimodules are always constructed by the process of completion. The setup we usually have is as follows. Let $A_0$ be a dense $C^{*}$-subalgebra of $A$ and $B_0$ be a dense $C^{*}$-subalgebra of $B$ where $A$ and $B$ are $C^{*}$-algebras. 
   \begin{dfn}
   \label{pre}
   A  pre $A_0$-$B_0$ imprimitivity bimodule is a vector space $E_0$ which is a $A_0$-$B_0$ bimodule with an $A_0$-valued and a $B_0$-valued semi-definite inner products such that 
   \begin{enumerate}
   \item[(1)] the $A_0$-valued inner product is linear in the first variable and conjugate linear in the second variable,
   \item[(2)] the $B_0$-valued inner product is linear in the second variable and conjugate linear in the first variable,
   \item[(3)] for $x,y \in E_0$ and $a \in A_0$, $\langle ax|y \rangle_{B_0}=\langle x|a^{*}y \rangle_{B_0}$,
   \item[(4)] for $x,y \in E_0$ and $b \in B_0$, $\langle xb| y \rangle_{A_0}=\langle x|yb^{*} \rangle_{A_0}$, 
   \item[(5)] for $x,y,z \in E_0$, $\langle x|y\rangle_{A_0}z=x \langle y|z \rangle_{B_0}$, 
   \item[(6)] for $x \in E_0$, $a \in A_0$ and $b \in B_0$, $\langle ax|ax \rangle_{B_0} \leq ||a||^{2} \langle x|x \rangle_{B_0}$ and $\langle xb|xb \rangle_{A_0} \leq ||b||^{2} \langle x| x \rangle_{A_0}$, and
   \item[(7)] the set $\{\langle x|y\rangle_{A_0}:x,y \in E_0\}$ and the $\{\langle x|y \rangle_{B_0}:x,y \in E_0\}$ span dense ideals in $A$ and $B$ respectively. 
      \end{enumerate}
   \end{dfn}
   \begin{ppsn}
   Let $E_0$ be a $A_0$-$B_0$ imprimitivity bimodule. For $x \in E$, \[||\langle x| x \rangle_{A_0}||=||\langle x|x \rangle_{B_0}||.\] 
   
   \end{ppsn}
   \textit{Proof.} Let $x \in E_0$ be given. Let $a=\langle x|x \rangle_{A_0}$. Calculate as follows to observe that 
   \begin{align*}
   ||a||^{2} \langle x|x \rangle_{B_0} &\geq \langle ax|ax \rangle_{B_0} \\
   & \geq \langle \langle x|x \rangle_{A_0}x|\langle x|x \rangle_{A_0}x \rangle_{B_0} \\
      & \geq \langle x \langle x|x \rangle_{B_0}|x \langle x|x \rangle_{B_0} \rangle_{B_0}  \\
      & \geq \langle x|x \rangle_{B_0}^{3}.
   \end{align*}
   Taking norms and cancelling $||\langle x|x \rangle_{B_0}||$, we get $|| \langle x|x \rangle_{A_0} \geq ||\langle x|x \rangle_{B_0}||$. A similar argument yields $||\langle x|x \rangle_{B_0}|| \geq ||\langle x|x \rangle_{A_0}||$. This completes the proof. \hfill $\Box$
   
   \begin{rmrk}
   Suppose $E_0$ is a pre $A_0$-$B_0$ imprimitivity bimodule, first we mod out the null vectors and then complete to obtain a genuine $A$-$B$ imprimitivity bimodule. The previous proposition implies that the null vectors of $E_0$ are the same whether we give the norm induced by the $A_0$-valued inner product or the $B_0$-valued inner product. 
      \end{rmrk}
      
      We proceed towards proving Theorem \ref{Rieffel-Mackey} in the discrete setting. For the rest of this section, assume that $G$ is a discrete countable group and $H \subset G$ be a subgroup. Denote the set of left cosets of $H$ by $G/H$. The group $G$ acts on $G/H$ by left translations. Let $\alpha:=\{\alpha_{s}\}_{s \in G}$ be the action of $G$ on $C_{0}(G/H)$ induced by the left translation of $G$ on $G/H$. 
      
      Let us fix notation which will be used throughout. 
      Let $A:=C_{0}(G/H) \rtimes G$ and $B:=C^{*}(H)$. Denote the generating unitaries of $B$ by $\{v_{t}:t \in H\}$.  For $a \in C_{0}(G/H)$ and $s \in G$, let  $a \otimes\delta_{s} \in C_{c}(G,C_0(G/H))$ be the function whose value at $s$ is $a$ and vanishes elsewhere.  For $s \in G$, let $e_{sH} \in C_{c}(G/H)$ be the characteristic function at $sH$. Note that $\alpha_{s}(e_{tH})=e_{stH}$ for $s,t \in G$. Let $A_0$ be the linear span of $\{e_{sH} \otimes \delta_{t}: s \in G, t \in G\}$ and $B_0$ be the linear span of $\{v_{t}:t \in H\}$. Then $A_0$ and $B_0$ are dense $*$-subalgebras of $A$ and $B$ respectively. Also note that $\{v_{t}:t \in H\}$ and $\{e_{rH} \otimes \delta_{s}: r,s \in G\}$ form a basis for $B_0$ and $A_0$ respectively. 
      
      Let $E_0:=C_{c}(G)$ and let $\{\epsilon_{s}:s \in G\}$ be the standard basis for $E_0$. Define a left $A_0$ action and a right $B_0$ action on $E_0$ by
      \begin{align*}
      (e_{rH} \otimes \delta_{s}).\epsilon_t:&=1_{rH}(st)\epsilon_{st} \\
      \epsilon_{t}.v_{s}&=\epsilon_{ts}.
      \end{align*}
   Define an $A_0$-valued sesquilinear form (by extending linearly in the first variable) and a $B_0$-valued sesquilinear form (by extending linearly in the second variable) by
   \begin{align*}
   \langle \epsilon_{s}|\epsilon_{t}\rangle_{B_0}&=1_{H}(s^{-1}t)v_{s^{-1}t} \\
   \langle \epsilon_{s}|\epsilon_{t} \rangle_{A_0}&=e_{sH} \otimes \delta_{st^{-1}}.
      \end{align*}
      Theorem \ref{Rieffel-Mackey}, in the discrete case, follows from the next theorem. 
      \begin{thm}
      With the foregoing notation, $E_0$ is a pre $A_0$-$B_0$ imprimitivity bimodule. 
            \end{thm}
      \textit{Proof.} First we show that the sesquilinear forms defined are indeed positive semi-definite. 
            Let us first deal with the $B_0$-valued sesquilinear form. 
      Let $x:=\sum_{s \in G}a_{s}\epsilon_{s} \in E_0$ be given. Let $F:=\{s \in G: a_s \neq 0\}$. Then $F$ is a finite subset of $G$. Define an equivalence relation on $F$ by
      for $s_1,s_2 \in F$, $s_1 \sim s_2$ if and only if $s_1H=s_2H$. For $s \in F$, let $[s]$ be the equivalence class containing $s$. List the equivalence classes as 
      $[s_1], [s_2],\cdots,[s_m]$. Write $[s_i]=\{s_i h_{ij}:j = 1,2,\cdots,k_i\}$ for every $i=1,2,\cdots,m$.  
      
      Calculate as follows to observe that 
      \begin{align*}
      \langle x|x \rangle_{B_0}&=\sum_{i=1}^{m}\Big(\sum_{r=1}^{k_i}\overline{a_{s_ih_{ir}}}a_{s_ih_{is}}v_{h_{ir}^{-1}h_{is}}\Big) \\
      &=\sum_{i=1}^{m}\Big((\sum_{j=1}^{k_i}a_{s_ih_{ij}}v_{h_{ij}})^{*}(\sum_{j=1}^{k_i}a_{s_ih_{ij}}v_{h_{ij}})\Big)\\
      & \geq 0.
           \end{align*}
      This shows that the $B_0$-valued sesquilinear form is a semidefinite inner product. 
      
      Let $x:=\sum_{s \in G}a_{s}\epsilon_s \in E_0$ be given. Fix a non-degenerate representation of the crossed product $C_{0}(G/H) \rtimes G$. In other words, fix a covariant representation $(\pi,U)$ of the $C^{*}$-dynamical system $(C_0(G/H),G,\alpha)$. 
      Calculate as follows to observe that 
      \begin{align*}
      (\pi \rtimes U)(\langle x|x \rangle_{A_0})&=\sum_{s,t \in G}a_s\overline{a_t}\pi(e_{sH})U_{st^{-1}}\\
      &=\sum_{s,t \in G}a_{s}\overline{a_t}U_s\pi(e_H)U_{t}^{*} ~~(\textrm{since $U_{s}\pi(e_H)U_s^{*}=\pi(e_{sH}$}) )\\
      &= \big(\sum_{s \in G}a_{s}U_{s}\pi(e_H)\big)\big(\sum_{s\in G}a_sU_s\pi(e_H)\big)^{*}\\
      & \geq 0.
            \end{align*}
      Since $(\pi \rtimes U)(\langle x|x \rangle_{A_0}) \geq 0$ for every covariant representation $(\pi,U)$, it follows that $\langle x|x \rangle_{A_0} \geq 0$ in $A$. 
            The verifications of the axioms, except Axiom $(6)$, of Defn. \ref{pre} are routine and we leave the verification to the reader. 
      
      View $E_0$ as a pre Hilbert $B_0$-module. Mod out the null vectors and complete to obtain a Hilbert $B$-module $E$. For $s \in G$ and $x:=\sum_{t \in G}a_{t}\epsilon_{t}$, let $U_{s}(x)=\sum_{t \in G}a_{t}\epsilon_{st}$. Note that for $x,y \in E_0$, 
      \begin{align*}
      \langle U_{s}x|U_{s}x \rangle_{B_0}&=\langle x|x \rangle_{B_0} \\
      \langle U_{s}x|y \rangle_{B_0}&=\langle x|U_{s^{-1}}y \rangle_{B_0}.
            \end{align*}
      Thus there exists a unique adjointable operator on $E$, which again denote by $U_{s}$, such that $U_{s}\epsilon_{t}=\epsilon_{st}$. 
      
      For $x=\sum_{t \in G}a_{t}\epsilon_{t} \in E_0$, define $Px=\sum_{t \in G}1_{H}(t)a_{t}\epsilon_{t}$. It is clear that $P^{2}x=Px$ and $\langle Px|y \rangle=\langle x|Py\rangle$ for $x,y \in E_0$. For $x \in E_0$, calculate as follows to observe that 
      \begin{align*}
      \langle x|x \rangle_{B_0}&=\langle (1-P)x+Px|(1-P)x+x \rangle_{B_0} \\
      &= \langle (1-P)x|(1-P)x \rangle_{B_0}+ \langle Px|Px \rangle_{B_0}\\
      & \geq \langle Px|Px \rangle_{B_0}.
       \end{align*}
      The above inequality implies that there exists a unique adjointable operator, again denoted $P$, such that $P\epsilon_{t}=1_{H}(t)\epsilon_t$. For $s \in G$, set $P_{sH}=U_{s}PU_{s}^{*}$. Note that $P_{sH}$ is a projection. Then \[P_{sH}(\epsilon_t)=1_{sH}(t)\epsilon_t.\]
      Hence it follows that $P_{sH}P_{tH}=1_{H}(s^{-1}t)$. Making use of Proposition \ref{orthogonal projections universal}, we conclude that there exists a unique $^*$-homomorphism $\pi:C_{0}(G/H) \to \mathcal{L}_{B}(E)$ such that $\pi(e_{sH})=P_{sH}$. It is routine to verify that $(\pi, U)$ is a covariant representation of the dynamical system $(C_{0}(G/H),G,\alpha)$. Also for $a \in A_0$ and $x \in E_0$, $a.x =(\pi \rtimes U)(a)x$. Calculate as follows to observe that for $a \in A_0$ and $x \in E_0$, 
      \begin{align*}
      \langle ax|ax \rangle_{B_0} & = \langle (\pi \rtimes U)(a)x|(\pi \rtimes U)(a)x \rangle_{B} \\
      & \leq ||a||^{2} \langle x|x \rangle_{B} \\
      & \leq ||a||^{2} \langle x|x \rangle_{B_0}.
            \end{align*}
      The verification of the second half of Axiom $(6)$ is similar and therefore relegated to an exercise. \hfill $\Box$
      
      \begin{xrcs}
      Verify the second half of Axiom $(6)$ and complete the proof of the previous Theorem. 
            \end{xrcs}
       \begin{rmrk}     
            For examples and applications of Mackey's imprimitivity theorem, we recommend Tyrone Crisp's notes  available online at www.math.ru.nl/{~}tcrisp.
      \end{rmrk}
  
  

  
  \chapter{K-theory}
  \section{$K_0$ of a $C^{*}$-algebra}
  In this chapter, we  give a basic introduction to the subject of K-theory. We will not give complete proofs of many results and give only rough sketches at many places.
  The reader interested in a detailed development should consult \cite{Bla}, \cite{Rordam} or \cite{Olsen}.  
  
  Let $A$ be a unital algebra over $\mathbb{C}$. Let $E$ be a right $A$-module. We say that $E$ is finitely generated and projective if there exists a right $A$-module $F$ and a positive integer $n \geq 1$ such that 
  $E \oplus F \cong A^{n}$.   Denote the set of isomorphism classes of finitely generated projective right $A$-modules \footnote{We only consider right $A$-modules.} by $\mathcal{V}(A)$. Then, $\mathcal{V}(A)$ is an abelian semigroup with identity. 
  First we obtain a better description of $\mathcal{V}(A)$ in terms of idempotents.  We always think of elements of $A^{n}$ as column vectors. 
  
  \begin{xrcs}
  Let $m,n \geq 1$ be given. For $x \in M_{m \times n}(A)$, let $T_{x}:A^{n} \to A^{m}$ be defined by $T_{x}(v)=xv$. Show that the map 
  \[
  M_{m \times n}(A) \ni x \to T_{x} \in \mathcal{L}_{A}(A^{n},A^{m})\]
  is an isomorphism. Here, $\mathcal{L}_{A}(A^{n},A^{m})$ denotes the abelian group of $A$-linear maps from $A^{n}$ to $A^{m}$. 
  \end{xrcs}
  
  \begin{ppsn}
  \label{idempotents}
  We have the following.
  \begin{enumerate}
  \item[(1)] For an idempotent $e \in M_{n}(A)$, $eA^{n}$ is a finitely generated projective $A$-module. 
  \item[(2)] Let $M$ be a finitely generated projective $A$-module. Then, there exists a natural number $n$ and an idempotent $e \in M_{n}(A)$ such that $M$ is isomorphic to $eA^{n}$.
  \item[(3)] Let $e \in M_{m}(A)$ and $f \in M_{n}(A)$ be such that $e$ and $f$ are idempotents. Then, $eA^{m}$ and $fA^{n}$ are isomorphic as $A$-modules if and only if there exist $x \in M_{m \times n}(A)$ and $y \in M_{n \times m}(A)$ such that 
  $xy=e$ and $yx=f$. 
 \end{enumerate}
  \end{ppsn}
  \textit{Proof.} Let $e \in M_{n}(A)$ be an idempotent. Clearly, $eA^{n}\oplus (1-e)A^{n}=A^{n}$. Therefore $eA^{n}$ is a finitely generated projective $A$-module. This proves $(1)$. Let $M$ be a finitely generated projective $A$-module. Choose an $A$-module $N$ such that $M \oplus N=A^{n}$ for some $n$. Let $T:A^{n} \to A^{n}$ be the map defined by $T(m\oplus n)=m$. Then, $T$ is clearly an idempotent. Hence there exists $e \in M_{n}(A)$ such that $e$ is an idempotent and $T$ is given by left multiplication by $e$. Note that $M=eA^{n}$. This proves $(2)$. 
  
  Let $e \in M_{m}(A)$ and $f \in M_{n}(A)$ be such that $e$ and $f$ are idempotents. Suppose that $eA^{n}$ and $fA^{m}$ are isomorphic. Let $T:fA^{n} \to eA^{m}$ be an isomorphism and let $S$ be the inverse of $T$. Decompose $A^{n}$ as $A^{n}=fA^{n}\oplus (1-f)A^{n}$ and $A^{m}$ as $A^{m}=eA^{m} \oplus (1-e)A^{m}$. Let $X:A^{n} \to A^{m}$ be defined by $X(u,v)=(Tu,0)$ and $Y:A^{m} \to A^{n}$ be defined by $Y(u,v)=(Su,0)$. Then $XY$ is given by left multiplication by $e$ and $YX$ is given by left multiplication by $f$. Let $x \in M_{m \times n}(A)$ be the matrix corresponding to $X$ and $y \in M_{n \times m}(A)$ be the matrix corresponding to $Y$. Then, $xy=e$ and $yx=f$. This proves the ``only if" part. 
  
  Suppose there exists   $x \in M_{m \times n}(A)$ and $y \in M_{n \times m}(A)$ such that $xy=e$ and $yx=f$. Replacing $x$ by $exf$ and $y$ by $fye$, we can assume that $exf=x$ and $fye=y$. Let $X:A^{n} \to A^{m}$ and $Y:A^{m} \to A^{n}$ be the $A$-linear maps that correspond to $x$ and $y$ respectively. Then, $X$ maps $fA^{n}$ into $eA^{m}$ and $Y$ maps $eA^{m}$ into $fA^{n}$. Clearly, when restricted to $fA^{n}$ and $eA^{m}$, $X$ and $Y$ are inverses of each other. This proves the ``if part". This completes the proof. \hfill $\Box$
  
  In view of Prop. \ref{idempotents}, the semigroup $\mathcal{V}(A)$ can be described as follows. For $n \geq 1$, let $E_{n}(A)$ be the set of idempotents in $M_{n}(A)$. 
  Set 
  \begin{align*}
  M_{\infty}(A):&=\bigcup_{ n = 1}^{\infty}M_{n}(A) \\
  E_{\infty}(A):&=\bigcup_{n=1}^{\infty}E_{n}(A).
    \end{align*}
    Define an equivalence relation on $E_{\infty}(A)$ as follows: For $e \in E_{m}(A)$ and $f \in E_{n}(A)$, we say $e \sim f$ if there exist $x \in M_{m \times n}(A)$ and $y \in M_{n \times m}(A)$ such that $xy=e$ and $yx=f$. Then, $\mathcal{V}(A)=E_{\infty}(A)/\sim$. Moreover the addition operation is as follows. For $e,f \in E_{\infty}(A)$, 
    \[
    e \oplus f= \begin{bmatrix}
    e & 0 \\
    0 & f 
    \end{bmatrix}.\]
    
    In a $C^{*}$-algebra, we can replace idempotents by projections and the equivalence relation then becomes Murray-von Neumann equivalence. Let $A$ be a unital $C^{*}$-algebra. For $n \geq 1$, let $P_{n}(A)$ be the set of projections in $M_{n}(A)$. Set 
    \[
    P_{\infty}(A):=\bigcup_{ n=1}^{\infty}P_{n}(A).\]
  Let $p \in P_{m}(A)$ and $q \in P_n(A)$ be given. We say that $p$ and $q$ are \emph{Murray-von Neumann} equivalent if there exists a partial isometry\footnote{An element $u \in M_{m \times n}(A)$ is said to be a partial isometry if $u^*u$, or equivalently $uu^*$, is a projection. Another equivalent criteria is that $u^*uu^*=u^*$, or equivalently $uu^*u=u$.} $u \in M_{n \times m}(A)$ such that $u^{*}u=p$ and $uu^{*}=q$. 
  
  \begin{xrcs}
  \label{adding zeros}
  Let $p \in P_{m}(A)$ be given. Show that $p$ is Murray-von Neumann equivalent to $\begin{bmatrix}
  p & 0 \\
  0 & 0_{n} 
  \end{bmatrix}$. 
  \end{xrcs}
  \textit{Hint:} Consider the `column vector' $\begin{bmatrix}
   p \\
   0
   \end{bmatrix}$.
  
  \begin{ppsn}
  Let $A$ be a unital $C^{*}$-algebra. 
  \begin{enumerate}
  \item[(1)] Let $e \in E_{\infty}(A)$ be given. Then, there exists $p \in P_{\infty}(A)$ such that $e \sim p$. 
  \item[(2)] Let $p,q \in P_{\infty}(A)$ be given. Then, $p \sim q$ if and only if $p$ and $q$ are Murray-von Neumann equivalent. 
   \end{enumerate}
    \end{ppsn}
  \textit{Proof.} Let $e \in E_{n}(A)$ be given. Without loss of generality, we can assume that $n=1$. Represent $A$ faithfully as bounded operators on a Hilbert space $\clh$ in a unital fashion. Let $p$ be the orthogonal projection onto $Ran(e)=Ker(1-e)$. 
  Decompose $\clh$ as $\clh:=Ran(p) \oplus Ker(p)$. With respect to this decomposition, $e$ has the following matrix form
  \[
  e= \begin{bmatrix}
        1 & x \\
        0 & 0
          \end{bmatrix}.\]
  Set $z:=1+(e-e^{*})(e^{*}-e)$. Then, $z$ is invertible in $A$. A simple matrix calculation implies that \[
  ee^{*}z^{-1}=\begin{bmatrix}
         1+xx^{*} & 0 \\
         0 & 0 
         \end{bmatrix} \begin{bmatrix}
                                 (1+xx^{*})^{-1} & 0 \\
                                 0 & (1+x^{*}x)^{-1}
                                 \end{bmatrix} = \begin{bmatrix}
                                 1 & 0 \\
                                 0 & 0. 
                                 \end{bmatrix}
                                \]
  Hence $p=ee^{*}z^{-1}$. This implies in particular that $p \in A$. Let $x=e$ and $y=p$. Again a direct matrix calculation implies that $xy=p$ and $yx=e$. Hence $e \sim p$ and the proof of $(1)$ is complete. 
  
  Let $p,q \in P_{\infty}(A)$ be given. Suppose $p \sim q$. By adding zeros along the diagonal, we can assume that $p$ and $q$ are of the same size. Again without loss of generality, we can assume $p,q \in A$. Let $x,y \in A$ be such that $xy=p$ and $yx=q$. Replacing $x$ by $pxq$ and $y$ by $qyp$, if necessary, we can assume that $pxq=x$ and $qyp=y$. Represent $A$ faithfully on a Hilbert space in a unital fashion.  Note that $y:p\clh \to q\clh$ and $x:q\clh \to p\clh$ are inverses of each other. Also, $y^{*}$ maps $q\clh$ to $p\clh$ and $x^{*}$ maps $p\clh$ to $q\clh$. Hence, $y^{*}y:p\clh \to p\clh$ is invertible. Moreover $y^{*}y \in pAp$. Hence there exists $r \in pAp$ such that $(y^{*}y)^{\frac{1}{2}}r=r(y^{*}y)^{\frac{1}{2}}=p$. Set $u:=yr$. Note that $y=u(y^{*}y)^{\frac{1}{2}}$. Clearly, $u^{*}u =p$ and $uu^{*} \leq q$. Calculate as follows to observe that 
  \begin{align*}
  q&=qq^{*}\\
    &=(yx)(x^{*}y) \\
    & \leq ||x||^{2}yy^{*}\\
    &= ||x||^{2}u (y^{*}y)^{\frac{1}{2}}(y^{*}y)^{\frac{1}{2}}y^{*} \\
    & \leq ||x||^{2}||y||^{2} uu^{*}.
    \end{align*} 
  Therefore $q \leq uu^{*}$. Consequently, $uu^{*}=q$. This implies that $p$ and $q$ are Murray-von Neumann equivalent. This completes the proof. \hfill $\Box$. 
  
 \textbf{Grothendieck construction:}  The Grothendieck construction allows us to construct an abelian group out of the semigroup $\mathcal{V}(A)$. Let $(R,+)$ be an abelian semigroup with identity $0$. 
 Define an equivalence on $R \times R$ as follows: for $(a,b), (c,d) \in R \times R$, we say $(a,b) \sim (c,d)$ if there exists $e \in R$ such that $a+d+e=b+c+e$. Then, $\sim$ is an equivalence relation on $R \times R$. 
 Denote the set of equivalence classes by $G(R)$. Then, $G(R)$ becomes an abelian group with the addition defined as 
 \[
 [(a,b)]+[(c,d)]=[(a+c,b+d)].\]
 For any $a \in R$, $[(a,a)]$ represents the identity element and the inverse of $[(a,b)]$ is $[(b,a)]$. For $a \in R$, let $[a]:=[(a,0)]$. With this notation, 
 \[
 G(R)=\{[a]-[b]: a,b \in R\}.\]
 Note that $[a]=[b]$ if and only if $a+c=b+c$ for some $c \in R$. 
 
 Let $A$ be a unital algebra over $\mathbb{C}$. The Grothendieck group $G(\mathcal{V}(A))$ is denoted $K_{00}(A)$. Note that 
 \[
 K_{00}(A)=\{[p]-[q]: p,q \in E_{\infty}(A)\}.\]
 Also $[p]=[q]$ if and only if there exists $r \in E_{\infty}(A)$ such that $\begin{bmatrix}
 p & 0 \\
 0 & r 
 \end{bmatrix} \sim \begin{bmatrix}
 q & 0 \\
 0 & r 
 \end{bmatrix}$. 
 
  Let $\phi:A \to B$ be a unital homomorphism. For $n \geq 1$, let $\phi^{(n)}:M_{n}(A) \to M_{n}(B)$ be the amplification of $\phi$, i.e. 
 \[
 \phi^{(n)}((a_{ij}))=(\phi(a_{ij})).\]
 To save notation, we denote $\phi^{(n)}$ again by $\phi$. A moment's reflection with definitions reveal that there exists a unique homomorphism denoted $K_{00}(\phi):K_{00}(A) \to K_{00}(B)$ such that 
 \[
 K_{00}(\phi)([p]-[q])=[\phi(p)]-[\phi(q)].\]
 In short, $K_{00}$ is a covariant functor from the category of unital algebras to the category of abelian groups.

 \begin{rmrk}
 Let $A$ be a unital $C^{*}$-algebra. 
 \begin{enumerate}
 \item[(1)] Let $p,q \in A$ be projections such that $pq=0$. Then $[p+q]=[p]+[q]$. For if we set $u:=(p,q)$ then $u^{*}u=\begin{bmatrix}
  p & 0 \\
  0 & q 
  \end{bmatrix}$ and $uu^{*}=p+q$. 
 \item[(2)] Let $p,q \in P_{n}(A)$ be given. Then $[p]=[q]$ in $K_{00}(A)$ if and only if there exists $m \geq 1$ such that $\begin{bmatrix}
 p & 0 \\
 0 & 1_{m}
 \end{bmatrix}$ is Murray-von Neumann equivalent to $\begin{bmatrix}
 q & 0 \\
 0 & 1_{m}
 \end{bmatrix}$

 The ``if part" is clear. For the "only if" part, suppose $[p]=[q]$ in $K_{00}(A)$. Then there exists $r \in P_{m}(A)$ such that $p \oplus r \sim q \oplus r$. Note that 
 \begin{align*}
 p \oplus 1_{m} & \sim p \oplus (r+1-r)\\
 & \sim p \oplus (r \oplus (1-r)) \\
 &\sim (q \oplus r)\oplus (1-r) \\
& \sim q \oplus (r+1-r) \\
& \sim q \oplus 1_m.
  \end{align*}
  \end{enumerate}
  \end{rmrk}
  
  \begin{xrcs}
  \begin{enumerate}
 \item[(1)] Show that $K_{00}(\mathbb{C})$ is isomorphic to $\mathbb{Z}$ and $[1]$ forms a $\mathbb{Z}$-basis for $K_{00}(\mathbb{C})$. 
 \item[(2)] Let $n \geq 1$. Show that $K_{00}(M_n(\mathbb{C}))=\mathbb{Z}$. Let $p$ be a minimal projection in $M_{n}(\mathbb{C})$. Show that $[p]$ is a $\mathbb{Z}$-basis for $K_{00}(M_n(\mathbb{C}))$.  
  \item[(3)] Let $\clh$ be an infinite dimensional  Hilbert space. Prove that $K_{00}(B(\clh))=0$. 
   \end{enumerate}
    \end{xrcs}
  
  \begin{xrcs}
  \label{direct sum k00}
  Let $A_1$ and $A_2$ be unital algebras and set $A:=A_1\oplus A_2$. Show that the map $K_{00}(\pi_1) \oplus K_{00}(\pi_2): K_{00}(A) \to K_{00}(A_1) \oplus K_{00}(A_2)$ is an isomorphism.

  \end{xrcs}

 Next, we define $K_0$ for a $C^{*}$-algebra. Let $A$ be a $C^{*}$-algebra (unital or non-unital). Set $A^{+}:=\{(a,\lambda): a \in A, \lambda \in \mathbb{C}\}$. Recall that  the addition and scalar multiplication on $A^{+}$ are defined coordinate wise. 
 The multiplication rule is given by 
 \[
 (a,\lambda)(b,\mu)=(ab+\lambda b+\mu a, \lambda \mu).\]
 Let $\epsilon: A^{+} \to \mathbb{C}$ be the map defined by $\epsilon(a,\lambda)=\lambda$. Let $s:A^{+} \to A^{+}$ be defined by $s(a,\lambda)=(0,\lambda)$. The map $s$ is called the ``scalar map" as it  remembers only the scalar part.  We denote the amplifications of $\epsilon$ and $s$ by $\epsilon$ and $s$ itself. 
 
 Define \[
 K_{0}(A):=Ker(K_{00}(\epsilon)).\]
 
  \begin{ppsn}[The standard picture]
  \label{the standard picture}
  Let $A$ be a $C^{*}$-algebra. Then 
  \[
  K_{0}(A)=\{[p]-[s(p)]: p \in P_{\infty}(A^{+})\}.\]
    \end{ppsn}
  \textit{Proof.} It is clear that for $p \in P_{\infty}(A^{+})$, $[p]-[s(p)] \in K_0(A)$. Let $x \in K_0(A)$ be given. Write $x=[p]-[q]$ with $p,q$ projections of same size say of size $n$. The fact that $x \in K_{00}(\epsilon)$ implies that $\epsilon(p)$ and $\epsilon(q)$ are of same rank. Choose a  unitary $u \in M_n(\bbc)$ such that $\epsilon(p)=u\epsilon(q)u^{*}$. Replacing $q$ by $uqu^{*}$, we can assume that $x=[p]-[q]$ with $\epsilon(p)=\epsilon(q)$. Set 
  $
  e:=\begin{bmatrix}
      p & 0 \\
      0 & 1-q
      \end{bmatrix}
      $ and $f:=\begin{bmatrix}
                        0 & 0 \\
                        0 & 1
                        \end{bmatrix}$. 
  Then, $x=[e]-[f]$. Also, note that $[s(e)]=[s(p)]+[1_n-s(q)]=[1_n]=[f]$. Therefore $x=[e]-[s(e)]$. This completes the proof. \hfill $\Box$
 
 \begin{ppsn}
 Let $A$ be a unital $C^{*}$-algebra. Then, $K_0(A)$ is isomorphic to $K_{00}(A)$. 
  \end{ppsn}
 \textit{Proof.} Note that  the map $A^{+} \ni (a,\lambda) \to (a+\lambda 1_{A},\lambda) \in A \oplus \mathbb{C}$ is an isomorphism. With respect to this isomorphism, the map $\epsilon$ becomes the second projection. The result follows immediately from Exercise \ref{direct sum k00}. \hfill $\Box$. 
 
  \textbf{$K_0$ as a functor:} Let $\phi:A \to B$ be a $^*$-algebra homomorphism. The map $\phi$ induces a map $\phi^{+}:A^{+} \to B^{+}$ which is defined as 
  \[
  \phi^{+}((a,\lambda))=(\phi(a),\lambda).\]
  Note that $\epsilon_{B} \circ \phi^{+}=\epsilon_{A}$. Hence $K_{00}(\epsilon_{B})\circ K_{00}(\phi^{+})=K_{00}(\epsilon_{A})$. Therefore, $K_{00}(\phi^{+})$ maps $K_{0}(A)$ to $K_{0}(B)$. We denote the restriction of $K_{00}(\phi^{+})$ to $K_0(A)$ by $K_0(\phi)$. 
  Thus $K_0$ is a functor from the category of $C^{*}$-algebras to the category of abelian groups. 
     The functor $K_0$ is stable, homotopy invariant, half-exact and split-exact. We explain this in what follows. 
  
  \textbf{Stability:} Let $A$ be a $C^{*}$-algebra and $p$ be a minimal projection of $M_{n}(\mathbb{C})$. Let $\omega:A \to M_{n}(A)=A \otimes M_n(\mathbb{C})$ be defined by
  \[
  \omega(a):=a \otimes p.\]
  Then $K_{0}(\omega):K_{0}(A) \to K_{0}(M_n(A))$ is an isomorphism. The reason is a matrix with entries being matrices over $A$ is again a matrix with entries in $A$. We omit the proof and refer the reader to \cite{Rordam}. 
  
  \textbf{Homotopy invariance:} Let $A$ and $B$ be $C^{*}$-algebras and let $\phi, \psi: A \to B$ be $^*$-homomorphisms. We say that $\phi$ and $\psi$ are homotopy equivalent if there exists a family of $^*$-homomorphisms $\phi_{t}:A \to B$ for $t \in [0,1]$ such that 
  \begin{enumerate}
  \item[(1)] for $a \in A$, the map $[0,1] \ni t \to \phi_t(a) \in B$ is norm continuous, and
  \item[(2)] $\phi_0=\phi$ and $\phi_1=\psi$. 
    \end{enumerate}
  The homotopy invariance of $K_0$ means that if $\phi$ and $\psi$ are two homotopy equivalent $^*$-homomorphisms, then $K_0(\phi)=K_0(\psi)$. A moment's thought reveals that this amounts to proving the next lemma. 
  
  \begin{lmma}
  \label{key to homotopy}
  Let $e,f \in A$ be such that $e$ and $f$ are projections. Suppose that $||e-f||<1$. Then, $e$ and $f$ are Murray-von Neumann equivalent. 
    \end{lmma}
   \textit{Proof.} Let $x:=ef$. Note that $||x^{*}x-f||=||f(e-f)f||<1$. Hence $x^{*}x$ is invertible in $fAf$. Choose $r \in fAf$ such that $r(x^{*}x)^{\frac{1}{2}}=(x^{*}x)^{\frac{1}{2}}r=f$. 
   Set $u:=xr$. Then, $u^{*}u=f$. Since $eu =u$, it follows that $uu^{*} \leq e$.  Represent $A$ faithfully on a Hilbert space, say on $\clh$.
   Suppose that $uu^{*}$ is a proper subprojection of $e$. Then, there exists $\xi \in \clh$ such that $e\xi=\xi \neq 0$ but $u^{*}\xi=0$. Hence $rx^{*}\xi=0$. Note that $x^{*}\xi \in Ran(f)$ and $r$ is $1$-$1$ on the range space of $f$. 
   Hence $x^{*}\xi=0$, i.e. $fe\xi=0$. Calculate as follows to observe that 
   \begin{align*}
   ||\xi||&=||e^{2}\xi-fe\xi||\\ 
    &=||(e-f)e\xi|| \\
    &< ||e\xi||=||\xi||
      \end{align*}
   which is a contradiction. Hence, $uu^{*}=e$. This completes the proof. \hfill $\Box$
   
   Let $A$ and $B$ be $C^{*}$-algebras. The $C^{*}$-algebras $A$ and $B$ are said to be \emph{homotopy equivalent} if there exists $^*$-homomorphisms $\phi:A \to B$ and $\psi:B \to A$ such that $\phi \circ \psi$ and $\psi \circ \phi$ are homotopy equivalent to the identity homomorphisms. 
   As an example, consider $A:=C[0,1]$ and $B:=\mathbb{C}$. Define $\epsilon:A \to B$ by $\epsilon(f)=f(0)$ and $\sigma:B \to A$ by $\sigma(\lambda)=\lambda$. Then, clearly $\epsilon \circ \sigma$ is identity and $\sigma \circ \epsilon$ is homotopy equivalent to the identity. 
   
   \begin{xrcs}
   Suppose $A$ and $B$ are homotopy equivalent. Show that $K_0(A)$ and $K_0(B)$ are isomorphic. Conclude that $K_{0}(C(X))=\mathbb{Z}$ for a compact contractible space $X$. 
      \end{xrcs}
   
   The next important property of $K_0$ is that it is half-exact and sends split exact sequences to split exact sequences. 
   \begin{ppsn}
   \label{half exact}
   Let $0 \longrightarrow I \longrightarrow A \stackrel{\pi} \longrightarrow B \longrightarrow 0$ be a short exact sequence of $C^{*}$-algebras. Then, the sequence 
   \[
   K_0(I) \longrightarrow K_0(A) \stackrel{K_0(\pi)} \longrightarrow K_0(B)\]
   is exact in the middle. 
      \end{ppsn}
      \textit{Proof.} Let $x:=[p]-[s(p)] \in K_0(I)$ be given. Since $\pi^{+}(p)=\pi^{+}(s(p))$, it follows that $x \in Ker(K_0(\pi))$. Let $x \in Ker(K_0(\pi))$ be given. Write $x:=[p]-[s(p)]$ with $p$ a projection in $M_n(A^{+})$. Replacing $p$ by $\begin{bmatrix}
      p & 0 \\
      0 & 1_{m}
      \end{bmatrix}$ for large $m$, we can assume that $\pi^{+}(p)$ and $s(\pi^{+}(p))$ are Murray-von Neumann equivalent. 
      
      Let $v$ be a partial isometry in $M_{n}(B^{+})$ be such that $v^{*}v=\pi^{+}(p)$ and $vv^{*}=s(\pi^{+}(p))$. Let $U:=\begin{bmatrix}
      v & 0 \\
      0 & v^{*}
      \end{bmatrix}$. Note that $U\begin{bmatrix}
                                                    \pi^{+}(p) & 0 \\
                                                     0 & 0 
                                                     \end{bmatrix}U^{*}=\begin{bmatrix}
                                                                                     s(\pi^{+}(p)) & 0 \\
                                                                                     0 & 0 
                                                                                     \end{bmatrix}$. Thus, by replacing $p$ by $\begin{bmatrix}
                                                                                      p & 0 \\
                                                                                      0 & 0 
                                                                                      \end{bmatrix}$, we can assume that $\pi^{+}(p)$ And $s(\pi^{+}(p))$ are unitarily equivalent.  
               Let $a \in M_{n}(A^{+})$ be a contraction such that $\pi^{+}(a)=U$. Set 
       \[
       V:=\begin{bmatrix}
              a & \sqrt{1-aa^{*}} \\
              -\sqrt{1-a^{*}a} & a^{*}.\end{bmatrix}\]
              Then $V$ is a unitary and $\pi^{+}(V)=\begin{bmatrix}
              U & 0 \\
              0 & U^{*}
              \end{bmatrix}$. 
      
      Note that $\pi^{+}(V \begin{bmatrix}
       p & 0 \\
       0 & 0 \end{bmatrix}V^{*})$ is a scalar matrix. This implies in particular that $q:= V\begin{bmatrix}
       p & 0 \\
       0 & 0 \end{bmatrix}V^{*}$ lies in $M_{2n}(I^{+}))$.  
       Also the scalar part of $q$  is $\begin{bmatrix}
                                               s(p) & 0  \\
                                               0 & 0
                                               \end{bmatrix}$. Consequently, $x=[q]-[s(q)] \in Im(K_0(i))$ where $i:I \to A$ denotes the inclusion. This completes the proof. \hfill $\Box$                                                                     
   
   Next we show that $K_0$ is split exact. Let \[0 \longrightarrow I \longrightarrow A  \stackrel{\pi} \longrightarrow B \longrightarrow 0\] be a short exact sequence of $C^{*}$-algebras . We say that it is split exact if there exists a $^*$-homomorphism $\mu:B \to A$ such that $\pi \circ \mu=id_{B}$. The map $\mu$ will then be called a splitting. 
   
 \begin{ppsn}
Let $0 \longrightarrow I \longrightarrow A  \stackrel{\pi} \longrightarrow B \longrightarrow 0$ be a split exact sequence of $C^{*}$-algebras with the splitting given by $\mu:B \to A$. Then the sequence 
\[
0 \longrightarrow K_0(I) \longrightarrow K_0(A) \stackrel{K_0(\pi)} \longrightarrow K_0(B) \longrightarrow 0\] is a split exact sequence of abelian groups with the splitting given by $K_0(\mu)$.  
  \end{ppsn}
  \textit{Proof.} Let $i:I \to A$ be the inclusion. We have already shown the exactness at $K_0(A)$. Since $K_0(\pi) \circ K_0(\mu)=Id$, it follows that $K_0(\pi)$ is onto. The only thing that requires proof is that $K_0(i)$ is injective. To that effect, let $x:=[p]-[s(p)] \in K_0(I)$ be such that $x \in Ker(K_0(i))$.
  
  Arguing as in Prop. \ref{half exact}, we can assume that there exists a unitary $u \in M_{n}(A^{+})$ such that $upu^{*}=s(p)$. Set $w:=(\mu^{+} \circ \pi^{+}(u^{*}))u$. Note that $\pi^{+}(w)$ is a scalar. Hence $w \in M_{n}(I^{+})$. Calculate as follows to observe that 
  \begin{align*}
  wpw^{*}&:=(\mu^{+} \circ \pi^{+})(u^{*})upu^{*}(\mu^{+} \circ \pi^{+})(u) \\
  &=(\mu^{+} \circ \pi^{+})(u^{*}s(p)u)\\
  &=(\mu^{+} \circ \pi^{+})(p) \\
  &=s(p) ~~(\textrm{since $p \in I^{+}$}).
  \end{align*}
 This proves that $p$ and $s(p)$ are Murray-von Neumann equivalent in $M_{n}(I^{+})$. Hence $x=0$. This completes the proof. \hfill $\Box$
  
  \begin{xrcs}
  Let $A_1$ and $A_2$ be $C^{*}$-algebras and $A:=A_1 \oplus A_2$. Denote the projection of $A$ onto $A_i$ by $\pi_i$. Show that the map $K_{0}(\pi_1) \oplus K_0(\pi_2): K_0(A) \to K_{0}(A_1) \oplus K_0(A_2)$ is an isomorphism.
   \end{xrcs}
   
  \section{$K_1$ of a $C^{*}$-algebra}
  In this section, we define another functor, denoted $K_1$, from the category of $C^{*}$-algebras to abelian groups. It shares the same functorial properties with $K_0$. This is not a coincidence as we will see later that $K_1$ can indeed be defined in terms of $K_0$. To define $K_1$, we work  with invertible elements or  unitaries. 
  
  Let $A$ be a unital Banach algebra. Denote the set of invertible elements of $M_n(A)$ by $GL_n(A)$. Note that $GL_n(A)$ is a topological group. Denote the connected  component of $1_n$ by $GL_n^{0}(A)$. Then $GL_n^{0}(A)$ is a normal subgroup of $GL_n(A)$. Consider the quotient group $GL_n(A)/GL_n^{(0)}(A)$. There is a natural map from $GL_n(A)/GL_{n}^{0}(A) \to GL_{n+1}(A)/GL_{n+1}^{0}(A)$ given by
  \[
  x \to \begin{bmatrix}
      x & 0 \\
      0 & 1
      \end{bmatrix}.\]
  The group $\widetilde{K}_1(A)$ is defined as the inductive limit  $\displaystyle \lim_{n} GL_n(A)/GL_n^{0}(A)$. 
  
  \begin{xrcs}
  \label{connectedness}
  Show that $GL_n(\mathbb{C})$ is connected. Conclude that $\widetilde{K}_{1}(\mathbb{C})=0$. 
   \end{xrcs}
  Use the previous exercise to show that  for $x \in GL_n(A)$, the elements $\begin{bmatrix}
   x & 0 \\
   0 & 1
   \end{bmatrix}$ and $\begin{bmatrix}
            1 & 0 \\
            0 & x
            \end{bmatrix}$ represent the same element in $GL_{n+1}(A)/GL_{n+1}^{0}(A)$. 
  
  \begin{ppsn}
  Let $A$ be a unital Banach algebra. 
  \begin{enumerate}
  \item[(1)] We have $\widetilde{K}_{1}(A)=\{[x]: x \in GL_n(A), n \geq 1\}$. 
  \item[(2)] For $x,y \in GL_n(A)$, $[x]=[y]$ if and only if there exists  $m$ and a path of invertibles in $GL_{n+m}(A)$ connecting $\begin{bmatrix}
          x & 0 \\
          0 & 1_m
          \end{bmatrix}$ and $\begin{bmatrix}
           y & 0 \\
            0 & 1_m
            \end{bmatrix}$. 
  \item[(3)] The group operation on $\widetilde{K}_{1}(A)$ is given by 
  $[x] \oplus [y]:=\begin{bmatrix}
  x & 0 \\
  0 & y 
  \end{bmatrix}$. 
  \item[(4)] The group $\widetilde{K}_{1}(A)$ is abelian. 
    \end{enumerate}
    \end{ppsn}
  \textit{Proof.} $(1)$ and $(2)$ are just rephrasing the definition of the inductive limit. Statements $(3)$ and $(4)$ follows from Exercise \ref{connectedness}. \hfill $\Box$
  
  \begin{rmrk}
  Let $A$ be a unital $C^{*}$-algebra. Suppose $a \in A$ is invertible. Then $u:=a|a|^{-1}$ is a unitary. Note that $(a|a|^{-t})_{t \in [0,1]}$ is a path of invertible elements connecting $a$ to $u$. Using this it is 
  routine to see that in the definition of $\widetilde{K}_{1}(A)$, we could have taken unitaries in place of invertible elements. 
  We usually work with unitaries in the case of $C^{*}$-algebras. 
  
  We denote the set of unitaries in $M_n(A)$ by $\mathcal{U}_{n}(A)$ and the connected component of $1_n$ by $\mathcal{U}_{n}^{0}(A)$. For unitaries $u,v \in \mathcal{U}_{n}(A)$, we write $u \sim v$ if $u$ and $v$ represent the same element in $\mathcal{U}_{n}(A)/\mathcal{U}_{n}^{0}(A)$. 
    \end{rmrk}
  It is clear that $A \to \widetilde{K}_{1}(A)$ is a functor from the category of unital $C^{*}$-algebras to the category of abelian groups. 
  \begin{xrcs}
  \label{additivity}
  Let $A_1$ and $A_2$ be unital $C^{*}$-algebras and $A:=A_1 \oplus A_2$. Denote the projection of $A$ onto $A_i$ by $\pi_i$. Show that the map 
  $\widetilde{K}_{1}(\pi_1) \oplus \widetilde{K}_{1}(\pi_2): \widetilde{K}_{1}(A) \to \widetilde{K}_{1}(A_1) \oplus \widetilde{K}_{1}(A_2)$ is an isomorphism. 
    \end{xrcs}
  
  For any $C^{*}$-algebra $A$, define $K_1(A):=\widetilde{K}_{1}(A^{+})$. For unital $C^{*}$-algebras, we have $A^{+}=A \oplus \mathbb{C}$. Since $\widetilde{K}_{1}(\mathbb{C})=0$, it follows that $K_1(A)=\widetilde{K}_{1}(A)$. 
  Also $K_1$ is a functor. If $\phi:A \to B$ is a $^*$-homomorphism then there exists a unique group homomorphism $K_1(\phi):K_1(A) \to K_1(B)$ such that 
  \[
  K_1(\phi)([u])=[\phi^{+}(u)].\]
  Next we discuss the functorial properties of $K_1$. 
  
  \textbf{Stability:} Let $A$ be a $C^{*}$-algebra and $p$ be a minimal projection in $M_n(\mathbb{C})$. Let $\omega:A \to A \otimes M_{n}(\mathbb{C})=M_n(A)$ be defined by
  \[
  \omega(a):=a \otimes p.\]
  Then $K_1(\omega):K_1(A) \to K_1(M_n(A))$ is an isomorphism. As with $K_0$, we omit its proof and refer the reader to \cite{Rordam}. 
  
  \textbf{Homotopy invariance:} Let $A$ and $B$ be $C^{*}$-algebras. Suppose $\phi:A \to B$ and $\psi:A \to B$ are $^*$-homomorphisms that are homotopy equivalent. Then $K_1(\phi)=K_1(\psi)$. This is obvious since homotopy invariance is built in the definition of $K_1$. 
  
  \begin{lmma}
  \label{exponential form}
  Let $A$ be a unital $C^{*}$-algebra and $u \in A$ be a unitary. Then $u \in \mathcal{U}^{0}(A)$ if and only if there exists $a_1,a_2,\cdots,a_n \in A$ such that $a_i$'s are self-adjoint and
  $u=e^{ia_1}e^{ia_2}\cdots e^{ia_n}$. 
    \end{lmma}
    \textit{Proof.} For a self-adjoint element $a$, $(e^{ita})_{t \in [0,1]}$ is a path of unitaries connecting $1$ to $e^{ia}$. Thus the ``if part" is clear. Suppose $u \sim 1$. Let $(u_t)_{t \in [0,1]}$ be a path of unitaries such that $u_0=1$ and $u_1=u$. By uniform continuity, there exists a partition $0=t_0<t_1<t_2 < \cdots <t_n=1$ such that $||u_{t_i}-u_{t_{i-1}}||<1$. Set $u_i:=u_{t_i}$. 
    
    We claim $u_1$ is of the required form.   Since $||u_1-1||<1$, it follows that $-1 \notin \sigma(u_1)$. Define $a:=-ilog(u_1)$. Then $u_1=e^{ia}$. This proves the claim.
    Note that $||u_1^{*}u_2-1||=||u_1-u_2||<1$. Applying the above argument, we conclude that $u_1^{*}u_2$ is of the required form. But $u_2=u_1u_1^{*}u_2$. This proves that $u_2$ is the required form.
    Proceeding this way, we see that $u_n=u$ is of the form required. This completes the proof. \hfill $\Box$
  
  \begin{ppsn}
   \label{half exact for K1}
   Let $0 \longrightarrow I \longrightarrow A \stackrel{\pi} \longrightarrow B \longrightarrow 0$ be a short exact sequence of $C^{*}$-algebras. Then the sequence 
   \[
   K_1(I) \longrightarrow K_1(A) \stackrel{K_1(\pi)} \longrightarrow K_1(B)\]
   is exact in the middle. 
  \end{ppsn}
  \textit{Proof.} Let $i:I \to A$ be the inclusion. Let $[u] \in K_0(I)$ be given. Then $\pi^{+}\circ i^{+}(u)$ is a scalar matrix. Consequently, $[\pi^{+} \circ i^{+}(u)]=[1]$. Hence $Im(K_1(i)) \subset Ker(K_1(\pi))$.
  Let $u \in \mathcal{U}_{n}(A^{+})$ be such that $[\pi^{+}(u)]=[1_n]$.  Replacing $u$ by $\begin{bmatrix}
  u & 0 \\
  0 & 1_m
  \end{bmatrix}$ for $m$ sufficiently large, we can assume that $\pi^{+}(u) \sim 1_n$. Choose self-adjoint elements $b_1,b_2,\cdots,b_r \in M_{n}(B^{+})$ such that 
  \[
  \pi^{+}(u)=e^{ib_1}e^{ib_2}\cdots e^{ib_r}.\]
  
  Choose $a_i \in M_{n}(A^{+})$ such that $a_i$ is self-adjoint and $\pi^{+}(a_i)=b_{i}$. Set $v:=e^{ia_1}e^{ia_2}\cdots e^{ia_r}$. Then $\pi^{+}(uv^{*})=1$. This implies that there exists $w \in \mathcal{U}_{n}(I^{+})$ such that 
  $uv^{*}=i^{+}(w)$. Since $[v]=1$, it follows that $[u]=[uv^{*}]=K_{1}(i)([w])$. Hence $Im(K_1(i))=Ker(K_1(\pi))$. This completes the proof. \hfill $\Box$
  
  \begin{ppsn}
  Let $0 \longrightarrow I \longrightarrow A  \stackrel{\pi} \longrightarrow B \longrightarrow 0$ be a split exact sequence of $C^{*}$-algebras with the splitting given by $\mu:B \to A$. Then the sequence 
\[
0 \longrightarrow K_1(I) \longrightarrow K_1(A) \stackrel{K_1(\pi)} \longrightarrow K_1(B) \longrightarrow 0\] is a split exact sequence of abelian groups with the splitting given by $K_1(\mu)$.  
  \end{ppsn}
  \textit{Proof.} Let $i:I \to A$ be the inclusion. We have already shown the exactness at $K_1(A)$. Since $K_1(\pi) \circ K_1(\mu)=Id$, it follows that $K_1(\pi)$ is onto. The only thing that requires proof is that $K_1(i)$ is injective. 
  
  Let $u \in \mathcal{U}_{n}(I^{+})$ be such that $K_1(i)([u])=[1_n]$. By ``amplifying $u$", if necessary, we can assume that $i^{+}(u) \sim 1_n$. Let $(w_{t})_{t \in [0,1]}$ be a path of unitaries in $M_{n}(A^{+})$ such that 
  $w_0=1_n$ and $w_1=i^{+}(u)$. Set $v_{t}:=(\mu^{+} \circ \lambda^{+})(w_t^{*})w_t$. Then $\pi^{+}(v_t)=1$. Hence there exists $u_{t} \in \mathcal{U}_{n}(I^{+})$ such that $i^{+}(u_t)=w_t$. Note that $(u_{t})_{t \in [0,1]}$ is a path 
  of unitaries in $M_{n}(I^{+})$ connecting $1_n$ to $xu$ where $x$ is a scalar matrix. Hence $[u]=[1_n]$. Therefore, $K_1(i)$ is injective. This completes the proof. \hfill $\Box$  

  \begin{xrcs}
  Let $A_1$ and $A_2$ be $C^{*}$-algebras and $A:=A_1 \oplus A_2$. Denote the projection of $A$ onto $A_i$ by $\pi_i$. Show that the map $K_{1}(\pi_1) \oplus K_1(\pi_2): K_1(A) \to K_{1}(A_1) \oplus K_1(A_2)$ is an isomorphism.
   \end{xrcs}
  
  \section{Inductive limits and $K$-theory}
  An important property of $K$-theory that allows to compute the $K$-groups for a large class of $C^{*}$-algebras, called AF-algebras, is that it preserves direct limits. The purpose of this section is to explain this. The data that we require to define the inductive limit of $C^{*}$-algebras is as follows. 
  
  Let $(A_n)_{n \geq 1}$ be a sequence of $C^{*}$-algebras and $\phi_n:A_n \to A_{n+1}$ be a $^*$-homomorphism. The above data is usually given pictorially as follows:
  \[
  A_1 \stackrel{\phi_1} \longrightarrow A_2 \stackrel{\phi_2} \longrightarrow A_3 \longrightarrow \cdots\]
  
  For $m < n$, let $\phi_{n,m}:A_{m} \to A_n$ be defined by $\phi_{m,n}:=\phi_{n-1} \circ \phi_{n-2} \circ \cdots \circ \phi_m$. For $m=n$, set $\phi_{m,n}=Id$. Note that for $\ell \leq m \leq n$, 
  \[
  \phi_{n,\ell}=\phi_{n,m} \circ \phi_{m,\ell}.\]
  
  Let $\mathcal{B}:=\{(a,n): a \in A_n, n \geq 1\}.$ Define an equivalence relation on $\mathcal{B}$ as follows. For $(a,m), (b,n) \in \mathcal{B}$, we say $(a,m) \sim (b,n)$ if there exists $r,s$ such that $m+r=s+n$ and $\phi_{m+r,m}(a)=\phi_{n+s,n}(b)$. Denote the set of equivalence classes by $\mathcal{A}_{\infty}$. The set $\mathcal{A}_{\infty}$ has a $^*$-algebra structure where addition, scalar multiplication, multiplication and the $*$-structure are as follows.
  \begin{align*}
  [(a,m)]+[(b,n)]&=[(\phi_{m+n,m}(a)+\phi_{m+n,n}(b),m+n)]\\
  \lambda [(a,m)]&=[(\lambda a, m)] \\
  [(a,m)][(b,n)]&=[(\phi_{m+n,m}(a)\phi_{m+n,n}(b),m+n)]\\
  [(a,m)]^{*}&=[(a^{*},m)].
   \end{align*}
  On $\mathcal{A}_{\infty}$, define a $C^{*}$-seminorm as follows.
  \[
  ||[(a,m)]||:=\lim_{n \to \infty}||\phi_{m+n,m}(a)||.\]
  Mod out the null vectors and complete to get a genuine $C^{*}$-algebra which we denote by $A_{\infty}$. Also $A_{\infty}$ is called the inductive limit of $(A_n,\phi_n)$. 
  
  Let $i_n: A_n \to A_{\infty}$ be defined by $i_n(a):=[(a,n)]$. Note that $i_{n} \circ \phi_{n,m}=i_{m}$. This implies in particular that $i_{n}(A_n)$ is an increasing sequence of $C^{*}$-algebras. Moreover, the union $\bigcup_{n \geq 1}i_{n}(A_n)$ is dense in $A_{\infty}$. 
  
  \begin{ppsn}[The universal property]
  Keep the foregoing notation. Suppose $B$ is a $C^{*}$-algebra and there exists $^*$-homomorphisms $j_n:A_n \to B$ such that $j_n \circ \phi_{n,m}=j_m$. Then there exists a unique $^*$-algebra homomorphism $\phi:A_{\infty} \to B$ such that 
  \[
  \phi(i_n(a))=j_n(a)\]
  for $a \in A_n$. Moreover $A_{\infty}$ is characterised by this property. 
   \end{ppsn}
  \textit{Proof.} Left to the reader. 
  
  \begin{rmrk}
  Inductive limits of systems indexed by a general directed set can be defined. We have chosen to work with sequences for simplicity. 
   \end{rmrk}
  
  \begin{xrcs}
  Discuss inductive limits in the category of abelian groups. Formulate and prove a universal property in this context. 
  \end{xrcs}
  
  The main theorem about inductive limits and $K$-theory is the following. 
  
  \begin{thm}
  Let $(A_n,\phi_n)$ be a directed system of $C^{*}$-algebras and let $A_\infty$ be the direct limit. Then 
  \[
  K_{i}(A_\infty)=\lim_{n \to \infty}(K_{i}(A_n),K_i(\phi_n))\]
  for $i=0,1$. 
   \end{thm}
  The soul of the proof of the above theorem relies in the following  two propositions. The reader should convince herself that  it is indeed so. 
  
  \begin{ppsn}
  \label{K0 and inductive limit}
  Let $A$ be a $C^{*}$-algebra. Suppose $A_n$ is an increasing sequence of $C^{*}$-subalgebras of $A$ such that $\bigcup_{n=1}A_n$ is dense in $A$. 
  \begin{enumerate}
  \item[(1)] Let $e \in A$ be a projection. Then there exists a projection $f \in A_m$ for some $m$ such that $e \sim f$. 
  \item[(2)] Let $e,f \in A_m$ be projections. Suppose $e \sim f$ in $A$. Then there exists $n$ large such that $e \sim f$ in $A_{m+n}$. 
  \end{enumerate}
  \end{ppsn}
  
  \begin{lmma}
  Suppose $A$ is a $C^{*}$-algebra. Let $U$ be a non-empty open subset of $\mathbb{C}$. Then $E:=\{a \in A: spec(a) \subset U\}$ is an open subset of $A$. 
    \end{lmma}
    \textit{Proof.} Let $C$ be the complement of $U$ and $F$ be the complement of $E$. We show that $F$ is closed. Let $a_n$ be a sequence in $F$ such that $a_n \to a$. Then there exists $\lambda_n \in C$ such that $
    \lambda_n \in spec(a_n)$. Since $(||a_n||)$ is bounded, it follows that $\lambda_n$ is bounded. By passing to a subsequence, we can assume that $\lambda_n$ converges. Let $\lambda:=\lim_n \lambda_n$. 
    Since $C$ is closed, $\lambda \in C$. Suppose $a-\lambda$ is invertible. Since $a_n-\lambda_n \to a-\lambda$, it follows that $a_n-\lambda_n$ is invertible for large $n$ which is a contradiction. This forces that $\lambda \in spec(a)$. 
    Hence $a \in F$. This proves that $F$ is closed and hence the proof. \hfill $\Box$
  
  \textit{Proof of Prop. \ref{K0 and inductive limit}.} Let $\mathcal{B}:=\bigcup_{n \geq 1}A_n$. Suppose $e \in A$ is a projection. Since $\mathcal{B}$ is dense, there exists $a \in \mathcal{B}$ such that $a=a^{*}$, $||a^{2}-a||<\frac{1}{4}$, $||a-e||<\frac{1}{2}$ and 
  $spec(a) \subset U:=(-\frac{1}{4},\frac{1}{4})\cup (\frac{3}{4},\frac{5}{4})$. Choose $m$ such that $a \in A_m$. Let $h:U \to \mathbb{R}$ be defined by 
   \begin{equation} h(t):=\begin{cases}
 0  & \mbox{ if
} t \in (\frac{-1}{4},\frac{1}{4}),\cr
   &\cr
    1 &  \mbox{ if } t \in (\frac{3}{4},\frac{5}{4}).
         \end{cases}
\end{equation}
Set $f:=h(a)$. Clearly $f$ is a projection in $A_m$. Note that $||a-h(a)||<\frac{1}{2}$. Hence $||e-f|| \leq ||e-a||+||a-h(a)||<1$. By Lemma \ref{key to homotopy}, it follows that $e$ and $f$ are Murray-von Neumann equivalent. This proves $(1)$. 

Let $e,f \in A_m$ be projections. Suppose that $e \sim f$ in $A$. Let $u \in A$ be such that $u^{*}u=e$ and $uu^{*}=f$. Choose a sequence $u_n \in \mathcal{B}$ such that $u_n \to u$. Set $v_n:=fu_ne$. Then $v_{n}^{*}v_n \to e$ and $v_nv_{n}^{*} \to f$. Note that $v_n \in \mathcal{B}$. Thus, there exists $v \in \mathcal{B}$ such that $||v^{*}v-e||<1$, $||vv^{*}-f||<1$ and $fv=ve=v$. Let $n>m$ be such that $v \in A_n$. 

Note that $v^{*}v$ is invertible in $eA_ne$ and $vv^{*}$ is invertible in $fA_nf$. Let $r \in eA_ne$ and $s \in fA_nf$ be such that $r(v^{*}v)^{\frac{1}{2}}=e$ and $s(vv^{*})^{\frac{1}{2}}=f$. Set $w:=vr$. Then $w^{*}w=e$. We claim $w=sv$. 
To see this, note that $v(v^{*}v)^{\frac{1}{2}}=(vv^{*})^{\frac{1}{2}}v$. Multiply by $r$ on the right to deduce that $v=ve=(vv^{*})^{\frac{1}{2}}vr$. Multiply on the left by $s$ to deduce that $sv=s(vv^{*})^{\frac{1}{2}}vr=fvr=vr$. 
This proves the claim. 

Calculate as follows to observe that 
\begin{align*}
ww^{*}&=svv^{*}s \\
           &=s(vv^{*})^{\frac{1}{2}}(vv^{*})^{\frac{1}{2}}s\\
           &=f.
\end{align*}
This proves that $e$ and $f$ are Murray-von Neumann equivalent in $A_n$. This completes the proof. \hfill $\Box$

  \begin{ppsn}
  Let $A$ be a unital $C^{*}$-algebra. Suppose $A_n$ is an increasing sequence of unital $C^{*}$-subalgebras of $A$ such that $\bigcup_{n=1}A_n$ is dense in $A$. 
  \begin{enumerate}
  \item[(1)] Let $u \in A$ be a unitary element. Then there exists a unitary $v \in A_m$ for some $m$ such that $u \sim v$. 
  \item[(2)] Let $u,v \in A_m$ be unitaries. Suppose $u \sim v$ in $A$. Then there exists $n$ large such that $u \sim v$ in $A_{m+n}$. 
    \end{enumerate}
    \end{ppsn}
    \textit{Proof.} Let $\mathcal{B}:=\bigcup_{n=1}^{\infty}A_n$. Suppose $u \in A$ is a unitary element. Since $\mathcal{B}$ is dense, there exists $a \in \mathcal{B}$ such that $||u-a||<1$. Choose $m$ such that $a \in A_m$. Note that  $||1-u^{*}a||<1$. 
    Hence $b:=u^{*}a$ is invertible. Moreover $(-\infty,0]$ is disjoint from $spec(b)$. Choose a holomorphic branch, say $\ell$, of the logarithm defined on $\mathbb{C}\backslash(-\infty,0]$. Set $c:=\ell(a)$ where $c$ is defined using the holomorphic functional calculus. Then $(e^{tc})_{t \in [0,1]}$ is a path of invertibles connecting $1$ to $b$. This implies in particular that $u \sim a$. Let $v:=a|a|^{-1}$. Then $u \sim v$ and $v$ is a unitary in $A_m$. This proves $(1)$. 
   
   Let $u,v$ be unitaries in $A_m$. Suppose $\widetilde{w}:=(w_t)$ is a path of unitaries in $A$ such that $w_0=u$ and $w_1=v$. 
   Set \begin{align*}
   \widetilde{A}:&=C([0,1],A)  \\
   \widetilde{A_n}:&=C([0,1],A_n).
   \end{align*}
    We can view $\widetilde{A}_{n}$ as a unital subalgebra of $\widetilde{A}$. Note that $\bigcup_{n=1}^{\infty}\widetilde{A}_{n}$ is dense in $\widetilde{A}$. Think of $\widetilde{w}$ as an element in $\widetilde{A}$. As in  part $(1)$, extract a path  $(a_{t})_{t \in [0,1]}$ of invertibles in $A_n$ for some $n$ with $n>m$ such that $||w_t-a_t||<1$ for every $t \in [0,1]$. Arguing as in $(1)$ in $A_n$, we see that $u=w_0 \sim a_0$ in $A_n$ and $v=w_1 \sim a_1$ in $A_n$. But $a_0 \sim a_1$ in $A_n$. Therefore $u \sim v$ in $A_n$. This completes the proof. \hfill $\Box$
    
Let $A$ be a $C^{*}$-algebra. We say that $A$ is \emph{approximately finite dimensional}, also called an AF algebra,  if there exists a sequence $(A_n)$ of finite dimensional $C^{*}$-subalgebras of $A$ such that $\bigcup_{n=1}^{\infty}A_n$ is dense in $A$. 
Let $i_{n}:A_n \to A_{n+1}$ be the inclusion. Then $A:=\lim_{n \to \infty}(A_n,i_n)$. 
Note that if $A$ is a finite dimensional algebra then $A$ is isomorphic to $M_{n_1}(\mathbb{C})\oplus M_{n_2}(\mathbb{C})\oplus \cdots M_{n_r}(\mathbb{C})$. Consequently, $K_0(A)=\mathbb{Z}^{r}$ and $K_1(A)=0$. Since $K_i$ preserves inductive limits, in principle, it is possible to compute the $K$-groups of an AF-algebra. In particular, $K_1(A)=0$ for any AF-algebra. The reader should do the following $K$-group computation. 

\begin{enumerate}
\item[(1)] Let $\clh$ be an infinite dimensional separable Hilbert space. Denote the algebra of compact operators by $\mathcal{K}(\clh)$. Then $\mathcal{K}(\clh)$ is AF and $K_0(\mathcal{K}(\clh))=\mathbb{Z}$. Moreover, if $p$ is a minimal projection in $\mathcal{K}(\clh)$ then $[p]$ is a $\mathbb{Z}$-basis for $K_{0}(\mathcal{K}(\clh))$. 
\item[(2)] Set $A_{n}:=M_{2^{n}}(\mathbb{C})$. Let $\phi_{n}:A_n \to A_{n+1}$ be defined by
$\phi_{n}(A):=\begin{bmatrix}
   A & 0 \\
   0 &  A 
   \end{bmatrix}$. The inductive limit $A_{\infty}:=\lim_{n \to \infty}(A_n,\phi_n)$ is called the CAR algebra. Then  \[K_0(A_\infty)=\mathbb{Z}[\frac{1}{2}]:=\Big\{\frac{m}{2^{n}}: m \in \mathbb{Z}, n \in \mathbb{N} \cup \{0\}\Big\}.\] 
\item[(3)] Let $X:=\{0,1\}^{\mathbb{N}}$ be the Cantor set. Then $C(X)$ is an AF-algebra. Its $K$-group is given by $\displaystyle K_0(C(X))=\bigoplus_{ n \in \mathbb{N}}\mathbb{Z}$. 

\end{enumerate}

\begin{xrcs}
Prove the following stability result for $K$-theory. 
Let $\mathcal{K}$ be the $C^{*}$-algebra of compact operators on an infinite dimensional separable Hilbert space. Suppose $p$ is a minimal projection in $\mathcal{K}$ and $A$ is a $C^{*}$-algebra. Let $\omega:A \to \mathcal{K} \otimes A$ be defined by 
\[\omega(a):=p \otimes a.\]
Prove that $K_{i}(\omega)$ is an isomorphism.

\end{xrcs}

\begin{rmrk}
One of the first significant results in the subject is the classification of AF-algebras in terms of its $K$-theory. This was due to Elliot. Elliot's theorem asserts roughly that two AF-algebras are isomorphic if and only
their $K$-theoretic invariants are the same. For a precise statement, we refer the reader  to \cite{Davidson_Ken}.
\end{rmrk}  
  \section{Six term exact sequence}
  An important computational tool that enables us to calculate the $K$-groups explicitly is the six term exact sequence. We omit the proof altogether and merely explain the consequences.  Let 
  \[
  0 \longrightarrow I \longrightarrow A \stackrel{\pi} \longrightarrow B \longrightarrow 0\]
  be an exact sequence of $C^{*}$-algebras. Then there exists maps $\partial:K_1(B) \to K_0(I)$, called the \emph{index map}, and $\sigma:K_0(B) \to K_1(I)$ which makes the following six term sequence exact.
    \begin{equation*}
\def\labelstyle{\scriptstyle}
\xymatrix@C=25pt@R=20pt{
K_0(I)\ar[r]& K_0(A)\ar[r]& K_0(B) \ar[d]_{\sigma} \\
K_1(B)\ar[u]^{\partial}&
K_1(A)\ar[l] & K_1(I)\ar[l] }
\end{equation*}
Moreover the maps $\partial$ and $\sigma$ are ``natural".

The construction of the index map $\partial$, though tedious, is not that difficult. It is explicitly described below. Let $[u] \in K_1(B)$ be given where $u$ is a unitary in $M_{n}(B^{+})$. Choose a unitary $V \in M_{2n}(A^{+})$ such that $\pi^{+}(V)=\begin{bmatrix}
u & 0 \\
0 & u^{*}
\end{bmatrix}$. The justification for the existence of such a unitary is given in Prop. \ref{half exact}. Then 
\[
\partial([u])=\Big[ V \begin{pmatrix}
                                1_{n} & 0 \\
                                0 & 0 
                                \end{pmatrix} V^{*} \Big] - \Big[\begin{pmatrix}
                                             1_n & 0 \\
                                             0 & 0 
                                             \end{pmatrix} \Big].
\]
The map $\partial$ defined above is well defined, i.e. it is independent of the various choices made and makes the diagram exact at $K_1(B)$ and $K_0(I)$. The construction of $\sigma$ is more difficult and requires Bott periodicity. 

The following is often used in applications. 
\begin{ppsn}
\label{used in applications}
Let $0 \longrightarrow I \longrightarrow A \stackrel{\pi} \longrightarrow B \longrightarrow 0$ be a short exact sequence of $C^{*}$-algebras. Let $u \in M_{n}(B^{+})$ be a unitary. Suppose there exists a partial isometry $v \in M_{n}(A^{+})$ such that 
$\pi^{+}(v)=u$. Then $\partial([u])=[1_n-v^{*}v]-[1_n-vv^{*}]$.
\end{ppsn}
\textit{Proof.} Let $V:=\begin{bmatrix}
 v & 1-vv^{*} \\
 1-v^{*}v v^{*}& v^{*} 
 \end{bmatrix}$. Then $V$ is a unitary ``lift" of $\begin{bmatrix}
   u & 0 \\
   0 & u^{*}
   \end{bmatrix}$.
 Clearly, $V\begin{bmatrix}
                   1_n & 0 \\
                    0 & 0 
                    \end{bmatrix}V^{*}=\begin{bmatrix}
                                                   vv^{*} & 0 \\
                                                    0 & 1-v^{*}v
                                                    \end{bmatrix}$. Calculate as follows to observe that 
  \begin{align*}
  \partial([u])&=\Big[\begin{pmatrix}
                                                          vv^{*} & 0 \\
                                                          0 & 1-v^{*}v 
                                                          \end{pmatrix}\Big]-\Big[\begin{pmatrix}
                                                                                              1-vv^{*}+vv^{*} & 0 \\
                                                                                              0 & 0 
                                                                                              \end{pmatrix}\Big] \\
                                                                                              &=  [vv^{*}]+[1-v^{*}v]-[1-vv^{*}]-[vv^{*}]\\
                                                                                              &=[1-v^{*}v]-[1-vv^{*}].
                                          \end{align*}        
                                          This completes the proof. \hfill $\Box$

As a first application, we deduce that $K_1$ can be defined in terms of $K_{0}$ and $K_0$ can be defined in terms of $K_1$. We say a $C^{*}$-algebra $B$ is \emph{contractible} if the identity homomorphism is homotopy equivalent to the zero map. 
If $B$ is contractible then $K_0(B)=K_1(B)=0$. 

Let $A$ be a $C^{*}$-algebra. Denote the $C^{*}$-algebra of continuous $A$-valued functions  on $[0,1]$ by $C([0,1],A)$. The norm here is the supremum norm. Set 
\begin{align*}
CA:&=\{f \in C([0,1],A): f(0)=0\}\\
SA:&=\{f \in CA: f(1)=0\}.
\end{align*}
The $C^{*}$-algebra $CA$ is called the \emph{cone over $A$} and $SA$ is called the \emph{suspension over $A$}. Note that $A \to CA$ and $A \to SA$ are functors from the category of $C^{*}$-algebras to $C^{*}$-algebras. 

\begin{lmma}
The cone $CA$ is contractible. Hence $K_0(CA)=0$ and $K_1(CA)=0$. 
\end{lmma}
\textit{Proof.} For $t \in [0,1]$, let $\epsilon_{t}:CA \to CA$ be defined by $\epsilon_{t}(f)(s)=f(st)$. Then $(\epsilon_{t})_{t \in [0,1]}$ is a homotopy of $^*$-homomorphisms connecting the zero map with the identity map. This completes the proof. \hfill $\Box$.

\begin{crlre}
For any $C^{*}$-algebra $A$, $K_1(A) \cong K_0(SA)$ and $K_0(A) \cong K_1(SA)$. 
\end{crlre}
\textit{Proof.} Let $\epsilon:CA \to A$ be defined by $\epsilon(f)=f(1)$. Then we have the following exact sequence 
\[
0 \longrightarrow SA \longrightarrow CA \stackrel{\epsilon} \longrightarrow A \longrightarrow 0 .\]
The conclusion is immediate if we apply the six term exact sequence and the fact that $K_0(CA)=K_1(CA)=0$. \hfill $\Box$

\begin{rmrk}
Actually, we have cheated a lot. In fact, the isomorphism $K_0(A) \cong K_1(SA)$ is needed apriori to define the map $\sigma$ in the six term sequence. After first proving this, $\sigma$ is defined as the composite of the following maps 
\[
K_0(B) \cong K_1(SB) \stackrel{\partial} \longrightarrow  K_0(SI) \cong K_1(I).\]
The isomorphism $K_0(A) \cong K_1(SA)$ is called the Bott periodicity in $K$-theory. We will discuss Cuntz' proof of it in the next two sections. 
\end{rmrk}

\begin{xrcs}
\begin{enumerate}
\item[(1)] Note that for $A=\mathbb{C}$, $SA:=C_{0}(\mathbb{R})$. Conclude that $K_0(C_0(\mathbb{R}))=0$ and $K_1(C_0(\mathbb{R}))=\mathbb{Z}$. 
\item[(2)] Let $\epsilon:C(\mathbb{T}) \to \mathbb{C}$ be the evaluation map at $1$. Then the short exact sequence 
\[
0 \longrightarrow C_{0}(\mathbb{R}) \longrightarrow C(\mathbb{T}) \stackrel{\epsilon} \longrightarrow \mathbb{C} \longrightarrow 0\]
is split exact. Conclude that $K_0(C(\mathbb{T})) \cong \mathbb{Z}$ and $K_1(C(\mathbb{T}))\cong \mathbb{Z}$. Show that $[1]$ forms a $\mathbb{Z}$-basis for $C(\mathbb{T})$. 
\item[(3)] Compute the $K$-groups of $C(S^{2})$ where $S^{2}$ is the unit sphere in $\mathbb{R}^{3}$. 
\end{enumerate}
\end{xrcs}

As an application of the six term sequence, we compute the $K$-groups of the Toeplitz algebra. Recall that the Toeplitz algebra $\mathcal{T}$ is the universal $C^{*}$-algebra generated by an isometry $v$. 
Let $z \in C(\mathbb{T})$ be the generating unitary. Then there exists a unique surjective $^*$-homomorphism $\pi:\mathcal{T} \to C(\mathbb{T})$ such that $\pi(v)=z$. Also, the kernel of $\pi$ is isomorphic 
to the $C^{*}$-algebra of compact operators, denoted $\mathcal{K}$, on a separable infinite dimensional Hilbert space. Let $i:\mathcal{K} \to \mathcal{T}$ be the inclusion. 
We apply the six term sequence to the following exact sequence.
\[
0 \longrightarrow \mathcal{K} \longrightarrow \mathcal{T} \stackrel{\pi} \longrightarrow C(\mathbb{T}) \longrightarrow 0\]

Consider the six term exact sequence    \begin{equation*}
\def\labelstyle{\scriptstyle}
\xymatrix@C=25pt@R=20pt{
K_0(\mathcal{K})\ar[r]& K_0(\mathcal{T})\ar[r]& K_0(C(\mathbb{T})) \ar[d]_{\sigma} \\
K_1(C(\mathbb{T}))\ar[u]^{\partial}&
K_1(\mathcal{T})\ar[l] & K_1(\mathcal{K})\ar[l] }
\end{equation*}

We claim $\partial([z])=-[p]$ where $p$ is a rank one projection in $\mathcal{K}$. Since $\pi(v)=z$ and $v$ is an isometry, it follows from Prop. \ref{used in applications} that $\partial([z])=[1-v^{*}v]-[1-vv^{*}]=-[p]$. Note that $[p]$ is a $\mathbb{Z}$-basis for $K_0(\mathcal{K})$. Also, we know that $K_1(C(\mathbb{T}))\cong \mathbb{Z}$. Hence $[z]$ is a $\mathbb{Z}$-basis for $K_1(\mathbb{C}(\mathbb{T}))$. 
Therefore, $\partial$ is an isomorphism. Consequently, $K_0(i)=0$. This implies that $Ker(K_0(\pi))=0$. Note that $K_1(\mathcal{K})=0$. This implies that $\sigma$ is the zero map. Hence $Im(K_0(\pi))=K_0(C(\mathbb{T}))$. Consequently, $K_0(\pi)$ is an isomorphism. Therefore $K_0(\mathcal{T})$ is isomorphic to $\mathbb{Z}$ and $[1]$ is a $\mathbb{Z}$-basis for $K_0(\mathcal{T})$. 
Note that $K_1(\mathcal{K})=0$ and $Im(K_1(\pi))=Ker(\partial)=0$. Thus we have the short exact sequence \[0 \longrightarrow K_1(\mathcal{T}) \longrightarrow 0.\] As a consequence, $K_1(\mathcal{T})=0$. 

  \section{Quasi-homomorphisms}
  We conclude these notes by discussing Cuntz' proof of Bott periodicity. An important technical tool that we need is the notion of quasi-homomorphisms. This also offers a first glimpse of KK-theory. The notion of quasi-homomorphisms is due to Cuntz. We know that homomorphisms between $C^{*}$-algebras induce maps at the $K$-theory level. The important observation due to Cuntz is that quasi-homorphisms, a sort of a generalised morphism between $C^{*}$-algebras, too induce maps at the $K$-theory level. 
  
  Let $A$ and $J$ be $C^{*}$-algebras. By a quasi-homomorphism from $A \to J$, we mean the following data. There exists a $C^{*}$-algebra $\mathcal{E}$ which contains $J$ as an ideal and two $^*$-homomorphisms $\phi_{+},\phi_{-}:A \to \mathcal{E}$ such that for $a \in A$, $\phi_{+}(a)-\phi_{-}(a) \in J$.  We simply say that \[\phi:=(\phi_{+},\phi_{-}):A \to \mathcal{E} \trianglerighteq J\] is a quasi-homomorphism from $A$ to $J$ to mean the above data. Strictly speaking, there exists an embedding $i:J \to \mathcal{E}$ such that $i(J)$ is an ideal in $\mathcal{E}$. As usual, we suppress the embedding to be economical with notation. 
  
  \begin{xmpl}
  Suppose $\sigma: A \to J$ is a homomorphism. Then $(\sigma,0):A \to J \trianglerighteq J$ is a quasi-homomorphism. More generally, suppose $\sigma_1,\sigma_2:A \to J$ are homomorphisms. Then $\sigma:=(\sigma_1,\sigma_2):A \to J \trianglerighteq J$ is a quasi-homomorphism. 
  \end{xmpl}
  
  Let $\phi:=(\phi_{+},\phi_{-}):A \to \mathcal{E} \trianglerighteq J$ be a quasi-homomorphism. Set 
  \[
  A_{\phi}:=\{(a,x) \in A \oplus \mathcal{E}: \phi_{+}(a) \equiv x\mod J\}.\]
  Let $\pi:A_{\phi} \to A$ be defined by $\pi(a,x)=a$. Then $Ker(\pi):=\{(0,x): x \in J\}$ which we identify with $J$. Let $j:J \to A_{\phi}$ be the embedding $j(x)=(0,x)$. 
  Define $\widetilde{\phi_{+}}:A \to A_{\phi}$ and $\widetilde{\phi_{-}}:A \to A_{\phi}$ by 
  \begin{align*}
  \widetilde{\phi_{+}}(a)&=(a,\phi_{+}(a))\\
  \widetilde{\phi_{-}}(a)&=(a,\phi_{-}(a)).
    \end{align*} 
  Then we have the following split exact sequence of $C^{*}$-algebras with the splitting given by either $\widetilde{\phi_{+}}$ or $\widetilde{\phi_{-}}$. 
  \[
  0 \longrightarrow J \longrightarrow A_{\phi} \stackrel{\pi} \longrightarrow A \longrightarrow 0\]
  
  Hence $K_{i}(j)$ is injective. Note that $K_{i}(\widetilde{\phi_+})-K_{i}(\widetilde{\phi_{-}}) \in Ker(K_i(\pi))=Im(K_i(j))$. For $i=0,1$, we define $\widehat{K_i}(\phi):K_i(A) \to K_{i}(J)$ by the formula
  \[
  \widehat{K_i}(\phi):=K_{i}(j)^{-1}\circ \Big(K_i(\widetilde{\phi_+})-K_i(\widetilde{\phi_{-}})\Big).\]
  
  Next we derive a few basic properties about quasi-homomorphisms.
  \begin{ppsn}
  Let $\sigma:=(\sigma_1,\sigma_2):A \to J \trianglerighteq J$ be a quasi-homomorphism. Then $\widehat{K_i}(\sigma)=K_{i}(\sigma_1)-K_i(\sigma_2)$. 
    \end{ppsn}
    \textit{Proof.} Note that $A_{\sigma}=A \oplus J$. Then we can identify $K_{i}(A_{\sigma})$ with $K_i(A) \oplus K_{i}(J)$. Once this identification is made, $K_{i}(j)^{-1}$ on $Im(K_{i}(j))$ is nothing but $K_{i}(pr_2)$ where $pr_2:A \oplus J \to J$ is the second projection. The conclusion is now obvious. \hfill $\Box$
    
    In view of the above proposition, for a quasi-homomorphism $\phi$, we simply denote $\widehat{K}_{i}(\phi)$ by $K_i(\phi)$. Next we discuss how to precompose a quasi-homomorphism with a homomorphism. Let $\phi:=(\phi_{+},\phi_{-}):A \to \mathcal{E} \trianglerighteq J$ be a quasi-homomorphism. Suppose $\epsilon:B \to A$ is a homomorphism. Then $\psi:=(\psi_{+},\psi_{-}):B \to \mathcal{E} \trianglerighteq J$ is a quasi-homomorphism where $\psi_{+}=\phi_{+} \circ \epsilon$ and $\psi_{-}=\phi_{-} \circ \epsilon$. 
    
    \begin{ppsn}
    With the foregoing notation, we have $K_i(\psi)=K_{i}(\phi) \circ K_{i}(\epsilon)$. 
    \end{ppsn}
    \textit{Proof.} Let $j_{B}: J \to B_{\psi}$ and $j_{A}: J \to A_{\phi}$ be the embeddings. Define $\eta:B_{\psi} \to A_{\phi}$ by 
    \[
    \eta(b,x):=(\epsilon(b),x).\]
    Note that $\eta \circ j_{B}=j_{A}$, $\eta \circ \widetilde{\psi_{+}}=\widetilde{\phi_{+}} \circ \epsilon$ and $\eta \circ \widetilde{\psi_{-}}=\widetilde{\phi_{-}} \circ \epsilon$. Let $y \in K_{i}(B)$ given. Choose $x \in K_{i}(J)$ such that $K_{i}(j_B)x=(K_{i}(\widetilde{\psi}_{+})-K_i(\widetilde{\psi}_{-}))y$. Then $K_i(\psi)y=x$. 
    
    To show that $K_{i}(\phi)\circ K_{i}(\epsilon)y=x$, it suffices to show that $K_{i}(j_A)x=(K_{i}(\widetilde{\phi}_{+})-K_{i}(\widetilde{\phi}_{-}))K_i(\epsilon)y$. Calculate as follows to observe that 
    \begin{align*}
    K_{i}(j_A)(x)&=K_{i}(\eta)K_{i}(j_B)x \\
     &= K_{i}(\eta)K_{i}(\widetilde{\psi}_+)y-K_{i}(\eta)K_{i}(\widetilde{\psi}_{-})y\\
     &=K_{i}(\widetilde{\phi}_+ \circ \epsilon)y-K_{i}(\widetilde{\phi}_{-}\circ \epsilon)y  \\
     &=(K_{i}(\widetilde{\phi}_+)-K_{i}(\widetilde{\phi}_{-}))K_{i}(\epsilon)y.  
    \end{align*}
    This completes the proof. \hfill $\Box$
    
    Post composing a quasi-homomorphism with a homomorphism is a bit tricky. The data we need is the following. Suppose $\phi:=(\phi_{+},\phi_{-}):A \to \mathcal{E} \trianglerighteq J$ is a quasi-homomorphism. Let $\epsilon^{'}:J \to J^{'}$ be a $^*$-homomorphism. To define $\epsilon^{'} \circ \phi$, we need an extension of $\epsilon^{'}$. Suppose there exists a $C^{*}$-algebra $\mathcal{E}^{'}$ containing $J^{'}$ as an ideal such that $\epsilon^{'}$ extends to a map from $\mathcal{E} \to \mathcal{E}^{'}$. Denote an extension again by $\epsilon^{'}$. Set $\psi_{+}:=\epsilon^{'} \circ \phi_{+}$ and $\psi_{-}:=\epsilon^{'} \circ \phi_{-}$. Then $\psi:=(\psi_{+},\psi_{-}):A \to \mathcal{E}^{'} \trianglerighteq J^{'}$ is a quasi-homomorphism.
    
    \begin{ppsn}
    \label{post composition}
    With the foregoing notation, we have $K_{i}(\psi)=K_{i}(\epsilon^{'})\circ K_{i}(\phi)$. 
        \end{ppsn}
        \textit{Proof.} Let $j:J \to A_{\phi}$ and $j^{'}:J^{'} \to A_{\psi}$ be the embeddings. Let $\eta:A_{\phi} \to A_{\psi}$ be defined by $\eta(a,x):=(a,\epsilon^{'}(x))$. Then $\eta \circ \widetilde{\phi_+}=\widetilde{\psi_{+}}$ and $\eta \circ \widetilde{\phi_{-}}=\widetilde{\psi_{-}}$. Also $\eta \circ j=j^{'} \circ \epsilon^{'}$. 
        
        Let $y \in K_i(A)$ be given. Choose $x \in K_{i}(J)$ such that $K_{i}(\widetilde{\phi}_{+})y-K_{i}(\widetilde{\phi}_{-})y=K_{i}(j)x$. To prove $K_{i}(\psi)y=K_{i}(\epsilon^{'})K_{i}(\phi)y$, it suffices to show that $K_{i}(j^{'})K_{i}(\epsilon^{'})x=K_{i}(\widetilde{\psi_+})y-K_{i}(\widetilde{\psi_{-}})y$. Calculate as follows to observe that 
        \begin{align*}
        K_{i}(j^{'})K_{i}(\epsilon^{'})x&=K_{i}(\eta)K_{i}(j)x \\
        &=K_{i}(\eta)(K_{i}(\widetilde{\phi_+})y-K_{i}(\eta)K_{i}(\widetilde{\phi_{-}})y\\
        &=K_{i}(\widetilde{\psi_+})y-K_{i}(\widetilde{\psi_{-}})y.
           \end{align*}
           This completes the proof. \hfill $\Box$
  
  For $t \in [0,1]$, let $\phi^{t}:=(\phi_{+}^{t},\phi_{-}^{t}):A \to \mathcal{E} \trianglerighteq J$ be a family of quasi-homomorphisms. We say that $(\phi^{t})_{t \in [0,1]}$ is a homotopy if $(\phi_{+}^{t})_{t \in [0,1]}$ and $(\phi_{-}^{t})_{t \in [0,1]}$ are homotopy of $^*$-homomorphisms. 
  \begin{ppsn}
  Let $\phi^{t}:=(\phi_{+}^{t},\phi_{-}^{t}):A \to \mathcal{E} \trianglerighteq J$ be a homotopy of $^*$-homomorphisms. Then $K_{i}(\phi^{t})$ is independent of $t$. 
   \end{ppsn}
   \textit{Proof.} Let $\widetilde{\mathcal{E}}:=C([0,1],\mathcal{E})$ and $\widetilde{J}:=C([0,1],J)$. Define $\widetilde{\phi_+}:A \to \widetilde{\mathcal{E}}$ by the formula
   \[
   \widetilde{\phi_+}(a)(t):=\phi_{+}^{t}(a).\]
  Similarly define $\widetilde{\phi_{-}}$. 
  Then $\widetilde{\phi}:=(\widetilde{\phi_+},\widetilde{\phi_{-}}):A \to \widetilde{\mathcal{E}} \trianglerighteq \widetilde{J}$ is a quasi-homomorphism. For $t \in [0,1]$, let $\epsilon_{t}:\widetilde{\mathcal{E}} \to \mathcal{E}$ be the evaluation at $t$. By Prop. \ref{post composition}, $K_{i}(\phi^{t})=K_{i}(\epsilon_t) \circ K_{i}(\widetilde{\phi})$. However, $K_i(\epsilon_t)$ is constant by the homotopy invariance of $K$-theory. Hence the proof. \hfill $\Box$
  
  Next we discuss the additive property of $K$-theory. Let $\phi,\psi:A \to B$ be homomorphisms. We say $\phi$ and $\psi$ are othogonal and write $\phi \perp \psi$ if for $x,y \in A$, $\phi(x)\psi(y)=0$. Note that if $\phi \perp \psi$ then $\phi+\psi:A \to B$ is a $^*$-homomorphism. 
  Let $\phi:=(\phi_{+},\phi_{-}):A \to \mathcal{E} \trianglerighteq J$ and $\psi:=(\psi_{+},\psi_{-}):A \to \mathcal{E} \trianglerighteq J$ be two quasi-homomorphisms. We say that $\phi$ and $\psi$ are orthogonal if $\phi_{+} \perp \psi_{+}$ and $\phi_{-} \perp \psi_{-}$. 
  If $\phi$ and $\psi$ are orthogonal, then clearly $\phi + \psi:=(\phi_{+}+\psi_{+},\phi_{-}\circ \psi_{-}):A \to \mathcal{E} \trianglerighteq J$ is a quasi-homomorphism. 
  
  \begin{ppsn}
  \label{additivity}
    Suppose $\phi:=(\phi_{+},\phi_{-}):A \to \mathcal{E} \trianglerighteq J$ and $\psi:=(\psi_{+},\psi_{-}):A \to \mathcal{E} \trianglerighteq J$ are orthogonal quasi-homomorphisms. Then $K_{i}(\phi+\psi)=K_i(\phi)+K_{i}(\psi)$. 
    \end{ppsn}
  
  \begin{lmma}
  Let $A$ be a $C^{*}$-algebra. Let $i_1,i_2 \to A \to A \oplus A$ be defined by $i_1(a)=(a,0)$ and $i_2(a)=(0,a)$. 
  Then $K_i(i_1+i_2)=K_{i}(i_1)+K_{i}(i_2)$. 
    \end{lmma}
  \textit{Proof.} Let $\pi_1,\pi_2:A \oplus A \to A$ be the first and the second projections respectively. We know that $K_i(\pi_1)\oplus K_i(\pi_2):K_{i}(A \oplus A) \to K_{i}(A) \oplus K_i(A)$ is an isomorphism. To verify the equality $K_{i}(i_1+i_2)=K_i(i_1)+K_i(i_2)$, it suffices to verify the following two equalities.  
  \begin{align*}
  K_{i}(\pi_1) \circ K_{i}(i_1+i_2)&=K_{i}(\pi_1)\circ K_i(i_1) + K_{i}(\pi_1) \circ K_{i}(i_2)\\
  K_{i}(\pi_2)\circ K_{i}(i_1+i_2)&=K_{i}(\pi_2)\circ K_{i}(i_1)+K_{i}(\pi_2)\circ K_i(i_2)
  \end{align*} 
  This verification is obvious. \hfill $\Box$.
  
  \textit{Proof of Prop. \ref{additivity}:} Let $\Sigma_{+}: A \oplus A \to \mathcal{E} \trianglerighteq J$ be defined by $\Sigma_{+}(a,b)=\phi_{+}(a)+\psi_{+}(b)$. Similarly define $\Sigma_{-}$. Then $\Sigma:=(\Sigma_{+},\Sigma_{-}):A \to \mathcal{E} \trianglerighteq J$ is a quasi-homomorphism. Let $\Delta:A \to A \oplus A$ be defined by $\Delta(a)=(a,a)$. Then $\Delta=i_1+i_2$. 
    Note that $\phi+\psi:= \Sigma \circ \Delta$. Clearly $\phi=\Sigma \circ i_1$ and $\psi=\Sigma \circ i_2$. Calculate as follows to observe that 
    \begin{align*}
    K_{i}(\phi+\psi)&=K_{i}(\Sigma) \circ K_{i}(\Delta) \\
     & = K_{i}(\Sigma) \circ (K_i(i_1)+K_i(i_2)) \\
     &=K_{i}(\Sigma)\circ K_{i}(i_1)+ K_{i}(\Sigma)\circ K_i(i_2) \\
     &=K_{i}(\Sigma \circ i_1)+K_i(\Sigma \circ i_2)\\
     &= K_{i}(\phi)+K_i(\psi).
     \end{align*}
  This completes the proof. \hfill $\Box$

  \section{Bott periodicity}
  In this section, we discuss Cuntz' proof of Bott periodicity. The main result in Cuntz' proof is to first compute the $K$-theory of the Toeplitz algebra. We had already computed the $K$-groups of the Toeplitz algebra assuming Bott periodicity. Here we compute it without this assumption.   Recall that the Toeplitz algebra $\mathcal{T}$ is the universal $C^{*}$-algebra generated by a single isometry $v$.

  We need to use tensor products of $C^{*}$-algebras, a delicate topic, in what follows. The reader should consult  \cite{Ozawa-Brown} for a detailed treatment. We ask the reader to accept the statements made here in good faith. 
    Let $A_1$ and $A_2$ be $C^{*}$-algebras. Consider the algebraic tensor product $A_1 \otimes_{alg} A_2$. Then $A_1 \otimes_{alg} A_2$ is a $^*$-algebra where the multiplication and $*$-structure are given by 
\begin{align*}
(a_1 \otimes a_2)(b_1 \otimes b_2)&=a_1b_1 \otimes a_2b_2 \\
(a \otimes b)^{*}&=a^{*} \otimes b^{*}.
\end{align*}
Let $||~||$ be a $C^{*}$-norm on $A_1 \otimes_{alg} A_2$. The norm $||~||$ is said to be a cross-norm on $A_1 \otimes_{alg} A_2$ if $||a \otimes b||=||a||||b||$. It is true that there exists $C^{*}$-algebras $A_1$ and $A_2$ such that $A_1 \otimes_{alg} A_2$ admits more than one $C^{*}$ cross norm. 

\begin{dfn}
A $C^{*}$-algebra $A$ is called nuclear if the following holds. For  every $C^{*}$-algebra $B$, there is only one $C^{*}$ cross norm on the algebraic tensor product $A \otimes_{alg} B$.
 If $A$ is nuclear then $A \otimes B$ denotes the completion of $A \otimes_{alg} B$ with respect to any $C^{*}$ cross norm  
\end{dfn}

\begin{xrcs}
Show that $M_{n}(\mathbb{C})$ is nuclear. 
\end{xrcs}

\textbf{Spatial tensor product:} 
It is always possible to define a $C^{*}$-cross norm as follows. Let $A_1$ and $A_2$ be two $C^{*}$-algebras. Let $\pi_1:A_1 \to B(\clh_1)$ and $\pi_2: A_2 \to B(\clh_2)$ be faithful representations. Define $\pi_1 \otimes \pi_2 : A_1 \otimes_{alg} A_2 \to B(\clh_1 \otimes \clh_2)$ by the equation \[(\pi_1 \otimes \pi_2)(a_1 \otimes a_2):=\pi_1(a_1) \otimes \pi_2(a_2).\]  Then $\pi_1 \otimes \pi_2$ is a $^*$-homomorphism and is injective. 
For $x \in A_1 \otimes_{alg} A_2$, let $||x||:=||\pi_1 \otimes \pi_2(x)||$. Then $||~||$ is a norm on $A_1 \otimes_{alg} A_2$. It is a non-trivial fact that $||~||$ is independent of the chosen faithful representations $\pi_1$ and $\pi_2$. This norm on $A_1 \otimes_{alg} A_2$ is called the spatial norm and the completion of $A_1 \otimes_{alg} A_2$ is called \emph{the spatial tensor product}.  
The reader can assume that the tensor product of $C^{*}$-algebras that we consider is always the spatial one without much loss. 

\begin{xrcs}
Let $D$ be a $C^{*}$-algebra. Show that the map $A \to A \otimes D$ is a functor from the category of $C^{*}$-algebras to the category of $C^{*}$-algebras. Here the tensor product is the spatial one.

\end{xrcs}

We need the following facts regarding nuclear $C^{*}$-algebras and tensor products. 

\begin{enumerate}
\item Commutative $C^{*}$-algebras are nuclear.
\item Inductive limits of nuclear $C^{*}$-algebras are nuclear. 
\item Let $0 \rightarrow I \rightarrow A \rightarrow B \rightarrow 0$ be a short exact sequence. If $I$ and $B$ are nuclear then $A$ is nuclear. 
\item Let $0 \rightarrow I \rightarrow A \rightarrow B \rightarrow 0$ be a short exact sequence of nuclear $C^{*}$-algebras. If $D$ is a $C^{*}$-algebra then the sequence \[
0 \rightarrow D \otimes I \rightarrow D \otimes A \rightarrow D \otimes B \rightarrow 0\] is exact. 
\end{enumerate}

\begin{xrcs}
Use the above facts to conclude that the Toeplitz algebra is nuclear. 
\end{xrcs}

\begin{xrcs}
Let $X$ be locally compact Hausdorff topological space and $A$ be a $C^{*}$-algebra. Assume that $C_{0}(X)$ is nuclear. Use this assumption to
show that $C_{0}(X) \otimes A \cong C_{0}(X,A)$. 
\end{xrcs}
\textit{Hint:} The map $C_{0}(X) \otimes_{alg} A \ni f \otimes a \to f.a  \in C_{0}(X,A)$ is an embedding. Here $f.a$ stands for the map which sends $x$ to $f(x)a$.

Let us return to the discussion on Bott periodicity. Let $q:\mathcal{T} \to \mathbb{C}$ be defined by $q(v)=1$ and $j:\mathbb{C} \to \mathcal{T}$ be defined by $j(1)=1$. The map $q$ exists by Coburn's theorem. 
We claim that $K_{i}(q):K_{i}(\mathcal{T}) \to K_{i}(\mathbb{C})$ is an isomorphism with inverse given by $K_{i}(j)$. Since $q \circ j=id$, it follows that $K_i(q) \circ K_{i}(j)=id$. 

Let $p=1-vv^{*}$. Let $\omega:\mathcal{T} \to \mathcal{K} \otimes \mathcal{T}$ be defined by $\omega(x)=p\otimes x$. Since $K_i(\omega)$ is an isomorphism, to show that $K_i(j) \circ K_{i}(q)=Id$, it suffices to show that $K_{i}(\omega) \circ K_i(j\circ q)=K_{i}(\omega)$. Let $\sigma_1:=\omega \circ j \circ q$ and $\sigma_2=\omega$. Then 
\begin{align*}
\sigma_1(v)&=p \otimes 1\\
\sigma_1(v)&=p \otimes v.
\end{align*}
We need to show that $K_i(\sigma_1)=K_i(\sigma_2)$. 

\begin{xrcs}
Prove the following version of Coburn's theorem. Let $A$ be a $C^{*}$-algebra. Suppose $w$ is a partial isometry in $A$ such that $ww^{*} \leq w^{*}w$. Then there exists a unique $^*$-homomorphism $\sigma:\mathcal{T} \to A$ such that $\sigma(v)=w$. 
\end{xrcs}  
 \textit{Hint:} Set $p:=w^{*}w$ and consider $pAp$.  
 
 Keep the notation preceeding the above exercise. 
 
 \begin{thm}
We have $K_i(\sigma_1)=K_i(\sigma_2)$. Therefore, the  map $K_i(q):K_i(\mathcal{T}) \to K_i(\mathbb{C})$ is an isomorphism. 
  \end{thm}
 \textit{Proof.}  By Coburn's theorem, there exists a $^*$-homomorphism $\epsilon:\mathcal{T} \to \mathcal{T} \otimes \mathcal{T}$ such that $\epsilon(v)=v(1-p)\otimes 1$. Note that \[\sigma:=(\sigma_1,\sigma_2):\mathcal{T} \to \mathcal{T} \otimes \mathcal{T} \trianglerighteq \mathcal{K} \otimes \mathcal{T}\] and \[\widetilde{\epsilon}:=(\epsilon,\epsilon):\mathcal{T} \to \mathcal{T} \otimes \mathcal{T} \trianglerighteq \mathcal{K} \otimes \mathcal{T}\] are quasi-homomorphisms. Moreover $\sigma \perp \widetilde{\epsilon}$. By the properties of quasi-homomorphisms discussed in the previous section, it follows that $K_{i}(\sigma+\widetilde{\epsilon})=K_{i}(\sigma)+K_{i}(\widetilde{\epsilon})=K_i(\sigma_1)-K_i(\sigma_2)$. We will be done if we show that $K_i(\sigma+\widetilde{\epsilon})=0$. 
  
  Let \begin{align*}
  v_{t}:&=cos(\frac{\pi}{2}t)(p \otimes 1)+sin(\frac{\pi}{2}t)(vp\otimes 1) + v(1-p) \otimes 1 \\
  w_{t}:&=cos(\frac{\pi}{2}t)(p \otimes v)+sin(\frac{\pi}{2}t)(vp\otimes 1) + v(1-p) \otimes 1
   \end{align*}
  Note that $v_t$ and $w_t$ are isometries in $\mathcal{T} \otimes \mathcal{T}$. For every $t \in [0,1]$, by  Coburn's theorem, there exists   $^*$-homomorphisms $\sigma^{(t)}_{+}:\mathcal{T} \to \mathcal{T} \otimes \mathcal{T}$ and $\sigma^{(t)}_{-}:\mathcal{T} \to \mathcal{T} \otimes \mathcal{T}$ such that 
  \begin{align*}
  \sigma^{(t)}_{+}(v)&=v_t \\
  \sigma^{(t)}_{-}(v)&=w_t.
   \end{align*}
   Clearly $\sigma^{(t)}:=(\sigma^{(t)}_{+},\sigma^{(t)}_{-}):\mathcal{T} \to \mathcal{T} \otimes \mathcal{T} \trianglerighteq \mathcal{K} \otimes \mathcal{T}$ is a homotopy of quasi-homorphisms. Note that  $\sigma^{(0)}:=\sigma+\widetilde{\epsilon}$. By the homotopy invariance, we have
   \[
   K_{i}(\sigma+\widetilde{\epsilon})=K_{i}(\sigma^{(1)}_{+},\sigma^{(1)}_{-}).\]
   But $\sigma^{(1)}_{+}=\sigma^{(1)}_{-}$. Hence $K_i(\sigma+\widetilde{\epsilon})=0$. Consequently, we have $K_i(\sigma_1)=K_i(\sigma_2)$. This completes the proof. \hfill $\Box$.

  \begin{crlre}
  Let $\mathcal{T}_0:=Ker(q)$. Then $K_i(\mathcal{T}_0)=0$. 
    \end{crlre}
    \textit{Proof.} Note that the short exact exact sequence 
    \[
    0 \longrightarrow \mathcal{T}_0 \longrightarrow \mathcal{T} \stackrel{q} \longrightarrow \mathbb{C} \longrightarrow 0\]
    is split exact with the splitting given by $j$. Since $K_i(q)$ is an isomorphism, the conclusion follows. \hfill $\Box$
    
    The next step in the proof of Bott periodicity is to establish that $K_i(\mathcal{T}_0 \otimes B)=0$ for every $C^{*}$-algebra $B$. Actually, we do not have to do anything. If we go through the proofs once again, we realise that all we need to know about the functor $K_i$ is that it is stable, homotopy invariant and sends split exact sequences to split exact sequences. 
         The proof is applicable for any functor from the category of (nuclear) $C^{*}$-algebras to the category of abelian groups which is split exact, homotopy invariant and is stable. 
    
        Fix a $C^{*}$-algebra $B$. Let $F$ be the functor from the category of nuclear $C^{*}$-algebras to the category of abelian groups defined by $F(A)=K_i(A \otimes B)$. Then $F$ is split exact, stable and homotopy invariant. Therefore $F(\mathcal{T}_0)=0$, i.e. $K_i(\mathcal{T}_0 \otimes B)=0$. With this in hand, we can complete the proof of Bott periodicity. 
    
    \begin{thm}[Bott periodicity]
    For any $C^{*}$-algebra $B$, we have \[K_0(B) \cong K_1(SB)=K_1(C_0(\mathbb{R})\otimes B).\] 
    \end{thm}
    \textit{Proof.} Let $\sigma: \mathcal{T} \to C(\mathbb{T})$ be the map that sends $v$ to $z$. Denote by $ev_1$, the evaluation map from $C(\mathbb{T}) \to \mathbb{C}$ at $1$. Note that $q=ev_1 \circ \sigma$. We can identify $C_{0}(\mathbb{R})$ with 
    $Ker(ev_1)$. Consequently, we have the following short exact sequence
    \[
    0 \longrightarrow \mathcal{K} \longrightarrow \mathcal{T}_{0} \longrightarrow C_{0}(\mathbb{R}) \longrightarrow 0.\] 
    Tensor the above short exact sequence to obtain the following. 
    \[
    0 \longrightarrow \mathcal{K} \otimes B  \longrightarrow \mathcal{T}_{0} \otimes B  \longrightarrow C_{0}(\mathbb{R}) \otimes B \longrightarrow 0.\]    
    Since $K_0(\mathcal{T}_0 \otimes B)=K_1(\mathcal{T}_0\otimes B)=0$, it follows that the index map \[\partial: K_1(C_0(\mathbb{R})\otimes B) \to K_0(\mathcal{K} \otimes B) \cong K_0(B)\] is an isomorphism. This completes the proof. \hfill $\Box$
    
    \begin{rmrk}
    In fact, the statement of Bott periodicity is a bit more. Bott periodicity gives an explicit map from $K_0(B) \to K_1(SB)$. We explain this for unital $C^{*}$-algebras. Let $B$ be a unital $C^{*}$-algebra. Then
    \[
    M_{n}((SB)^{+}):=\{f:\mathbb{T} \to M_n(B): f(1)_{ij} \in \mathbb{C}1_B\}.\]
    For a projection $p \in M_n(B)$, let $f_p: \mathbb{T} \to M_n(B)$ be defined by
    \[
    f_p(z)=zp+1_n-p.\]
    Then $f_p$ is a unitary in $M_{n}((SB)^{+})$. 
    
    Bott periodicity asserts that there exists a unique map $\beta:K_0(B) \to K_1(SB)$, called the Bott map, which is an isomorphism such that 
    \[
    \beta([p])=[f_p].\]
        If we carefully work through the proofs and unwrap all the identifications, we can prove that the Bott map is indeed an isomorphism. The reader should carry out this verification. 
    
    \end{rmrk}
    
   \begin{rmrk}
   Much of the material on $K$-theory is based on the lectures given by Cuntz during a conference held at Oberwolfach in 2014.    
   \end{rmrk}

\nocite{Arveson_invitation}
\nocite{Raeburn_Williams}
\nocite{Arveson_spectral}
\nocite{Williams_Dana}
\nocite{Cuntz_Meyer}
\nocite{Rordam}
\nocite{Olsen}
\nocite{Faraut_Lie}
\nocite{Davidson_Ken}
\nocite{Murphy_book}

\bibliography{references}
 \bibliographystyle{amsplain}

\end{document}